\pdfoutput=1

\documentclass{book}[a4paper,11pt,openright]
\usepackage{thesispreamble}

\begin{document}

\pagenumbering{gobble}
\singlespacing

~

\begin{center}
    \vspace{\fill}

    \textbf{\Huge \thesistitle}
    
    \vspace{1cm}
    
    \textbf{\LARGE \thesissubtitle}
    
    \vspace{4cm}
    
    \textbf{\Large Joe Allen}

    \
    
    {}

\end{center}

\vspace{\fill}
\noindent\copyright~2021 CC BY-NC-ND 4.0\hspace{\fill}approx. 64~000 words

\newpage
\pagenumbering{roman}
\setcounter{page}{1}
\addcontentsline{toc}{chapter}{Frontmatter}

\section*{Abstract}
\addcontentsline{toc}{section}{Abstract}
We present discrete models of special biserial (SB) algebras and their string modules, drawing inspiration from cellular automata, and cast new light on patterns among syzygies. We explore applications of our models to open questions in homological algebra regarding certain triangulated subcategories of derived categories, with implications for the finitistic dimension conjectures.

More pertinently, our models provide the inner workings for a new, original \GAP\ package called \sbstrips, written and implemented by the author. Its source code is freely available online and its documentation is included as an appendix. The package calculates syzygies $\syzygy^k X$ of string modules $X$ (and much more besides) using specialised methods much more efficient than the general methods currently employed by the \QPA\ package.

\section*{About this document}
This is an arXiv-friendly version of a PhD dissertation submitted to the University of Bristol. The prose is the same but the line-spacing is smaller (and the figures have been moved around to accommodate), so you should print this version but cite the original. At the time of writing, the original may be found at 
\begin{center}
\noindent\href{
https://research-information.bris.ac.uk/en/studentTheses/syzygies-of-string-modules-for-special-biserial-algebras}{
\small
\mono{https://research-information.bris.ac.uk/en/studentTheses/ syzygies-of-string-modules-for-special-biserial-algebras}.}    
\end{center}
Alternatively, the author suggests you look online for "Bristol student theses" and then use his name (Joe Allen) and/or this thesis' title (\thesistitle).

\newpage
\section*{Dedication}
\addcontentsline{toc}{section}{Dedication and acknowledgements}

I dedicate this to my grandparents, Sheila and John, Bessie and Sydney, with love. I miss you every day.
\begin{center}
\textit{\noindent While three chambers of my heart\\
beat true and strong with love for another,\\
the fourth,\\
the fourth is yours forever.}
\end{center}

% \vspace{\fill}

\section*{Acknowledgements}

I never understood why others prefaced their work with such mawkish, grandstanding acknowledgements. But then I wrote a thesis. Now I see.

My foremost debt of gratitude goes to my supervisor Jeremy Rickard, whose extensive knowledge and advice I have heartily welcomed for almost five years. You have left me a great deal wiser. I hope you are not long haunted by the image of my endless coloured pens, highlighters and other stationery.

I thank Peter Green, John Mackay, Fatemeh Mohammadi, Mark Hagen, Fionnuala Hill, Alicia Midgeley, Christina Rozon and the rest of the admin team for their behind-the-scenes support of my studies. 

I recognise the funding from EPSRC that enabled this research, as well as the NPIF allocation which paid for my internship at the Archive Trust for Research. I am grateful to Michael Wright and Edmund Davies at the Trust for hosting me for two months and to Aby Sankaran, formerly of the Bristol Doctoral College, who organised the opportunity.

Without doubt, the greatest privilege of my PhD has been sharing the experience with such inspiring, creative and generous peers. It fills me with tremendous hope to know that the future is being written by such great minds and kind souls as you. Any list I attempt will be incomplete but, nonetheless, I will name and acclaim Emma Bailey, Hazel Brickhill, Paisley Carter, Jenny Chakravarty, Kieran Child, Sam Colvin, Elisa Covato, Charley Cummings, Stamatis Dimopoulous, Scott Harper, Ayesha Hussain, Mai Thanatkrit Kaewtem, Luke Kershaw, Kathryn Leeming, Chris Lutsko, Dan Saattrup Nielsen, Simon Peacock, Felipe P\'erez, Wrenna Robson, Nick Rome, Sophie Stevens and Sophie Watson.

My time in Bristol was not entirely spent with my nose in a book (sorry Jeremy); I tried my hand at much else that the University, the Students' Union and the city had to offer. In particular, I would like to give thanks to Angie Edwards and the jazztinas at AE Dance + Fitness, who taught me joy, community and fierceness through movement, to Robiu Salisu, Caitlin Flint and Ann O'Malley, who helped me in post as Faculty of Science PGR Representative, and to Anna M{\"u}ller, for many uplifting lunchtime chats at French Fridays. I additionally thank my various flatmates during my time at Bristol, including Lizzie, Luke, Hannah and Serena, for responding to my burnout with compassion and patience.

Constant throughout this journey has been the love and support from my nearest and dearest. Trite as it sounds, it is no exaggeration that I could not have managed this without you. I thank my parents, Ruth and Steve, and my brother, Dan, for their unflinching belief in me. I thank my darling Kate, with whom it has been a delight but no surprise to fall in love during my studies, as well as her parents Andrea and Richard, who opened up their home to me and made me feel so welcome. I thank my closest and most enduring friends, Sam, Jim, Liam and Monty: brothers to me, then, now and always. And I thank Dan, Ana, Roopesh, Charley, Abigail and Irena for their honest friendship and encouragement, from one university to the next.

Finally, crucially, I would like to recognise that this thesis has been obtained at great cost to my mental health and wellbeing. The causes are diverse. One is the traumatic experience of surviving a global epidemic (not yet over). Another, more structural, is the intense pressure and unsure footing of a research degree. 

It should be hard, but not this hard.

I encourage all of my fellow PhD students~-- and any human being in general~-- to ask for the help they need when they need it. I did. I would do it again. And, in fact, I wish I had done it sooner. To this end, I thank Debbie Price at the University of Brstiol Student Wellbeing Service and Rachel Grimstead at the Student Counselling Service, and Chaim Golker at Hertfordshire Partnership University NHS Foundation Trust.

\newpage
\setcounter{tocdepth}{2}
\tableofcontents
\listoffigures

\newpage
\pagenumbering{arabic}
\setcounter{page}{1}

\chapter{Introduction}

Special biserial (SB) algebras arose in the twentieth century in the modular representation theory of groups. History has shown they are easy to define and nontrivial to study, yet still prove a rich source of interesting behaviour.

They have been a frequent testing ground for the latest developments through the years in the representation theory of finite-dimensional algebras. An early example: shortly after Drozd proved his famed trichotomy result about representation type \cite{Dro80}, SB algebras were duly proven to be tame algebras by Wald and Waschb{\"u}ch~\cite{WW85} using the techniques of Gel'fand and Ponomarev~\cite{GP68}. More recently, Huisgen-Zimmermann~\cite{HZ16} demonstrated that these representation-tame algebras nonetheless display behaviour that is not "homologically-tame", this meaning that there exist SB algebras on which the difference between big and little finitistic dimensions is arbitrarily large. Schr{\"o}er~\cite{Sch00} gives an excellent historic overview of additional uses of SB algebras.

SB algebras are perhaps the best-studied class of algebras for which two homological questions remain open. These questions concern the satisfaction of properties nicknamed \emph{injectives generate} (IG) and \emph{projectives cogenerate} (PC) and are both phrased in terms of triangulated subcategories of the derived category \cite{Ric19}. The initial aim of this thesis was to prove IG and PC for SB algebras; the final aim was the following conjecture.
\begin{quote}
    \noindent\textbf{Conjecture.\hspace{1ex}} \textit{Over SB algebras, projective modules have finite cosyzygy type and injective modules have finite syzygy type.}
\end{quote}
These two conditions respectively imply IG and PC (but neither is a necessary condition). Extensive calculations on the author's part have unearthed no counterexamples.

These two conjectured conditions are equivalent to one another when dealing with finite-dimensional modules (as we mostly do). It suffices to focus on only one; we choose the latter for the superficial reason that syzygies seem to appear in the literature more frequently than cosyzygies. That being said, we do not believe the previous exploration of syzygy behaviour over SB algebras to be truly thorough. It appears to this author that researchers in the field currently lack the requisite conceptual and computational tools.

To fill this absence, this thesis provides a novel, original approach to string modules and their syzygies over SB algebras that draws inspiration from cellular automata. This means introducing a model for this process and~-- more philosophically but no less importantly~-- giving a language with which to recognise and describe syzygy patterns and phenomena. Our model formalises and disambiguates the informal graphical approaches that have historically sufficed. We strongly recommend the reader consult Subsection \ref{subsec:preview-of-strips} for an extended, illustrative example of our model in action.

This formality has enabled us to implement the key parts of our model in the \GAP\ programming language as the \sbstrips\ package, for which the source code and download instructions are freely available online \cite{All21}. We highlight that the specialised algorithms underlying our package dramatically outperform the current syzygy techniques implemented for general finite-dimensional algebras by the \QPA\ package for \GAP. Our package also implements other useful tools for string modules in SB algebras.

The main results of this thesis are as follows. (For reference, the \define{pin graph} $\pingraph$ is a quiver we associate to a SB algebra $A$, whose vertex set equals that of the defining quiver $Q$ of $A$.) Firstly, we prove the above conjecture for SB algebras $A$ under certain conditions.
\begin{quote}
    \noindent\textbf{Theorem \ref{thm:explicit-bounds-on-syzygy-type}\ref{subthm:syzygy-finite-strings-when-pin-graph-acyclic}.\hspace{1ex}} \textit{The conjecture holds if $A$ has acyclic pin graph.}

    \noindent\textbf{Proposition \ref{thm:ig-and-pc-if-A-has-2-simples}.\hspace{1ex}} \textit{The conjecture holds if $A$ has at most two simple modules.}
\end{quote}
Secondly, we prove a result very specific to our framework: in a technical sense, the symbols with which we denote injective (string) modules appear inconsistent with those modules having infinite syzygy type.
\begin{quote}
    \noindent\textbf{Proposition \ref{prop:descent-orbits-of-inj-sylls-do-not-have-interior-tails}.} \textit{No tail of a $\descent$-orbit of an injective syllable contains only interior syllables.}
\end{quote}
For the final result, let us mention again that the IG and PC properties are phrased in terms of particular triangulated subcategories of the derived category. We show that certain simple modules belong to these subcategories under mild assumptions.
\begin{quote}
    \noindent\textbf{Theorem \ref{thm:acyclic-pin-simples-for-2-reg-quiv-in-subcats}.} \textit{If $A=\kk Q/I$ is special biserial, if $Q$ is $2$-regular and if $\Psi$ is an acyclic connected component of the pin graph of $A$ having vertex set $\Psi_0$, then $S_i \in \Coloc(\Projcat A)$ and $S_i\in \Loc(\Injcat A)$ for all $i\in \Psi_0$.}
\end{quote}

In Chapter \ref{chap:background} we give relevant background material. Nothing here is original. Our purpose is to contextualise this thesis and establish notation.

In Chapter \ref{chap:perm-data-and-syllables}, we provide a combinatorial model for SB algebras in terms of posets of paths on certain 1-regular quivers. Our approach makes clear the connection with combinatorial tools used in the study of Nakayama algebras.

In Chapter \ref{chap:patches-strips-fabric}, we provide a model for string modules and their syzygies over SB algebras. The syzygies of a SB algebra are arranged as (rigid) arrays of symbols. Strips~-- our formulation of string graphs/modules~-- are the rows of this array. The columns are, thanks to our formalism, largely just the orbits of a single partial operation on a finite set.

This reduction to finiteness is employed in Chapter \ref{chap:patterns-and-applications}, where we prove the main results mentioned above as well as some related results that appear stated or implied in the literature but whose proofs are nowhere to be found (until now).

One of our results concerns SB algebras with at most 2 simples. This requires a brute force calculation that encompasses the analysis of dozens of cases. The gory details are relegated to Appendix \ref{app:sb-algs-with-few-verts}. Appendix \ref{app:sbstrips-documentation} is then the documentation of the most recent version of our \sbstrips\ package. Aside from the requisite information befitting a technical manual, it also features an extended chapter of worked examples. These demonstrate the capabilities of \sbstrips~-- which extend notably beyond just the syzygies treated in the thesis~-- as well as the integration of \sbstrips\ with the \QPA\ package.

\chapter{Background}\label{chap:background}

Our ambient set theory is ZFC (in case it matters). We fix once and for all an algebraically closed field $\kk$ (of arbitrary characteristic).

We write \define{$ 2^X $} for the power set of $X$, we write \define{$X \setminus Y$} for the set difference $\{x \in X \colon x \notin Y \}$ and we write \define{$\omega$} for the cardinality of the natural numbers $\mdefine{\naturals} \ceq \{0, 1, 2, \dots\}$.

Where possible we write $fg$ for the composition of functions $\ilarrow{X}{f}{Y}$ and $\ilarrow{Y}{g}{Z}$ and $xf$ for the image in $f$ of an element $x$; this makes $xfg$ the image in the composition. However we defer to the more traditional notations $f(x)$ when convenient or when firmly entrenched in the literature, an example being $\source( \alpha )$ and $\target( \alpha )$ for the sources and target of an arrow $\alpha$ in a quiver. We accordingly write $g \circ f$ for the composition and $(g \circ f) ( x )$ for the image. We try to feature as few hybrids of these conventions as we can; when they do appear, careful use of brackets will always make the order of composition unambiguous.

Our categories are locally small; that is, the class $\categoryfont{C}(X,Y)$ of morphisms of any two objects $X,Y$ of $\categoryfont{C}$ is a set.

\section{Miscellaneous mathematical prerequisites}

\subsection{Partial functions}

\npass Partial functions can be described in two synonymous ways, each with its own jargon and notation, both encapsulating the notion of being ``functions not necessarily defined on all domain elements''.

Since we use them both interchangeably in this thesis, we outline both approaches below.

To be more specific we present two categories, the morphisms of either of which we believe deserve to be called partial functions. We demonstrate the well-known equivalence between them.

Both of these categories are underpinned by the familiar category $\setcat$ of sets. Recall that in $\setcat$ we can write the composition ``$f$ then $g$'' as $fg$ and the image of $x$ in $f$ as $xf$, or we can write that composition as $g \circ f$ and that image as $f(x)$. We prefer the first notation.

\nlpass{Support approach to partial functions} Let $\partsetcat$ be the category with the following data.
\begin{enumerate}
    \item The objects are sets.
    
    \item A morphism $\ilpart{X}{f}{Y}$ is a quadruple $(X,Y,\supp f, F)$ of sets comprising the \emph{domain} $X$, the \emph{codomain} $Y$, the \define{support} $\supp f \subseteq X$ and the \emph{graph set} $F \subseteq ( \supp f ) \times Y$. The graph set has the property that for every $x \in \supp f$ there is a unique $y \in Y$ with $(x,y) \in F$. We call this $y$ the \define{image} of $x$ and denote it $xf$.
    
    The \define{supported} elements of the domain of $f$ are those in $\supp f$.
    
    The identity morphisms $\ilpart{ X }{ \id }{ X }$ have support $X$ and graph set $\{(x,x) \colon x \in X\}$.
    
    \item The composition of two morphisms $\ilpart{X}{f}{Y}$ and $\ilpart{Y}{g}{Z}$, having data $f = (X,Y,\supp f, F) $ and $ g = (Y, Z, \supp g, G)$ respectively, is the morphism $\ilpart{X}{fg}{Z}$ with domain $X$, codomain $Z$, support
    $
    \supp fg
    \ceq
    \{
        x \in X
        \colon
        x \in \supp f \text{ and } xf \in \supp g
    \}
    $
    and graph
    $$
    \big\{
      \big(x,(xf)g\big) \in (\supp fg ) \times Z \colon x \in \supp fg
    \big\}\tstop
    $$
\end{enumerate}

\nlpass{Basepoint approach to partial functions} Let $\basedsetcat$ be the category with the following data.
\begin{enumerate}
    \item The objects are pairs $(X,x_\star)$ comprising a set $X$ and a member $x_\star \in X$. 
    
    We call the pair a based set and $x_\star$ the \define{basepoint} or \emph{distinguished element}.
    
    \item A morphism $\ilarrow{ (X,x_\star) }{ f }{ ( Y,y_\star) }$ is a function $f$ (ie, a morphism of $\setcat$) such that $x_\star f = y_\star$.
    
    The \define{support $\supp f$} of a morphism $\ilarrow{ ( X, x_0 )}{ f }{ ( Y, y_0 ) }$ is $\{ x \in X \colon xf \neq y_\star \}$.
    
    The identity morphisms $\ilarrow{ ( X, x_\star ) }{}{ ( X, x_\star ) }$ are those of $\setcat$.
    
    \item The composition of $\ilarrow{ (X,x_\star) }{ f }{ (Y, y_\star) }$ and $\ilarrow{ ( Y, y_\star ) }{ g }{ ( Z, Z_\star ) }$ is then $\ilarrow{ (X,x_\star) }{ fg }{ (Z, Z_\star) }$; ie, $\partsetcat$ inherits its composition law from $\setcat$.
\end{enumerate}

\nlpass{Equivalence of $\partsetcat$ and $\basedsetcat$} Any set $X$ gives rise to a based set $( X \cup \{x_\star\}, x_\star )$ by freely adjoining a new element $x_\star \notin X$, for example $x_\star \ceq X$ itself, to ensure that the union is disjoint. (That $X \notin X$ is a consequence of the foundation axiom for sets.) Any morphism $\ilpart{X}{f}{Y}$ of $\partsetcat$ gives rise to a morphism $\ilarrow{ ( X \cup \{x_\star\}, x_\star ) }{ f' }{ ( Y \cup \{y_\star\}, y_\star ) }$  in $\basedsetcat$ by defining $ xf \ceq y_\star $ for all $x \notin \supp f$. The assignments induce a functor $\ilarrow{ \partsetcat  }{}{ \basedsetcat }$.

Conversely, to a based set $(X,x_\star)$ we assign the set $X \setminus \{ x_\star \}$, and to a morphism $\ilarrow{ ( X, x_\star ) }{ f }{ ( Y, y_\star ) }$ of $\basedsetcat$ we assign the morphism
$\big(
  X \setminus \{ x_\star \},
  Y \setminus \{ y_\star \},
  \supp f,
  \big\{ (x, xf) \colon x \in \supp f \big\}
\big)
$
of $\partsetcat$. These too give rise to a functor $ \ilarrow{ \basedsetcat }{ }{ \partsetcat } $.

These functors are mutually inverse equivalences of categories \cite[Example 1.5.6]{Rie14}.

\nlpass{Defining a partial function} Now that we are acquainted with  $ \partsetcat $ and $ \basedsetcat $, we can stipulate that by \define{partial function} we formally mean a morphism in either category. Partial functions modelled one way can be translated into one modelled the other way using the equivalences outlined above. We do this freely.

\nlpass{Particular partial functions} To obtain ``partial'' variants of familiar concepts like surjectivity, we describe supported elements rather than domain elements.

A partial function $\ilpart{X}{f}{Y}$ is a \define{partial injection} if $xf\neq x'f$ for any $x,x' \in \supp X$ with $x \neq x'$. It is a \define{partial surjection} if for any $y \in Y$ there exists $x \in \supp f$ with $y = xf$.

A partial function $\ilarrow{ X }{ g }{ X }$ is a \define{partial involution} if for all $x\in \supp g$ we have  $xg\in \supp g$ and $x g^2 = x$. It is \define{fixpointfree} if further there is no $x \in \supp g$ with $x g = x$.

\subsection{Order theory}

\nlpass{Orders} A \define{partial order $\leq$}, or simply \define{order}, on a set $X$ is a binary relation such that for all elements $x,x',x'' \in X$ we have
\begin{enumerate}
    \item $ x \leq x $ (\emph{reflexivity}),
    \item $x \leq x'$ and $x' \leq x$ implies $x = x'$ (\emph{antisymmetry}), and
    \item $ x \leq x' $ and $x' \leq x'' $ implies $x \leq x'' $ (\emph{transitivity})
\end{enumerate}
and, if additionally
\begin{enumerate}
    \setcounter{enumi}{3}
    \item $x \leq x'$ or $x' \leq x$ (\emph{totality})
\end{enumerate}
for all $x,x' \in X$, then it is a \define{total order}. A (\define{totally-} or \define{partially-}) \define{ordered set} is a pair $(X, \leq)$ comprising a set $X$ and a (total or partial) order $\leq$ on it. We sometimes abbreviate \emph{totally-} and \emph{partially-ordered set} to \define{toset} and \define{poset}.

(We normally use the terms \emph{partial order} and \emph{partially-ordered set} to emphasise that the order is not total and use \emph{order} and \emph{ordered set} when the distinction is unimportant. Any relation denoted $\leq$ is an order.)

Being no more perfect than our forebears we indulge in common abuses, such as writing $x \geq x'$ to mean $x' \leq x$ when convenient, writing $x < x'$ to mean that $x \leq x'$ and $ x \neq x' $ (and similarly $x > x'$), and referring just to the set $X$ of a pair $(X,\leq)$ as the poset.

Two elements $x,x'$ of a poset are \define{comparable} if $x \leq x'$ or $x' \leq x$, and \define{incomparable} otherwise. Comparability is a symmetric relation; the equivalence classes it generates are the \define{comparability components} of the ordered set, or simply \define{components} for short.

Subsets of a poset inherit its order in an obvious fashion. 

\nlpass{Poset elements and their relations} An element $x$ of a poset $X$ is
\begin{enumerate}
    \item
    \begin{enumerate}
        \item \define{minimal} if there is no $x' \in X$ with $x' < x$,
        \item \define{maximal} if there is no $x' \in X$ with $x < x'$,
        \item \define{extremal} if it is minimal or maximal,
        
    \end{enumerate}
    
    \item\label{item:poset-extrema}
    \begin{enumerate}
        \item a \define{minimum} if $x \leq x'$ for all $x'\in X$, in which case we denote it \define{$\min X$},
        \item a \define{maximum} if $x' \leq x$ for all $x'\in X$, in which case we denote it \define{$\max X$},
        \item an \define{extremum} if it is a minimum or a maximum.
    \end{enumerate}
    
    (The plural forms of these terms are \emph{minima}, \emph{maxima} and \emph{extrema}.)
\end{enumerate}
We see that any poset has at most one minimum element and at most one maximum element. 

If $x < x'$ and there is no $y \in X$ with $x < y < x'$, then we say that \define{$ x $ is covered by $x'$} and, reciprocally, \define{$x'$ covers $x$}.

We will consider several posets that readily lend themselves to illustration, such as the set of paths in a given quiver. These can be visualised in a \emph{Hasse diagram}: a quiver (oriented graph) with members of the poset for vertices and with arrows $\ilarrow{x}{}{y}$ when $x$ is covered by $y$. Figure \ref{fig:hasse-diagram-of-paths-in-ZZ} gives an example.

\nlpass{Bounds and divergence} Suppose $(X, \leq)$ is an ordered set and $Y \subseteq X$. We say $x \in X$ is an \define{upper bound for $Y$} if $y \leq x$ for all $y \in Y$. If $Y$ has an upper bound, then we say $Y$ is \define{bounded above} (by any of its upper bounds). The \define{supremum $\sup(Y)$} of $Y$ is $\min\{ x \in X \colon x \text{ is an upper bound for } Y \}$ if this minimum exists. We define \define{lower bounds}, \define{boundedness below} and the \define{infimum $\inf Y$} of $Y$ dually, and we say that $Y$ is \define{bounded} if it is bounded above and below.

Suppose that $X$ is the toset $\naturals$. We say that a sequence $(x_k)_{k \geq 0}$ of positive integers is \define{bounded} if $\{x_k \geq 0 \colon k \geq 0\}$ is bounded, and we say the sequence \define{diverges to $+\infty$} if for all $N \geq 0$ there exists $K > 0$ such that $ x_k \geq N $ for all $k \geq K$. 

\nlpass{Chains and antichains} A subset $Y \subseteq X$ of a poset $X$ is a \define{chain} if all elements of $Y$ are comparable (in the inherited order) and an \define{antichain} if no two distinct elements are comparable.

\nlpass{Order ideals} An \define{(order) ideal} $I \subseteq X$ of a poset is a subset such that if $x,x' \in X$ satisfy $x' \in I$ and $x \leq x'$, then $ x \in I $.

\nlpass{Monotone functions} A (set-theoretic) function $\ilarrow{ X }{ f }{ Y }$ between ordered sets is (\define{weakly}) \define{monotone} if $x \leq y$ implies $xf \leq yf$; it is \define{strictly} so when $x < y$ implies $xf < yf $.

Again, we only use the term \emph{weakly} when underscoring that the monotone function is not strictly monotone.

\subsection{Graph theory}

We define graphs and associated notions. Note that our framework permits loops and multiple edges in a graph.

\nlpass{Graphs} A graph $\Gamma$ is a triple $\Gamma = (\Gamma_0, \Gamma_1, \ee )$ comprising a set $\Gamma_0$ of \define{vertices}, a set $\Gamma_1$ of \define{edges} and an \emph{incidence function $\ee$}; that is, a function $\ee \colon \ilarrow{ \Gamma_1 }{}{2^{\Gamma_0} }$ such that $1 \leq |\ee(e)| \leq 2$ for all edges $e \in \Gamma_1$.

The graph is \define{finite} if both $\Gamma_0$ and $\Gamma_1$ are finite.

\nlpass{Incidence, adjacency and loops} We say that an edge $e \in \Gamma_1$ is \emph{incident} to $v \in \Gamma_0$ if $v \in \ee(e)$. A \define{loop} is an edge with $|\ee(e)|=1$.

Two vertices $v,v' \in \Gamma_0$ are \define{adjacent} if there is some edge $e \in \Gamma_1$ with $\ee(e)=\{v, v'\}$.

In the standard formulation, the \define{valency} of a vertex $v$ is $2a+b$, where $a$ is the number of loops and $b$ the number of nonloop edges incident to $v$.

\nlpass{Subgraphs} A subgraph $\Delta$ of $\Gamma$ is a graph $(\Delta_0, \Delta_1, \mathrm{d})$ such that $\Delta_0 \subseteq \Gamma_0$ and $\Delta_1 \subseteq \Gamma_1$ and such that   $\mathrm{d}(e) \ceq \mathrm{e}(e)$ is a well-defined function $\mathrm{d}=\mathrm{e}|_{\Delta_1} \colon \ilarrow{ \Delta_1 }{}{ 2^{\Delta_0} }$.

The \emph{full subgraph specified by $V \subseteq \Gamma_0$} is the subgraph $\Delta$ where $ \Delta_0 = V $ and where $ \Delta_1 $ comprises exactly those edges incident only to members of $V$. In general, a subgraph is \define{full} if there is some subset of vertices that specify it in the above sense.

\nlpass{Graph homomorphisms} A \define{graph homomorphism} $f=(f_0, f_1) \colon \ilarrow{ \Gamma }{}{ \Delta }$ is a pair comprising (set-theoretic) functions $ f_0 \colon \ilarrow{ \Gamma_0 }{}{ \Delta_0 } $ and $ f_1 \colon \ilarrow{ \Gamma_1 }{}{ \Delta_1 } $ respecting incidence. This means that the square
$$
\begin{tikzcd}
\Gamma_1 \ar[r, "\ee"] \ar[d, "f_1"]
    \& 2^{ \Gamma_0 } \ar[d, "f_0"]
        \\
\Delta_1 \ar[r, "\mathrm{d}"]
    \& 2^{ \Delta_0 }
\end{tikzcd}
$$
commutes, where we abuse notation and also denote by $f_0 \colon \ilarrow{ 2^{ \Gamma_0 } }{ }{ 2^{ \Gamma_1 } }$ the map induced by $ f_0 \colon \ilarrow{ \Gamma_0 }{}{ \Delta_0 } $.

\nlpass{Paths and cycles} A \emph{path} of length $l \geq 0$ in $\Gamma$ is a sequence $(v_0,e_1,v_1,\dots,e_\ell,v_\ell)$ alternately comprising vertices $v_j$ and edges $e_j$ and such that $\ee(e_j) = \{v_{j-1},v_j\}$ for each $1 \leq j \leq \ell$. (We identify any such path with its reverse $(v_\ell,e_\ell,v_{\ell-1},\dots,e_1,v_0)$, since the graph is undirected.) It is a \emph{simple path} if the vertices $\{v_0, \dots, v_\ell\}$ visited are all distinct.

A \emph{cycle} of length $\ell > 0$ is a path $(v_0,e_1,v_1,\dots,e_\ell,v_\ell)$ whose extremal vertices $v_0$ and $v_\ell$ are the same. It is a \emph{simple cycle} if all of the visited vertices are distinct, save for $v_0=v_\ell$. A graph that has no cycles is \emph{acyclic}.

\nlpass{The integers $\integers$ as a graph} We model the graph
$\big(
\begin{tikzcd}[cramped, sep=small]
\cdots \ar[r, no head]
  \& \syllableplaceholder \ar[r, no head]
    \& \syllableplaceholder \ar[r, no head]
      \& \syllableplaceholder \ar[r, no head]
        \& \syllableplaceholder \ar[r, no head]
          \& \syllableplaceholder \ar[r, no head]
            \& \cdots
\end{tikzcd}
\big)$
using the integers
$\big(
\begin{tikzcd}[cramped, sep=small]
\cdots \ar[r, no head]
  \& (-2) \ar[r, no head]
    \& (-1) \ar[r, no head]
      \& 0 \ar[r, no head]
        \& 1 \ar[r, no head]
          \& 2 \ar[r, no head]
            \& \cdots
\end{tikzcd}
\big)$.

\nlpass{Connectedness for graphs} Two vertices of a graph $\Gamma$ are \define{connected} if they are the extremal vertices $v_0, v_l$ of some path $(v_0,e_1,v_1,\dots,e_l,v_l)$. Connectedness of vertices induces an equivalence relation on the vertex set $\Gamma_0$. The classes of this equivalence relation specify full subgraphs of $\Gamma$ called \define{connected components}. If there is only one equivalence class (and hence only one connected component) then we call the graph \define{connected}.

\nlpass{Rooted trees} A \define{tree} is a graph that is connected and acyclic. Equivalently, it is a graph in which any two vertices are connected by a unique path.

A \define{rooted tree} is a pair $(T,r)$ comprising a tree $T$ and a vertex $r \in T_0$ called the \define{root}. When the root is clear from context then it is a standard abuse of notation to denote $(T,r)$ just by $T$. A homomorphism $\ilarrow{(T,r)}{}{(T',r')}$ of rooted trees is a homomorphism $\ilarrow{T}{}{T'}$ of graphs that sends $\ilmapsto{r}{}{r'}$.

A vertex $v \in T_0$ is in the \emph{$l$th level} of the tree if the unique path connecting it to the root has length $l$. This means that the root is the unique vertex in level $0$.

The \define{children}, or \define{child vertices}, of a vertex in level $l \geq 0$ are those vertices adjacent to it in level $l+1$; the transitive closure of the child relation gives the \define{descendant} relation. The \define{parent}, or \define{parent vertex}, of a vertex in level $l > 0$ is the (unique) vertex adjacent to it in level $l-1$; the transitive closure of the parent relation gives the \define{ancestor} relation. We see that the root is the ancestor of all vertices in a rooted tree.

The \define{valency of a vertex} is the number child vertices it has. The \define{valency of a rooted tree} is the maximal valency of any of its vertices. A rooted tree is \define{regular} if all of its vertices have the same valency.

\nlpass{The regular rooted tree $\integers^*$ of countable valency} Let $\mdefine{\integers^*} \ceq \bigcup_{l \geq 0} \integers^l$ denote the set of finite sequences of integers, which includes the empty sequence $\emptyseq \ceq ()$. We create a rooted tree with root $\emptyseq$ and vertex set $\integers^*$ by specifying that the children of $(x_0, x_1, \dots, x_{l-1})$ are $(x_0, x_1, \dots, x_{l-1}, k)$, for $k \in \integers$.

We easily verify that $\integers^*$ is regular, having valency $\omega$. The length of a sequence gives its level in the tree.

Any rooted tree of finite (or indeed countable) valency embeds as a subgraph of $\integers^*$.

\subsection{Quiver theory}

\nlpass{Quivers} Formally, a \define{quiver} is a tuple $Q=(Q_0, Q_1, \source, \target)$ comprising sets $Q_0$ of \define{vertices} and $Q_1$ of \define{arrows} together with \define{source} and \define{target} functions $\source,\target \colon \ilarrow{Q_1}{}{Q_0}$. It is \define{finite} if $Q_0$ and $Q_1$ are both finite, and \define{locally finite} if every vertex is the source and target of finitely many vertices.

The \define{opposite quiver $Q^\op$} of $(Q_0,Q_1,\source,\target)$ is $(Q_0, Q_1, \target, \source)$. That $(Q^\op)^\op=Q$ is immediate.

\nlpass{Quiver homomorphisms} A \define{quiver homomorphism} $\phi \colon \ilarrow{ (Q_0, Q_1, \source, \target) }{}{ (Q'_0, Q'_1, \source', \target') }$ is a pair of (set-theoretic) maps $\phi_k \colon \ilarrow{ Q_k }{}{ Q'_k }$ ($k \in \{0,1\}$) that commute the squares
$$
  \begin{tikzcd}
  Q_1 \ar[r,"\source"] \ar[d,"\phi_1"]
    \& Q_0 \ar[d,"\phi_0"]
      \\
  Q'_1 \ar[r,"\source'"]
    \& Q'_0
  \end{tikzcd}
  \text{\hspace{1cm}and\hspace{1cm}}
  \begin{tikzcd}
  Q_1 \ar[r,"\target"] \ar[d,"\phi_1"]
    \& Q_0 \ar[d,"\phi_0"]
      \\
  Q'_1 \ar[r,"\target'"]
    \& Q'_0
  \end{tikzcd}\tstop
$$

We call $\phi$ injective or an \define{inclusion} if both $\phi_k$ are injective, in which case we say $Q$ \define{includes into} $Q'$, we write $\phi \colon \ilmono{Q}{}{Q'}$ and call $Q$ a \define{subquiver} of $Q'$. Dually, we call $\phi$ surjective or a \define{projection} if both $\phi_k$ are surjective, in which case we say $Q$ \define{projects onto} $Q'$, write $\phi \colon \ilepi{Q}{}{Q'}$ and call $Q'$ a \define{quotient} of $Q$.

We call $Q=(Q_0, Q_1, \source, \target)$ a \define{full subquiver} of $Q' = (Q'_0, Q'_1, \source', \target')$ if there is some $X \subseteq Q_0$ such that
$$
Q_0 = X\tcomma
\hspace{0.5cm}
Q_1 = \{\alpha\in Q'_1\colon \source(\alpha) \in X \text{ and } \target(\alpha) \in X\}\tcomma
\hspace{0.5cm}
\source = \source'|_{Q_1}
\text{\hspace{0.5cm}and\hspace{1cm}}
\target=\target'|_{Q_1}\tcomma
$$
in which case we call it the full subquiver of $Q'$ \define{specified by $X$}.

\nlpass{Underlying graph} The \define{underlying graph} of $Q=(Q_0,Q_1,\source,\target)$ is $\big(Q_0, Q_1, \ilmapsto{ \alpha }{  }{ \big\{ \source( \alpha ), \target( \alpha ) \big\} }\big)$.

\nlpass{Connectedness for quivers} A \define{connected component} of a quiver $Q$ is a full subquiver of $Q$ specified by the vertices in a connected component of the underlying graph of $Q$. A quiver is \define{connected} if it has one connected component (iff its underlying graph does).

\nlpass{Paths} Informally, a path in a quiver is a list of compatibly oriented arrows starting and ending at some vertices.

Formally, let us extend $\source,\target$ to maps $\source,\target \colon \ilarrow{Q_0\sqcup Q_1}{}{Q_0}$ by specifying $\source(i)\ceq i$ and $\target(i)\ceq i$ for all vertices $i$. Then we formally define a \define{path (in $Q$)} as a finite, ordered tuple $p \ceq (u_0; u_1,\dots,u_\ell; u_{\ell+1})$ with $\ell\geq 0$ such that: $u_0$ and $u_{\ell+1}$ are vertices; all other $u_r$ are arrows; and for all $0\leq r \leq \ell$ we have $\target(u_r)=\source(u_{r+1})$. In this case, the \define{length} of the path is $\ell$, the \define{source} of the path is $u_0$ and the \define{target} of the path is $u_{\ell+1}$. We write \define{$\length p$} for the length of $p$ and, by abuse of notation, \define{$\source(p)$} for \define{$\target(p)$} for the source and target of $p$. We say that $p$ \define{passes through} the vertices $u_0, u_{\ell+1}$ and $\source(u_r),\target(u_r)$  ($0 \leq r \leq \ell $).

We call a path \define{stationary} or \define{trivial} iff it has length $0$. Outside of this preliminary section, stationary paths $(i,i)$ will be denoted just by $e_i$ and nonstationary paths $(u_0;u_1,u_2,\dots,u_\ell;u_{\ell+1})$ by $u_1u_2\cdots u_{\ell}$. Paths should be read left-to-right.

The set of all paths in $Q$ is $\Subpath(Q)$.

\nlpass{Concatenation of paths} Let $p\ceq(u_0;u_1,\dots,u_{\ell};u_{\ell+1})$ and $q\ceq(v_0;v_1,\dots,v_m;v_{m+1})$ be paths. If $\target(p)=\source(q)$, which is exactly to say if $u_{\ell+1}=v_0$, then we define the \define{concatenation of $p$ and $q$} to be the path $pq \ceq (u_0;u_1,\dots,u_{\ell},u_{\ell+1},v_1,\dots,v_m;v_{m+1})$. This has $\source(p)$ as its first entry and $\target(q)$ as its last entry, and combines (in order) the lists of constituent arrows of $p$ and $q$.

When $\target(p)\neq \source(q)$, we leave $pq$ undefined. Observe that $e_iu$ is defined iff $i=\source(u)$, in which case $e_iu=u$. Similarly, $ue_j$ is defined iff $j=\target(u)$, in which case $ue_j=u$. 

\nlpass{Path category}\label{psg:path-categories} The paths in any quiver $Q$ form a category. Specifically, the \define{path category of $Q$} is the (small) category with object set $Q_0$ and morphism set $\Subpath(Q)$. The domain and codomain of any path $p$ are $\source(p)$ and $\target(p)$ respectively. The identity morphisms are the trivial paths $e_i$. Concatenation of paths gives composition.

We comment that the full subcategory on any subset $V \subseteq Q_0$ of vertices contains those paths whose source and target both lie in $V$.

Suppose quivers $Q$ and $Q'$ respectively have path categories $\categoryfont{Q}$ and $\categoryfont{Q'}$. Any morphism $\ilarrow{Q}{}{Q'}$ of quivers induces a functor $\ilarrow{ \categoryfont{Q} }{}{ \categoryfont{ Q' } }$, however not all functors between path categories arise in this way. A trivial example, provided that $Q$ has at least one arrow and $Q'$ at least one vertex, is any functor $\ilarrow{ \categoryfont{Q} }{}{ \categoryfont{Q'} }$ that maps all $Q$-paths to some stationary $Q'$-path. For a slightly less trivial example, let $Q \ceq \big(\begin{tikzcd}
1 \ar[r, "\alpha"] \& 2
\end{tikzcd}\big)$ and $Q' \ceq \big(\begin{tikzcd}
1 \ar[r, "\beta_1"] \& 2 \ar[r, "\beta_2"] \& 3
\end{tikzcd}\big)$, and consider the path category functor sending $\ilmapsto{ \alpha }{}{ \beta_1\beta_2 }$ (and therefore $\ilmapsto{e_1}{}{e_1}$ and $\ilmapsto{e_2}{}{e_3}$).

\nlpass{Subpaths, prefixes and suffixes} In general, when paths $p,q,u,v$ satisfy $q=upv$ then $p$ is a \define{subpath} of $q$ and $q$ is a \define{superpath} of $p$. Some particular subpaths of $q$ are especially important.

The \define{prefixes} of a path $q$ are those paths $u$ for which there exists a path $v$ with $uv=q$, in which case we call $v$ the \define{prefix complement} of $u$ in $q$. \define{Suffixes} and their \define{(suffix) complements} are defined dually.

\nlpass{Subpath order on $\Subpath(Q)$} The relation ``is a subpath of'' is a partial order on the $\Subpath(Q)$. The minimal elements of this poset are the stationary paths, and this poset is finite iff $Q$ is acyclic and finite.

%In the subpath poset of a locally-finite quiver, any path $p$ is covered by finitely many paths. These will be $\alpha p$ and $p\beta$, for $\alpha,\beta$ satisfying $\target(\alpha)  = \source(p)$ and  $\source(\beta) = \target(p)$.
We may grade $\Subpath(Q)$ by path length; that is, $p \leq q$ implies $\len p \leq \len q$.

\centresubsec{Regular and subregular quivers}

\nlpass{Regularity and subregularity} Let $m$ be a positive integer. A quiver is \define{sub-$m$-regular} iff
$
\text{for all }
i \in Q_0 
\text{ we have }
\big|
    \source^{-1}\big(
        \{i\}
    \big)
\big|
\leq
m
\text{ and }
\big|
    \target^{-1}\big(
        \{i\}
    \big)
\big|
\leq m
$,
and it is \define{$m$-regular} iff the above inequalities may both be replaced by equalities. Here, $\source^{-1}$ and $\target^{-1}$ denote the preimage set.

\nlpass{Lemma (augmenting sub-regular to regular)} \textit{Any finite, connected sub-$m$-regular quiver $Q$ may be augmented to a connected $m$-regular quiver by adding only arrows.}

\proof Since $Q$ is connected and sub-$m$-regular, $Q_1$ satisfies $ |Q_0|-1 \leq |Q_1| \leq m|Q_0|$. We then define $r_Q \ceq m|Q_0|-|Q_1|$; this nonnegative integer is $0$ iff $Q$ is $m$-regular.

If $r>0$, then a counting argument establishes the existence of vertices $i$ and $j$ with $|s^{-1}(i)|< m$ and $|t^{-1}(j)| < m$. Choose such an $i$ and $j$ and augment $Q_1$ with an extra arrow $\ilarrow{i}{}{j}$. Denoting the augmented arrow $Q'$, we find that
$
r_{Q'} = m|Q'_0|-|Q'_1|=m|Q_0|-(|Q_1|+1) = r_Q-1
$.
Finitely many such augmentations yield an $m$-regular quiver of which $Q$ is a subquiver, as required. \qed

\nlpass{Remark (regularity and opposite quivers)} A quiver $Q$ is $m$-regular or sub-$m$-regular iff its opposite $Q^\op$ is. If $Q'$ is a $m$-regular augmentation of $Q$, then $(Q')^\op$ is one of $Q^\op$.

\nlpass{Example (almost-uniqueness of regular augmentation)}\label{psg:almost-uniqueness-of-regular-augmentations} Consider the quiver $\big(\begin{tikzcd}1 \ar[r, very thick] \& 2\end{tikzcd}\big)$. It has a unique $1$-regular augmentation $\big(\begin{tikzcd}1 \ar[r, very thick, shift left] \& 2 \ar[l, shift left]\end{tikzcd}\big)$. It has two $2$-regular augmentations:
\begin{center}
    $
    \begin{tikzcd}[]
    1 \ar[r, very thick, shift left] \ar[loop left]
     \& 2 \ar[l, shift left] \ar[loop right]
    \end{tikzcd}
    $
    \hspace{1cm}
    and
    \hspace{1cm}
    $
    \begin{tikzcd}[]
    1 \ar[r, very thick, shift left] \ar[r, bend left=50]
     \& 2 \ar[l, shift left] \ar[l, bend left=50]
    \end{tikzcd}
    $\tstop
\end{center}
In fact, it is a subquiver of every (connected) $m$-regular quiver on $2$ vertices. One easily shows that these have adjacency matrices of the form $\left[\begin{smallmatrix}(m-r) & r\\ r & (m-r)\end{smallmatrix}\right]$ for $1\leq r \leq m$ and, therefore, that there are $m$ of them.

This exemplifies the general truth that a sub-$m$-regular quiver does not necessarily have a unique $m$-regular augmentation. Nonetheless, there are only finitely many. % We call this finite divergence from unicity \define{almost uniqueness}.

\centresubsec{Sub-$1$-regular quivers}

\nlpass{Type $\AAA$ and $\widetilde{\AAA}$ quivers} Fix some positive integer $n>0$. There are two isomorphism classes of connected sub-$1$-regular quivers having $n$ vertices, namely those represented by
$$
\begin{tikzcd}[sep=small]
    1 \ar[r, leftarrow]
      \& 2 \ar[r, leftarrow]\ar[rr, phantom, ""{coordinate, name=C}]
        \& \cdots \ar[r, leftarrow]
          \& n-1 \ar[r, leftarrow]
            \& n \ar[llll, leftarrow, rounded corners, to path={
                    -- ([xshift=1.5ex]\tikztostart.east)
                    |- ([yshift=-2ex]C)[pos=1]\tikztonodes
                    -| ([xshift=-1.5ex]\tikztotarget.west)
                    -- (\tikztotarget.west)
                    }]
\end{tikzcd}
\text{\hspace{1cm}and\hspace{1cm}}
\begin{tikzcd}[sep=small]
1 \ar[r, leftarrow]
  \& 2 \ar[r, leftarrow]\ar[rr, phantom, ""{coordinate, name=C}]
    \& \cdots \ar[r, leftarrow]
      \& n-1 \ar[r, leftarrow]
        \& n
\end{tikzcd}\tstop
$$
In the literature these are respectively called the equioriented type $\widetilde{\AAA}_{n-1}$ and type $\AAA_n$ quivers. The former is $1$-regular and the latter is not. Indeed, the former is the unique $1$-regular augmentation of the latter.

% In this thesis we will frequently deal with sub-$1$-regular quivers. We will represent them as above, using integers for vertices and drawing arrows $\ilarrow{i}{}{(i-1)}$ (vertices considered modulo some $n$ when appropriate).

\nlpass{Alternative characterisation of $1$-regularity} An equivalent definition of the $1$-regularity of $Q$ is that its source and target maps $\source$ and $\target$ be bijections $\ilarrow{Q_1}{}{Q_0}$. When so, they have respective inverses $\source^{-1}$ and $\target^{-1}$ and we thus obtain permutations $ (\target \circ \source^{-1}) \colon \ilarrow{ Q_0 }{}{ Q_0 } $ and $ ( \source^{-1} \circ \target ) \colon \ilarrow{ Q_1 }{}{ Q_1 } $.

\nlpass{The quiver $\integers$} The quintessential $1$-regular quiver is
$
\big( \begin{tikzcd}[sep=small]
\cdots
  \& -2 \ar[l]
    \& -1 \ar[l]
      \& 0 \ar[l]
        \& 1 \ar[l]
          \& 2 \ar[l]
            \& \cdots \ar[l]
\end{tikzcd} \big)
$, which we naughtily call $\integers$.

Powers of $(\target \circ \source^{-1})$ or $(\source^{-1} \circ \target)$ give $\integers$-actions on the vertices and arrows, moving things ``with the flow of the quiver''. We write these actions additively and on the right so that, for $k\in\integers$, $\mdefine{i-k} \ceq (\target \circ \source^{-1})^k(i)$ and $\mdefine{\alpha-k} \ceq (\source^{-1} \circ \target)^k(\alpha)$. The pictures to keep in mind are the following:
$$\begin{tikzcd}[row sep=tiny]
    \cdots 
     \& (i-2) \ar[l]
      \& (i-1) \ar[l]
       \& i \ar[l]
        \& (i+1) \ar[l]
         \& (i+2) \ar[l]
          \& \cdots \ar[l]
  \end{tikzcd}$$
and
  $$\begin{tikzcd}[row sep=tiny, column sep=large]
    \cdots \ar[from=r, "\alpha-2"']
     \& \cdot \ar[from=r, "\alpha-1"']
      \& \cdot \ar[from=r, "\alpha"']
       \& \cdot \ar[from=r, "\alpha+1"']
        \& \cdot \ar[from=r, "\alpha+2"']
        \& \cdots
  \end{tikzcd}\tstop
  $$
This notation for the action is chosen so as to be compatible with any sensible identification of the arrows of vertices of the quiver with $\integers$.

\nlpass{Paths in $\integers$} Three pieces of information associated to each $p \in \Subpath(\integers)$ are its source $\source(p)$, target $\target(p)$ and length $\len p$. Any two of these determine the third since $\source(p) - l = \target(p)$.

We write \define{$\big(\ilslpath{i}{{l}}\big)$} for the path in $\integers$ with source $i$ and length $l$. (The brackets are for clarity when written inline rather than a strict part of the notation.) In the much rarer events that we specify a path in $\integers$ by its target and length or its source and target, we will use the notations \define{$\big(\illtpath{l}{i}\big)$} and \define{$\big(\ilstpath{i}{j}\big)$}.

The $\integers$-actions on the vertices of $\integers$ induce one on $\Subpath(\integers)$. An obvious notation would be to write $p-k$ for $\big( \ilstpath{(\source(p)-k)}{(\target(p)-k)}\big)$, where $k\in\integers$, however this is unnecessary for we only mean to mention the existence of the action and show how it follows from the action on vertices.

\nlpass{Subpath poset of $\integers$} The poset $\Subpath(\integers)$ of paths in $\integers$ (ordered by the subpath order) is very transparent. Since $\target^{-1}(\source(p))$ and $\source^{-1}(\target(p))$ are uniquely defined for any path $p$, the whole story is told by ``fundamental diamonds''
$$\begin{tikzcd}[sep=small]
 \& \alpha p \beta
  \&
  \\
\alpha p \ar[ur]
 \&
  \& p \beta \ar[ul]
  \\
 \& p \ar[ul] \ar[ur]
\end{tikzcd}$$
in the Hasse diagram. These fit together as in Figure \ref{fig:hasse-diagram-of-paths-in-ZZ}.
\begin{thesisfigure}
\centering

$$
\begin{tikzcd}[column sep=20]
\&
  \& \vdots
    \&
      \& \vdots
        \&
          \& \vdots
            \&
              \&  \vdots
                \& ~
                  \\
\& ~
  \&
    \& ~
      \& 
        \& ~
          \& 
            \& ~
              \&
                \& ~
                \\
\cdots
\&
    \& \cdot \ar[ur,dotted]\ar[ul, dotted]
      \&
        \& \cdot \ar[ur,dotted]\ar[ul,dotted]
          \&
            \& \cdot \ar[ur,dotted]\ar[ul,dotted]
              \&
                \& \cdot \ar[ul,dotted] \ar[ur, dotted]
                  \& 
                    \& \cdots
                  \\
\& ~ \ar[ur, dotted]
  \&
    \& \cdot \ar[ur]\ar[ul]
      \&
        \& \cdot \ar[ur]\ar[ul]
          \&
            \& \cdot \ar[ur]\ar[ul]
              \&
                \& ~ \ar[ul, dotted]
                  \\
\cdots
\&
  \&
    \cdot \ar[ur] \ar[ul, dotted]
      \&
        \& \cdot \ar[ur]\ar[ul]
          \&
            \& \cdot \ar[ur]\ar[ul]
              \&
                \& \cdot \ar[ul] \ar[ur, dotted]
                  \&
                    \& \cdots
                  \\
\& ~ \ar[ur, dotted]
  \&
    \& e_{i-1} \ar[ul]\ar[ur]
      \&
        \& e_i \ar[ul]\ar[ur]
          \&
            \& e_{i+1} \ar[ul]\ar[ur]
                \&
                  \& ~\ar[ul, dotted]
\end{tikzcd}
$$

\caption[%
  The Hasse diagram of $\Subpath(\integers)$
]{%
  \label{fig:hasse-diagram-of-paths-in-ZZ}
  \textit{The Hasse diagram of $\Subpath(\integers)$.}
}

\end{thesisfigure}

\nlpass{Notation for other connected $1$-regular quivers} Any finite, connected $1$-regular quiver $Q$ is a quotient of $\integers$. This means we can represent vertices, arrows or paths of $Q$ using those in $\integers$. Context will always clarify whether the source and target vertices of $\big(\ilstpath{i}{j}\big)$ are to be understood modulo $n$ and, when so, what the value of $n$ is.

\nlpass{Notation for general $1$-regular quivers} An arbitrary $1$-regular quiver $Q$ is the quotient of the disjoint union of copies of $\integers$, one for each connected component.

It will be necessary to discuss several paths in such a quiver $Q$ at once, however those paths will all lie in the same connected component and that component should be clear from context. Therefore, in keeping with the above convention, each path of interest in that component can be represented using paths in $\integers$.

\nlpass{Other posets of paths in $\integers$} If we partially order $\integers$-paths using the prefix relation, then two paths are comparable iff they have they same source, and two comparable paths satisfy $p \leq q$ iff $\len p \leq \len q$. This means that the prefix-poset of $\integers$-paths is a disjoint union of chains
$$
\big(\ilslpath{i}{0}\big)
<
\big(\ilslpath{i}{1}\big)
<
\big(\ilslpath{i}{2}\big)
<
\big(\ilslpath{i}{3})\big)
<
\cdots\tcomma
$$
one copy for each $i \in \integers$, each chain being order-isomorphic to $\naturals$. Essentially the same statements hold for the suffix poset; we need only compare paths by their targets rather than sources.

\section{Algebraic prerequisites}

\subsection{Quiver algebras}

We assume throughout that $\kk$ is some fixed field.

\nlpass{Algebras and opposites} By an \define{algebra $A$}, we mean an associative and unital (but not necessarily commutative) ring with a compatible $\kk$-vector space structure.

Any algebra $A$ has an opposite algebra $A^\op$. The underlying vector space is the same but the multiplication is reversed.

\nlpass{Path algebras and their elements} Suppose $Q$ is a finite quiver: that is, a directed graph with finitely many vertices and finitely many arrows, where loops and/or multiple edges are permitted. The paths of $Q$, also called \define{$Q$-paths} to emphasise their parent quiver (and including the stationary paths $e_i$ at each vertex $i$), form the basis of a vector space. Multiplication can be defined on basis vectors $p$ and $q$ by "concatenation extended by zero"; more precisely, $p \cdot q = pq$ ("$p$ then $q$") if $pq$ is a path in $Q$, and $p \cdot q = 0$ otherwise. This defines the \define{path algebra} $KQ$. Its multiplicative unit is the sum $\sum_{i \in Q_0} e_i $ of stationary paths. It has finite $\kk$-dimension iff $Q$ contains no (directed) cycles.

An arbitrary element of $\kk Q$ is thus a finite $\kk$-linear sum $\sum_{k=1}^m \lambda_k p_k$ of $Q$-paths $p_k$. The \define{components} of such an element are those paths $p_k$ with $\lambda_k\neq 0$.

\nlpass{Arrow ideal} Let $J \ideal KQ$ denote the \define{arrow ideal} of $KQ$: the smallest two-sided ideal of $KQ$ containing the arrows of $Q$. This is the Jacobson radical of $A$; that is, the intersection of all maximal ideals of $A$.

\nlpass{Ideals of path algebras} An ideal $I \ideal KQ$ is \define{admissible} iff there is an integer $N \geq 2$ with $J^N \subseteq I \subseteq J^2$.

If $\rho \subseteq \kk Q$ is a collection of elements, then we denote the smallest two-sided ideal of $\kk Q$ to contain $\rho$ by \define{$\langle \rho \rangle$}; when $\rho$ is finite, say $\rho=\{r_1, \cdots, r_s\}$, then we alternatively denote it by \define{$\langle r_1, \cdots, r_s \rangle$}. Any admissible ideal $I$ of a path algebra is of the form $I = \langle \rho \rangle$ for some finite $\rho$ \cite[Lem II.2.8]{ASS06}; without loss of generality, we may therefore assume the generating set $\rho$ to be minimal.

The elements of a minimal $\rho$ are called \define{relations}. The components of any relation $\sum_{k=1}^m \lambda_k p_k$ all have the same source and target without loss of generality, since this is true of $e_i(\sum_{k=1}^m \lambda_k p_k)e_j$ (for primitive idempotents $e_i,e_j$) and
$$
\sum_{k=1}^m \lambda_k p_k = \left(\sum_{i \in Q_0} e_i \right) \left( \sum_{k=1}^m \lambda_k p_k \right) \left( \sum_{j \in Q_0} e_j \right) = \sum_{i,j\in Q_0} e_i \left( \sum_{k=1}^m \lambda_k p_k \right)e_j\tstop
$$
A relation $\sum_{k=1}^m \lambda_k p_k$ is called a \define{monomial relation} if it has $1$ component and a \define{(skew) commutativity relation} if it has $2$.

\nlpass{Quiver algebras} By a (\define{bound}) \define{quiver algebra}, we mean a quotient $KQ/I$ of a path algebra $KQ$ by an admissible ideal $I$. We will refer to $Q$ as the \define{ground quiver} or \define{presenting quiver} of the algebra. 

Any quiver algebra is a finite-dimensional \cite[Prop II.2.6]{ASS06} algebra. Conversely, any basic, connected algebra is isomorphic to a quiver algebra \cite[Thm II.3.7]{ASS06}. To understand why it suffices to consider algebras that are basic and connected (terms defined in see \cite[Def I.6.1 and Sec I.4]{ASS06}) note that any finite-dimensional algebra is a direct product of connected ones (trivially) and the module category of any connected finite-dimensional algebra is equivalent to that of a basic one \cite[Cor I.6.10]{ASS06}. 

\nlpass{Paths and residues}\label{psg:A-paths-and-residues} Retaining the notation of the previous passage, we define the \define{$A$-residue} of a $Q$-path $p$ to be the element $p + I$. (Others have called this the canonical image of $p$ in $A$.)

We use the term \define{$A$-path} to mean the nonzero $A$-residue of some $Q$-path.

\subsection{Representation theory}

\nlpass{Modules} A \define{representation} of $A$ is a homomorphism of algebras $ \phi \colon A \to \End X$ whose target is the endomorphism algebra $\End X$ of a $\kk$-vector space $X$; for convenience, write $\phi_a$ for the image of $a$ in $\phi$. In this case, we call $X$ a \define{(right) $A$-module}, with associated action $xa \ceq x\phi_a$. The module is \define{finite dimensional} iff \emph{X} is. 

\nlpass{Representations of quivers} Since $A=\kk Q/I$, we can work in terms of \define{(bound) representations of quivers}. These are assignments of a vector space $X_i$ to each vertex $i$ of $Q$ and of a linear map $\theta_\alpha \colon X_{\source(\alpha)} \to X_{\target(\alpha)}$ to each arrow $\alpha$ of $Q$ such that for all $\kk$-linear sum of paths $ \sum_{k=1}^m \lambda_k p_k \in I$, the associated map $\sum_{k=1}^m \lambda_k \theta_{p_k}  \colon \bigoplus_i X_i \to \bigoplus_i X_i$ is zero. Here, $\theta_{p_k}=\theta_{\alpha_1}\cdots\theta_{\alpha_r}$ for a decomposition of a nonstationary $p_k$ into a product $\alpha_1\cdots\alpha_r$ of arrows, and $\theta_{i}=\id_{X_i}$ for any stationary path at $i$. 

\nlpass{Equivalence} As is well-known, representations of quivers are equivalent to modules \cite[III.1]{ASS06}. More specifically, the categories $\Repcat(Q,I)$ of bound quiver representations and $\Modcat A$ of $A$-modules are equivalent, and this equivalence restricts to their respective full subcategories $\repcat(Q,I)$ and $\modcat A$ of finite-dimensional objects. We use the two formulations interchangeably throughout.

We note in particular that all of the categories in the previous paragraph are
abelian: thus, we can speak of the direct sum of modules (denoted with $\oplus$).

\nlpass{Indecomposability and unique decomposition} We call a nonzero module $X$ \define{indecomposable} if $X = Y \oplus Z$ implies $Y$ or $Z$ is zero.

If $X$ is a finite-dimensional module (or bound quiver representation) and there are indecomposable modules $Y_r, Y'_r$ satisfying $X \isom Y_1 \oplus \cdots \oplus Y_m \isom Y'_1 \oplus \cdots \oplus Y'_{m'} $, then $m=m'$ and there exists some permutation $\sigma$ of $\{1,\dots,m\}$ such that $Y_r \isom Y'_{r \sigma}$. 

\nlpass{Additive generation} If $\mathcal{X}$ is any set of $A$-modules, we define the \define{additive closure of $\mathcal{X}$} as the full subcategory of $\modcat A$ whose objects are isomorphic to direct summands of finite direct sums of members of $\mathcal{X}$. We denote this subcategory by \define{$\add \mathcal{X}$} or, when $\mathcal{X} = \{ X \}$ is a singleton, simply \define{$\add X$}.

If a category $\categoryfont{C}$ has $\categoryfont{C} = \add X$, then we say $X$ \define{additively generates}, or is the \define{additive generator of}\!, $\categoryfont{C}$.

\nlpass{Representation type} One can seek to classify the isomorphism classes of indecomposable
(finite-dimensional) modules of an algebra. A deep theorem of Drozd \cite{Dro80} establishes that all finite-dimensional algebras fall into exactly one of
three representation types.

In increasing order of complexity, the options are \emph{representation finite}, \emph{tame} or \emph{wild}. The first type simply means the algebra has only finitely many isoclasses of indecomposables. The second type means, informally speaking, that in each dimension all but finitely many modules lie in one of finitely many classes, each parameterised by the field. Speaking even more informally, the third type are those for which the classification problem is intractible in a very strong way.

Discussion and formal definitions of representation type can be found in \cite[Sec 4.4]{Ben95}. We mention representation type only to contextualise the classification of indecomposables for SB algebras.

\nlpass{Submodules and quotients}\label{psg:submods-and-quots} An $A$-\define{submodule} $Y \leq X$ is a vector subspace $Y$ of $X$ such that $ya \in Y$ for all $y \in Y$ and $a\in A$. For any submodule $Y \leq X$, the vector-space quotient $X/Y$ has a well-defined, natural $A$-module structure. Therefore we also denote by $X/Y$ the \define{quotient} $A$-module.

We call $X$ \define{simple} if it has no submodules beside $0$ and $X$. The (isomorphism classes of) simple $A$-modules correspond to the vertices of the ground quiver $Q$ of $A$. We write \define{$S_i$} for the simple module associated to $i \in Q_0$ comprising $\kk$ at vertex $i$, $0$ at all other vertices and all maps zero. A module is \define{semisimple} if it is a direct sum of simple modules.

The largest semisimple submodule of $X$ is the \define{socle $\soc X$} of $X$. The \define{radical $\rad X$} is the the smallest submodule $Y \leq X$ such that $X/Y$ is semisimple; the \define{head $\head X$} of $X$ is that quotient $X/\rad X$.

We call $X$ \define{uniserial} if its poset of submodules is a chain.

% A \define{composition series} for a module $X$ is a strictly ascending chain of submodules 
% $$
% 0=X_0 < X_1 < X_2 < \cdots < X_{l-1} < X_l=X
% $$
% of $X$ such that each consecutive quotient $X_{k+1}/X_k$ is simple.

\nlpass{Projective and injective modules} We define two special classes of $A$-modules. Each class has several equivalent definitions, some of which emphasise the duality between the definitions, others which are more useful in practice.
\begin{thesisfigure}
\centering

\subcaptionbox{%
  \label{subfig:projective-definition-diag}
  \textit{Projective modules $X$.} In particular note that $\ilarrow{U}{}{V}$ is surjective.
}[%
  0.40\textwidth
]{%
  $$\begin{tikzcd}
          \& X \ar[d] \ar[dl, dashed]
              \\
        U \ar[r]
          \& V \ar[r]
            \& 0
        \end{tikzcd}$$
}
\hspace{0.1\textwidth}
\subcaptionbox{%
  \label{subfig:injective-definition-diag}
  \textit{Injective modules $X$.} In particular note that $\ilarrow{U}{}{V}$ is injective.
}[%
  0.40\textwidth
]{%
  $$
  \begin{tikzcd}
    0 \ar[r]
      \& U \ar[r]\ar[d]
        \& V\ar[dl, dashed]
          \\
      \& X
  \end{tikzcd}
  $$
}

\caption[%
  Defining diagrams for projectives and injectives
]{%
  \label{fig:defining-diags-for-projs-and-injs}
  \textit{Defining diagrams for projectives and injectives.}
}
\end{thesisfigure}

\begin{enumerate}
    \item The following three conditions are equivalent. Any module satisfying them is \define{projective}.
    \begin{enumerate}
        \item The covariant functor $\Hom_A(X,-)\colon \Modcat A \to \Modcat \kk$ is exact
        
        \item For any solid diagram with exact rows as in Figure \ref{subfig:projective-definition-diag}, there exists a homomorphism $\ilarrow{X}{}{U}$ commuting it.
        
        \item The module $X$ is a direct summand of some direct power of $A$.
    \end{enumerate}
    
    \item The following two conditions are equivalent. Any module satisfying them is \define{injective}.
    \begin{enumerate}
        \item The contravariant functor $\Hom_A(-,X) \colon \Modcat A \to \Modcat \kk$ is exact.
        
        \item For any solid diagram with exact rows as in Figure \ref{subfig:injective-definition-diag}, there exists a homomorphism $\ilarrow{U}{}{X}$ commuting it.
    \end{enumerate}
\end{enumerate}
We write \define{$\Projcat A$} and \define{$\Injcat A$} for the full subcategories of $\Modcat A$ comprising the projective and injective modules, respectively, and we write \define{$\projcat A$} and \define{$\injcat A$} for the corresponding full subcategories of $\modcat A$.

As with the simples, there corresponds to each vertex $i \in Q_0$ an indecomposable projective module $\mdefine{P_i}\ceq e_i A$ and an indecomposable injective module denoted \define{$E_i$}. It follows that $\bigoplus_{i \in Q_0} P_i$ and $\bigoplus_{i \in Q_0} E_i$ are additive generators for $\projcat A$ and $\injcat A$ respectively.

\nlpass{Dualities}\label{psg:dualities} We write \define{$\dual$} for the vector-space duality $\Hom_\kk(-, \kk) \colon \modcat A \to \modcat A^\op $ of modules; by \define{duality}, we mean that this functor gives a contravariant equivalence of categories $\ilarrow{ \modcat A}{}{\modcat A^\op}$ \cite[I.2.9]{ASS06}. It restricts to dualities $\ilarrow{\projcat A}{}{\injcat A^\op}$ and $\ilarrow{\injcat A}{}{ \projcat A^\op }$~ \cite[Thm I.5.13]{ASS06} that we also denote by $\dual$.

Additionally, we also write \define{$(~)^*$} (ie, a superscript asterisk) for the functor $$\Hom_A(-, A) \colon \ilarrow{\modcat A}{}{\modcat A^\op}$$ which, importantly, restricts to a duality $\ilarrow{ \projcat A}{}{\projcat A^\op}$ (see first paragraph of \cite[\S IV.2]{ASS06}).

\subsection{Special biserial algebras}

\nlpass{SB algebras}\label{def:sb-algs} A \define{special biserial (SB) algebra} is a quiver algebra $\kk Q/I$ where
\begin{enumerate}
    \item $Q$ is sub-$2$-regular and
    \item\label{subdef:arrow-residues-for-sb-algs} for any arrow $\alpha \in Q_1$,
    \begin{enumerate*}
        \item\label{subdef:at-most-one-next-arrow} there exists at most one arrow $\beta \in Q_1$ with $\alpha \beta \notin I$ and
        \item\label{subdef:at-most-one-prev-arrow} there exists at most one arrow $\gamma \in Q_1$ with $\gamma \beta \notin I$. 
    \end{enumerate*}
\end{enumerate}
We remark in passing that $A$ is special biserial iff $A^\op$ is.

\npass Previous authors have restricted attention to particular classes of SB algebras. These include string algebras (for which the defining ideal $I$ is generated by paths), gentle algebras ($I$ is generated by paths of length $2$) and symmetric or representation-finite or radical-cube-zero SB algebras. We make no such restriction.

In the following proposition, we survey some of the most significant theorems about (unrestricted) SB algebras.

\begin{proposition}\label{prop:properties-of-sb-algebras} Let $A=\kk Q/I$ be a SB algebra.

(We introduce projective, injective, nonuniserial modules and string modules in greater detail later in this subsection. Band modules are defined in Subsection \ref{subsec:band-modules}. For the definition of the representation dimension of an algebra, see \cite[\S1]{EHIS04}. For discussion of the big and little finitistic dimensions $\Findim A$ and $\findim A$ of an algebra $A$, as well as the big and little finitistic dimension conjectures related to them, refer forward to Passage \ref{psg:finitistic-dimension-conjectures}.)
\begin{enumerate}
    \item\label{subprop:sb-algs-rels-monomial-or-commutativity} \emph{\cite[Corollary to Lem 1]{SW83} \cite[Prop 1.3]{WW85}}~ The defining admissible $I$ is generated by monomial relations and (skew) commutativity relations.

    \item\label{subprop:sb-algs-characterized-indecs} \emph{\cite[Prop 2.3]{WW85}}~ The indecomposable finite-dimensional modules of an $A$ come in three distinct flavours.
    \begin{enumerate}
        \item \emph{Projective, injective, nonuniserial modules:} there are only finitely many of these, possibly zero.
        \item\label{subprop:fin-many-string-mods-of-given-dim} \emph{String modules:} there are only finitely many of these with any given $\kk$-dimension.
        \item \emph{Band modules:} these come in families parameterised by elements of the ground field $\kk$.
    \end{enumerate}
    Consequently $A$ is representation-finite or tame.
    
    \item \label{subprop:rejection-lemma} \emph{\cite[\S2]{WW85} \cite[II.1.3]{Erd90}~} Decompose the regular module of $A$ as $P \oplus P'$, where all direct summands of $P$ are projective, injective and nonuniserial but no summands of $P'$ are. Then $\soc P$ is in fact a two-sided ideal of $A$, and the only indecomposable $A$-modules that are not $A/(\soc P)$-modules are the indecomposable direct summands of $P$.
    
    \item \label{subprop:sb-algs-repdim} \emph{\cite[Cor 1.3]{EHIS04}}~ The representation dimension of $A$ is at most $3$ and so $A$ satisfies the little finitistic dimension conjecture.
    
    \item \label{subprop:hz16-sb-algs}
    \emph{\cite[Thm 3.1]{HZ16}~} For any $r>0$, there exists a SB algebra $A_r$ that satisfies $\findim A_r = r+1$ and $\Findim A_r = 2r+1$. 
\end{enumerate}
\end{proposition}

\proof We direct the reader to the references given with each statement. We comment that the finiteness statement of \ref{subprop:sb-algs-characterized-indecs}\ref{subprop:fin-many-string-mods-of-given-dim} is implicit.

The proof of part \ref{subprop:rejection-lemma} amounts to \cite[Corollary to Lem 1]{SW83} but it bears repeating here. We assume part \ref{subprop:sb-algs-rels-monomial-or-commutativity}, and thus suppose $I = \langle \rho \rangle$ for some minimal set of monomial and commutativity relations. Suppose $\lambda p - \mu q$ is a commutativity relation (and so $\lambda \neq 0 \neq \mu$). Both components $p$ and $q$ are clearly paths of length at least $2$ since $I$ is admissible; say they factor as $p = \alpha_1 p' \alpha_2$ and $q=\beta_1 q' \beta_2$ for arrows $\alpha_1,\alpha_2,\beta_1,\beta_2$ and paths $p',q'$. Then, by definition,
\begin{equation}\label{eqn:cpts-of-commu-rel-are-cong}
 p+I = (\alpha_1 p' \alpha_2) + I = (\lambda^{-1}\mu) (\beta_1 q' \beta_2) + I = (\lambda^{-1}\mu)(q + I)\tstop
\end{equation}
Evidently, left-multiplication by $e_{\source(p)}=e_{\source(q)}$ or right-multiplication by $e_{\target(p)}=e_{\target(q)}$ fixes $p+I$ and $q+I$, while left- or right-multiplication by any other idempotent annihilates them. For any arrow $\gamma \in Q_1$, consider right-multiplying (\ref{eqn:cpts-of-commu-rel-are-cong}) by $\gamma$ to give
$$(p + I)\gamma  = p\gamma+I = \alpha_1 p' \alpha_2 \gamma + I = (\lambda^{-1}\mu) \beta_1 q' \beta_2 \gamma + I = (\lambda^{-1}\mu)(q + I) \gamma\tstop
$$
If this is nonzero then $\alpha_2\gamma \notin I$ and $\beta_2 \gamma \notin I$, which contradicts Definition \ref{def:sb-algs}\ref{subdef:arrow-residues-for-sb-algs}\ref{subdef:at-most-one-prev-arrow}; thus arrows act on the right on $p+I$ and $q+I$ as zero. One similarly shows they act as zero on the left. 

We know that $\soc P$ is spanned by
$\{ p + I \colon p \text{ is a component of a commutativity relation in } \rho \}
$.
Evidently $\soc P$ is a right submodule of $P$ and hence of $A$ by transitivity of inclusion, and therefore \emph{a fortiori} a vector subspace of $A$. The preceding paragraph shows that this subspace is closed under the action of $A$ on both the left and right. Thus, $\soc P$ is a two-sided ideal of $A$, as claimed. It follows that $A/(\soc P) = \kk Q/\langle \rho' \rangle$, where
$\rho' = \rho \cup \left\{ p \colon p \text{ is a component of a commutativity relation in }\rho \right\}$. \qed

\nrunningexample\label{example:basis-of-projectives} Throughout this thesis, we will consistently work with the SB algebra
$$
A \ceq \kk\big(
\begin{tikzcd}
  1 \ar[loop left, "\alpha"] \ar[r, "\beta", shift left]
    \& 2 \ar[l, "\gamma", shift left] \ar[loop right, "\delta"]
\end{tikzcd}
\big)
/
\big\langle
  \alpha^2,
  \beta\delta,
  \gamma\beta,
  \delta\gamma,
  \alpha\beta\gamma-\beta\gamma\alpha,
  \gamma\alpha\beta,
  \delta^4
\big\rangle\tstop
$$
The $A$-paths of length $2$ are $\alpha\beta, \beta\gamma, \gamma\alpha, \delta^2$. The paths give the basis of the regular $A$-module as depicted in Figure \ref{fig:structure-of-running-example-algebra}.

Our graphical notation follows that in, for example, \cite{ZH91}. The vertices of the graphs are $A$-paths, which form a basis. An edge labelled by some $\alpha' \in Q_1$ indicates that $\alpha'$ maps the higher vector to some (nonzero) scalar multiple of the lower vector, while all other arrows annihilate the higher node.
\begin{thesisfigure}
\centering 
\begin{tikzpicture}[scale=0.6]
      % Give "anchoring" nodes for nonzero summands
      \draw (-7,0) node[](1){};
      \draw (0,0) node[](2){};
      
      % Draw nodes of lefthand projective
      \draw (1) ++(2,6) node[](1e1){$e_1$};
      \draw (1) ++(4,4) node[](1b){$\beta$};
      \draw (1) ++(4,2) node[](1bc){$\beta\gamma$};
      \draw (1) ++(2,0) node[](1bca){$\alpha\beta\gamma$};
      \draw (1) ++(0,4) node[](1a){$\alpha$};
      \draw (1) ++(0,2) node[](1ab){$\alpha\beta$};
      
      % Draw nodes of righthand projective
      \draw (2) ++(2,6) node[](2e2){$e_2$};
      \draw (2) ++(4,4) node[](2d){$\delta$};
      \draw (2) ++(4,2) node[](2dd){$\delta^2$};
      \draw (2) ++(4,0) node[](2ddd){$\delta^3$};
      \draw (2) ++(0,4) node[](2c){$\gamma$};
      \draw (2) ++(0,2) node[](2ca){$\gamma\alpha$};
      
      % Draw action lines on lefthand projective
      \draw (1e1) --node[pos=0.5,right](){$\beta$} (1b);
      \draw (1b) --node[pos=0.5,right](){$\gamma$} (1bc);
      \draw (1bc) --node[pos=0.5,right](){$\alpha$} (1bca);
      \draw (1e1) --node[pos=0.5,right](){$\alpha$} (1a);
      \draw (1a) --node[pos=0.5,right](){$\beta$} (1ab);
      \draw (1ab) --node[pos=0.5,right](){$\gamma$} (1bca);
      
      % Draw action lines on righthand projective
      \draw (2e2) --node[pos=0.5,right](){$\delta$} (2d);
      \draw (2d) --node[pos=0.5,right](){$\delta$} (2dd);
      \draw (2dd) --node[pos=0.5,right](){$\delta$} (2ddd);
      \draw (2e2) --node[pos=0.5,right](){$\gamma$} (2c);
      \draw (2c) --node[pos=0.5,right](){$\alpha$} (2ca);
    \end{tikzpicture}
    
  \caption[
    Structure of running example algebra
  ]{
    \label{fig:structure-of-running-example-algebra}
    \textit{Structure of the regular module of running example algebra $A$.} The lefthand graph illustrates the structure of $P_1 = e_1 A$; the righthand graph, that of $P_2 = e_2 A$.
  }
\end{thesisfigure}

This algebra is built into the \sbstrips\ package as \mono{SBStripsExampleAlgebra( 1 )}.
\begin{Verbatim}[commandchars=!@|,fontsize=\small,frame=single,baselinestretch=1.2]
  !gapprompt@gap>| !gapinput@SBStripsExampleAlgebra( 1 );|
  <Rationals[<quiver with 2 vertices and 4 arrows>]/
  <two-sided ideal in <Rationals[<quiver with 2 vertices and 4 arrows>]>,
    (7 generators)>>
\end{Verbatim}

\nlpass{Special biserial algebas are biserial}\label{psg:sbas-are-biserial} A finite-dimensional algebra $A$ is \define{biserial} if for every indecomposable projective $A$-module $P$ there exists uniserial $A$-modules $U,V$ such that $\rad P = U+V$ and $U \cap V$ is simple or zero.

The SB algebra axioms ensure that the principal right ideal $\alpha A$ generated by an arrow $\alpha$ is uniserial, being spanned by the $A$-paths of the form $\alpha p$ ($p$ some $Q$-path). We find that $\rad P_i = \sum \alpha A$, where the sum runs over all arrows $\alpha$ with $\source(\alpha)=i$, of which there are at most two. We deduce that all SB algebras are biserial. The converse is not necessarily true \cite[\S1]{SW83}.

% \nlpass{Indecomposable projectives with simple socle} Suppose $e_i A$ is an indecomposable, projective, injective, nonuniserial module. There exists a skew commutativity relation $\lambda \alpha p' - \mu \beta q'$, where $\lambda,\mu\in \kk$, where $\alpha,\beta\in Q_1$ are distinct paths common target $i$ and where $p'$ and $q'$ are distinct paths with common target. We find that $\alpha A + \beta A = \rad (e_i A)$ as above, and we see that $\alpha A \cap \beta A $ is the $1$-dimensional space (hence, simple module) for which either $\{\alpha p' + I\}$ or $\{\beta q' + I\}$ is a basis. (Here, $I$ is the defining ideal of the SB algebra.)

% Conversely, suppose that $e_i A$ is an indecomposable projective module such that $\rad e_i A = U + V$ with $U \cap V$ simple, as above. If $e_i A$ is pin, then in particular it is indecomposable injective and so has simple socle, and those $U \cap V$ is simple. 

% Any indecomposable projective module has simple head. By the classification in Proposition \ref{prop:properties-of-sb-algebras}\ref{subprop:sb-algs-characterized-indecs} of indecomposables for a SB algebra $A$, 

\centresubsec{Pin modules}

\nlpass{Pin modules} We abbreviate \emph{projective, injective, nonuniserial module} to \define{pin module}. We further insist that all pin modules are indecomposable.

\nrunningexample In our running example algebra from Passage \ref{example:basis-of-projectives}, $P_1=e_1 A$ is a pin module. Its structure graph, as in Figure \ref{fig:structure-of-running-example-algebra}, has a lozenge shape. Reading down either side of the lozenge, the labels of the edges spell out $\alpha\beta\gamma$ and $\beta\gamma\alpha$; that the $A$-residues of these two paths are linearly dependent on one another is a consequence of the defining relation $\alpha\beta\gamma-\beta\gamma\alpha$ for $A$.

\nrmks\label{psg:socles-of-projectives} \ilitem{} If an indecomposable module is projective, injective and nonuniserial then, in particular, it is an indecomposable injective and so has simple socle. (This is dual, using the functor $\dual$, to the fact that indecomposable projectives have simple head \cite[Cor I.5.17]{ASS06}.)

\ilitem{} The pin modules for a SB algebra correspond to commutativity relations. This general fact seems to be well-known, appearing, at least implicitly, in the literature in (eg) \cite[\S2, para 5]{LM04} and throughout \cite{HZ16}. We remarked on the particular case of our running example algebra in the previous passage.

\ilitem{} Continuing the notation of the previous passage, let the pair $(U,V)$ of submodules of $e_i A$ be: $(0,0)$ if $i$ is the source of no arrows; $(\alpha A, 0)$ if $i$ is the source of one arrow $\alpha$ or; $(\alpha A, \beta A)$ if $i$ is the source of two arrows $\alpha \neq \beta$.

To have $U \cap V \neq 0$, we must not only have that $i$ be the source of two arrows (since otherwise $V=0$), but also that some $A$-path is represented by (a linear multiple of) a $Q$-path with first arrow $\alpha$ and one with first arrow $\beta$. This occurs iff there is a commutativity relation with source $i$ iff $P_i$ is a projective, injective, nonuniserial module. When this does occur, the simple intersection is isomorphic to $S_j$, where $j$ is the target of the commutativity relation that has source $i$. This is easy to see, since a basis vector for this simple socle is given by either of the components of the commutativity relation that has source $i$. Confer with our example basis in Figure \ref{fig:structure-of-running-example-algebra} or, more generally, the basis for projectives of SB algebras given in \cite[\S2, para 5]{LM04}

% \nrmk Any indecomposable pin module is isomorphic to a principal right ideal of the form $e_i A$, for $A$ the SB algebra and $e_i$ the primitive idempotent associated to the source $i$ of a commutativity relation $p-\lambda q$ ($\lambda \neq 0$). There is a convenient basis of $e_i A$ corresponding to right subideals of $e_i A$.

% For $u$ a path with source $i$, let $uA \ideal A$ be the principal (right) ideal. A vertex or arrow $x$ of the algebra acts on the right ideals by multiplication $\ilmapsto{uA}{x}{uxA}$; thus, each ideal is annihilated by every primitive idempotent save one which fixes it, and each arrow sends a primitive ideal to another one or to $0$. With respect to inclusion, the nonzero ideals form a lozenge shape: that is, besides a maximal element ($e_i A$) and a minimal element ($pA = qA=\soc(e_i A)$), there are two disjoint chains. We hope the example below is illuminating.

% These ideals $uA$ form a basis of $e_i A$. Additionally, each arrow of $A$ sends a basis element to zero or to a scalar multiple of another basis element.
% % acts as a matrix of a special form. Each row contains at most one nonzero entry; with the exception of a single $\lambda$, those nonzero entries are $1$.

% We will always work with this basis for pin modules unless otherwise noted.

\centresubsec{String modules}

\nlpass{String graphs}\label{defn:string-graph} A \define{string graph (for $A$)} is a quiver homomorphism $w \colon \ilarrow{G}{}{Q} $ such that:
\begin{enumerate}
    \item the underlying (undirected) graph of $G$ is linear, each connected component being either finite
    $\big(
      \begin{tikzcd}[sep=small, cramped]
        \syllableplaceholder \ar[r, no head]
          \& \syllableplaceholder \ar[r, no head]
            \& \cdots \ar[r, no head]
              \& \syllableplaceholder
      \end{tikzcd}
    \big)$, unbounded in one direction $\big(
      \begin{tikzcd}[sep=small, cramped]
        \syllableplaceholder \ar[r, no head]
          \& \syllableplaceholder \ar[r, no head]
            \& \syllableplaceholder \ar[r, no head]
              \& \cdots
      \end{tikzcd}
    \big)$ or unbounded in both directions $\big(
      \begin{tikzcd}[sep=small, cramped]
        \cdots \ar[r, no head]
          \& \syllableplaceholder \ar[r, no head]
            \& \syllableplaceholder \ar[r, no head]
              \& \cdots
      \end{tikzcd}
    \big)$;
    \item
    \begin{enumerate}
        \item for any subgraph $\big(\begin{tikzcd}[sep=small,cramped]
    \syllableplaceholder \& i \ar[l,"x"'] \ar[r, "y"] \& \syllableplaceholder
    \end{tikzcd}\big)$ of $G$ featuring a source vertex $i$ and its two outgoing arrows, we have $w(x)\neq w(y)$;
        \item for any subgraph
    $\big(\begin{tikzcd}[sep=small,cramped]
    \syllableplaceholder \ar[r, "x"] \& i \& \ar[l,"y"'] \syllableplaceholder
    \end{tikzcd}\big)$ of $G$ featuring a sink vertex $i$ and its two incoming arrows, we have $w(x)\neq w(y)$;
    \end{enumerate}
    \item\label{subdefn:linear-subgraphs-linearly-independent} for any path
    $\big(\begin{tikzcd}[sep=small]
    \syllableplaceholder \ar[r, "x_1"]
      \& \syllableplaceholder \ar[r, "x_2"]
        \& \cdots \ar[r, "x_\ell"]
          \& \syllableplaceholder
    \end{tikzcd}\big)$
    in $G$, the $A$-path $p\ceq w(x_1 x_2 \cdots x_l)$ is linearly independent of any other $A$-path.
\end{enumerate}
We follow the commonplace convention of depicting a string graph as a labelled graph, each vertex $v$ or arrow $x$ being labelled by $w(v)$ or $w(x)$.

We do not assume that $G$ is connected. When $G$ is connected, we call $w$ \define{indecomposable}. The restriction of $w$ to any connected component of $G$ gives an indecomposable string graph.

\nrmks \ilitem{} Condition \ref{defn:string-graph}\ref{subdefn:linear-subgraphs-linearly-independent} is just a formal way of saying that $p\notin I$ and $p$ is not a component of a commutativity relation.

\ilitem{} The empty graph is a string graph.

\nlpass{String modules} Let $w$ be a string graph for $A$. The \define{string module $\Str(w)$} associated to $w$ is the $A$ module constructed as follows. Its basis is the vertex set of $w$. Any basis vector $v$ is fixed by the idempotent $e_{w(v)}$ and annihilated by all others. For each arrow
$\begin{tikzcd}[sep=small]
v \ar[r, "x"] \& v'
\end{tikzcd}$ of the string graph, the $A$-arrow $w(x)$ sends $v$ to $v'$. Otherwise, each $A$-arrow acts as zero on basis vectors. Thanks to condition \ref{defn:string-graph}\ref{subdefn:linear-subgraphs-linearly-independent} and Proposition \ref{prop:properties-of-sb-algebras}\ref{subprop:sb-algs-rels-monomial-or-commutativity}, this does indeed give an $A$-action.

In general, a string module is a module isomorphic to some $\Str w$. However we will only ever work with these representative examples.

\nexamples \ilitem{} Any simple module is a string module: the corresponding string graph is a single vertex. More generally, any linear graph  $\big(\begin{tikzcd}[sep=small]
    \ell \ar[r, "x_1"]
      \& (\ell-1) \ar[r, "x_2"]
        \& \cdots \ar[r, "x_\ell"]
          \& 1
\big)\end{tikzcd}$ gives a uniserial module having a number of vertices equal to the composition length of $X$. In fact, if the composition series of $X$ is $0=X_0 < X_1 < \cdots <X_{\ell-1} < X_\ell = X$, then the composition factor $X_{k}/X_{k-1}$ is $S_{w(k)}$.

\ilitem{} The zero module is the string module associated to the empty graph.

\ilitem{} The indecomposable direct summands of $\Str(w)$ correspond to the restrictions of $w$ to its connected components.

\nlpass{Rescaling string graphs} Suppose that in addition to a string graph $w$ we are given a nonzero scalar $\lambda_x$ for each arrow $\ilarrow{v}{x}{v'}$ of $w$. We can then define an $A$-module as in the above construction of string graphs, except now we stipulate that the $A$-arrow $w(x)$ sends $v$ to $\lambda_x v'$ instead of $v'$.

This gives a construction of a string module that seems more general than the one given above, but is not. The modules given by either construction are isomorphic; the isomorphism may be realised by rescaling the basis vectors independently.

Going forward, we will take such rescaling for granted when recognising string modules.

\npass String modules can be defined over any finite-dimensional algebra but they have proved especially attractive over SB algebras. The following propositions collect some highlights.

\begin{proposition}\label{prop:functors-send-strings-to-strings}
Let $X \ceq \Str(w)$ be a string module over a SB algebra.

(See \cite[\S IV.2]{ASS06} for the definitions and discussion of tranpose $\Tr$ and the Auslander-Reiten translate $\artranslate$ and inverse translate $\arinversetranslate$. Refer forward to Definition \ref{def:proj-covers-and-syzs} for syzygies $\syzygy=\syzygy^1$. See \cite[\S1.1]{Gel20} for discussion of suspension $\suspension$.)
\begin{enumerate}
    \item\label{subprop:d-of-string-is-string} \emph{\cite[Lem 3.1(1)]{WW85}\hspace{1ex}} The vector-space dual $\dual X \ceq \Hom_\kk(X, \kk)$ of $X$ is a string module (for $A^\op$).
    
    \item\label{subprop:tr-of-string-is-string} \emph{\cite[Lem 3.2(2)]{WW85}\hspace{1ex}}The transpose $\Tr X $ of $X$ is a string module (for $A^\op$).
    
    \item\label{subprop:syzygy-of-string-is-string} \emph{\cite[Prop 2.2]{LM04}\hspace{1ex}} The syzygy $\syzygy^1 X$ of $X$ is a string module.
    
    \item The Auslander-Reiten translate $\artranslate X \ceq (\dual \circ \Tr)(X)$ and inverse translate $\arinversetranslate X \ceq (\Tr \circ \dual)(X)$ of $X$ are string modules.
    
    \item The suspension $\suspension X \ceq (\Tr \circ \syzygy \circ \Tr)(X)$ of $X$ is a string module. 
\end{enumerate}
\end{proposition}

\proof For parts \ref{subprop:d-of-string-is-string}, \ref{subprop:tr-of-string-is-string} and \ref{subprop:syzygy-of-string-is-string}, we direct the reader to the references given with each statement. The remaining parts follow by composition.\qed

\begin{proposition}
Any indecomposable projective or injective $A$-module that is not a pin module is a string module.
\end{proposition}

\proof For projective modules, this is mentioned explicitly in the final sentence of the "String modules" subsection of \cite[\S1]{HZ16} and implicitly in the fifth paragraph (beginning "Furthermore") of \cite[\S2]{LM04}. 

The statement about injective $A$-modules can be obtained by applying $\dual$ to the injective $A^\op$-modules since: $\dual$ sends strings to strings (by the previous proposition); exchanges injective and projective modules (by Passage \ref{psg:dualities}) and; respects (non)uniseriality (this follows from \cite[Thm 5.13(a)(c)]{ASS06}) \qed

\nsbstrips These constructions on string modules in Proposition \ref{prop:functors-send-strings-to-strings} are implemented in the \sbstrips\ package. We refer the reader to the \sbstrips\ documentation provided in Appendix \ref{app:sbstrips-documentation}. 

\subsection{The derived category}

In the following, $A$ denotes an arbitrary finite-dimensional algebra over $\kk$.

The purpose of this subsection is twofold. First we must define the derived category $\derivedcat(A)$ of $A$. This category, initially considered in the 1960s by the circle around Grothendieck and formalised in the PhD thesis of his student Verdier \cite{Ver67,Ver77,Ver96} (see also \cite{Har66} for a contemporaneous account), is a natural setting in which to investigate the homological properties of $A$. It encompasses $\Modcat A$ (and so also $\modcat A$) but, unlike them, has the structure of a triangulated category rather than an abelian category. We will race through the subtle minutiae of its construction and defer often to the literature for details.

Many homological questions about $A$ can be framed in terms of $\derivedcat(A)$ and its (triangulated) subcategories. Achieving this is our second goal.

\centresubsec{Defining the derived category}

Our prevailing reference is \cite[Chapters 1, 10]{Wei94}.

\nlpass{Complexes and $\complexcat(A)$} The \define{complex category of $A$} is the category with the following data. We denote it \define{$\complexcat(A)$}, an abbreviation for the more correct $\complexcat(\Modcat A)$

Its objects are diagrams $
X^\bullet \ceq (
\begin{tikzcd}
\cdots \ar[r]
    \& X^{-1} \ar[r,"d_X^{-1}"]
        \& X^0 \ar[r, "d_X^{0}"]
            \& X^1 \ar[r, "d_X^{1}"]
                \& X^2 \ar[r]
                    \& \cdots
\end{tikzcd})
$
in $\Modcat A$ where the boundary morphisms $d_X^k$ ($k \in \integers$) satisfy $d_X^k d_X^{k+1} = 0$. This implies $\im( d_X^k ) \leq \ker( d_X^{k+1} )$.

A morphism $f^\bullet \colon \ilarrow{ X^\bullet }{}{ Y^\bullet }$ is a family $(f^k)_{k \in \integers}$  of module homomorphisms commuting the diagram
$$
\begin{tikzcd}
\cdots \ar[r]
    \& X^{-1} \ar[r,"d_X^{-1}"] \ar[d, "f^{-1}"]
        \& X^0 \ar[r, "d_X^{0}"] \ar[d, "f^0"]
            \& X^1 \ar[r, "d_X^{1}"] \ar[d, "f^1"]
                \& X^2 \ar[r] \ar[d, "f^2"]
                    \& \cdots
                        \\
\cdots \ar[r]
    \& Y^{-1} \ar[r,"d_Y^{-1}"]
        \& Y^0 \ar[r, "d_Y^{0}"]
            \& Y^1 \ar[r, "d_Y^{1}"]
                \& Y^2 \ar[r]
                    \& \cdots\tstop
\end{tikzcd}
$$
Morphisms compose termwise using the composition law of $\Modcat A$; that is, $f^\bullet g^\bullet = (f^kg^k)_{k\in\integers}$.

Other obvious termwise constructions carry the addition of morphisms over to $\complexcat(A)$ from $\Modcat A$, as well as kernels, cokernels, products and coproducts, making $\complexcat(A)$ an abelian category.

The \define{support} of a complex $X^\bullet$ is $\mdefine{\supp X^\bullet} \ceq \{ k \in \integers \colon X^k \neq 0 \}$. We say that $X^\bullet$ is \define{bounded above}, \define{bounded below} or \define{bounded} if $\supp X^\bullet \subseteq \integers$ is too, in the sense of ordered sets.

Any module can be viewed as a bounded complex supported in a single term, usually the $0$th, and so $\Modcat A$ embeds into $\complexcat(A)$, as does its full subcategory $\modcat A$.

\nlpass{Shift functors} For $r \in \integers$, the shift functor $[r] \colon \ilarrow{ \complexcat(A) }{}{ \complexcat(A) }$ is defined on objects by
$$
X^\bullet \ceq \big(
\begin{tikzcd}[sep=small]
\cdots \ar[r]
  \& X^0 \ar[r,"d^0"]
    \& X^1 \ar[r]
      \& \cdots
\end{tikzcd}
\big)
\hspace{0.5cm}
\ilmapsto{}{}{}
\hspace{0.5cm}
X^\bullet[r] \ceq \big(
\begin{tikzcd}[sep=small]
\cdots \ar[r]
  \& X^r \ar[rrr,"(-1)^r d^r"]
    \&
      \&
        \& X^{r+1} \ar[r]
          \& \cdots
\end{tikzcd}
\big)\tstop
$$
We may think of it as moving all terms $r$ step leftwards (or more accurately one step ``against the flow of arrows'') and introducing signs $(-1)^r$ to all boundary maps. It is defined on morphisms as $\ilmapsto{ (f^k)_k }{}{ ( f^{k+r} )_k }$. It clearly respects composition and identities and,
evidently, $[r][s]=[r+s]$ for all $r, s \in\integers$.

\nlpass{Homotopy and the homotopy category $\homotcat(A)$} A morphism $f^\bullet \colon \ilarrow{ X^\bullet }{}{ Y^\bullet }$ of complexes is \define{nullhomotopic} if there exists a family $(s^k \colon \ilarrow{X^k}{}{Y^{k-1}})_{k \in \integers}$ of module homomorphisms with the property that $d_X^{k+1}s^{k+2} + s^{k+1}d_Y^k = f^{k+1}$ for all $k \in\integers$. This family $(s^k)_k$ is called a \define{contraction}.
% They fit into the following visually-pleasing diagram:
% $$
% \begin{tikzcd}
% \cdots \ar[r]
%     \& X^{-1} \ar[r,"d_X^{-1}"] \ar[dl]
%         \& X^0 \ar[r, "d_X^{0}"] \ar[dl, "s^0"']
%             \& X^1 \ar[r, "d_X^{1}"] \ar[dl, "s^1"']
%                 \& X^2 \ar[r] \ar[dl, "s^2"']
%                     \& \cdots \ar[dl]
%                         \\
% \cdots \ar[r]
%     \& Y^{-1} \ar[r,"d_Y^{-1}"']
%         \& Y^0 \ar[r, "d_Y^{0}"']
%             \& Y^1 \ar[r, "d_Y^{1}"']
%                 \& Y^2 \ar[r]
%                     \& \cdots
% \end{tikzcd}
% $$.

The nullhomotopic morphisms form an ideal of $\complexcat(A)$. The resulting quotient category is the \define{homotopy category $\homotcat(A)$} of $\Modcat A$. (As before we use $\homotcat(A)$ as a convenient shorthand for $\homotcat(\Modcat A)$.)

For clarity, we underline that $\homotcat(A)$ has the same objects as $\complexcat(A)$ and that morphisms $f^\bullet,g^\bullet \colon \ilarrow{ X^\bullet }{}{ Y^\bullet }$ of complexes represent the same $\homotcat(A)$-morphism if their difference $f^\bullet-g^\bullet$ is nullhomotopic.

The homotopy category $\homotcat(A)$ is additive but not (typically) abelian.

The shift functors on $\complexcat(A)$ induce shift functors on $\homotcat(A)$.

\nlpass{Homology} Fix a complex $X^\bullet$. For any $k \in\integers$, $ \ker d_X^{k}/\im d_X^{k+1}$ is a well-defined object of $\Modcat A$ thanks to the condition $d_X^{k-1}d_X^k = 0$. We call it the \define{$k$th homology  $H^k( X^\bullet )$ module of $X^\bullet$}. A routine diagram chase shows that any term $f^k \colon \ilarrow{ X^k }{}{ Y^k }$ satisfies $f^k(\ker d_X^k) \subseteq \ker d_Y^k$ and $f^k(\im d_X^{k-1}) \subseteq \im d_Y^{k-1}$. We find that $f^k$ induces a well-defined map $\ilarrow{ H^k(X^\bullet) }{}{ H^k(Y^\bullet) }$ which respects identities and composition.

Consequently we obtain the \define{$k$th homology functor} $H^k \colon \ilarrow{ \complexcat(A) }{}{ \Modcat A }$ for each $k \in \integers$. It is straightforward to verify that this factors through the canonical quotient functor $\ilarrow{ \complexcat(A) }{}{ \homotcat(A) }$, and therefore we are justified in treating homology as a functor family $H^k \colon \ilarrow{ \homotcat(A)}{}{ \Modcat A }$ ($k \in \integers$).

Homology is compatible with shifting in an obvious way: $H^k(X^\bullet[r]) \isom H^{k+r}(X^\bullet)$.

\nlpass{Quasiisomorphisms} We call a $\complexcat(A)$- or $\homotcat(A)$-morphism $f^\bullet$ a \define{quasiisomorphism} if all $H^k(f^\bullet)$ are isomorphisms in $\Modcat A$. 

The property of two complexes $X^\bullet, Y^\bullet$ being connected by a quasiisomorphism (in either direction $\ilarrow{X^\bullet}{}{Y^\bullet}$ or $\ilarrow{Y^\bullet}{}{X^\bullet}$) generates an equivalence relation; we call two complexes \define{quasiisomorphic} if they belong to the same class under this equivalence relation.

Standard functor properties prove that any isomorphism, and in particular any identity morphism, is a quasiisomorphism.

\nlpass{Resolutions give quasiisomorphisms} We know that any $A$-module $M$ admits a projective resolution and an injective resolution. These are respectively exact sequences
\begin{equation}\label{diag:usual-resolutions}
\begin{tikzcd}[sep=small]
\cdots \ar[r]
  \& P^{-2} \ar[r]
    \& P^{-1} \ar[r]
      \& P^0 \ar[r, "\pi"]
        \& M \ar[r]
          \& 0
\end{tikzcd}
\text{\hspace{0.75cm}and\hspace{0.75cm}}
\begin{tikzcd}[sep=small]
0 \ar[r]
  \& M \ar[r, "\iota"]
    \& E^0 \ar[r]
      \& E^1 \ar[r]
        \& E^2 \ar[r]
          \& \cdots
\end{tikzcd}
\end{equation}
with each $P^k$ projective and each $E^k$ injective. Removing $M$ from these sequences gives complexes $P^\bullet$ and $E^\bullet$  over $\Modcat A$, as shown belown. The module homomorphisms $\pi$ and $\iota$ induce morphisms $\pi^\bullet$ and $\iota^\bullet$ in $\complexcat(A)$, and subsequently $\homotcat(A)$.
$$
\begin{tikzcd}[column sep=small]
P^\bullet \ar[d, "\pi^\bullet"]
  \&
    \& (\cdots \ar[r]
      \& P^{-2} \ar[r] \arrow[d,gray]
        \& P^{-1} \ar[r] \arrow[d,gray]
          \& P^0 \arrow[r,gray] \ar[d, "\pi"]
            \& 0 \arrow[r,gray] \arrow[d,gray]
              \& 0 \arrow[r,gray] \arrow[d,gray]
                \& \cdots)
                  \\
M^\bullet \ar[d, "\iota^\bullet"]
  \&
    \& (\cdots \arrow[r,gray]
      \& 0 \arrow[r,gray] \arrow[d,gray]
        \& 0 \arrow[r,gray] \arrow[d,gray]
          \& M \arrow[r,gray] \ar[d, "\iota"]
            \& 0 \arrow[r,gray] \arrow[d,gray]
              \& 0 \arrow[r,gray] \arrow[d,gray]
                \& \cdots)
                  \\
E^\bullet
  \&
    \& (\cdots \arrow[r,gray]
      \& 0 \arrow[r,gray]
        \& 0 \ar[r,gray]
          \& E^0 \ar[r]
            \& E^1 \ar[r]
              \& E^2 \ar[r]
                \& \cdots)
\end{tikzcd}
$$
It follows from the exactness of the sequences in (\ref{diag:usual-resolutions}) that all three complexes have zero homology outside of the $0$th term where $H^0(P^\bullet)=H^0(E^\bullet)=M$. We deduce that $\pi^\bullet$ and $\iota^\bullet$ are quasiisomorphisms.

\nlpass{Mapping cones}\label{psg:mapping-cones} The \define{mapping cone $\cone^\bullet(f^\bullet)$} of a $\complexcat(A)$-morphism $f^\bullet \colon \ilarrow{ X^\bullet }{}{ Y^\bullet }$ is the complex
$$
\begin{tikzcd}[]
\cdots \ar[r]
  \& X^k \oplus Y^{k-1}\ar[r,"d^{k-1}"]
    \& X^{k+1} \oplus Y^k\ar[r, "d^k"]
      \& X^{k+2} \oplus Y^{k+1} \ar[r]
        \& \cdots
\end{tikzcd}
$$
with $k$th term $X^{k+1} \oplus Y^k$ and $k$th boundary map
$d^k 
\ceq
\left[
  \begin{smallmatrix}
    -d_X^{k+1}
      & -f^{k+1} \\
    0 & d_Y^k
  \end{smallmatrix}
\right]$.

One verifies that
$\big(
\begin{tikzcd}[sep=small]
0 \ar[r]
  \& Y^\bullet \ar[rr, "u^\bullet"]
    \& \& \cone^\bullet(f^\bullet) \ar[rr, "v^\bullet"]
      \& \& X^\bullet[1] \ar[r]
        \& 0
\end{tikzcd}\big)
$
is a short exact sequence, for maps given termwise by $u^k \ceq \left[ 0~~ 1 \right] \colon \ilarrow{ Y^k }{}{ X^{k+1} \oplus Y^k }$ and $v^k \ceq \left[
\begin{smallmatrix}
  -1 \\ 0
\end{smallmatrix}
\right] \colon \ilarrow{ X^{k+1} \oplus Y^k }{}{ X^{k+1} }$. In the associated long exact sequence 
$$
\begin{tikzcd}[sep=small]
\cdots \ar[r]
  \& H^k( Y^\bullet ) \ar[r]
    \& H^k( \cone^\bullet( f^\bullet ) ) \ar[r]
      \& H^{k}( X^\bullet[1] ) \ar[r, "\partial"]
        \& H^{k+1}( Y^\bullet ) \ar[r]
          \& \cdots
\end{tikzcd}
$$
of homology modules, the connecting homomorphism $ \big( \ilarrow{ H^{k+1}( X^\bullet ) = \big) H^{k}( X^\bullet[1] ) }{ \partial }{ H^{k+1}( Y^\bullet ) } $ is just the map induced by $f^\bullet$. The mapping construction thus ensures that an arbitrary morphism of complexes can be placed into a long exact sequence, and implies that $f^\bullet$ is a quasiisomorphism iff $\cone^\bullet( f^\bullet )$ is exact.

If $
\begin{tikzcd}[sep=small]
0 \ar[r]
 \& X^\bullet \ar[rr, "f^\bullet"]
   \&\& Y^\bullet \ar[rr, "g^\bullet"]
     \&\& Z^\bullet \ar[r]
       \& 0
\end{tikzcd}
$
is a short exact sequence in $\complexcat(A)$, then there is a quasiisomorphism $\varphi^\bullet \colon \ilarrow{ \cone^\bullet( f^\bullet ) }{}{ Z^\bullet }$ given termwise by $\varphi^k \ceq \left[
\begin{smallmatrix}
  0 \\ -g^k
\end{smallmatrix}
\right] \colon \ilarrow{ X^{k+1} \oplus Y^k }{}{ Z^k }$. (This is not typically an isomorphism.)

\nlpass{Inverses of quasiisomorphisms and $\derivedcat(A)$}\label{psg:inverses-of-quasiisomorphisms} The \define{derived category $\derivedcat(A)$} is the category obtained from $\homotcat(A)$ by adjoining inverses to quasiisomorphisms. 

The introduction of formal inverses should remind the reader of the localisation $S^{-1}R$ of a ring $R$ at a multiplicative subset $S\subseteq R$. This is the correct intuition. What is happening is a localisation of $\homotcat(A)$ at the class of quasiisomorphisms, much like with rings at a multiplicative subset.

There exists a concrete description of $\derivedcat(A)$-morphisms $\ilarrow{ X^\bullet }{}{ Y^\bullet }$ as so-called fractions of $\homotcat(A)$-morphisms $
\begin{tikzcd}
X^\bullet
  \& Z^\bullet \ar[l, squiggly, "q^\bullet"'] \ar[r, "f^\bullet"]
    \& Y^\bullet
\end{tikzcd}$, with numerator $f^\bullet$ and denominator $q^\bullet$ (connecting via some intermediate complex $Z^\bullet$). The denominator is a quasiisomorphism, a fact we signify using a squiggly arrow. The composition law is given by some contrived ``calculus of fractions''.

This is, once again, analogous to how elements of $S^{-1}R$ are constructed and multiplied. However such constructions are not to be performed and employed carelessly. Morphisms in $\derivedcat(A)$ involve some mediating object ($Z^\bullet$ in the above notation) of which there is potentially a proper class. It is no longer clear \emph{a priori} that each morphism class $\Hom_{\derivedcat(A)}(X^\bullet,Y^\bullet)$ is actually a set, and so unclear whether $\derivedcat(A)$ is a category as we understand the term.

Other notions of derived categories fall prey to such problems. The derived category $\derivedcat(A)$ of an algebra does not \cite[\S10.3]{Wei94}. Consequently we may proceed safe in the assumption that $\derivedcat(A)$ is a category in our universe, knowing that what forays we make into investigating it are well-founded.

Indeed our interactions with $\derivedcat(A)$ will be so unsophisticated that we do not even need to articulate any fine points about the calculus of fractions. We instead focus only on what pertains directly to our exposition.

% Such foundational concerns can be safely put to rest for $\derivedcat(A)$. The finer points do not concern us here, nor do the intricacies of the calculus of fractions~-- see \cite[\S10.3]{Wei94} for information about both. Suffice it to say that we need not and will not worry about the details: our 

\nlpass{Some key facts about $\derivedcat(A)$} What matters in this thesis is that $\derivedcat(A)$ has the same objects as $\homotcat(A)$ and $\complexcat(A)$~-- namely complexes of $A$-modules~-- and that two complexes which are quasiisomorphic in $\homotcat(A)$ are isomorphic in $\derivedcat(A)$. In particular, any module is isomorphic in $\derivedcat(A)$ to its projective or injective resolutions, viewed as complexes as above.

Additionally, the shift functors $[r]$ on $\complexcat(A)$ and $\homotcat(A)$ induce shift functors $[r]$ on $\derivedcat(A)$ ($r \in \integers$).

\centresubsec{Triangulated categories and subcategories}

\nlpass{Sketch of triangulation} Although $\homotcat(A)$ and $\derivedcat(A)$ are additive, neither is typically abelian and in neither of them it is generally sensible to talk about short exact sequences. Instead one endows $\homotcat(A)$ and $\derivedcat(A)$ with a related structure~-- called exact triangles~-- that makes them into \define{triangulated categories}.

For a given category like $\homotcat(A)$ or $\derivedcat(A)$, an \define{exact triangle} is any triangle isomorphic to a strict triangle. Here: by \define{triangle} we mean a triple $(f^\bullet, g^\bullet, h^\bullet)$ of morphisms $\begin{tikzcd}
X^\bullet \ar[r, "f^\bullet"]
  \& Y^\bullet \ar[r, "g^\bullet"]
    \& Z^\bullet \ar[r, "h^\bullet"]
      \& X^\bullet [1]
\end{tikzcd}
$; by a \define{strict triangle} we mean a triangle belonging to some distinguished class, and; by an isomorphism of triangles $\ilarrow{ (r^\bullet, s^\bullet, t^\bullet) }{}{ (f^\bullet, g^\bullet, h^\bullet) }$ we mean isomorphisms forming the vertical arrows in the following commuting diagram. (The righthand isomorphism is the $[1]$-shift of the lefthand isomorphism.)
$$
\begin{tikzcd}
R^\bullet \ar[r, "r^\bullet"] \ar[d, "\isom"]
  \& S^\bullet \ar[r, "s^\bullet"] \ar[d, "\isom"]
    \& T^\bullet \ar[r, "t^\bullet"] \ar[d, "\isom"]
      \& R^\bullet [1] \ar[d, "\isom"]
        \\
X^\bullet \ar[r, "f^\bullet"]
  \& Y^\bullet \ar[r, "g^\bullet"]
    \& Z^\bullet \ar[r, "h^\bullet" ]
      \& X^\bullet[1]
\end{tikzcd}
$$

The strict triangles of $\homotcat(A)$ are the triangles
$\begin{tikzcd}
X^\bullet \ar[r, "f^\bullet"]
  \& Y^\bullet \ar[r, "u^\bullet"]
    \& \cone^\bullet(f^\bullet) \ar[r, "v^\bullet" ]
      \& X^\bullet[1]
\end{tikzcd}
$ arising from mapping cones (see passage \ref{psg:mapping-cones}). We cannot describe the strict triangles of $\derivedcat(A)$ without expounding further on the calculus of fractions. For this reason, we wave our hands and claim that the strict triangles of $\derivedcat(A)$ are triangles whose numerators essentially form a strict triangle in $\homotcat(A)$.

These triangles must satisfy certain axioms, the enunciation of which is surplus to the needs of this thesis. Details can be found in \cite[\S10.2]{Wei94} and \cite[Chapter 1]{Nee01}.

\nlpass{Exact triangles from $\complexcat(A)$} Any short exact sequence
$
\begin{tikzcd}[sep=small]
0 \ar[r]
 \& X^\bullet \ar[rr, "f^\bullet"]
   \&\& Y^\bullet \ar[rr, "g^\bullet"]
     \&\& Z^\bullet \ar[r]
       \& 0
\end{tikzcd}
$
in $\complexcat(A)$ yields an exact triangle in $\derivedcat(A)$.

Specifically, it yields one isomorphic to the strict triangle $
\begin{tikzcd}[sep=small]
X^\bullet \ar[r]
  \& Y^\bullet \ar[r]
    \& \cone^\bullet(f^\bullet) \ar[r]
      \& X^\bullet[1]
\end{tikzcd}
$, whose morphisms (in $\derivedcat(A)$) are the ``fractions'' 
$$
\begin{tikzcd}
X^\bullet
  \& X^\bullet \ar[l, squiggly, "1^\bullet"'] \ar[r, "f^\bullet"]
    \& Y^\bullet
\end{tikzcd}
\text{,\hspace{0.5cm}}
\begin{tikzcd}
Y^\bullet
  \& Y^\bullet \ar[l, squiggly, "1^\bullet"'] \ar[r, "u^\bullet"]
    \& \cone^\bullet(f^\bullet)
\end{tikzcd}
\text{,\hspace{0.5cm}}
\begin{tikzcd}
\cone^\bullet(f^\bullet)
  \& Z \ar[l, squiggly, "\varphi^\bullet"'] \ar[r, "g^\bullet"]
    \& X^\bullet[1]
\end{tikzcd}\tstop
$$
of $\homotcat(A)$-morphisms, in the notation of passage \ref{psg:mapping-cones}. %We comment in passing that this relies on the existence of the quasiisomorphism $\varphi^\bullet$.

\nlpass{Triangulated subcategories} We will only examine triangulated subcategories of $\derivedcat(A)$, therefore we will only define them in this generality.

\newcommand{\triangsubcat}{\categoryfont{T}}

A full subcategory $\triangsubcat \subseteq \derivedcat(A)$ is a \define{triangulated subcategory} if the inclusion functor $\ilarrow{ \triangsubcat }{}{ \derivedcat(A) }$ is additive, commutes with the shift functors $[r]$ and sends exact triangles in $\triangsubcat$ to exact triangles in $\derivedcat(A)$, and additionally if every exact triangle of $\derivedcat(A)$ is exact in $\triangsubcat$.

\nlpass{Localising and colocalising subcategories} A triangulated subcategory is \define{localising} if it is closed under arbitrary (set-indexed) coproducts. A triangulated subcategory is \define{colocalising} if it is closed under arbitrary (set-indexed) products.

We denote by \define{$\Loc(\mathcal{X})$} or \define{$\Coloc(\mathcal{X})$} the smallest localising or colocalising subcategory containing some class $\mathcal{X}$ of $\derivedcat(A)$-objects.

\nlpass{Key properties of localising or colocalising subcategories} There are several, well-established properties of localising subcategories. Some of these properties were collected into a very convenient list in an article of Rickard. We reproduce that convenient list below, with a small specialisation: the claims that Rickard makes about the derived category of an arbitrary ring, we make only for a finite-dimensional algebra over a field.

In Proposition \ref{prop:properties-of-coloc-subcats} we prove a corresponding list of properties of colocalising subcategories. Some aspects of Rickard's standard arguments carry over directly, but some stages require more work, in a manner which does require the specialisation to finite-dimensional algebras.

\newcommand{\locsubcat}{\categoryfont{L}}

\begin{proposition}\label{lem:properties-of-loc-subcats}
\emph{\cite[Prop 2.1]{Ric19}} Let $\locsubcat$ be a localising subcategory of $\derivedcat(A)$.
\begin{enumerate}
    \item
    \begin{enumerate}
        \item If
        $\begin{tikzcd}[sep=small]
        0 \ar[r]
          \& X^\bullet \ar[r]
            \& Y^\bullet \ar[r]
              \& Z^\bullet \ar[r]
                \& 0
        \end{tikzcd}$
        is a short exact sequence of complexes and two of the three objects $X^\bullet,Y^\bullet,Z^\bullet$ are in $\locsubcat$, then so is the third.
        
        \item If a complex $X^\bullet$ is in $\locsubcat$, then so is $X^\bullet[r]$ for every $r\in\integers$.
        
        \item If $X^\bullet$ and $Y^\bullet$ are quasiisomorphic complexes and $X^\bullet$ is in $\locsubcat$, then so is $Y^\bullet$.
        
        \item If $\{X_i^\bullet \colon i \in I \}$ is a set of objects of $\locsubcat$, then $\bigoplus_{i \in I}X_i$ is in $\locsubcat$.
    \end{enumerate}
    
    \item If $X^\bullet \oplus Y^\bullet $ is in $\locsubcat$, then so are $X^\bullet$ and $Y^\bullet$.
    
    \item If $X^\bullet$ is a bounded complex, where the module $X^k$ is in $\locsubcat$ for every $k \in \integers$, then $X^\bullet$ is in $\locsubcat$.
    
    \item If
    $\big(
    \begin{tikzcd}
    X_0^\bullet \ar[r,"f_0^\bullet"] 
      \& X_1^\bullet \ar[r,"f_1^\bullet"] 
        \& X_2^\bullet \ar[r,"f_2^\bullet"] 
          \& X_3^\bullet \ar[r]
            \& \cdots
    \end{tikzcd}
    \big)$
    is a sequence of $\complexcat(A)$-morphisms between complexes, with $X^\bullet_m$ in $\locsubcat$ for all $m \geq 0$, then $\varinjlim X_m^\bullet $ is in $\locsubcat$.
    
    \item If $X^\bullet$ is a bounded above complex where the module $X^k$ is in $\locsubcat$ for every $k \in \integers $, then $X^\bullet$ is in $\locsubcat$.
    
\end{enumerate}
\end{proposition}

\subsection{Homological questions}

We introduce homological questions of interest. Two questions concern finite projective dimensions of modules; another concerns triangulated subcategories of $\derivedcat(A)$.

\nlpass{Projective covers and syzygies}\label{def:proj-covers-and-syzs} For any finite-dimensional module $X$ there is up to isomorphism a unique smallest (in vector-space dimension) projective module $\projcover X$ that maps onto it, say by the map $\pi \colon \ilepi{ \projcover X}{}{X}$. We call this module the \define{projective cover $\projcover X$} of $X$.

The kernel of this map is the \define{(first) syzygy $\syzygy^1 X$} of $X$; it is uniquely determined up to isomorphism. We inductively define the $k$th syzygy as $\mdefine{\syzygy^{k+1} M} = \syzygy^1( \syzygy^k M )$ for $k \geq 0$ and, by convention, we set $\syzygy^0 M$ to be $X/P$, for $P$ the largest projective direct summand of $X$. 

\nlpass{Finitistic dimension conjectures}\label{psg:finitistic-dimension-conjectures} For an $A$-module $X$, we define the \define{projective dimension} of $X$ as $\mdefine{\projdim X} \ceq \min \{r\geq 0 \colon \syzygy^r X\text{ is projective}\}$, or $+ \infty$ if no such $r$ exists. We can then respectively define the \define{big} and \define{little finitistic dimensions} of $A$ as
$$
\begin{array}{lcr}
  \mdefine{\Findim A}
    & \ceq
      & \sup\{\projdim X \geq 0 \colon X \in \Modcat A\}\tcomma
        \\
  \mdefine{\findim A}
    & \ceq
      & \sup \{ \projdim X \geq 0 \colon X \in \modcat A\}\tstop
\end{array}
$$
An algebra $A$ satisfies the \define{big} or \define{little finitistic dimension conjectures} when
$$
\Findim A < +\infty
\text{\hspace{1cm}or\hspace{1cm}}
\findim A < +\infty\tcomma 
$$
respectively. (We clearly have $\findim A \leq \Findim A$ always, and so the big implies the little.) These conjectures originated in a paper by Bass \cite{Bas60}, who attributed them to Rosenberg and Zalinsky.

\nlpass{Status of the finitistic dimension conjectures} The ``social'' status of the finitistic dimension conjectures is great, since they sit at the top of a network of interconnected homological properties. We direct the interested reader to \cite[\S1]{GPS18}, where this network is drawn in part and further references provided. 

In terms of academic status, the finitistic dimension conjectures are open in general but known to hold in select cases, the majority of work focussing on the little conjecture. What follows is a (not necessarily exhaustive) list; we  direct the reader to the cited references for pertinent definitions.

The little finitistic dimension conjecture holds for algebras $A$ that:
\begin{enumerate}
    \item are local, since modules over local algebras have projective dimension either $0$ or $\infty$ (stated in \cite[\S8]{Ric19} but see also \cite[proof of Prop 2.1(b)]{AR91}) and thus $\findim A = 0 < \infty$, and also since every idempotent ideal $\langle e \rangle $ of a local algebra $A$ is isomorphic to $A$ itself and hence is projective, implying $\findim A \leq 1$ by \cite[Prop 1]{Mer98};
    \item have radical-cubed zero \cite{GZH91,ZH93, IT05};
    \item have representation dimension at most $3$ \cite{IT05}, this class including that of SB algebras \cite{EHIS04});
    \item\label{subprop:syz-fin-algs} are syzygy-finite, since this means there exists a $k \geq 0$ and a finite number of indecomposable modules $Y_1, \dots, Y_m$ with $\{\syzygy^r X \in \modcat A \colon r \geq k, X \in \modcat A\} \subseteq \add (Y_1 \oplus \cdots \oplus Y_m)$, in which case $\findim A \leq k + \max \{\projdim Y \geq 0 \colon Y\in\{Y_1,\dots,Y_m\} \text{ has } \projdim Y < \infty\}$;
    \item\label{subprop:rep-fin-algs} are representation-finite, since these are syzygy-finite as above with $k=0$;
    \item\label{subprop:monomial-algs} are monomial \cite{ZH91,ZH92};
    \item\label{subprop:sym-selfinj-gor-algs} are symmetric, self-injective or Gorenstein (see next sentence).
\end{enumerate}
Indeed, algebras in classes \ref{subprop:syz-fin-algs}, \ref{subprop:rep-fin-algs}, \ref{subprop:monomial-algs} or \ref{subprop:sym-selfinj-gor-algs} satisfy the big finitistic dimension conjecture as well, as a consequence of \cite[Thms 8.1 and 4.3]{Ric19}.

\nlpass{Differences of finitistic dimensions} In the aforementioned paper, Bass also wondered whether any finite-dimensional algebra satisfies the strict inequality $\findim A < \Findim A$. The earliest answer came in the affirmative three decades later, when Zimmermann-Huisgen \cite{ZH92} found monomial algebras $A_r$ ($r\geq 2$) with
$$
(\findim A_r, \Findim A_r)=(r,r+1)\tcomma
$$
having proven in an earlier paper \cite{ZH91} that the discrepancy of $1$ was the maximal achievable among monomial algebras. Later, Smal{\o} \cite{Sma98} found a family of algebras $\Lambda_r$ ($r \geq 1$) with
$$
(\findim \Lambda_r, \Findim \Lambda_r)=(1,r)\tcomma
$$
his examples being radical-cubed zero algebras with wild representation type. Huisgen-Zimmermann \cite{HZ16} then found SB algebras $\Gamma_r$ ($r \geq 1$) with
$$
(\findim \Gamma_r, \Findim \Gamma_r) = (r+1, 2r+1)\tcomma
$$
demonstrating that the difference between big and little finitistic dimensions can be arbitrarily large even among tame algebras.

\nlpass{Generation properties of the derived category} Rickard \cite{Ric19} publicised homological properties, one of which he credits to Keller. We adopt his (Rickard's) notation in our work.

We respectively say that \define{injectives generate for $A$} or \define{projectives cogenerate for $A$} when $$\Loc(\Injcat A)=\derivedcat(A)
\text{\hspace{1cm}or\hspace{1cm}}
\Coloc(\Projcat A)= \derivedcat(A)$$ (Here, as usual, we identify $\Injcat A$ and $\Projcat A$ with subcategories of $\derivedcat(A)$ in the standard way.)

\nlpass{Status of injective generation and projective cogeneration} Rickard showed that if injectives generate for $A$, or if projectives cogenerate for $A^\op$, then $A$ satisfies the big finitistic dimension conjecture (and consequently the small one); indeed, he gives a characterisation of the big finitistic dimension conjecture using the derived category \cite[Thm 4.4]{Ric19}. From this it follows that injective generation and projective generation are stronger than the finitistic dimension conjectures.

It is shown in the aforementioned paper that injectives generate for syzygy-finite and Gorenstein algebras. The former class encompasses algebras that are representation-finite, radical-squared zero or monomial; the latter, algebras that are self-injective or symmetric (such as group algebras $\kk G$) or that have finite global dimension.

Several questions regarding these homological notions of (co)generation remain open.

The first avenue of inquiry is to prove directly that injectives generate and/or projectives cogenerate for other algebras. The most obvious obvious candidates for study are those for which the finitistic dimension conjectures hold. This thesis was born of a desire to address SB algebras. Our partial results in this direction appear in \ref{chap:patterns-and-applications}. 

A second approach is to explore how various ring constructions impact these generation properties. This is the path taken by Cummings \cite{Cum20}, who has fruitfully explored various relations between rings and their module or derived categories.

It also remains open whether the converses to Rickard's implications hold under any general circumstance, and whether there is a connection between injective generation of an algebra and projective cogeneration (of its opposite, if need be).

\section{Computational tools for algebra}

\GAP\ \cite{GAP} is a computer algebra system for computational discrete algebra. Initially developed in 1986 at Lehrstuhl D f{\"u}r Mathematik in RWTH Aachen, it is now jointly coordinated at centres in Aachen, Braunschweig and Kaiserslautern in Germany, Fort Collins in the US and St Andrews in the UK. It offers a library of functions that implement algebraic algorithms written in a programming language also called \GAP, in addition to data libraries of algebraic objects.

User-supplied programmes called packages extend its core functionality. One of these is \QPA\ \cite{QPA}, a package that implements quivers and path algebras and quotients thereof. This in turn uses Gr{\"o}bner basis machinery provided by the \textsf{GBNP} package.

The \sbstrips\ package adds further functionality to \GAP, building on \QPA. Worked examples in this thesis featuring the \sbstrips\ package are presented as follows, in a fashion emulating the on-screen output at a terminal.

\begin{Verbatim}[commandchars=!@|,fontsize=\small,frame=single,baselinestretch=1.2]
 !gapprompt@gap>| !gapinput@LoadPackage( "sbstrips" );|
 ----------------------------------------------------------------------
 Loading  SBStrips v0.6.5 (for syzygies of string modules over special biserial algebras)
 by Joe Allen (https://research-information.bris.ac.uk/en/persons/joe-allen).
 Homepage: https://jw-allen.github.io/sbstrips/
 Report issues at https://github.com/jw-allen/sbstrips/issues
 ----------------------------------------------------------------------
 true
\end{Verbatim}

\chapter[Permissible data and syllables]{Permissible data and syllables of a SB algebra}\label{chap:perm-data-and-syllables}

The first step to effectively calculate and investigate syzygies of string modules is represent them in a convenient way. Ultimately, we will represent string graphs as words, and indecomposable projective modules as grids, made up of symbols from an alphabet; the symbols will be called \emph{syllables}, the grids \emph{patches} and the words \emph{strips}. These three collections all depend on a discrete model of the SB algebra itself, which is to say certain collections of paths in a certain $1$-regular graph obtained from the algebra. In this chapter, we present that model.

Our model forms the foundation of the \sbstrips\ package for syzygies of strings. For this reason, once we have constructed our model of SB algebra, we associate to it a finite set of symbols called \emph{syllables}, together with a complementation operation on them called \emph{descent}. (This choice of name will be clarified when we treat arrays later.)

\section{Permissible data}

\subsection{Definition of permissible data}

The \sbstrips\ package treats SB algebras using combinatorial data. In this subsection we outline that combinatorial data.

Throughout, $A = \kk Q/\langle \rho \rangle$ is a SB algebra with ground quiver $Q$ and minimal relation set $\rho$ comprising paths and commutativity relations. By definition, $Q$ is a sub-$2$-regular quiver.

\newcommand{\aug}{\widetilde{Q}}
\newcommand{\preoq}{\widehat{\OQ}}

\npass{} The following handful of passages are technical in nature. A concrete example follows in Passage \ref{example:overquiver-and-dagger-of-running-example}, which we advise the reader to consult to help make sense of the coming definitions.

\nlpass{Overquiver $\OQ$ and vertex exchange map $\dagger$} Choose some $2$-regular augmentation $\aug$ of $Q$ from the finitely many options (recall passage \ref{psg:almost-uniqueness-of-regular-augmentations}) and identify $Q$ with its image in $\aug$.

The $A$-residue of a $Q$-arrow has already been defined (in Passage \ref{psg:A-paths-and-residues}). Those $\aug$-arrows not arising from $Q$-arrows, we define as having $A$-residue $0\in A$.

Strictly speaking, we define an \define{overquiver} to be a quiver homomorphism $\ilepi{\OQ}{}{\aug}$ such that
\begin{enumerate}
    \item $\OQ$ is $1$-regular and has $2|Q_0|$ vertices (although it need not be connected),
    \item the induced map $\ilepi{\OQ_1}{}{\aug_1}$ on arrows is a bijection,
    \item the induced map $\ilepi{\OQ_0}{}{\aug_0}$ on vertices is $2$-to-$1$ (ie, the preimage set of any $i \in Q_0$ has cardinality exactly $2$),
    \item all $\aug$-paths with nonzero residue in $A$ are images of $\OQ$-paths in (the functor on path categories induced by) this homomorphism.
\end{enumerate}
Of course, the homomorphism $\ilepi{\OQ}{}{\aug}$ is unambiguous once the correspondence between arrows of $\OQ$ and those of $\aug$ is given. Context and our notation will supply this correspondence; hence we will suppress explicit mention of the homomorphism \emph{per se} and abuse notation to just call its domain $\OQ$ an \define{overquiver}.

Evidently if $\OQ$ is an overquiver of $A$ then $\OQ^\op$ is an overquiver of $A^\op$ in a natural way.

The \define{vertex exchange map $\dagger$} of an overquiver $\OQ$ is the fixpointfree involution $\dagger \colon \ilarrow{\OQ_0}{}{\OQ_0}$ that exchanges $\OQ$-vertices with common image in $Q$.

\nlpass{$A$-residues of $\OQ$-paths} Recall that $\aug_0=Q_0$, and so the $A$-residue of any $\aug$-vertex (or, rather, of the associated stationary) is already known. Recall also that the $A$-residue of an $\aug$-arrow $\alpha$ was defined in the previous passage either $\alpha + \langle \rho \rangle$ for $\alpha \in Q_1$ or $0+\langle \rho \rangle$ for $\alpha \in \aug_1 \setminus Q_1$. We can extend this by composition to define the $A$-residue of a nonstationary $\aug$-path $\alpha_1\alpha_2\cdots\alpha_r$ to be the product (in order) of the $A$-residues of the $\alpha_k$. Thus, any $\aug$-path has a meaningful $A$-residue.

The overquiver homomorphism $\ilarrow{\OQ}{}{\aug}$ gives rise to a functor between the path categories of $\OQ$ and of $\aug$, as discussed in Passage \ref{psg:path-categories}. We define the $A$-residue of any $\OQ$-path to be the $\aug$-path this functor maps it to.

\nlpass{Existence of overquivers} Any SB algebra has at least one overquiver; the argument for existence essentially comprises the discussions of tracks in \cite[\S1]{WW85} and how tracks be combined or augmented in \cite[proof of Thm 1.4]{WW85}. We sketch it below.

Initially, let $\pi \colon \ilarrow{\aug_1}{}{\aug_1}$ be the partial function sending any $\alpha \in \aug_1$ to the arrow $\beta \in \aug_1$ such that the $A$-residue of $\alpha\beta$ is nonzero if it exists (or, rather, such that the product in this order of the $A$-residues of $\alpha$ and $\beta$ is nonzero). Clearly if $\alpha \in \supp \pi$ then $\target(\alpha)=\source(\alpha\pi)$ and, as a consequence of Definition \ref{def:sb-algs}\ref{subdef:arrow-residues-for-sb-algs}), this partial function is a partial injection. Since $\aug_0$ is $2$-regular, we have that $|X_i|=|Y_i|$ for any $i \in \aug_0$, where
$$\begin{array}{rcl}
X_i
  & \ceq
    & \big\{
    \alpha \in \aug_1
    \colon
    \target(\alpha)=i \text{ and } \alpha\notin \supp \pi
\big\}\tcomma
  \vspace{\parskip}
    \\
Y_i
  & \ceq
    & \big\{
    \alpha\in\aug_1
    \colon
    \source(\alpha)=i\text{ and } \forall \beta\in \aug_1~ (\beta\pi\neq \alpha)
    \big\}\tstop
\end{array}
$$
Accordingly, for each $i$ arbitrarily choose any bijection $\pi_i \colon \ilarrow{X_i}{}{Y_i}$, and use these bijections to extend $\pi$ to a bijection $\pi' \colon \ilarrow{\aug_1}{}{\aug_1}$ as follows:
$$\alpha \pi' \ceq \left\{
\begin{array}{ll}
  \alpha \pi
    & \text{if }\alpha \in \supp \pi\tcomma 
      \\
  \alpha \pi_{\target(\alpha)}
    & \text{if } \alpha\notin \supp\pi\tstop
\end{array}
\right.$$
This permutation $\pi'$ of $\aug_1$ has cycles; for each cycle $\sigma$ let $\Gamma_\sigma$ be the subquiver of $\aug_1$ whose arrows are the arrows in $\sigma$ and whose vertices are sources and targets of those arrows. The disjoint union of the inclusion maps $\ilmono{\Gamma_\sigma}{}{\aug}$ yields a function $\ilarrow{\bigsqcup_\sigma \Gamma_\sigma}{}{\aug}$ from the disjoint union $\bigsqcup_\sigma \Gamma_\sigma$ of these subquivers to $\aug$, a function which satisfies the defining properties of an overquiver by construction.

Our constructed overquiver $\bigsqcup_\sigma \Gamma_s$ depended on the choices of the bijections $\pi_i \colon \ilarrow{ X_i }{}{ Y_i }$ but, of course, there is only one possibility for $\pi_i$ whenever $|X_i|=|Y_i|\in \{0,1\}$. It follows that we are faced with a nontrivial choice for $\pi_i$ (and, even then, just a binary choice) only when $|X_i|=|Y_i|=2$, which occurs iff there are no $\aug$-paths of length $2$ that simultaneously have nonzero $A$-residue and as well as $e_i$ as a strict subpath.

\nrex\label{example:overquiver-and-dagger-of-running-example} The ground quiver $Q=\big( \begin{tikzcd}
1 \ar[loop left, "\alpha"] \ar[r, shift left, "\beta"]\& 2 \ar[l, shift left, "\gamma"] \ar[loop right, "\delta"]
\end{tikzcd} \big)$ of our running example algebra (see Passage \ref{example:basis-of-projectives}) is itself $2$-regular and therefore equals its $2$-regular augmentation $\aug$. This algebra has a unique overquiver $\OQ$, namely
$$\begin{tikzcd}
\source(\beta) \ar[r, "\beta"] \ar[d, dashed, no head]
  \& \source(\gamma) \ar[dl, "\gamma"] \ar[d, dashed, no head]
   \\
\source(\alpha) \ar[u, bend left, "\alpha"]
  \& \source(\delta) \ar[loop right, "\delta"]
    \\
\end{tikzcd}\tstop
$$
This overquiver is a finite graph which, on occasion, it will be convenient to represent by the acyclic cover
$
\big(
\begin{tikzcd}[sep=tiny]
\cdots \ar[r]
  \& \source(\alpha) \ar[r, "\alpha"]
    \& \source(\beta) \ar[r, "\beta"]
      \& \source(\gamma) \ar[r, "\gamma"]
        \& \source(\alpha) \ar[r, "\alpha"]
          \& \source(\beta) \ar[r, "\beta"]
            \& \source(\gamma) \ar[r, "\gamma"]
              \& \source(\alpha) \ar[r]
                \& \cdots
\end{tikzcd}
\big)
\sqcup
\big(
\begin{tikzcd}[sep=tiny]
\cdots \ar[r]
  \& \source(\delta) \ar[r, "\delta"]
    \& \source(\delta) \ar[r, "\delta"]
      \& \source(\delta) \ar[r]
        \& \cdots
\end{tikzcd}
\big)
$. As mentioned in the previous passage, the naming of $\OQ$-arrows makes the homomorphism $\ilepi{\OQ}{}{\aug}\text{(}=Q\text{)}$ unambiguous.

All $A$-paths (recall the definition in Passage \ref{example:basis-of-projectives}) are not just represented by $Q$-paths, but also by $\OQ$-paths. In particular, all of the $A$-paths represented by $Q$-paths of length $2$ (namely $\alpha\beta$, $\beta\gamma$, $\gamma\alpha$ and $\delta^2$) are represented by $\OQ$-paths.

The vertex exchange map $\dagger$ exchanges $\source(\alpha)$ with $\source(\beta)$ and $\source(\gamma)$ with $\source(\delta)$. We depict this fact using dashed lines in the quiver above.

\nsbstrips In \sbstrips, the overquiver of an SB algebra is stored in the attribute \mono{OverquiverOfSBAlg}.

\begin{Verbatim}[commandchars=!@|,fontsize=\small,frame=single,baselinestretch=1.2]
  !gapprompt@gap>| !gapinput@alg1 := SBStripsExampleAlgebra( 1 );;|
  !gapprompt@gap>| !gapinput@oq1 := OverquiverOfSBAlg( alg1 );|
  Quiver( ["v1","v2","v3","v4"], [["v3","v1","a_over"],["v1","v2","b_over"\
  ],["v2","v3","c_over"],["v4","v4","d_over"]] )  
\end{Verbatim}
In this example, the names we have been using for the $\OQ$-vertices (eg, $\source(\alpha)$) and those that \sbstrips\ uses (eg, \mono{v1}) correspond as follows:
$$
\ilcorresp{\source(\alpha)}{}{\mono{v3}}\tcomma
\hspace{0.5cm}
\ilcorresp{\source(\beta)}{}{\mono{v1}}\tcomma
\hspace{0.5cm}
\ilcorresp{\source(\gamma)}{}{\mono{v2}}\tcomma
\hspace{0.5cm}
\ilcorresp{\source(\delta)}{}{\mono{v4}}\tstop
$$
There is no connection in general between the name that \sbstrips\ assigns to an $\OQ$-vertex and the name of the $Q$-vertex that that $\OQ$-vertex represents. This can be seen in this particular example, because $\source(\alpha),\source(\beta)\in\OQ_0$ represent the vertex $1\in Q_0$, yet \sbstrips\ calls them $\mono{v3}$ and $\mono{v1}$. 

\nlpass{Nonzero paths $N \ideal \Subpath(\OQ) $} The $\OQ$-paths with nonzero $A$-residue form an order ideal of $\Subpath(\OQ)$; that is, if $p \leq q$ in the subpath order and $q$ has nonzero residue, then $p$ has nonzero residue. We denote this by order ideal by \define{$N$} and we call the paths in it \define{nonzero}.

By construction, every stationary path in $\OQ$ is nonzero. An $\OQ$-arrow lies \emph{outside} of $N$ (ie, has zero residue) iff it represents an arrow of $\aug_0 \setminus Q_0$.

\nlpass{Components $C \subseteq N$} Amongst the nonzero paths, there are some whose $A$-residues depend linearly on the residues of another nonzero path. We call these \define{components} and write $\mdefine{C} \subseteq N$ for the set of them.

Recall that the components $p,q$ of a commutativity relation $p-q \in \rho$ have common source and common target but distinct first arrows and distinct last arrows. When we lift the $Q$-paths $p,q$ to $\OQ$-paths $p,q$, first and last arrows are distinct as before. The sources and targets of the lifted paths are now distinct, but exchanged by $\dagger$.

There is an obvious involution on components which exchanges $p \in C$ with the unique other component $q \in C$ such that the $A$-residues of $p$ and $q$ are linearly dependent. We can denote this \define{component exchange map} by \mdefine{$\dagger$} also, since it is compatible with the source and target maps of $\OQ$ and the vertex exchange map $\dagger \colon \ilarrow{ \OQ_0 }{}{ \OQ_0 }$ via the following commutative diagram.
$$
\begin{tikzcd}
\OQ_0 \ar[d, "\dagger"]
  \& C \ar[l, "\source"'] \ar[r, "\target"] \ar[d, "\dagger"]
    \& \OQ_0 \ar[d, "\dagger"]
      \\
\OQ_0
  \& C \ar[l, "\source"'] \ar[r, "\target"]
    \& \OQ_0
\end{tikzcd}
$$

We finally mention that any component is a maximal element of $N$ and that the components form an antichain. These facts follow from the axioms of a special biserial algebra.

\nrex\label{example:nonzero-paths-and-components-in-running-example} Figure \ref{fig:nonzero-and-component-paths-in-OQ} illustrates the sets $N$ and $C$ of nonzero paths and components. We have 
$$C = \{ \alpha\beta\gamma, \beta\gamma\alpha \} \subset \{e_{\source(\alpha)}, e_{\source(\beta)}, e_{\source(\gamma)}, e_{\source(\delta)}, \alpha, \beta, \gamma, \delta, \alpha\beta, \beta\gamma, \gamma\alpha, \delta^2, \alpha\beta\gamma, \beta\gamma\alpha,\delta^3\} = N
\tstop
$$
\begin{thesisfigure}
\centering

\begin{tikzpicture}[scale=0.65]
    \newcommand{\fadedcolour}{black!30}
    \newcommand{\nodeHspacing}[1]{(#1,0)}
    \newcommand{\nodeVspacing}[1]{(0,#1*1.5)}
    \newcommand{\nonzeronode}[4]{%
      \draw (o) ++\nodeHspacing{#1} ++\nodeVspacing{#2} node[black](#3){$ #4 $};
    }
    \newcommand{\componentnode}[4]{%
      \draw (o) ++\nodeHspacing{#1} ++\nodeVspacing{#2} node[circle, fill=promptColor!10](#3){$ #4 $};
    }
    \newcommand{\nodeplaceholder}{\cdot}
    \newcommand{\zeronodeplaceholder}{\circ}
    \newcommand{\zeronodeframework}[5]{
      \draw (o) ++\nodeHspacing{#1} ++\nodeVspacing{#2} node[rectangle, rounded corners, inner sep=5, #5](#3){$ #4 $};
    }
    \newcommand{\zeronode}[4]{%
      \zeronodeframework{#1}{#2}{#3}{#4}{\fadedcolour}
    }
    \newcommand{\hasseconnector}[2]{%
      \draw[\fadedcolour] (#1) -- (#2);
    }
    \newcommand{\cylinderboundary}[2]{%
      \draw[dotted] (#1) -- (#2);
    }
    
    %%  LEFTHAND COMPONENT 
    
    %%  Place "lattice" of nodes
    % 0th (bottom) row
    \draw (0,0) node[](o){};
    \nonzeronode{0}{0}{sa0}{\source(\alpha)};
    \nonzeronode{2}{0}{sb0}{\source(\beta)};
    \nonzeronode{4}{0}{sc0}{\source(\gamma)};
    \nonzeronode{6}{0}{sa1}{\source(\alpha)};
    % 1st row
    \nonzeronode{1}{1}{a0}{\alpha};
    \nonzeronode{3}{1}{b0}{\beta};
    \nonzeronode{5}{1}{c0}{\gamma};
    % 2nd row
    \nonzeronode{0}{2}{ca0}{\gamma\alpha};
    \nonzeronode{2}{2}{ab0}{\alpha\beta};
    \nonzeronode{4}{2}{bc0}{\beta\gamma};
    \nonzeronode{6}{2}{ca1}{\gamma\alpha};
    % 3rd row
    \zeronode{1}{3}{cab0}{\gamma\alpha\beta};
    \componentnode{3}{3}{abc0}{\alpha\beta\gamma};
    \componentnode{5}{3}{bca0}{\beta\gamma\alpha};
    % 4th row
    \zeronode{0}{4}{bcab0}{\beta\gamma\alpha\beta};
    \zeronode{2}{4}{cabc0}{\gamma\alpha\beta\gamma};
    \zeronode{4}{4}{abca0}{\alpha\beta\gamma\alpha};
    \zeronode{6}{4}{bcab1}{\beta\gamma\alpha\beta};
    % 5th row
    \zeronode{1}{5}{bcabc0}{\beta\gamma\alpha\beta\gamma};
    \zeronode{3}{5}{cabca0}{\gamma\alpha\beta\gamma\alpha};
    \zeronode{5}{5}{abcab0}{\alpha\beta\gamma\alpha\beta};
    % 6th row
    \zeronode{0}{6}{abcabc0}{};
    \zeronode{2}{6}{bcabca0}{};
    \zeronode{4}{6}{cabcab0}{};
    \zeronode{6}{6}{abcabc1}{};

    %%  Draw connecting lines in Hasse diagram
        % 0th row to 1st row
    \hasseconnector{sa0}{a0};
    \hasseconnector{sb0}{b0};
    \hasseconnector{sc0}{c0};
    \hasseconnector{sb0}{a0};
    \hasseconnector{sc0}{b0};
    \hasseconnector{sa1}{c0};
    % 1st row to 2nd row
    \hasseconnector{a0}{ca0};
    \hasseconnector{b0}{ab0};
    \hasseconnector{c0}{bc0};
    \hasseconnector{a0}{ab0};
    \hasseconnector{b0}{bc0};
    \hasseconnector{c0}{ca1};
    % 2nd row to 3rd row
    \hasseconnector{ca0}{cab0};
    \hasseconnector{ab0}{abc0};
    \hasseconnector{bc0}{bca0};
    \hasseconnector{ab0}{cab0};
    \hasseconnector{bc0}{abc0};
    \hasseconnector{ca1}{bca0};
    % 3rd row to 4th row
    \hasseconnector{cab0}{bcab0};
    \hasseconnector{abc0}{cabc0};
    \hasseconnector{bca0}{abca0};
    \hasseconnector{cab0}{cabc0};
    \hasseconnector{abc0}{abca0};
    \hasseconnector{bca0}{bcab1};
    % 4th row to 5th row
    \hasseconnector{bcab0}{bcabc0};
    \hasseconnector{cabc0}{cabca0};
    \hasseconnector{abca0}{abcab0};
    \hasseconnector{cabc0}{bcabc0};
    \hasseconnector{abca0}{cabca0};
    \hasseconnector{bcab1}{abcab0};
    % 5th row to 6th row
    \hasseconnector{bcabc0}{abcabc0};
    \hasseconnector{cabca0}{bcabca0};
    \hasseconnector{abcab0}{cabcab0};
    \hasseconnector{bcabc0}{bcabca0};
    \hasseconnector{cabca0}{cabcab0};
    \hasseconnector{abcab0}{abcabc1};
    
    % Draw dotted lines up the side, identifying cylinder
    \cylinderboundary{sa0}{ca0};
    \cylinderboundary{sa1}{ca1};
    \cylinderboundary{ca0}{bcab0};
    \cylinderboundary{ca1}{bcab1};
    \cylinderboundary{bcab0}{abcabc0};
    \cylinderboundary{bcab1}{abcabc1};
    
    %%  RIGHTHAND COMPONENTS
    
    %%  Place "lattice" of nodes
    % 0th (bottom) row
    \nonzeronode{10}{0}{sd0}{\source(\delta)};
    \nonzeronode{12}{0}{sd1}{\source(\delta)};
    % 1st row
    \nonzeronode{11}{1}{d0}{\delta};
    % 2nd row
    \nonzeronode{10}{2}{dd0}{\delta^2};
    \nonzeronode{12}{2}{dd1}{\delta^2};
    % 3rd row
    \nonzeronode{11}{3}{ddd0}{\delta^3};
    % 4th row
    \zeronode{10}{4}{dddd0}{\delta^4};
    \zeronode{12}{4}{dddd1}{\delta^4};
    % 5th row
    \zeronode{11}{5}{ddddd0}{\delta^5};
    % 6th row
    \zeronode{10}{6}{dddddd0}{};
    \zeronode{12}{6}{dddddd1}{};
    
    %%  Draw connecting lines in Hasse diagram
    % 0th row to 1st row
    \hasseconnector{sd0}{d0};
    \hasseconnector{sd1}{d0};
    % 1st row to 2nd row
    \hasseconnector{d0}{dd0};
    \hasseconnector{d0}{dd1};
    % 2nd row to 3rd row
    \hasseconnector{dd0}{ddd0};
    \hasseconnector{dd1}{ddd0};
    % 3rd row to 4th row
    \hasseconnector{ddd0}{dddd0};
    \hasseconnector{ddd0}{dddd1};
    % 4th row to 5th row
    \hasseconnector{dddd0}{ddddd0};
    \hasseconnector{dddd1}{ddddd0};
    % 5th row to 6th row
    \hasseconnector{ddddd0}{dddddd0};
    \hasseconnector{ddddd0}{dddddd1};
    
    % Draw dotted lines up the side, identifying cylinder
    \cylinderboundary{sd0}{dd0};
    \cylinderboundary{sd1}{dd1};
    \cylinderboundary{dd0}{dddd0};
    \cylinderboundary{dd1}{dddd1};
    \cylinderboundary{dddd0}{dddddd0};
    \cylinderboundary{dddd1}{dddddd1};
  \end{tikzpicture}

\caption[%
  Nonzero and component paths within $\Subpath(\OQ)$
]{%
  \label{fig:nonzero-and-component-paths-in-OQ}
  \textit{Nonzero and component paths within $\Subpath(\OQ)$.} The subpath poset of $\OQ$ forms two cylinders, one for each connected component. The nonzero paths are printed in black, the others in grey. Amongst the nonzero paths are the components, highlighted.
}
\end{thesisfigure}

\nsbstrips The nonzero paths and components of an SB algebra can be accessed in \sbstrips\ by the attributes \mono{NonzeroPathsOfSBAlg} and \mono{ComponentsOfSBAlg}.

\begin{Verbatim}[commandchars=!@|,fontsize=\small,frame=single,baselinestretch=1.2]
  !gapprompt@gap>| !gapinput@NonzeroPathsOfSBAlg( alg1 );|
  [ v1, v2, v3, v4, a_over, b_over, c_over, d_over, a_over*b_over,
  b_over*c_over, c_over*a_over, d_over^2, a_over*b_over*c_over,
  b_over*c_over*a_over, d_over^3 ]
  !gapprompt@gap>| !gapinput@ComponentsOfSBAlg( alg1 );|
  [ a_over*b_over*c_over, b_over*c_over*a_over ]
\end{Verbatim}

\nlpass{Permissible data} The \define{permissible data} of a SB algebra $A$ is the tuple $(\OQ, N, C, \dagger)$, where $\OQ$ is an overquiver for $A$, $C \subseteq N \ideal \Subpath(\OQ)$ are the collections of components and nonzero paths, and $\dagger$ is the compatible pair of involutions $\ilarrow{ C }{}{ C }$ and $\ilarrow{ \OQ_0 }{}{ \OQ_0 }$, as above.

\nsbstrips The \sbstrips\ attribute \mono{PermDataOfSBAlg} stores the set $N \setminus C$ of paths with linearly independent residues and the set of components $C$, obtained from $(\OQ,N,C,\dagger)$.

\begin{Verbatim}[commandchars=!@|,fontsize=\small,frame=single,baselinestretch=1.2]
  !gapprompt@gap>| !gapinput@PermDataOfSBAlg( alg1 );|
  [ [ v1, v2, v3, v4, a_over, b_over, c_over, d_over, a_over*b_over,
      b_over*c_over, c_over*a_over, d_over^2, d_over^3 ],
  [ a_over*b_over*c_over, b_over*c_over*a_over ] ]
\end{Verbatim}

\subsection{Encodings of permissible data}

We encode the sets $C \subseteq N$ of permissible data $(\OQ, N, C, \dagger)$ numerically, in terms of an integer sequence and a $\{0,1\}$-sequence, both of which are indexed by $\OQ_0$. This is convenient for our calculations and for highlighting connections between SB algebras and Nakayama algebras. We will see shortly how it generalises the notion of admissible sequences used in the study of the latter. Many of the ideas and results in this section generalise those presented in \cite{Ful68} for Nakayama algebras; indeed, we have used the adjective \emph{permissible} as though it were the next evolution of \emph{admissible}, which Fuller employed in his notation.

Assume for the following that the permissible data $(\OQ,N,C,\dagger)$ of a SB algebra are fixed.

\nlpass{Dependent paths} Recall that any path $p \in C$ has $A$-residue depending linearly on another $A$-path (namely, that path represented by $p^\dagger$).

Write $N^\complement \ceq \Subpath(\OQ)\setminus N$. All $\OQ$-paths in $N^\complement$ have zero $A$-residue (by definition), which in a trivial way depends linearly on another $A$-path. 

Consequently, the $A$-residue of any path $p \in N^\complement \cup C$ depends linearly on the some other $A$-path. For this reason, we sometimes refer to paths $N^\complement \cup C$ as the set of \define{dependent paths}. In contrast, the $A$-residues of the remaining $\OQ$-paths (which are exactly those in $N \setminus C$) are linearly independent of all other $A$-paths.

Since $N$ is finite, all sufficiently long $\OQ$-paths are dependent paths. It follows that for any $i\in \OQ$ there exists a shortest dependent path with source $i$ and a shortest dependent path with target $i$.

\nlpass{Source encoding}\label{def:source-encoding} For $i \in \OQ_0$, let $\mdefine{a_i} \ceq \max \big\{ \len p \geq 0 \colon \source(p)=i \text{ and } p \in N \big\}$ and let
$$
    \mdefine{b_i}
    \ceq 
    \left\{
     \begin{array}{ll}
      0
       & \exists p \in C \text{ with } \source(p)=i\tcomma
        \\
      1
       & \text{otherwise.}
      \end{array}
      \right.
  $$
The \define{source encoding} of the permissible data is the pair $\big( (a_i)_{i}, (b_i)_{i} \big)$. Since $(a_i)_i$ is $\integers$-valued and $(b_i)_i$ is $\{0,1\}$-valued, we respectively call these the \define{integer sequence} and \define{bit sequence} of the encoding.\footnote{\ What we here call $(a_i)_i$ should be compared to the sequence $(c_i + 1)_i$ of \cite{Ful68} (which is Fuller's sequence $(c_i)_i$ with all terms incremented by $1$).}

For some $i \in \OQ_0 $, consider $p \ceq \big(\ilslpath{i}{a_i+b_i}\big)$. If $b_i=0$, then $p$ is the longest nonzero path having source $i$, a vertex which is the source of a component. Components are maximal among nonzero paths so, in fact, $p \in C$. Alternatively, if $b_i=1$, then $p$ covers the path $q \ceq \big(\ilslpath{i}{a_i}{}\big)$ in the prefix order. Since $q$ is the longest path in $N$ with source $i$, we deduce that $p$ is the shortest path $ N^\complement$ with source $i$.

Combining these two observations, we deduce that $a_i+b_i$ is the length of the shortest dependent path having source $i$. 

\nlpass{Target encoding} For $i \in \OQ_0$, let $\mdefine{c_i} \ceq \max \big\{ \len p \geq 0 \colon \target(p)=i \text{ and } p \in N \big\}$ and let
$$
    \mdefine{d_i}
    \ceq 
    \left\{
     \begin{array}{ll}
      0
       & \exists p \in C \text{ with } \target(p)=i\tcomma
        \\
      1
       & \text{otherwise.}
      \end{array}
      \right.
  $$
This integer sequence and bit sequence comprise the \define{target encoding} $\big( (c_i)_{i}, (d_i)_{i} \big)$ of the permissible data.\footnote{\ Similarly, what we here call $(b_i)_i$ should be compared to the incremented version $(d_i + 1)_i$ of the sequence Fuller calls $(d_i)_i$ \cite{Ful68}.}

We remark that, comparable to the previous passage, $c_i+d_i$ is the length of the shortest dependent path with target $i$.

\nrex\label{example:encodings-of-running-example} Recall our running example algebra, with the permissible data as constructed in Passages \ref{example:overquiver-and-dagger-of-running-example} and \ref{example:nonzero-paths-and-components-in-running-example}. It has the following encoding.
$$
\begin{array}{ccccc}
i
  & a_i
    & b_i
      & c_i
        & d_i
          \\
\hline
\source(\alpha)
  & 3
    & 0
      & 3
        & 0
          \\
\source(\beta)
  & 3
    & 0
      & 3
        & 0
          \\
\source(\gamma)
  & 2
    & 1
      & 2
        & 1
          \\
\source(\delta)
  & 3
    & 1
      & 3
        & 1
\end{array}
$$

\nsbstrips These encodings can be accessed using the following commands in \sbstrips.

\begin{Verbatim}[commandchars=!@|,fontsize=\small,frame=single,baselinestretch=1.2]
  !gapprompt@gap>| !gapinput@source_enc := SourceEncodingOfPermDataOfSBAlg( alg1 );|
  [ <vertex-indexed integer sequence>, <vertex-indexed bit sequence> ]
  !gapprompt@gap>| !gapinput@Display( source_enc[1] ); Display( source_enc[2] );|
  <integer sequence indexed by vertices of <quiver with 4 vertices and
  4 arrows>>
    v1 := 3,
    v2 := 2,
    v3 := 3,
    v4 := 3
  <bit sequence indexed by vertices of <quiver with 4 vertices and
  4 arrows>>
    v1 := 0,
    v2 := 1,
    v3 := 0,
    v4 := 1
  !gapprompt@gap>| !gapinput@TargetEncodingOfPermDataOfSBAlg( alg1 );|
  [ <vertex-indexed integer sequence>, <vertex-indexed bit sequence> ]
\end{Verbatim}

\nlpass{Retrieving $N$ and $C$ from encodings} These sequences encode $N$ and $C$ since
$$
\begin{array}{rcl}
  N
    & =
      & \big\{
          \big(\ilslpath{i}{\ell}\big) \in \Subpath(\OQ)
      \colon
          i \in \OQ_0,~ \ell \leq a_i
      \big\}
        \\
    & =
      & \big\{
          \big(\illtpath{\ell}{i}\big) \in \Subpath(\OQ)
      \colon
          i \in \OQ_0,~ \ell \leq c_i
      \big\}
\end{array}
$$
and
$$
\begin{array}{lcr}
  C
    & =
      & \big\{
        \big(\ilslpath{i}{a_i}\big) \in \Subpath(\OQ)
    \colon
        i \in \OQ_0 \text{ satisfies } b_i = 0  
    \big\}
      \\
    & =
      & \big\{
          \big(\illtpath{c_i}{i}\big) \in \Subpath(\OQ)
      \colon
          i \in \OQ_0 \text{ satisfies } d_i = 0  
      \big\}
      \tstop
\end{array}
$$

\nlpass{Compatibility with opposites} If the SB algebra $A$ can be described by permissible data $(\OQ, N, C, \dagger )$, then $A^\op$ can be described by ``opposite'' permissible data $(\OQ^\op, N^\op, C^\op, \dagger^\op)$, the components of which are defined in an obvious way. Then the source and target encodings of the permissible data for $A$ coincide with the target and source encoding of the permissible data $A^\op$ respectively.

\nrmk One can of course obtain the source encoding from the target encoding (or vice versa) via the intermediate step of $N$ and $C$. Alternatively, one can move between the two directly using the following formulas: 
\begin{equation}\label{eqn:ci-di-from-ai-bi}
c_i = \max\{r \geq 0 \colon r \leq a_{i+r}\}\tcomma
\text{\hspace{1cm}}
d_i = \left\{
    \begin{array}{ll}
    0   & \exists j \in \OQ_0~ (b_j = 0 \text{ and } j-a_j = i)\tcomma    \\
    1   & \text{otherwise}
    \end{array}
\right.
\end{equation}
and, dually,
$$
a_i = \max\{r \geq 0 \colon r \leq c_{i-r}\}\tcomma
\text{\hspace{1cm}}
b_i = \left\{
    \begin{array}{ll}
    0   & \exists j \in \OQ_0~ (d_j = 0 \text{ and } j+c_j = i)\tcomma    \\
    1   & \text{otherwise.}
    \end{array}
\right.
$$
The formula for $c_i$ is analogous to the formula in  \cite[Thm 2.2(b)]{Ful68}.

Since these pairs of formulas are dual, we will explain only the first pair: those in (\ref{eqn:ci-di-from-ai-bi}). Towards this goal, fix some $i \in \OQ_0$.

For any $r\geq 0$ we evidently have $\big(\illtpath{r}{i}\big)=\big(\ilslpath{i+r}{r}\big)$. By the definition of $a_{i+r}$, this path belongs to $N$ iff $r \leq a_{i+r} $. Because $c_i$ is the longest path in $N$ with target $i$, we are able to deduce that $c_i = \max\{r \geq 0 \colon r \leq a_{i+r}\}$, as claimed.

Additionally, $d_i = 0$ iff $i$ is the target of some component $p \in C$ (necessarily a unique one, since distinct components are incomparable in the subpath order and thus, in particular, have distinct targets). The source $j$ of this component satisfies $b_j=0$ and it must have length $a_j$ thus, when this $p$ exists, we have $p=\big(\ilslpath{j}{a_j}\big)=\big(\illtpath{a_j}{i}\big)$ whence we find $j-a_j = i$. This establishes the formula for $d_i$.

\newcommand{\OT}{\mathcal{T}}

\nlpass{Retrieving $A$ from permissible data} Let $\OT$ be the subquiver of $\OQ$ having vertex set $\OT_0 = \OQ_0$ and arrow set $\OT_1 = \{ \alpha \in \OQ_1 \colon \alpha \in N \}$, and for any component $p \in C$ we set $\kappa_{p, p^\dagger} \ceq 1\in \kk$.

We find that, in the notation of \cite[\S1]{WW85}, the connected components of $\OT$ are \emph{tracks} that identify $A$ up to the coefficients $\kappa_{p,p^\dagger}$. We can therefore essentially reconstruct $A$ from its permissible data.

Representing a SB algebra by permissible data does lose the coefficients in each commutativity relation $\lambda p - \mu q$ but, when it comes to calculating syzygies of string modules, all that matters is that that $\lambda$ and $\mu$ are both nonzero.\footnote{\ This point is implicit in Section 2, paragraph 2 of \cite{LM04}, the article where Liu and Morin prove that the syzygy of a string module over a SB algebra $A$ is a string module, because their assumption that all of the skew commutativity relations $\lambda p - \mu q$ in the definition of $A$ satisfy $\lambda = \mu$ does not play a important role in their proof.} Consequently, we are justified in presuming that they are both $1 \in \kk$ and otherwise forgetting about coefficients altogether.

\begin{lemma}
For any $i \in \OQ_0$ we have $a_i - 1 \leq a_{i-1}$ and, dually, $ c_i - 1 \leq c_{i+1} $. (Compare this with \cite[Thm 1.1(b)(c)]{Ful68}.)
\end{lemma}

\proof If $a_i=0$ then there is nothing to prove, so assume that $a_i > 0$. We find that $\big(\ilslpath{i}{a_i}{} \big) \in N$, and so its longest strict suffix $\big(\ilslpath{i-1}{a_i-1}{} \big)\in N$ also. The definition of the source encoding implies $a_{i-1} \geq a_i - 1$, as claimed. The statement about $c_i$ is proven dually.\qed

\begin{corollary}\label{cor:a_i-and-addition}
For any $i \in \OQ_0$ and $r \geq 0$ we have $a_i - r \leq a_{i-r}$ and, dually, $c_i - r \leq c_{i+r}$.
\end{corollary}

\proof The claim is a tautology when $r = 0$. Otherwise we iterate the proof of the previous lemma.\qed

% \npass A result of Fuller \cite{Ful68} characterises Nakayama algebras using integer sequences known as \emph{admissible sequences}. In Proposition \ref{prop:permissible-sequences-for-sb-algebras} we show how encodings of permissible data generalise this result to SB algebras. The proof requires slightly more notation, which we introduce below.

\nlpass{Decoupling and reduction of permissible data} We may relate a SB algebra $A$ to another $A'$ that has ``simpler'' permissible data; one of the two algebras will always be a quotient of the other.

We first remind the reader that if $p$ is a component of a commutativity relation then the one-arrow superpaths $\alpha p, p\alpha$ have zero $A$-residue for any arrow $\alpha\in Q_1$. We demonstrated this earlier, in our proof of Proposition \ref{prop:properties-of-sb-algebras}\ref{subprop:rejection-lemma}, where there was no assumption that the coefficients of the commutativity relation $\lambda p - \mu q$ satisfy $\lambda = 1 = \mu$. Crucially, this means that in the presence of any commutativity relation $\lambda p - \mu q$, the one-arrow superpaths of components may be considered as redundant relations for the algebra.

To \define{decouple} permissible data $(\OQ,N,C,\dagger)$ is to remove from $C$ a $\dagger$-pair of components; that is
$$
\ilmapsto{(\OQ, N, C, \dagger)}{}{ (\OQ, N, C', \dagger )}
$$
where $C' \ceq C \setminus\{p,p^\dagger\}$ for some $p \in C$. The algebra $A'$ with permissible data $(\OQ, N, C', \dagger)$ has the same presentation as $A$, except that the commutativity relation $r$ with components $p$ and $p^\dagger$ is replaced by the the shortest strict superpaths of those components, namely $\alpha p, p \beta, \gamma p^\dagger, p^\dagger \delta$ for appropriate $\alpha,\beta,\gamma,\delta \in Q_0$ whenever these arrows exist. It follows that $r$ spans a two-sided ideal of $A'$ and that $A'/\langle r \rangle = A$.

To \define{reduce} permissible data $(\OQ, N, C, \dagger)$ is to remove from $N \setminus C$ a maximal nonstationary element; ie,
$$\ilmapsto{ ( \OQ, N, C, \dagger ) }{}{ (\OQ, N', C, \dagger) }
$$
where $N' \ceq N - \left\{ p \right\}$ for some $p \in N\setminus C$, maximal among paths in $N$, satisfying $\len p > 0$. Writing $A'$ for the SB algebra with permissible data $(\OQ, N', C, \dagger)$, we find that $A' = A/\langle p \rangle$.

Decoupling and reduction diminish the finite sets $C$ and $N$ respectively. This is perhaps best seen through an example.

\nrex\label{example:decoupling-and-reducing-perm-data} We continue our running example, the permissible data of which appeared in Passage \ref{example:nonzero-paths-and-components-in-running-example}. Recall from there that
\begin{equation}\label{eqn:C-and-N-for-rex}
C = \{ \alpha\beta\gamma, \beta\gamma\alpha \} \subset \{e_{\source(\alpha)}, e_{\source(\beta)}, e_{\source(\gamma)}, e_{\source(\delta)}, \alpha, \beta, \gamma, \alpha\beta, \beta\gamma, \gamma\alpha, \alpha\beta\gamma, \beta\gamma\alpha,\delta,\delta^2,\delta^3\} = N
\tstop
\end{equation}
By first decoupling $\{\alpha\beta\gamma,\beta\gamma\alpha\}$ and then sequentially reducing $\delta^3$, $\delta^2$, $\delta$, $\beta\gamma\alpha$, $\alpha\beta\gamma$, $\gamma\alpha$, $\beta\gamma$, $\alpha\beta$, $\gamma$, $\beta$ and $\alpha$, we diminish $C$ to the empty set $\emptyset$ and we diminish $N$ to the set $\{e_{\source(\alpha)},e_{\source(\beta)},e_{\source(\gamma)},e_{\source(\delta)}\}$ comprising exactly the stationary paths. The decoupling turns $C$ into $\emptyset$ in a single step; each reduction removes from one element from $N$, starting from the furthest right (as written above in (\ref{eqn:C-and-N-for-rex})) and working leftwards.

\npass The effect that decoupling and reduction have on the encodings of permissible data is captured by the following lemma. The reader may also find it useful to see the effects in practice, as in Figure \ref{fig:dec-and-red-perm-data}.
\begin{thesisfigure}
\centering

\subcaptionbox{%
  \label{subfig:initial-encodings}
  Initially
}[%
  0.3\textwidth
]{
  $$
\begin{array}{ccccc}
i
  & a_i
    & b_i
      & c_i
        & d_i
          \\
\hline
\source(\alpha)
  & 3
    & 0
      & 3
        & 0
          \\
\source(\beta)
  & 3
    & 0
      & 3
        & 0
          \\
\source(\gamma)
  & 2
    & 1
      & 2
        & 1
          \\
\source(\delta)
  & 3
    & 1
      & 3
        & 1
\end{array}
$$
}
\subcaptionbox{%
  \label{subfig:after-decoupling}
  Decouple $\{\alpha\beta\gamma,\beta\gamma\alpha\}$
}[%
  0.3\textwidth
]{
  $$
\begin{array}{ccccc}
i
  & a_i
    & b_i
      & c_i
        & d_i
          \\
\hline
\source(\alpha)
  & 3
    & 1
      & 3
        & 1
          \\
\source(\beta)
  & 3
    & 1
      & 3
        & 1
          \\
\source(\gamma)
  & 2
    & 1
      & 2
        & 1
          \\
\source(\delta)
  & 3
    & 1
      & 3
        & 1
\end{array}
$$
}
\subcaptionbox{%
  \label{subfig:after-reducing-deltas}
  Reduce $\delta^3,\delta^2,\delta$
}[%
  0.3\textwidth
]{
  $$
\begin{array}{ccccc}
i
  & a_i
    & b_i
      & c_i
        & d_i
          \\
\hline
\source(\alpha)
  & 3
    & 1
      & 3
        & 1
          \\
\source(\beta)
  & 3
    & 1
      & 3
        & 1
          \\
\source(\gamma)
  & 2
    & 1
      & 2
        & 1
          \\
\source(\delta)
  & 0
    & 1
      & 0
        & 1
\end{array}
$$
}

\subcaptionbox{%
  \label{subfig:after-reducing-betagammaalpha}
  Reduce $\beta\gamma\alpha$
}[%
  0.3\textwidth
]{
  $$
\begin{array}{ccccc}
i
  & a_i
    & b_i
      & c_i
        & d_i
          \\
\hline
\source(\alpha)
  & 3
    & 1
      & 3
        & 1
          \\
\source(\beta)
  & 2
    & 1
      & 2
        & 1
          \\
\source(\gamma)
  & 2
    & 1
      & 2
        & 1
          \\
\source(\delta)
  & 0
    & 1
      & 0
        & 1
\end{array}
$$
}
\subcaptionbox{%
  \label{subfig:after-reducing-alphabetagamma-gammaalpha}
  Reduce $\alpha\beta\gamma,\gamma\alpha$
}[%
  0.3\textwidth
]{
  $$
\begin{array}{ccccc}
i
  & a_i
    & b_i
      & c_i
        & d_i
          \\
\hline
\source(\alpha)
  & 2
    & 1
      & 1
        & 1
          \\
\source(\beta)
  & 2
    & 1
      & 2
        & 1
          \\
\source(\gamma)
  & 1
    & 1
      & 2
        & 1
          \\
\source(\delta)
  & 0
    & 1
      & 0
        & 1
\end{array}
$$
}
\subcaptionbox{%
  \label{subfig:after-all-reductions}
  Remaining reductions
}[%
  0.3\textwidth
]{
  $$
\begin{array}{ccccc}
i
  & a_i
    & b_i
      & c_i
        & d_i
          \\
\hline
\source(\alpha)
  & 0
    & 1
      & 0
        & 1
          \\
\source(\beta)
  & 0
    & 1
      & 0
        & 1
          \\
\source(\gamma)
  & 0
    & 1
      & 0
        & 1
          \\
\source(\delta)
  & 0
    & 1
      & 0
        & 1
\end{array}
$$
}

\caption[%
  Effect of decoupling and reduction on encodings
]{%
  \label{fig:dec-and-red-perm-data}
  \textit{Effect of decoupling and reduction on encodings.} We show snapshots of the permissible data of our running example algebra as it undergoes the sequence of decoupling and reduction explained in Example \ref{example:decoupling-and-reducing-perm-data}. The subfigures are to be read in order, with the mentioned decouplings or reductions being cumulative.
}

\end{thesisfigure}

\begin{lemma}
Suppose $(\OQ, N, C, \dagger)$ are the permissible data of an algebra $A$, and their source and target encodings are respectively $\big( (a_i)_i, (b_i)_i \big)$ and $\big( (c_i)_i, (d_i)_i \big)$.
\begin{enumerate}
    \item\label{sublem:decoupling-and-encodings} If $p \in C$ then $b_{\source(p)}=d_{\target(p)}=b_{\source(p^\dagger)}=d_{\target(p^\dagger)}=0$. Decoupling by $\{p,p^\dagger\}$ increments all of these values by $1$ and leaves all other $b_i,d_i$ unchanged.
    
    \item\label{sublem:reduction-and-encodings} If $p \in N$ is maximal, then we have $a_{\source(p)}=c_{\target(p)} = \len p$. If furthermore $p \notin C$ then  reducing by $p$ leaves all $a_i,b_i,c_i,d_i$ unchanged except for $a_{\source(p)},c_{\target(p)}$, which are decremented by $1$.
    
    \item\label{sublem:sum-sequences-and-encodings} It follows that decoupling by $\{p,p^\dagger\}$ increments $$a_{\source(p)} + b_{\source(p)} = c_{\target(p)} + d_{\target(p)}
    \text{\hspace{1cm}and\hspace{1cm}}
    a_{\source(p^\dagger)} + b_{\source(p^\dagger)} = c_{\target(p^\dagger)} + d_{\target(p^\dagger)}$$ by $1$ but leaves all other $a_i+b_i, c_i+d_i$ unchanged, while reducing by $p$ decrements $$a_{\source(p)} + b_{\source(p)} = c_{\target(p)} + d_{\target(p)}$$ by $1$, but leaves all other $a_i+b_i,c_i+d_i$ are unchanged.
\end{enumerate}
\end{lemma}

\proof Parts \ref{sublem:decoupling-and-encodings} and \ref{sublem:reduction-and-encodings} follow directly from the definitions of the encodings of $(\OQ, N, C, \dagger)$. Part \ref{sublem:sum-sequences-and-encodings} is a consequence.\qed

\npass The preceding lemma performs all the heavy lifting for the following proposition, which shows a connecting between the source encoding of the permissible data and the target encoding. We consider the result below (and part \ref{sublem:ais-permutation-of-cis} of it especially) to be the analogue for SB algebras of the Fuller's \cite[Thm 2.2(a)]{Ful68} for Nakayama algebras. Our proof proceeds along similar very similar lines to his.

\begin{proposition}\label{prop:permissible-sequences-for-sb-algebras} Let $(\OQ, N, C, \dagger)$ be data and $\big( (a_i)_i, (b_i)_i \big),\big( (c_i)_i, (d_i)_i \big)$ its encodings as above.
\begin{enumerate}
    \item\label{sublem:ais-permutation-of-cis} There exists a permutation $\pi \colon \ilarrow{\OQ_0}{}{\OQ_0}$ such that $a_{i\pi}=c_{i}$.
    \item\label{lem:bit-seqs-of-encodings-are-mutual-permutations} There exists a permutation $\pi \colon \ilarrow{\OQ_0}{}{\OQ_0}$ such that $b_{i\pi}=d_{i}$.
    \item There exists a permutation $\pi \colon \ilarrow{\OQ_0}{}{\OQ_0}$ such that $a_{i\pi}+b_{i\pi}=c_{i}+d_{i}$.
\end{enumerate}
\end{proposition}

\proof We proceed by on induction on $m \ceq |N|+|C|$.

The base case is when $N=\{e_i \colon i \in \OQ_0\}$ and $C=\emptyset$. Here, $a_i=c_i=0$ and $b_i=d_i=1$ for all $i \in \OQ_0$ and the identity permutation suffices.

Otherwise assume $m > |\OQ_0|$, in which case it is possible to decouple or reduce the permissible data. From the previous lemma we know that either step has the same impact on at most two values of $(a_i)_i$ and $(b_i)_i$, two values of $(c_i)_i$ and $(d_i)_i$ or two values of $(a_i+b_i)_i$ and $(c_i+d_i)_i$, that impact being to change a common value by $1$ in each case. We find that the multiset of values taken by these sequences coincide pairwise for our permissible data iff they do for the decoupled or reduced data. The claim then follows from the inductive hypothesis.\qed

\section{Syllables}

In this section we introduce and explore symbols that we call \emph{syllables}. We consider these syllables to be a refinement of the identically-named symbols that other articles in the literature, among them \cite[\S2]{HZS05} and \cite[\S1]{HZ16}, have used in related topics. To best motivate some of our additional terminology, such as \emph{orientation}, we urge the reader to watch these symbols in action in Subsection \ref{subsec:preview-of-strips}.

We have some permissible data $(\OQ,N,C,\dagger)$ of $A$ fixed throughout.

\subsection{Definition of syllables}

In this section, we define the symbols with which we will formally describe string modules for SB algebras.

We assume that the permissible data $(\OQ,N,C,\dagger)$ of the SB algebra $A$ is fixed, and we write $\big((a_i)_i, (b_i)_i\big)$ and $\big((c_i)_i, (d_i)_i\big)$ for its source and target encodings respectively. For reference, we note that the permissible data of the running example algebra were found in Passages \ref{example:overquiver-and-dagger-of-running-example} and \ref{example:nonzero-paths-and-components-in-running-example} and that the encodings appear in \ref{example:encodings-of-running-example}.

\nlpass{Syllables}\label{def:syllables} Formally speaking, a \define{syllable} is a tuple $(p,\ep,s)$ that comprises a nonzero noncomponent path $p \in N \setminus C$ (which we call the \define{underlying path}), a bit $\ep \in \{0,1\}$ (which we call the \define{stability term}) and a sign $s \in \{+1,-1\}$ (which we call the \define{orientation}). The tuple must satisfy the property $0 < \len p + \ep$.

We write \define{$\Syll(A)$} for the \define{syllable set} of $A$ (with respect to the prescribed permissible data). The definition of the previous paragraph can be condensed into the following single, very terse, formula, that we use in the next subsection (and which appears in the source code of the \sbstrips\ package).
$$
\Syll(A) = \{(p, \ep, s) \in N \times \{0,1\} \times \{+1,-1\} \colon
0 < \len p + \ep < a_{\source(p)} + b_{\source(p)} + \ep
\}\tstop
$$

We will denote the syllable $(p,\ep,s)$ as $\big( \ilsyll{\syllableplaceholder}{ p }{\ep} \big)^s$. We omit the orientation $s$ when it is clear from context or is unimportant (which is almost all of the time). From time to time we will specify the underlying path $p$ by its length $\ell$ and either its source or target vertex $i$, in which case we will use the notation $\big( \ilsyll{ i }{ \ell }{ \ep } \big)^s$ or $\big( \injsyll{ \ell }{ i }{ \ep } \big)^s$. We can write these inline as $(i,\ell,\ep,s)$ or $(\ell,i,\ep,s)$.

A syllable $\big( \ilsyll{\syllableplaceholder}{p}{\ep} \big)^s$ inherits some adjectives from its constituent parts. It is \define{stationary} or \define{nonstationary} iff $p$ is and it is \define{positive} or \define{negative} iff $s$ is. We call the syllable an \define{interior} syllable if $\ep=0$ and a \define{boundary} syllable if $\ep = 1$. Observe that the condition $0 < \len p + \ep$ ensures that all stationary syllables are boundary.

The \define{compression} of a syllable $\bdp \ceq \big(\ilsyll{i}{\ell}{\ep}\big)$ is the path $q \ceq \big(\ilslpath{i}{\ell+\ep}\big)$. Note that generally this is different from the underlying path $\big(\ilslpath{i}{\ell}\big)$. Using the compression, we can define the \define{source} $\mdefine{\source(\bdp)}\ceq\source(q)$ and \define{target} $\mdefine{\target(\bdp)}\ceq \target(q)$ of a syllable.

\nlpass{Stationary syllables} Stationary syllables are denoted $ \mdefine{ \bde_i } \ceq \big( \ilsyll{\syllableplaceholder}{e_i}{1} \big)$, or $\bde_i^{\pm1}$ when we need to mention orientation.

\nrex The syllables for the running example algebra $A$ are listed in Figure \ref{fig:syllables-for-running-example-algebra}.
\begin{thesisfigure}
\centering

\begin{tabular}{cccc}
$\big(\ilsyll{\syllableplaceholder}{e_{\source(\alpha)}}{1}\big)^{\pm1}$,
 & $\big(\ilsyll{\syllableplaceholder}{e_{\source(\beta)}}{1}\big)^{\pm1}$,
   & $\big(\ilsyll{\syllableplaceholder}{e_{\source(\gamma)}}{1}\big)^{\pm1}$,
     & $\big(\ilsyll{\syllableplaceholder}{e_{\source(\delta)}}{1}\big)^{\pm1}$,
       \\
$\big(\ilsyll{\syllableplaceholder}{\alpha}{0}\big)^{\pm1}$,
 & $\big(\ilsyll{\syllableplaceholder}{\beta}{0}\big)^{\pm1}$,
   & $\big(\ilsyll{\syllableplaceholder}{\gamma}{0}\big)^{\pm1}$,
     & $\big(\ilsyll{\syllableplaceholder}{\delta}{0}\big)^{\pm1}$,
       \\
$\big(\ilsyll{\syllableplaceholder}{\alpha}{1}\big)^{\pm1}$,
 & $\big(\ilsyll{\syllableplaceholder}{\beta}{1}\big)^{\pm1}$,
   & $\big(\ilsyll{\syllableplaceholder}{\gamma}{1}\big)^{\pm1}$,
     & $\big(\ilsyll{\syllableplaceholder}{\delta}{1}\big)^{\pm1}$,
       \\

$\big(\ilsyll{\syllableplaceholder}{\alpha\beta}{0}\big)^{\pm1}$,
 & $\big(\ilsyll{\syllableplaceholder}{\beta\gamma}{0}\big)^{\pm1}$,
   & $\big(\ilsyll{\syllableplaceholder}{\gamma\alpha}{0}\big)^{\pm1}$,
     & $\big(\ilsyll{\syllableplaceholder}{\delta^2}{0}\big)^{\pm1}$,
       \\
$\big(\ilsyll{\syllableplaceholder}{\alpha\beta}{1}\big)^{\pm1}$,
 & $\big(\ilsyll{\syllableplaceholder}{\beta\gamma}{1}\big)^{\pm1}$,
   & $\big(\ilsyll{\syllableplaceholder}{\gamma\alpha}{1}\big)^{\pm1}$,
     & $\big(\ilsyll{\syllableplaceholder}{\delta^2}{1}\big)^{\pm1}$,
       \\
  & $\big( \ilsyll{\syllableplaceholder}{\delta^3}{0} \big)^{\pm1}$,
    & $\big( \ilsyll{\syllableplaceholder}{\delta^3}{1} \big)^{\pm1}$.
\end{tabular}

\caption[%
  Syllables for the running example algebra
]{%
  \label{fig:syllables-for-running-example-algebra}
  \textit{Syllables for the running example algebra.} The top row contains the stationary syllables.
}

\end{thesisfigure}

\nlpass{Blank syllable} In due course we will define some partial operations on $\Syll(A)$, which we can model as as total functions on a based version of the set. We will use the term \define{blank syllable} to mean the basepoint of $\Syll(A)$.

% \npass Next, we introduce some additional symbols for use within the formal calculation of the syzygy of a string module. These new symbols will be used to represent the occurrence of a pathological phenomenon within the syzygy calculation rather than representing any part of a string module. For this reason, we call these new symbols \emph{virtual syllables}.

\centresubsec{Virtual syllables}

\nlpass{Pathological syzygies} For any projective module $P$ there is the short exact sequence
$$\big(
\begin{tikzcd}[sep=small]
0 \ar[r]
  \& \soc P \ar[r]
    \& P \ar[r]
      \& P/\soc P\ar[r]
        \& 0
\end{tikzcd}\big)\tcomma
$$
from which we easily conclude $\syzygy^1(P/\soc P) = \soc P$. %In our framework for modelling string modules and their syzygies,
But when $P$ is (indecomposable and) pin, this syzygy short exact sequence of the string module $P/\soc P$ witnesses an unusual phenomenon. We will give the details later but, for now, suffice it to say that, upon calculating $\syzygy^1(P/\soc P)$ (when $P$ is indecomposable and pin), the two ``boundaries'' of $P/\soc P$ interact with one another before yielding the syzygy $\soc P$.

To mark this unusual interaction in our framework we need new, special symbols. These symbols are bespoke: they only signify this specific phenomenon, and they never represent any actual part of a string module. For this reason, we will call these new symbols \emph{virtual syllables}.

\nlpass{Virtual syllables} Recall from Remarks \ref{psg:socles-of-projectives} that the socle of a indecomposable pin module is a simple module at some $Q$-vertex $j$, which the target of some commutativity relation $\lambda p-\mu q$ for $A$. 
The components $p,q$ are represented by component $\OQ$-paths and $j$ is represented by some $\OQ$-paths $i$ which, being the targets of the component $\OQ$-paths, will satisfy $d_i=0$.

The \define{virtual syllables} at any such $i \in \OQ_0$ with $d_i =0$ are $(e_i,0,\pm1)=\big(\ilsyll{i}{0}{0}\big)^{\pm1}$.

\nrmks \ilitem{} Any reference to syllables excludes virtual syllables unless explicitly mentioned.

\ilitem{} It is unambiguous whether the symbol $\big( \ilsyll{\syllableplaceholder}{p}{\ep} \big)^{s}$ denotes a syllable or a virtual syllable, because $0 < \len p + \ep$ for syllables and $ 0 = \len p + \ep $ for virtual syllables.

\nrex The virtual syllables for the running example algebra are $\big( \ilsyll{\source(\alpha)}{0}{0} \big)^{\pm1}$ and $\big( \ilsyll{\source(\beta)}{0}{0} \big)^{\pm1}$.

\nsbstrips The syllable set and blank syllable of an algebra can be accessed using the commands \mono{SyllableSetOfSBAlg} and \mono{BlankSyllableOfSBAlg}. We mention that the set \mono{SyllableSetOfSBAlg} includes the blank syllables and the virtual syllables in addition to the ``actual'' syllables.

\begin{Verbatim}[commandchars=!@|,fontsize=\small,frame=single,baselinestretch=1.2]
  !gapprompt@gap>| !gapinput@SyllableSetOfSBAlg( alg1 );|
  [ ( ), ( v1, 0 ), ( v1, 1 ), ( v2, 1 ), ( v3, 0 ), ( v3, 1 ),
    ( v4, 1 ), ( b_over, 0 ), ( b_over, 1 ), ( c_over, 0 ),
    ( c_over, 1 ), ( a_over, 0 ), ( a_over, 1 ), ( d_over, 0 ),
    ( d_over, 1 ), ( b_over*c_over, 0 ), ( b_over*c_over, 1 ),
    ( c_over*a_over, 0 ), ( c_over*a_over, 1 ), ( a_over*b_over, 0 ),
    ( a_over*b_over, 1 ), ( d_over*d_over, 0 ), ( d_over*d_over, 1 ),
    ( d_over*d_over*d_over, 0 ), ( d_over*d_over*d_over, 1 ) ]
  !gapprompt@gap>| !gapinput@BlankSyllableOfSBAlg( alg1 );|
  ( )
\end{Verbatim}

\nlpass{Pin-boundary syllables} Virtual syllables turn up to describe the unusual syzygy behaviour of $\syzygy^1(P/\soc P)$. This means they are only ever seen near representations of socle-quotients $P/\soc P$ of pin modules.

The syllables arising in the representation of such socle-quotients need special attention. We call them \define{pin boundary syllables}. They are always of the form $\big(\ilsyll{i}{a_i+b_i-1}{1}\big)$ for $i \in \OQ_0$ satisfying $b_i=0$.

\subsection{Operations for syllables}

\npass We will formulate syzygy-taking as filling in an array. Excluding some boundary phenomena, the array is iteratively populated in two directions from some initial data, these two directions being downwards and sideways. Accordingly, the act of population is governed by two partial operations. The downwards-moving operation is called \emph{descent}; the sideways-moving operation, \emph{sidestep}. We define them both below.

Additionally, we define an operation used in describing the boundary phenomena.

\nlpass{Descent $\descent$} The following partial operation on syllables is denoted \define{$\descent$} and pronounced \define{descent}. For $\big(\ilsyll{i}{\ell}{\ep}\big)^s \in \Syll(A)$, we define
$$\big(\ilsyll{i}{\ell}{\ep}\big)^s \descent \ceq \big(\ilsyll{i-(\ell+\ep)}{a_i-(\ell+\ep)}{b_i}\big)^{-s}$$
whenever the righthand side is a syllable. The following lemma characterises $\supp \descent$.

\nrmk We highlight that if $\bdp$ has orientation $s$ (and if $\bdp\descent$ is defined) then $\bdp \descent$ has orientation $-s$, and that otherwise the orientation of $\bdp$ has no impact on $\bdp\descent$. Therefore we may safely omit mention of orientation of syllables.

\begin{lemma}
If $\bdp \ceq \big(\ilsyll{i}{\ell}{\ep}\big) \in \Syll(A)$, then $\bdp \notin \supp \descent$ iff $(\ell,\ep)=(a_i+b_i-1,1)$.
\end{lemma}

\proof Using the terse defining formula for $\Syll(A)$, we have that $\bdp\descent \in \Syll(A)$ iff
\begin{equation}\label{eqn:descent-support-inequality}
0 < a_i + b_i - (\ell+\ep) < a_{i-(\ell + \ep)} + b_{i - (\ell+\ep)} + b_i\tstop
\end{equation}
This is the conjunction of two strict inequalities. The above fails if either side does. We will find that the righthand inequality never fails, and that the lefthand inequality fails precisely when $(i, \ell, \ep)$ satisfies the given property.

In the first case, suppose that the righthand inequality of (\ref{eqn:descent-support-inequality}) fails. The negation of the righthand inequality is equivalent to
\begin{equation}\label{eqn:cancelled-negated-RH-support-ineq}
    a_{i-(\ell+\ep)} + b_{i-(\ell+\ep)} \leq a_i - (\ell+\ep)\tstop
\end{equation}
Using Corollary \ref{cor:a_i-and-addition} we know that $a_i - (\ell+\ep) \leq a_{i-(\ell+\ep)}$. Substituting this into (\ref{eqn:cancelled-negated-RH-support-ineq}) yields
\begin{equation}\label{eqn:case-i-syllable-sandwich-ineq}
a_{i-(\ell+\ep)} + b_{i-(\ell+\ep)} \leq a_i - (\ell+\ep) \leq a_{i-(\ell+\ep)}\tstop
\end{equation}

If $b_{i-(\ell+\ep)}=1$ then the contradiction $a_{i-(\ell+\ep)}+1 \leq a_{i-(\ell+\ep)}$ arises from (\ref{eqn:case-i-syllable-sandwich-ineq}). Else if $b_{i-(\ell+\ep)}=0$, then we deduce two facts: first, that $a_{i-(\ell+\ep)}=a_i-(\ell+\ep)$ from (\ref{eqn:case-i-syllable-sandwich-ineq}); second,  that $i-(\ell+\ep)$ is the source of the component $\big(\ilslpath{i-(\ell+\ep)}{a_{i-(\ell+\ep)}}\big)$ by definition. Combining these, we find that this component is a strict suffix of $q \ceq \big( \ilslpath{i}{a_i} \big)$. We have $q \in N$, by definition of $a_i$, and $q \notin N$ since it is a strict superpath of a component. We surmise that this case reduces to the absurd.

In the second case, suppose that the lefthand inequality of (\ref{eqn:descent-support-inequality}) fails. The negation of the lefthand inequality is equivalent to $a_i + b_i \leq \ell + \ep$. Since $(i,\ell,\ep)\in\Syll(A)$ we have $\ell+\ep < a_i + b_i + \ep$, from which we deduce
\begin{equation}\label{eqn:case-ii-syllable-sandwich-ineq}
a_i + b_i \leq \ell + \ep < a_i + b_i + \ep\tstop
\end{equation}
If $\ep=0$ then (\ref{eqn:case-ii-syllable-sandwich-ineq}) gives the absurdity $a_i + b_i < a_i + b_i$; this forces $\ep=1$. We deduce from (\ref{eqn:case-ii-syllable-sandwich-ineq}) that $\ell+\ep$ is sandwiched between consecutive integers. We infer $ a_i + b_i = \ell + \ep$, which is to say $\ell=a_i+b_i-1$. This proves the claim of the lemma.\qed

\begin{corollary}
If $\bdp$ is an interior syllable, then $\bdp \in \supp \descent$.
\end{corollary}

\proof Interior syllables have stability term $\ep=0 \neq 1$. We apply the previous lemma. \qed

\begin{corollary}
If $\bdp$ is a pin-boundary syllable, then $\bdp \notin \supp \descent$.
\end{corollary}

\proof Pin boundary syllables $\big(\ilsyll{i}{\ell}{\ep}\big)$ are precisely those with $(\ell,\ep)=(a_i + b_i - 1, 1)$ and $b_i=0$.\qed

\nlpass{Sidestep $\sidestep$} The operation on syllables denoted \define{$\sidestep$} and pronounced \define{sidestep} is defined a follows. For any syllable $\bdp \ceq \big( \ilsyll{i}{\ell}{\ep} \big)^s$, we set $\bdp \sidestep \ceq \bde_{\source(p)^\dagger}$.

This is the unique stationary syllable $\bde_j$ such that $\source(\bde_j)^\dagger = \source(\bdp)$. Clearly $\sidestep^3=\sidestep$.

\nlpass{Perturbation} The \define{perturbation} of any interior syllable $\big(\ilsyll{\syllableplaceholder}{p}{0}\big)^s$ is the corresponding boundary syllable $\big(\ilsyll{\syllableplaceholder}{p}{1}\big)^s$. This gives a partial operation on syllables, taking values on the boundary ones.

% What are peak and valley compatibility?

% What are $\descent$ and $\sidestep$?

\chapter{Patches, strips and the syzygy fabric}\label{chap:patches-strips-fabric}

In this chapter, we give the novel approach to syzygy calculations that underpins the \sbstrips\ package. Our approach involves representing string modules as words, and indecomposable projective modules as grids, made up of the syllables constructed in the previous chapter; the words are called \emph{strips} and the grids \emph{patches}. 

We first construct the patch set of a SB algebra and then demonstrate how to represent a string module as a strip and its projective cover as an associated line of patches. We verify that our method correctly calculates syzygies and then show how to record syzygies into an array. Following in the spirit of pin modules, string modules, patches and strips, this array also has a textile-based name: the \emph{syzygy fabric}. The purpose of the fabric is to facilitate the rigorous discussion of syzygy patterns given in the next chapter.

\section[Patches]{Patches of a SB algebra}

In this section, we introduce the patch set of a SB algebra $A$. This is a finite set whose members represent the indecomposable projective $A$-modules. We then verify important properties about our construction, to be used later when representing the syzygy operation.

\subsection{Construction of patches}\label{subsec:construction-of-patches}

In this subsection we show how to systematically construct the set of patches of a SB algebra in five steps. Throughout, we fix the permissible data $(\OQ,N,C,\dagger)$ of a SB algebra $A$ and we let $(a_i)_i, (b_i)_i$ and $(c_i)_i, (d_i)_i$ be the integer and bit sequences associated to the source and target encodings of the permissible data, as before.

As a preliminary remark, let us underscore that a patch will be a $(2 \times 2)$-grid populated with symbols.

\npass The orientations of any (nonblank) syllable in a patch are determined by its position in the patch, as shown in Figure \ref{fig:orientations-of-syllables-in-patches}. For this reason we omit mention of orientation in the forthcoming definitions. 
\begin{thesisfigure}
  \centering
  \begin{tikzpicture}[scale=\standardscale]
    \patchrep{(0,0)}{$-1$}{$+1$}{$+1$}{$-1$};
  \end{tikzpicture}
  
  \caption[%
    Orientations of syllables in patches
  ]{%
    \label{fig:orientations-of-syllables-in-patches} \textit{Orientations of syllables in patches.} %Any nonblank syllable in a patch has orientation befitting its location in the patch, as shown here. 
  }  
  
  \vspace{4ex}
  
    \begin{tikzpicture}[scale=\standardscale]
    \patchrep{(-3,0)}{ $\bdr_i$ }{ $\bdp$ }{}{ $\bdq$ }
    
    \patchrep{(3,0)}{ $\bdp$ }{ $\bdr_i$ }{ $\bdq$ }{}
  \end{tikzpicture}
  \caption[%
    Schematic of a patch with one pin-boundary syllable
  ]{%
    \label{fig:patch-with-one-pin-boundary-syllable}
    \textit{Schematic of a patch with one pin-boundary syllable.} For each $i$ with $b_i = 0$, a family of patches are created. We always have $\bdr_i \ceq (\ilsyll{i}{a_i-1}{1}  )$ the pin-boundary syllable at $i$, but $\bdp$ and $\bdq$ vary as described in Construction \ref{constr:patches-with-one-pin-boundary-syllable}.
  }
\end{thesisfigure}

\nsbstrips In our package, the five steps of patch construction are performed by the attribute \mono{PatchSetOfSBAlg}.

\centresubsec{Blank patch}

\nlpass{Construction} We construct a single patch with all cells blank. We call it the \define{blank patch} of the algebra. It is the basepoint of the set of patches of $A$.

\nsbstrips This blank patch of a SB algebra is returned by the attribute \mono{BlankPatchOfSBAlg}.

\centresubsec{Patches with no pin-boundary syllable in the top row}

\nlpass{Construction} For each pair $(\bdp,\bdp')$ of syllables such that $\source(\bdp)^\dagger = \source(\bdp')$ and neither $\bdp$ nor $\bdp'$ is pin boundary, we construct a patch. It has $\bdp$ above $\bdp\descent$ on one side and $\bdp'$ above $\bdp'\descent$ on the other.

From this original patch we then make additional patches if $\bdp$ is stationary or if $\bdp'$ is stationary, and we treat these two possibilities separately. If $\bdp$ is stationary, then we create a copy of the original patch except that $\bdp$ is replaced by the blank syllable (the remaining cells are unchanged). Independently, if $\bdp'$ is stationary, then we create a copy of the original patch except with $\bdp'$ replaced by the blank syllable (and the remaining cells unchanged). This means that if both $\bdp$ and $\bdp'$ are stationary, two copies are created. Each of these \define{amended patches} differs from the original patch in a single cell.

The symmetry in the conditions on $\bdp$ and $\bdp'$ mean that the reflection of any patch created at this step is also created at this step.

\centresubsec{Patches with one pin-boundary syllable in the top row}

\nlpass{Construction}\label{constr:patches-with-one-pin-boundary-syllable} Recall that $i\in \OQ_0$ has $b_i = 0$ iff it represents the source of a commutativity relation. When so, we write $\bdr_i$ for the pin boundary syllable at $i$.

For each such $i$ and each syllable $\bdp$ which is not pin boundary but does have $\source(\bdp)=\idag$, we let $\bdq$ be the perturbation of $\bdp\descent$ (ie, $\bdp\descent$ but with stability term $\ep=1$). We then create two patches. The first has $\bdr_i$ above a blank cell on one side and $\bdp$ above $\bdq$ on the other. The second is its reflection. These are illustrated in Figure \ref{fig:patch-with-one-pin-boundary-syllable}.

Similarly to the previous case, we make an \emph{amended} version of the patch if it features a stationary syllable in the top row. Specifically, if the top-row syllable $\bdp$ is stationary, we create a copy of the patches defined in the previous paragraph, except that the $\bdp$ is replaced by the blank syllable.

In terms of formulas, these patches are constructed for each triple comprising $i\in \OQ_0$, $\ell\geq 0$ and $\ep\in\{0,1\}$ such that $b_i=0$ and $0 < \ell + \ep < a_\idag$. Then we define $\bdp \ceq ( \ilsyll{\idag}{\ell}{\ep} )$ and $\bdq \ceq (\ilsyll{\idag-(\ell + \ep)}{a_\idag - (\ell + \ep)}{1})$.
\begin{thesisfigure}
  \subcaptionbox{%
    \label{subfig:running-example-opposite-pin-boundary-syllable-1} The four possible choices of $\bdp$ and $\bdq$ when $i=\source(\alpha)$.
  }[
    \textwidth
  ]{
    $$\begin{array}{l||lcl||l}
    \bdp
      & \ilsyll{\source(\beta)}{0}{1}
        & \hspace{1cm}
          & \bdp
            & \ilsyll{\source(\beta)}{1}{0}
        \\
    \bdq
      & \ilsyll{\source(\gamma)}{2}{1}
        & 
          & \bdq
            & \ilsyll{\source(\gamma)}{2}{1}
    \end{array}$$
    
    $$\begin{array}{l||lcl||l}
    \bdp
      & \ilsyll{\source(\beta)}{1}{1}
        & \hspace{1cm}
          & \bdp
            & \ilsyll{\source(\beta)}{2}{0}
              \\
    \bdq
      & \ilsyll{\source(\alpha)}{1}{1}
        & ~
          & \bdq
            & \ilsyll{\source(\alpha)}{1}{1}
    \end{array}$$
  }
  \subcaptionbox{%
    \label{subfig:running-example-opposite-pin-boundary-syllable-2} The four possible choice of $\bdp$ and $\bdq$ when $i=\source(\beta)$.
  }[%
    \textwidth
  ]{%
    $$\begin{array}{l||lcl||l}
    \bdp
      & \ilsyll{\source(\alpha)}{0}{1}
        & \hspace{1cm}
          & \bdp
            & \ilsyll{\source(\alpha)}{1}{0}
        \\
    \bdq
      & \ilsyll{\source(\beta)}{2}{1}
        & 
          & \bdq
            & \ilsyll{\source(\beta)}{2}{1}
    \end{array}$$
    
    $$\begin{array}{l||lcl||l}
    \bdp
      & \ilsyll{\source(\alpha)}{1}{1}
        & \hspace{1cm}
          & \bdp
            & \ilsyll{\source(\gamma)}{2}{0}
              \\
    \bdq
      & \ilsyll{\source(\alpha)}{1}{1}
        & ~
          & \bdq
            & \ilsyll{\source(\gamma)}{1}{1}
    \end{array}$$
  }
  \caption[
    Choices opposite $\bdr_i$, for the running example algebra $A$
  ]{
    \label{fig:example-choices-opposite-one-pin-boundary-syllable}\textit{Choices of $\bdp$ and $\bdq$ opposite $\bdr_i$, for the running example algebra $A$.} Here, each possible pair $(\bdp, \bdq)$ is given as a column in the table.
  }
\end{thesisfigure}

\nrex There are two vertices in $\OQ$ representing the source of a commutativity relation: one is $\source(\alpha)$ and one is $\source(\beta)$. (They are exchanged by $\dagger$.)

When $i = \source(\alpha)$, the possible pairs of $\bdp$ and $\bdq$ opposite $\bdr_i$ given in Figure \ref{fig:example-choices-opposite-one-pin-boundary-syllable}(\subref{subfig:running-example-opposite-pin-boundary-syllable-1}). When $i = \source(\beta) $, the possible pairs are as given in Figure \ref{fig:example-choices-opposite-one-pin-boundary-syllable}(\subref{subfig:running-example-opposite-pin-boundary-syllable-2}).

\centresubsec{Patches with two pin-boundary syllables in the top row}

\nlpass{Construction} Recall that a vertex $i\in \OQ_0$ has $b_i=0$ if and only if it represents the source of a commutativity relation.

For each such vertex $i$, we construct the patch with the data as shown in Figure \ref{fig:patch-with-two-pin-boundary syllables}. On the left side, the top entry is the pin-boundary syllable associated to $i$ while the bottom entry is the virtual syllable $\bdx_{i-a_i}$  at $i-a_i$. The right side is similar but for $\idag$.
\begin{thesisfigure}
  \centering 
  \begin{tikzpicture}[scale=\standardscale]
     \patchrep{(0,0)}{$\bdr_i$}{$\bdr_\idag$}{$\bdx_{i-a_i}$}{$\bdx_{\idag-a_\idag}$}
    \end{tikzpicture}
    
  \caption[
    Schematic of a patch with two pin-boundary syllables
  ]{
    \label{fig:patch-with-two-pin-boundary syllables} \textit{Schematic of a patch with two pin-boundary syllables.} One of these is created for each vertex in the set $\{i\in \OQ_0 \colon b_i=0\}$. Note that this set is closed with respect to the vertex exchange map $\dagger$, since $b_i=b_{\idag}$.
    
    Above, $\bdr_i \ceq ( \ilsyll{i}{a_i-1}{1} )$ the pin boundary syllable at $i$ (similarly $\bdr_\idag$) and $\bdx_{i-a_i} \ceq ( \ilsyll{ i-a_i }{ 0 }{ 0 } )$ is the virtual syllable at $i-a_i$ (and similarly $\bdx_{i-a_\idag}$.
  }

  \vspace{4ex}

  \begin{tikzpicture}[scale=0.8]
    \patchrep{(-3,0)}{$(\source(\alpha),2,1)$}{$(\source(\beta),2,1)$}{$(\source(\alpha),0,0)$}{$(\source(\beta),0,0)$}
    \patchrep{(3,0)}{$(\source(\beta),2,1)$}{$(\source(\alpha),2,1)$}{$(\source(\beta),0,0)$}{$(\source(\alpha),0,0)$}
  \end{tikzpicture}
  \caption[
    The patches with two-pin boundary syllables, for the running example algebra $A$.
  ]{
    \label{fig:example-of-patches-with-two-pin-boudnary-syllables}\textit{The patches with two-pin boundary syllables, for the running example algebra $A$.} These two patches are reflections of one another.
  }
  
  \vspace{4ex}
  
  \begin{tikzpicture}[scale=\standardscale]
    \patchrep{(-3,0)}{}{$\bdx_i$}{}{$\bde_i$}
    
    \patchrep{(3,0)}{$\bdx_i$}{}{$\bde_i$}{}
  \end{tikzpicture}
  \caption[%
    Schematic of virtual patches
  ]{%
    \label{fig:virtual-patches}\textit{Virtual patches.} These two patches are created for each $i$ with $d_i=0$. They are reflections of one another.
    
    Here, $\bdx_i \ceq ( \ilsyll{i}{0}{0} )$ is the virtual syllable and $\bde_i \ceq ( \ilsyll{i}{0}{1} )$ the stationary syllable at $i$. All unmarked cells are blank.
  }
\end{thesisfigure}

\npass Note that the reflection of the patch associated to $i$ is the patch associated to $i^\dag$. We know that $i$ is the source of a component iff $\idag$ is, so by varying $i$ around the family of vertices with $b_i=0$ we pick up its exchange partner $\idag$. This means we do not artificially need to add the reflections of these patches; we obtain them already.

\nrex The unique pin module for $A$ is $P_1=I_1$. The head vertex $1$ of this pin module is represented by two vertices in the overquiver $\OQ$, say $\source(\alpha)$ and $\source(\beta)$. We know that $a_{\source(\alpha)}=3$ and $a_{\source(\beta)}=3$ while $b_{\source(\alpha)}=b_{\source(\beta)}=0$. It so happens in this particular example that $ \source(\alpha) - a_{\source(\alpha)1} = \source(\alpha) $ and $ \source(\beta) - a_{\source(\beta)} = \source(\beta) $; this ``fixpoint'' behaviour is not to be expected in general. It follows that exactly two patches are created at this step. We show them in Figure \ref{fig:example-of-patches-with-two-pin-boudnary-syllables}.

\nexample Consider any Nakayama algebra. All of its indecomposable projectives are uniserial. This means it has no pin modules and therefore no patches of this type are created for it. The same is true more generally for any monomial SB algebra.

% \npass The purpose of these patches, and that of the virtual patches defined next, is technical and will likely only become clear later. They only arise when modelling the short exact sequence
% $$
% \begin{tikzcd}
% 0 \ar[r]
%   \& \soc P \ar[r]
%     \& P \ar[r]
%       \& P/\soc P \ar[r]
%         \& 0 
% \end{tikzcd}
% $$
% for $P$ pin; in words, this is precisely the syzygy short exact sequence of the socle-quotient of a pin module. As we will see, there is boundary behaviour at both the left and right hand sides with interfere with each other in an exceptional way.

\centresubsec{Virtual patches}

\nlpass{Construction} Recall that a vertex $i\in \OQ_0$ has $d_i=0$ if and only if it represents the target of a commutativity relation.

For each such vertex $i$, we construct the two patches shown in Figure \ref{fig:virtual-patches}. The first has a virtual syllable $(i,0,0)$ in the top-right cell and the corresponding stationary syllable $(i,0,1)$ underneath, but is otherwise blank. The second is its reflection.

\nsbstrips The virtual patches are those for which the property \mono{IsVirtualPatch} returns \true. Therefore, supposing the SB algebra in question is \mono{sba}, they are precisely the elements of the list returned by the following code. 
\begin{Verbatim}[commandchars=!@|,fontsize=\small,frame=single,baselinestretch=1.2]
  !gapprompt@gap>| !gapinput@Filtered( PatchSetOfSBAlg( sba ), IsVirtualPatch );|
\end{Verbatim}

\npass We emphasise that, as with the previous class of patches, virtual patches are only seen when calculating the syzygy of $P/\soc P$, for $P$ pin. 

\subsection{Properties of patches}\label{subsec:properties-of-patches}

In this subsection we make remarks and prove some properties about the patches of a SB algebra. There is no clever theoretical work to be seen here, only technical observations and their verification. These results are included in order to establish the validity of our model.

\begin{lemma}
The only patches containing virtual syllables are virtual patches (having one such, which lies in the top row) and patches with two pin-boundary syllables in the top row (having two such, in the bottom row). 
\end{lemma}

\proof This is true by inspection. We need only recall that, by definition, $\bdp \descent$ is never a virtual syllable for any $\bdp$. \qed

\begin{lemma}\label{lem:properties-of-patches} Let $X$ be a patch with entries $\bdp$ and $\bdp'$ in the top row, and with entry $\bdq$ under $\bdp$.
\begin{enumerate}
    \item\label{sublem:top-entries-of-patch-are-peak-compatible} If neither $\bdp$ nor $\bdp'$ is blank, then $\source(\bdp)^\dagger = \source( \bdp' )$.
    
    \item\label{sublem:columns-of-patches-are-usually-descent-images} If $\bdp$ is nonblank, then $\bdq = \bdp \descent$ unless
     \begin{enumerate}
         \item both top entries $\bdp,\bdp'$ of $X$ are pin boundary (in which case $\bdq$ is a virtual syllable), or
         
         \item the other top entry $\bdp'$ of $X$ is a pin boundary syllable (in which case $\bdq$ is the perturbed version of $\bdp\descent$), or
         
         \item $\bdp$ is a virtual syllable (in which case $\bdq$ is the corresponding stationary syllable).
     \end{enumerate}
    
    \item\label{sublem:entries-under-blank-entries} If $\bdp$ is blank, then $\bdq$ is either blank or is $\bde_i \descent $ for some stationary syllable $\bde_i$.
    
    \item\label{sublem:int-sylls-in-top-row-of-patches} If $\bdp = (\ilsyll{i}{\ell}{\ep})$ and $\ep=0$, then $\source(\bdq)=\target(\bdp)$.
\end{enumerate}
\end{lemma}

\proof The assumptions in part \ref{sublem:top-entries-of-patch-are-peak-compatible} exclude the blank patch, virtual patches and amended patches from consideration~-- in remaining cases \ref{sublem:top-entries-of-patch-are-peak-compatible} is seen to hold~-- but part \ref{sublem:entries-under-blank-entries} directs us exactly to those three cases. For the blank or virtual patches, any column featuring one blank entry is entirely blank as claimed in \ref{sublem:entries-under-blank-entries}. Amended patches only change the top syllable in a column comprising $\bde_i$ and $\bde_i\descent$. That top syllable $\bde_i$ is stationary, hence the claim of part \ref{sublem:entries-under-blank-entries} is valid in this case too.

Part \ref{sublem:columns-of-patches-are-usually-descent-images} holds by observation. We only remind the reader that $\descent$ maps a pin-boundary syllable to the blank syllable and the blank syllable to itself.

Towards part \ref{sublem:int-sylls-in-top-row-of-patches}, we can use part \ref{sublem:columns-of-patches-are-usually-descent-images}. Note that the assumption $\ep=0$ implies $\bdp$ is neither a pin-boundary syllable nor blank. If $\bdq=\bdp\descent$ then the result holds, since
$$
\source(\bdq) = \source(\bdp\descent) = \source(\ilsyll{i-\ell}{a_i+b_i-\ell}{b_i}) = i-l=\target(\bdp)\tstop
$$
If $\bdq$ is the perturbed version of $\bdp\descent$, then the result also holds because perturbation does not alter the source. If $\bdp$ is a virtual syllable $\big(\ilsyll{i}{0}{0}\big)$, then $\bdq$ is the corresponding stationary syllable $\bde_i = \big(\ilsyll{i}{0}{1}\big)$ and we clearly have
$ \source( \bde_i )=i=\target(\bdp) $.\qed

\begin{lemma}\label{lem:patches-have-different-tops}
No two patches have the same top row.
\end{lemma}

\proof We initially partition patches into two classes, depending on whether their top rows feature virtual syllables. Note that those who do are always virtual patches, each of which is uniquely distinguished by the single virtual syllable its contains and its orientation.

Amongst those whose top rows feature no virtual syllables, we further discriminate by the number of blank syllables in the top row. There is a single patch with two blank syllables: the blank patch. The patches with exactly one blank syllable are amended patches, each uniquely characterised by the nonblank syllable in the top row and its orientation.

Amongst those whose top rows feature no virtual syllables and no blank syllables, we can further subdivide by the number of pin-boundary syllables. These cases are treated in different parts of Subsection \ref{subsec:construction-of-patches}, from which the desired distinction follows immediately.\qed

\nlpass{Associating projectives to patches} \label{def:projective-of-a-patch} We define a function from the set of nonvirtual patches to the set containing the indecomposable projective $A$-modules with the prescribed bases as in Example \ref{example:basis-of-projectives} and the zero module (with basis $\emptyset$). Its purpose will become clear in Proposition \ref{prop:line-of-patches-is-projective-cover}. 

To the blank patch, we associate the zero module.

Otherwise, observe that any nonblank, nonvirtual patch $X$ has at least one nonblank syllable $\bdp$ in its top row. The source $\source(\bdp) $ represents a vertex $i$ of the ground quiver; if the top row features a second nonblank syllable $\bdq$ then $\source(\bdp)^\dagger=\source(\bdq)$ and so they represent the same vertex $i$. Therefore we associate to $X$ the projective $P_i$, for this well-chosen $i$.

\begin{thesisfigure}
\centering

\begin{tikzpicture}[scale=\standardscale]
  \patchrep{(0,0)}{$\bdp$}{$\bdp'$}{$\bdq$}{$\bdq'$}
\end{tikzpicture}

\caption[
  Notation for Definition \ref{def:projective-of-a-patch}
]{%
  \label{fig:notation-for-patches-repping-basis-of-projective}
  \textit{Notation for Definition \ref{def:projective-of-a-patch}.}
}
\end{thesisfigure}

\begin{lemma}\label{lem:bottom-row-of-pin-patch}
Suppose that the projective associated to the patch $X$ is a pin module.
\begin{enumerate}
    \item\label{sublem:bottom-row-of-socle-quotient-patch-is-virtual} If the top row of $X$ features two pin-boundary syllables, then the entries of the bottom row are both virtual syllables.
    
    \item\label{sublem:bottom-row-under-one-pin-boundary-is-boundary-valley} If the top row of $X$ features exactly one pin-boundary syllable, then one entry of the bottom row is a boundary syllable while the other is blank.
    
    \item\label{sublem:bottom-row-of-generic-pin-patch-is-interior-valley} If the top row of $X$ features no pin-boundary syllables, then the entries $\bdq$, $\bdq'$ of the bottom row satisfy $\bdq = \bdp \descent$ and $\bdq' = \bdp'\descent$ for some syllables $\bdp,\bdp'$ with $\source(\bdp)^\dagger = \source(\bdp')$. This implies $\target(\bdq)^\dagger=\target(\bdq')$.
\end{enumerate}
\end{lemma}

\proof For parts \ref{sublem:bottom-row-of-socle-quotient-patch-is-virtual} and \ref{sublem:bottom-row-under-one-pin-boundary-is-boundary-valley}, respectively refer to the construction of patches with two or one pin-boundary syllables in the top row. The claims hold by inspection.

For part \ref{sublem:bottom-row-of-generic-pin-patch-is-interior-valley}, let $\bdp$ and $\bdp'$ be the nonblank entries of the top row above $\bdq$ and $\bdq'$ respectively. If $X$ is an amended patch, instead let $\bdp,\bdp'$ be the corresponding nonblank entries of the top row of the original patch.

Because $X$ represents a pin module, we have $b_{\source(\bdp)}=b_{\source(\bdp')}=0$. Since $\bdp,\bdp'$ are nonblank syllables with this property, $\bdp\descent,\bdq\descent$ are defined and satisfy $\target(\bdp\descent)^\dagger = \target(\bdp'\descent)$. The construction of patches with no pin-boundary syllables in the top row implies $\bdq=\bdp\descent$ and $\bdq'=\bdp'\descent$. The result follows.\qed

\nrmk In addition to singling out a projective module $P$, a patch $X$ specifies two subsets of the standard basis: a top part $T$ and a bottom part $B$.

If the top row of $X$ features two pin-boundary syllables, then $B$ comprises only the vector spanning $\soc(P_i)$ and $T$ comprises the remaining vectors. Otherwise, suppose that the entries of $X$ are $\bdp,\bdp',\bdq,\bdq'$ as shown in Figure \ref{fig:notation-for-patches-repping-basis-of-projective}; if $X$ is an amended patch, let $\bdp,\bdp',\bdq,\bdq'$ be syllables of the original version. Write $p,p'$ for the compressions of $\bdp,\bdp'$ and $\leq$ for the prefix order. We define $T$ and $B$ via the following procedure.
\begin{enumerate}
    \item Define $T_p \ceq \{u \in \Subpath(\OQ) \colon u < p \}$, and define $T_{p'}$ similarly.
        
    The definition of $B_p$ depends on whether $\bdp$ is boundary or interior: if $\bdp$ is boundary then we define $B \ceq \{u \in \Subpath(\OQ) \colon p \leq u\}$; while if $\bdp$ is interior we define $ \{u \in \Subpath(\OQ) \colon p < u\}\subseteq B$ (in which case $p \notin B \cup T$). We define $B_{p'}$ similarly, again resorting to cases depending on whether $\bdp'$ is boundary or interior.
    
    \item Replace any path $u$ in $T_p$, $T_{p'}$, $B_p$ or $B_{p'}$ by the basis vector $u A$ represented by $u$.
    
    \item Set $T \ceq T_p \cup T_{p'}$ and $B \ceq B_p \cup B_{p'}$.
\end{enumerate}

We mention that if $\bdp$ is interior then (the basis vector represented by) $p$ does not belong to $T \cup B$, and that all basis elements not in $T \cup B$ are of this form.
\begin{thesisfigure}
  \centering
  
  \subcaptionbox{
    \label{subfig:patch-for-T-and-B}
    An example patch
  }[
    0.45\textwidth
  ]{
  \begin{tikzpicture}[scale=\standardscale]
    \begin{scope}[node distance=7]
      \patchrep{(0,0)}{}{}{}{};    
      \oddintsyll{nw}{}{$\alpha$}{};
      \evenbdysyll{ne}{}{$\beta$}{}{}{};
      \evenintsyll{sw}{}{$\beta\gamma$}{};
      \oddintsyll{se}{}{$\alpha$}{};
    \end{scope}
  \end{tikzpicture}
  }
  \subcaptionbox{
  \label{subfig:example-of-T-and-B}
  The associated sets $T$ (top) and $B$ (bottom)
  }[
    0.45\textwidth
  ]{
  \begin{tikzpicture}[scale=0.7*\standardscale]
    % Give "anchoring" nodes for nonzero summands
      \draw (0,0) node[](1){};
      
      % Draw nodes of lefthand projective
      \draw (1) ++(2,6) node[](1e1){$e_1$};
      \draw (1) ++(4,4) node[](1b){$\beta$};
      \draw (1) ++(4,2) node[](1bc){$\beta\gamma$};
      \draw (1) ++(2,0) node[](1bca){$\alpha\beta\gamma$};
      \draw (1) ++(0,4) node[](1a){$\alpha$};
      \draw (1) ++(0,2) node[](1ab){$\alpha\beta$};
      
      % Draw action lines on lefthand projective
      \draw (1e1) --node[pos=0.5,right](){$\beta$} (1b);
      \draw (1b) --node[pos=0.5,right](){$\gamma$} (1bc);
      \draw (1bc) --node[pos=0.5,right](){$\alpha$} (1bca);
      \draw (1e1) --node[pos=0.5,right](){$\alpha$} (1a);
      \draw (1a) --node[pos=0.5,right](){$\beta$} (1ab);
      \draw (1ab) --node[pos=0.5,right](){$\gamma$} (1bca);
      
      \begin{scope}[inner sep=0, minimum size=0]
        \newcommand{\tol}{0.5}
        \draw (1e1)+(\tol,\tol) node[](t1){};
        \draw (1e1)+(-\tol,\tol) node[](t2){};
        \draw (1e1)+(-\tol,-\tol) node[](t3){};
        \draw (1b)+(-\tol,-\tol) node[](t4){};
        \draw (1b)+(\tol,-\tol) node[](t5){};
        \draw (1b)+(\tol,\tol) node[](t6){};
        \draw[dotted] (t1) to[out=135,in=45] (t2) to[out=225,in=135] (t3) to[out=315,in=135] (t4) to[out=315,in=225] (t5) to[out=45,in=315] (t6) -- (t1);
        
        \draw (1ab)+(\tol,\tol) node(b0){};
        \draw (1ab)+(-\tol,\tol) node[](b1){};
        \draw (1ab)+(-\tol,-\tol) node[](b2){};
        \draw (1bca)+(-\tol,-\tol) node[](b3){};
        \draw (1bca)+(\tol,-\tol) node[](b4){};
        \draw (1bc)+(\tol,-\tol) node[](b5){};
        \draw (1bc)+(\tol,\tol) node[](b6){};
        \draw (1bc)+(-\tol,\tol) node[](b7){};
        \draw[dotted] (b0) to[out=135, in=45] (b1) to[out=225,in=135] (b2) to[out=315,in=135] (b3) to[out=335,in=225] (b4) to[out=45,in=225] (b5) to[out=45, in=315] (b6) to[out=135, in=45] (b7) to[out=225, in=315] (b0);
      \end{scope}
  \end{tikzpicture}
  }
  
  \caption[
    Basis subsets associated to a patch.
  ]{
    \label{fig:basis-subsets-associated-to-a-patch}
    \textit{Basis subsets associated to a patch}
  }
\end{thesisfigure}

\nrex In Figure \ref{fig:basis-subsets-associated-to-a-patch} we show a patch and the corresponding subsets of an associated patch. Notice that the top row contains the interior syllable $\big(\ilsyll{}{\alpha}{0}\big)$ and that, accordingly, the associated basis vector $\alpha A$ belongs to neither $T$ nor $B$.

\section{Strips}

In this section, we define strips. These are representations of string modules using the alphabet of syllables.

\setcounter{subsection}{-1}
\subsection{Preview of strips}
\label{subsec:preview-of-strips}

Let us explore our approach to string modules and their syzygies using worked examples.  We will start to see how to represent string graphs with no ambiguity, if a modicum of redundancy. The focus here is building intuition; formal definitions will come in due course. Recall the running example algebra
$$
A \ceq \kk\big(
\begin{tikzcd}
  1 \ar[loop left, "\alpha"] \ar[r, "\beta", shift left]
    \& 2 \ar[l, "\gamma", shift left] \ar[loop right, "\delta"]
\end{tikzcd}
\big)
/
\big\langle
  \alpha^2,
  \beta\delta,
  \gamma\beta,
  \delta\gamma,
  \alpha\beta\gamma-\beta\gamma\alpha,
  \gamma\alpha\beta,
  \delta^4
\big\rangle
$$
(whose uniserial $A$-modules can be identified with the strict prefixes of $\alpha\beta\gamma$, $\beta\gamma\alpha$, $\gamma\alpha\beta$ and $\delta^4$).

\nlpass{The big idea} We have a way of representing a string graph by entering symbols (\emph{syllables}) into a row of cells. These cells are alternately of two types, \emph{positive}
\iltikzpicture{
\evencell{(0,0)}{pos}{};
}
and \emph{negative}
\iltikzpicture{
\oddcell{(0,0)}{neg}{};
}, which means that pairs of cells alternately form \emph{peaks}
\iltikzpicture{
\peakrep{(0,0)}{pk1}{}{}
}
and \emph{valleys}
\iltikzpicture{
\valleyrep{(0,0)}{vl1}{}{}
}. Such a row of symbols is a \emph{strip}.

The contents of a peak in one row uniquely determine the contents of a valley in a row underneath it. One row may in general have several rows underneath it~-- a phenomenon we call \emph{branching}~-- but there is a fixed procedure for determining to which row the contents of a given valley are assigned. We thus create strips from a given strip by reading the latter a peak at a time.

As we demonstrate below, it is necessary that we operate ``peakwise'' (a peak upstairs specifies a valley downstairs) rather than ``cellwise'' (a positive or negative cell upstairs specifies a negative or positive cell downstairs) due to the existence of string modules related to socle-quotients of pin modules. These troublesome string modules are easy to spot; they are precisely those whose strip representation includes a certain kind of syllable, dubbed a \emph{pin-boundary} syllable. When a peak contains no pin-boundary syllables, the entry of each cell in the valley underneath it depends only on the individual cell directly above it.

There are several key features in our approach. We treat them all in their turn with examples.

\nlpass{Interiors, boundaries, peaks and valleys} Consider the string graph $v_1$, shown in Figure \ref{subfig:string-graph-v1}. It has four subgraphs that start at source vertices and end at sinks. The corresponding restrictions are
$\big(
  \begin{tikzpicture}[scale=0.75]
    \begin{scope}[minimum size=0, inner sep=1]
      \draw (0,0) node[](v1){$2$};
      \draw (1,0) node[](v2){$1$};
      \draw[->] (v2) --node[pos=0.5, above](){$\beta$} (v1);
    \end{scope}
  \end{tikzpicture}
\big)$,
$\big(
  \begin{tikzpicture}[scale=0.75]
    \begin{scope}[minimum size=0, inner sep=1]
      \draw (0,0) node[](v1){$1$};
      \draw (1,0) node[](v2){$1$};
      \draw[->] (v1) --node[pos=0.5, above](){$\alpha$} (v2);
    \end{scope}
  \end{tikzpicture}
\big)$,
$\big(
  \begin{tikzpicture}[scale=0.75]
    \begin{scope}[minimum size=0, inner sep=1]
      \draw (0,0) node[](v1){$1$};
      \draw (1,0) node[](v2){$2$};
      \draw[->] (v2) --node[pos=0.5, above](){$\gamma$} (v1);
    \end{scope}
  \end{tikzpicture}
\big)$
and
$\big(
  \begin{tikzpicture}[scale=0.75]
    \begin{scope}[minimum size=0, inner sep=1]
      \draw (0,0) node[](v1){$2$};
      \draw (1,0) node[](v2){$2$};
      \draw (2,0) node[](v3){$2$};
      \draw[->] (v1) --node[pos=0.5, above](){$\delta$} (v2);
      \draw[->] (v2) --node[pos=0.5, above](){$\delta$} (v3);
    \end{scope}
  \end{tikzpicture}
\big)$.
The first and last of these appear at the boundary of the string graph. Informally, this means that these subgraphs point outwards from the graph; more formally, that the sink vertices at which these end have indegree only $1$ (a property distinguishing these sink vertices from that at which the other two subgraphs end.) We can represent $v_1$ using the strip in Figure \ref{subfig:stripifying-v1}.
\begin{thesisfigure}
  \centering
  \subcaptionbox{%
      \label{subfig:string-graph-v1}
      \textit{String graph $v_1$}
    }[%
      0.95\textwidth
    ]{%
      \begin{tikzpicture}[scale=\standardscale]
        \newcommand{\firstexamplestring}[1]{%
          \begin{scope}[inner sep=1, minimum size=0]
            % Vertices of string graph
            \draw #1 ++(0,0) node[](1 1){$2$};
            \draw #1 ++(1,1) node[](1 2){$1$};
            \draw #1 ++(2,0) node[](1 3){$1$};
            \draw #1 ++(3,1) node[](1 4){$2$};
            \draw #1 ++(4,0) node[](1 5){$2$};
            \draw #1 ++(5,-1) node[](1 6){$2$}; 
            % Arrows of string graph
            \draw[->] (1 2) -- node[pos=0.5, below right](){$\beta$}  (1 1);
            \draw[->] (1 2) -- node[pos=0.5, above right](){$\alpha$} (1 3);
            \draw[->] (1 4) -- node[pos=0.5, below right](){ $\gamma$ } (1 3);
            \draw[->] (1 4) -- node[pos=0.5, above right](){ $\delta$ } (1 5);
            \draw[->] (1 5) -- node[pos=0.5, above right](){ $\delta$ } (1 6);
            
            % Name of string graph
            % \draw (1 3) ++(0.5,-3) node[](){$v_1$};
          \end{scope}
          }          
        \firstexamplestring{(0,0)};
      \end{tikzpicture}
    }
    
 \subcaptionbox{
    \label{subfig:stripifying-v1}
    \textit{Strip representing $v_1$.} The reflection of this strip also represents $v_1$.
  }[
    \textwidth
 ]{
    \begin{tikzpicture}[scale=0.7]
      % Cells
      \evenleftedgecell{(0,0)}{c0}{};
      \oddcell{(2,0)}{c1}{};
      \evencell{(4,0)}{c2}{};
      \oddcell{(6,0)}{c3}{};
      \evencell{(8,0)}{c4}{};
      \oddcell{(10,0)}{c5}{};
      \evencell{(12,0)}{c6}{};
      \oddcell{(14,0)}{c7}{};
      \evencell{(16,0)}{c8}{};
      \oddrightedgecell{(18,0)}{c9}{};

    % Entries
      \begin{scope}[inner sep=0, minimum size=0, node distance=7]
        \oddbdysyll{c3}{}{$\beta$}{}{}{};
        \evenintsyll{c4}{}{$\alpha$}{};
        \oddintsyll{c5}{}{$\gamma$}{};
        \evenbdysyll{c6}{}{$\delta^2$}{}{}{};
      \end{scope}
    \end{tikzpicture}
  }
  
  \caption[
    First example of representing a string graph with strips
  ]{
    \label{fig:v1-as-string-graph-and-strip}
    \textit{First example of representing a string graph with strips.}
  }
\end{thesisfigure}

Let us dissect the anatomy of this figure. To begin, note that it comprises an infinite ($\integers$-indexed) row of cells which alternately look like
\iltikzpicture{
  \oddcell{(0,0)}{oddcell}{};
}
and
\iltikzpicture{
  \evencell{(0,0)}{evencell}{};
}.
Therefore they alternately form peaks
\iltikzpicture{
  \peakrep{(0,0)}{pk1}{}{};
}
and valleys
\iltikzpicture{
  \valleyrep{(0,0)}{vl1}{}{};
}.

Next, notice that cells sufficiently far to the left and to the right are blank. Between these extremes is the \emph{support} of the strip: the concave interval of cells with nonblank entries. 

These entries are syllables that unambiguously represents the subgraph of the string graph. Recall that a syllable $\big( \ilsyll{\circ}{p}{\ep} \big)^{\pm1}$ contains three pieces of information: a path $p$ in the overquiver $\OQ$, a bit $\ep\in\{0,1\}$ which encodes whether the syllable is an interior ($\ep=0$) or boundary syllable ($\ep=1$) and an orientation which is either positive or negative ($+1$ or $-1$). Of these, orientation is the least interesting. All that really matters is that it alternates across the strip. In our illustrations, orientation of a (nonblank) syllable is a function of the kind of cell where it appears: entries in
\iltikzpicture{
  \oddcell{(0,0)}{oddcell}{};
}
cells have negative orientation while those in
\iltikzpicture{
  \evencell{(0,0)}{evencell}{};
}
cells have positive orientation. Orientation is also demonstrated by the diagonal along which the syllable is drawn.

%%% RETHINK
The underlying path of a syllable unambiguously represents a path in $A$, but we remind the reader of a technical point: underlying paths of a syllable are paths in the overquiver $\OQ$, not the ground quiver $Q$. In the interests of brevity this is not usually expressed in the notation. For instance, the rightmost syllable in Figure \ref{subfig:stripifying-v1} has underlying path $\delta^2$ (in the overquiver), but it represents the $A$-path that is also called $\delta^2$.

Boundary syllables, such as $\big( \ilsyll{\circ}{\beta}{1} \big) $ and $\big(\ilsyll{\circ}{\delta^2}{1}\big)$ in our example, only appear pointing outwards at the boundary of the strip. This means that if one entry in a valley is a boundary syllable then the other is blank, as may be verified in the example. 
% When a boundary syllable appears in a peak, its source is the $\dagger$-partner of the source of the syllable in the other cell of the peak. We see this in Figure \ref{subfig:stripifying-v1}. The peak on the left features the boundary syllable $\big( \ilsyll{ \circ }{ \beta }{ 1 } \big)$ next to $\big( \ilsyll{ \circ }{ \alpha }{ 0 } \big)$. Their sources are $\dagger$-partners, which is essentially to say they are distinct vertices of the overquiver $\OQ$ that represent the same vertex of the ground quiver.

In between the boundary syllables are the interior syllables, $\big( \ilsyll{\circ }{\alpha}{0} \big) $ and $\big(\ilsyll{\circ}{\gamma}{0}\big)$ here. Observe that if one syllable in a valley is an interior syllable then so too is the other and, moreover, their targets are $\dagger$-partners. The only example in Figure \ref{subfig:stripifying-v1} is the middle valley that contains $\big( \ilsyll{ \circ }{ \alpha }{ 0 } \big)$ and $\big( \ilsyll{ \circ }{\gamma }{ 0 } \big)$. More examples will appear when discussing string graph $v_2$, below.

Whether the syllables in a peak are boundary or interior, their sources are $\dagger$-partners. 

In our illustrations, boundary syllables have their stability term coloured, while interior syllables have theirs suppressed. This helps distinguish them visually.

\nlpass{Rounding-off} Now consider the string graph $v_2$ in Figure \ref{subfig:string-graph-v2}.%
\begin{thesisfigure}
  \centering

  \subcaptionbox{%
    \label{subfig:string-graph-v2}
    \textit{String graph $v_2$}
  }[%
    0.95\textwidth
  ]{%
    \begin{tikzpicture}[scale=\standardscale]
      \newcommand{\secondexamplestring}[1]{%
        \begin{scope}[inner sep=1, minimum size=0]
          % Vertices of string
          \draw #1 ++(0,0) node[](2 1){$1$};
          \draw #1 ++(1,-1) node[](2 2){$1$};
          \draw #1 ++(2,-2) node[](2 3){$2$};
          \draw #1 ++(3,-1) node[](2 4){$2$};
          \draw #1 ++(4,-2) node[](2 5){$1$};
          \draw #1 ++(5,-1) node[](2 6){$1$};
          \draw #1 ++(6,0) node[](2 7){$2$};
          % Arrows of string
          \draw[->] (2 1) -- node[pos=0.5,above right](){$\alpha$} (2 2);
          \draw[->] (2 2) -- node[pos=0.5,above right](){$\beta$} (2 3);
          \draw[->] (2 4) -- node[pos=0.5,below right](){$\delta$} (2 3);
          \draw[->] (2 4) -- node[pos=0.5,above right](){$\gamma$} (2 5);
          \draw[->] (2 6) -- node[pos=0.5,below right](){$\alpha$} (2 5);
          \draw[->] (2 7) -- node[pos=0.5,below right](){$\gamma$} (2 6);
          
          % Name of string graph
        %   \draw (2 4) ++(0,-3) node[](){$v_2$};
        \end{scope}
        }
      \secondexamplestring{(-2,1)};
    \end{tikzpicture}
  }
  
  \vspace{4ex}
  
  \subcaptionbox{
    \label{subfig:stripifying-v2}
    \textit{Strips representing $v_2$.} Here, $e$ and $e'$ are the stationary paths at $\source(\alpha)^\dagger$ and $\source(\gamma)^\dagger$, respectively. Since the stationary syllables at either end can be inferred from a neighbouring cell, they are both optional.
  }[
    0.95\textwidth
  ]{
    \begin{tikzpicture}[scale=0.7]
    % Cells of 1st row (top row)
    \oddleftedgecell{(0,2.5)}{c1 0}{};
    \evencell{(2,2.5)}{c1 1}{};
    \oddcell{(4,2.5)}{c1 2}{};
    \evencell{(6,2.5)}{c1 3}{};
    \oddcell{(8,2.5)}{c1 4}{};
    \evencell{(10,2.5)}{c1 5}{};
    \oddcell{(12,2.5)}{c1 6}{};
    \evencell{(14,2.5)}{c1 7}{};
    \oddcell{(16,2.5)}{c1 8}{};
    \evenrightedgecell{(18,2.5)}{c1 9}{};
    
    % Entries of 1st row (top row)
    \begin{scope}[inner sep=0, minimum size=0, node distance=7]
      \oddtrivbdysyll{c1 2}{}{$e$}{}{}{};
      \evenintsyll{c1 3}{}{$\alpha\beta$}{};
      \oddintsyll{c1 4}{}{$\delta$}{};
      \evenintsyll{c1 5}{}{$\gamma$}{};
      \oddintsyll{c1 6}{}{$\gamma\alpha$}{};
      \eventrivbdysyll{c1 7}{}{$e'$}{}{}{};
    \end{scope}
    
    % Cells of 2nd row
    \oddleftedgecell{(0,0)}{c1 0}{};
    \evencell{(2,0)}{c2 1}{};
    \oddcell{(4,0)}{c2 2}{};
    \evencell{(6,0)}{c2 3}{};
    \oddcell{(8,0)}{c2 4}{};
    \evencell{(10,0)}{c2 5}{};
    \oddcell{(12,0)}{c2 6}{};
    \evencell{(14,0)}{c2 7}{};
    \oddcell{(16,0)}{c2 8}{};
    \evenrightedgecell{(18,0)}{c2 9}{};
    
    % Entries of 2nd row
    \begin{scope}[inner sep=0, minimum size=0, node distance=7]
      \evenintsyll{c2 3}{}{$\alpha\beta$}{};
      \oddintsyll{c2 4}{}{$\delta$}{};
      \evenintsyll{c2 5}{}{$\gamma$}{};
      \oddintsyll{c2 6}{}{$\gamma\alpha$}{};
      \eventrivbdysyll{c2 7}{}{$e'$}{}{}{};
    \end{scope}
    
    % Cells of 3rd row
    \oddleftedgecell{(0,-2.5)}{c3 0}{};
    \evencell{(2,-2.5)}{c3 1}{};
    \oddcell{(4,-2.5)}{c3 2}{};
    \evencell{(6,-2.5)}{c3 3}{};
    \oddcell{(8,-2.5)}{c3 4}{};
    \evencell{(10,-2.5)}{c3 5}{};
    \oddcell{(12,-2.5)}{c3 6}{};
    \evencell{(14,-2.5)}{c3 7}{};
    \oddcell{(16,-2.5)}{c3 8}{};
    \evenrightedgecell{(18,-2.5)}{c3 9}{};
    
    % Entries of 3rd row
    \begin{scope}[inner sep=0, minimum size=0, node distance=7]
      \oddtrivbdysyll{c3 2}{}{$e$}{}{}{};
      \evenintsyll{c3 3}{}{$\alpha\beta$}{};
      \oddintsyll{c3 4}{}{$\delta$}{};
      \evenintsyll{c3 5}{}{$\gamma$}{};
      \oddintsyll{c3 6}{}{$\gamma\alpha$}{};
    \end{scope}
    
    % Cells of 4th row (bottom row)
    \oddleftedgecell{(0,-5)}{c4 0}{};
    \evencell{(2,-5)}{c4 1}{};
    \oddcell{(4,-5)}{c4 2}{};
    \evencell{(6,-5)}{c4 3}{};
    \oddcell{(8,-5)}{c4 4}{};
    \evencell{(10,-5)}{c4 5}{};
    \oddcell{(12,-5)}{c4 6}{};
    \evencell{(14,-5)}{c4 7}{};
    \oddcell{(16,-5)}{c4 8}{};
    \evenrightedgecell{(18,-5)}{c4 9}{};
    
    % Entries of 4th row (bottom row)
    \begin{scope}[inner sep=0, minimum size=0, node distance=7]
      \evenintsyll{c4 3}{}{$\alpha\beta$}{};
      \oddintsyll{c4 4}{}{$\delta$}{};
      \evenintsyll{c4 5}{}{$\gamma$}{};
      \oddintsyll{c4 6}{}{$\gamma\alpha$}{};
    \end{scope}
  \end{tikzpicture}
  }
  
  \caption[%
    Second example of representing a string graph with strips
  ]{%
    \label{fig:v2-as-string-graph-and-strip}
    \textit{Second example of representing a string graph with strips.}
  }
\end{thesisfigure}
What matters most about this example is that we can treat $v_2$ on an equal footing with $v_1$ despite their seeming difference in shape.

% while $v_2$ has a different shape to $v_1$ we have a robust manner of treating them as though they are comparably-shaped.

Informally speaking, $v_2$ resembles a ``generalised W''
\raisebox{1.2ex}{\iltikzpicture{%
  \begin{scope}[minimum size=0, inner sep=0]
    \draw (0,1) -- ++(1,-1) -- ++(1,1) -- ++(1,-1)  -- ++(1,1) node[](a){};
    \draw[dotted] (a)++(0,-0.5) -- ++(1,0);
    \draw (a)++(1,0) node[](b){};
    \draw (b)++(0,0) -- ++(1,-1) -- ++(1,1);
  \end{scope}
}} while $v_1$ resembles a ``generalised M''
\raisebox{1.2ex}{\iltikzpicture{%
  \begin{scope}[minimum size=0, inner sep=0]
    \draw (0,1) -- ++(1,1) -- ++(1,-1) -- ++(1,1) -- ++(1,-1) node[](a){};
    \draw[dotted] (a)++(0,0.5) -- ++(1,0);
    \draw (a)++(1,0) node[](b){};
    \draw (b)++(0,0) -- ++(1,1) -- ++(1,-1);
  \end{scope}
}}. It transpires that generalised Ms lend themselves more readily to syzygy calculations. There is an obvious trick to viewing a generalised W as a generalised M by making some trivial adjustment at the boundary:
\raisebox{1.2ex}{\iltikzpicture{%
  \begin{scope}[minimum size=0, inner sep=0]
    \draw (0,1)node[](0){} -- ++(1,-1) -- ++(1,1) -- ++(1,-1)  -- ++(1,1) node[](a){};
    \draw[dotted] (a)++(0,-0.5) -- ++(1,0);
    \draw (a)++(1,0) node[](b){};
    \draw (b)++(0,0) -- ++(1,-1) -- ++(1,1) node[](z){};
    \draw[densely dashed] (0) -- ++(-1,-1);
    \draw[densely dashed] (z) -- ++(1,-1);
  \end{scope}
}}. Our version of this trick is to round off the boundary of the strip with stationary syllables $\big( \ilsyll{ \circ }{e}{1} \big)$. 

(Similar tricks work for ``generalised Ns'' and ``generalised \reflectbox{N}s''. Also recall that by construction all stationary syllables used to represent string graphs are boundary. There are never stationary syllables in the interior of a strip.)

% When we represent $v_2$ as a strip, we use stationary syllables $\big( \ilsyll{ \circ }{e}{1} \big)$ to round off the data of the strip. (These must clearly be boundary syllables, hence the stability term $\ep=1$.)

See the uppermost strip in Figure \ref{subfig:stripifying-v2}. At each side, there is a unique stationary syllable whose source is the $\dagger$-partner of the source of the neighbouring syllable. That is the stationary syllable we use to round off the strip. To emphasise that the underlying path is trivial, we draw it as a dashed line in the illustration.

Since the appropriate stationary syllables at either end can be inferred from the remaining content of the strip, we will allow variations of the strip where one or both stationary syllables are absent. These are the other strips in Figure \ref{subfig:stripifying-v2}. This redundancy is a benefit, not a drawback: the framework we develop later will treat all of these variations the same.

In all cases, our comments from the previous passage about entries in peaks and valleys remain true. Boundary syllables, whenever they appear in a valley, do so next to a blank cell. The valley neighbour of an interior syllable is another interior syllable, and their targets are exchanged by $\dagger$. The sources of two (nonblank) syllables in a peak form a $\dagger$-pair, as do the targets of two interior syllables in a valley. However, whenever a stationary boundary syllable is only implied, interior syllables may form a peak with a blank cell. This final situation did not arise when treating $v_1$.

\nlpass{Simple modules} Examine the string graph $v_3$ in Figure \ref{subfig:string-graph-v3} that represents the simple module $S_1$.
\begin{thesisfigure}
\subcaptionbox{%
    \label{subfig:string-graph-v3}
    \textit{String graph $v_3$}
  }[%
    0.95\textwidth
  ]{%
    \begin{tikzpicture}[scale=\standardscale]
      \newcommand{\thirdexamplestring}[1]{%
        % Vertices of string graph
        \draw #1 ++(0,0) node[](3 1){$1$};
        
        % Name of string graph
        % \draw (3 1) ++(0,-3) node[](){$v_3$};
        }
        
      \thirdexamplestring{(7,0)};
    \end{tikzpicture}
  }

  \vspace{4ex}

  \subcaptionbox{
    \label{subfig:stripifying-v3}
    \textit{Strips representing $v_3$.} Here, $e$ and $e'$ are the two stationary paths at $\OQ$-vertices that represent the $Q$-vertex $1$. We can infer one of the stationary syllables from the other, but we cannot infer either of them out of nothing.
  }[
    0.95\textwidth
  ]{
    \begin{tikzpicture}[scale=0.7]
    % Cells of the 1st row (top row)
    \oddleftedgecell{(0,2.5)}{c0}{};
    \evencell{(2,2.5)}{c1}{};
    \oddcell{(4,2.5)}{c2}{};
    \evencell{(6,2.5)}{c3}{};
    \oddcell{(8,2.5)}{c4}{};
    \evencell{(10,2.5)}{c5}{};
    \oddcell{(12,2.5)}{c6}{};
    \evencell{(14,2.5)}{c7}{};
    \oddcell{(16,2.5)}{c8}{};
    \evenrightedgecell{(18,2.5)}{c9}{};
    
    % Entries of the 1st row (top row)
    \begin{scope}[inner sep=0, minimum size=0, node distance=7]
      \oddtrivbdysyll{c4}{}{$e$}{}{}{};
      \eventrivbdysyll{c5}{}{$e'$}{}{}{};
    \end{scope}
    
    % Cells of the 2nd row
    \oddleftedgecell{(0,0)}{c0}{};
    \evencell{(2,0)}{c1}{};
    \oddcell{(4,0)}{c2}{};
    \evencell{(6,0)}{c3}{};
    \oddcell{(8,0)}{c4}{};
    \evencell{(10,0)}{c5}{};
    \oddcell{(12,0)}{c6}{};
    \evencell{(14,0)}{c7}{};
    \oddcell{(16,0)}{c8}{};
    \evenrightedgecell{(18,0)}{c9}{};
    
    % Entries of the 2nd row
    \begin{scope}[inner sep=0, minimum size=0, node distance=7]
      \eventrivbdysyll{c5}{}{$e'$}{}{}{};
    \end{scope}
    
    % Cells of the 3rd row (bottom row)
    \oddleftedgecell{(0,-2.5)}{c0}{};
    \evencell{(2,-2.5)}{c1}{};
    \oddcell{(4,-2.5)}{c2}{};
    \evencell{(6,-2.5)}{c3}{};
    \oddcell{(8,-2.5)}{c4}{};
    \evencell{(10,-2.5)}{c5}{};
    \oddcell{(12,-2.5)}{c6}{};
    \evencell{(14,-2.5)}{c7}{};
    \oddcell{(16,-2.5)}{c8}{};
    \evenrightedgecell{(18,-2.5)}{c9}{};
    
    % Entries of the 3rd row (bottom row)
    \begin{scope}[inner sep=0, minimum size=0, node distance=7]
      \oddtrivbdysyll{c4}{}{$e$}{}{}{};
    \end{scope}
  \end{tikzpicture}
  }
  
  \caption[%
    Third example of representing a string graph with strips.
  ]{%
    \label{fig:v3-as-string-graph-and-strip}
    \textit{Third example of representing a string graph with strips.}
  }
  
  \vspace{4ex}
  
  \subcaptionbox{%
    \label{subfig:string-graph-v4}
    \textit{String graph $v_4$}
  }[%
    0.45\textwidth
  ]{%
    \begin{tikzpicture}[scale=\standardscale]
      \newcommand{\fourthexamplestring}[1]{%
        \begin{scope}[inner sep=1, minimum size=0]
          % Vertices of string graph
          \draw #1 ++(0,0) node[](4 1){$2$};
          \draw #1 ++(1,1) node[](4 2){$1$};
          \draw #1 ++(2,2) node[](4 3){$1$};
          \draw #1 ++(3,1) node[](4 4){$2$};
          \draw #1 ++(4,0) node[](4 5){$1$};
          
          % Arrows of string graph
          \draw[->] (4 2) -- node[pos=0.5, below right](){$\beta$} (4 1);
          \draw[->] (4 3) -- node[pos=0.5, below right](){$\alpha$} (4 2);
          \draw[->] (4 3) -- node[pos=0.5, above right](){$\beta$} (4 4);
          \draw[->] (4 4) -- node[pos=0.5, above right](){$\gamma$} (4 5);
          
        %   \draw (4 3) ++(0,-4) node[](){$v_4$};
        \end{scope}
      }
      \fourthexamplestring{(0,0)};
    \end{tikzpicture}
  }
  \subcaptionbox{
    \label{subfig:string-graph-v5}
    \textit{String graph $v_5$.}
  }[
    0.45\textwidth
  ]{
    \begin{tikzpicture}[scale=\standardscale]
      \newcommand{\fifthexamplestring}[1]{%
        \begin{scope}[inner sep=1, minimum size=0]
          % Vertices of string graph
          \draw #1 ++(0,0) node[](5 1){$2$};
          \draw #1 ++(1,1) node[](5 2){$1$};
          \draw #1 ++(2,2) node[](5 3){$1$};
          \draw #1 ++(3,1) node[](5 4){$2$};
          \draw #1 ++(4,0) node[](5 5){$1$};
          \draw #1 ++(5,1) node[](5 6){$1$};
          \draw #1 ++(6,0) node[](5 7){$2$};
          \draw #1 ++(7,-1) node[](5 8){$1$};
          
          % Arrows of string graph
          \draw[->] (5 2) -- node[pos=0.5, below right](){$\beta$} (5 1);
          \draw[->] (5 3) -- node[pos=0.5, below right](){$\alpha$} (5 2);
          \draw[->] (5 3) -- node[pos=0.5, above right](){$\beta$} (5 4);
          \draw[->] (5 4) -- node[pos=0.5, above right](){$\gamma$} (5 5);
          \draw[->] (5 6) -- node[pos=0.5, below right](){$\alpha$} (5 5);
          \draw[->] (5 6) -- node[pos=0.5, above right](){$\beta$} (5 7);
          \draw[->] (5 7) -- node[pos=0.5, above right](){$\gamma$} (5 8);
          
        %   \draw (4 3) ++(0,-4) node[](){$v_4$};
        \end{scope}
      }
      \fifthexamplestring{(0,0)};
    \end{tikzpicture}
  }
  
  \vspace{4ex}
  
  \subcaptionbox{
    \label{subfig:stripifying-v4}
    \textit{Strip representing $v_4$.}
  }[
    0.95\textwidth
  ]{
    \begin{tikzpicture}[scale=0.7]
    % Cells
    \oddleftedgecell{(0,0)}{c0}{};
    \evencell{(2,0)}{c1}{};
    \oddcell{(4,0)}{c2}{};
    \evencell{(6,0)}{c3}{};
    \oddcell{(8,0)}{c4}{};
    \evencell{(10,0)}{c5}{};
    \oddcell{(12,0)}{c6}{};
    \evencell{(14,0)}{c7}{};
    \oddcell{(16,0)}{c8}{};
    \evenrightedgecell{(18,0)}{c9}{};
    
    % Entries
    \begin{scope}[inner sep=0, minimum size=0, node distance=7]
      \oddbdysyll{c4}{}{$\alpha\beta$}{}{}{};
      \evenbdysyll{c5}{}{$\beta\gamma$}{}{}{};
    \end{scope}
  \end{tikzpicture}
  }
  
  \vspace{4ex}
  
  \subcaptionbox{
    \label{subfig:stripifying-v5}
    \textit{Strips representing $v_5$.} The boundary syllables in this strip also occur in the strip representing a socle-quotient of a pin module~-- compare (\subref{subfig:stripifying-v4}) above~-- but here they do not appear together in a single peak.
  }[
    0.95\textwidth
  ]{
    \begin{tikzpicture}[scale=0.7]
    % Cells
    \evenleftedgecell{(0,0)}{c0}{};
    \oddcell{(2,0)}{c1}{};
    \evencell{(4,0)}{c2}{};
    \oddcell{(6,0)}{c3}{};
    \evencell{(8,0)}{c4}{};
    \oddcell{(10,0)}{c5}{};
    \evencell{(12,0)}{c6}{};
    \oddcell{(14,0)}{c7}{};
    \evencell{(16,0)}{c8}{};
    \oddrightedgecell{(18,0)}{c9}{};
    
    % Entries
    \begin{scope}[inner sep=0, minimum size=0, node distance=7]
      \oddbdysyll{c3}{}{$\alpha\beta$}{}{}{};
      \evenintsyll{c4}{}{$\beta\gamma$}{};
      \oddintsyll{c5}{}{$\alpha$}{};
      \evenbdysyll{c6}{}{$\beta\gamma$}{}{}{};
    \end{scope}
  \end{tikzpicture}
  }

  \caption[%
    Fourth and fifth example of representing a string graph with strips.
  ]{%
    \label{fig:v4-and-v5-as-string-graphs-and-strips}
    \textit{Fourth and fifth examples of representing a string graph with strips.}
  }
\end{thesisfigure}%
To represent it as a strip we round off both sides with stationary syllables, giving the uppermost strip in Figure \ref{subfig:stripifying-v3} or its reflection.

Similarly to the previous example, we allow variations of the strip where either stationary syllable is omitted since one can be inferred from the other. These are the remaining strips in the figure. In contrast to the previous example, we do not allow them both to be omitted. Doing so would result in an entirely blank strip, which of course is only capable of representing the zero module.

\nlpass{Pin-boundary syllables and the need for operating peakwise} We turn to Figures \ref{subfig:string-graph-v4} and \ref{subfig:string-graph-v5} and the string graphs $v_4$ and $v_5$ they feature. Strips representing them are drawn in Figures \ref{subfig:stripifying-v4} and \ref{subfig:stripifying-v5}. Neither needs any rounding off. The compatibility conditions stated above regarding $\dagger$-partnership of sources or targets of peak- or valley-neighbouring syllables are readily seen to hold. In these senses, there is nothing new to say about these strips.

The noteworthy objects are the syllables at their boundaries; observe that these are the same for both strips. We call these \emph{pin-boundary} syllables since they arise in the strip representation of a socle-quotient of a pin module.

What proves noteworthy about these syllables is their unusual behaviour when put to their intended use in strips: calculating syzygies of string modules. To more clearly remark the unorthodoxy of $v_5$, and shortly afterwards $v_4$, we compare them to the very routine conduct of $v_3$ and $v_1$.

Consult Figure \ref{subfig:syzygy-of-v3}. It is a two-row picture, the first row representing $v_3$ and the second $\syzygy^1(\Str v_3)$.
Each valley of the second row is determined by the peaks above it and nothing else. Figures \ref{subfig:syzygy-of-v2} and \ref{subfig:syzygy-of-v5} are similar, representing $v_2$ and $v_5$ and their syzygies, with the syzygy in each case determined from the original strip according to a local rule.
\begin{thesisfigure}
  \centering
  \subcaptionbox{
    \label{subfig:syzygy-of-v3}
    \textit{Calculating $\syzygy^1(\Str v_3) = \Str\big(
      \begin{tikzpicture}[scale=0.75]
      \begin{scope}[minimum size=0, inner sep=1]
        % Vertices of string graph
        \draw (0,0) node[](v1){$ 1$};
        \draw (1,0) node[](v2){$ 2$};
        \draw (2,0) node[](v3){$ 1$};
        \draw (3,0) node[](v4){$ 1$};
        \draw (4,0) node[](v5){$ 2$};
        %
        % Arrows of string graph
        \draw[->] (v1) --node[pos=0.5,above](){$\beta$} (v2);
        \draw[->] (v2) --node[pos=0.5,above](){$\gamma$} (v3);
        \draw[->] (v4) --node[pos=0.5,above](){$\alpha$} (v3);
        \draw[->] (v5) --node[pos=0.5,above](){$\gamma$} (v4);
      \end{scope}
      \end{tikzpicture}
    \big)$ using strips.} % In this example, the location of the (implied) lefthand boundary in the lower row is column $2$, one column outward from that of the upper row in column $3$. (The same phenomenon occurs with the location of the implied righthand boundary, from column $4$ to column $5$.) This behaviour~-- for the boundary to move outwards~-- is typical.
  }[
    0.95\textwidth
  ]{
    \begin{tikzpicture}[scale=0.7]
      % Cells of the 1st row (top row)
      \oddleftedgecell{(0,1)}{1 0}{};
      \evencell{(2,1)}{1 1}{};
      \oddcell{(4,1)}{1 2}{};
      \evencell{(6,1)}{1 3}{};
      \oddcell{(8,1)}{1 4}{};
      \evencell{(10,1)}{1 5}{};
      \oddcell{(12,1)}{1 6}{};
      \evencell{(14,1)}{1 7}{};
      \oddcell{(16,1)}{1 8}{};
      \evenrightedgecell{(18,1)}{1 9}{};
      
      % Cells of the 2nd row (bottom row)
      \evenleftedgecell{(0,-1)}{2 0}{};
      \oddcell{(2,-1)}{2 1}{};
      \evencell{(4,-1)}{2 2}{};
      \oddcell{(6,-1)}{2 3}{};
      \evencell{(8,-1)}{2 4}{};
      \oddcell{(10,-1)}{2 5}{};
      \evencell{(12,-1)}{2 6}{};
      \oddcell{(14,-1)}{2 7}{};
      \evencell{(16,-1)}{2 8}{};
      \oddrightedgecell{(18,-1)}{2 9}{};
      
      \begin{scope}[inner sep=0, minimum size=0, node distance=7]
        % Entries of the 1st row (top row)
        \oddtrivbdysyll{1 4}{}{$e$}{}{}{};
        \eventrivbdysyll{1 5}{}{$e'$}{}{}{};
        
        % Entries of the 2nd row (bottom row)
        \evenintsyll{2 4}{}{$\beta\gamma$}{};
        \oddintsyll{2 5}{}{$\gamma\alpha$}{};
      \end{scope}
      
      % Write column indices
      \columnindexlabel{(1 1)}{0};
      \columnindexlabel{(1 2)}{1};
      \columnindexlabel{(1 3)}{2};
      \columnindexlabel{(1 4)}{3};
      \columnindexlabel{(1 5)}{4};
      \columnindexlabel{(1 6)}{5};
      \columnindexlabel{(1 7)}{6};
      \columnindexlabel{(1 8)}{7};
  \end{tikzpicture}
  }
  
  \vspace{4ex}
  
  \subcaptionbox{%
    \label{subfig:syzygy-of-v2}%
    \textit{Calculating $\syzygy^1(\Str v_2) =
    \Str\big(
      \begin{tikzpicture}[scale=0.75]
      \begin{scope}[minimum size=0, inner sep=1]
        % Vertices of string graph
        \draw (0,0) node[](v1){$ 2$};
        \draw (1,0) node[](v2){$ 1$};
        \draw (2,0) node[](v3){$ 1$};
        \draw (3,0) node[](v4){$ 2$};
        \draw (4,0) node[](v5){$ 2$};
        \draw (5,0) node[](v6){$ 2$};
        %
        % Arrows of string graph
        \draw[->] (v1) --node[pos=0.5,above](){$\gamma$} (v2);
        \draw[->] (v2) --node[pos=0.5,above](){$\alpha$} (v3);
        \draw[->] (v4) --node[pos=0.5,above](){$\gamma$} (v3);
        \draw[->] (v4) --node[pos=0.5,above](){$\delta$} (v5);
        \draw[->] (v5) --node[pos=0.5,above](){$\delta$} (v6);
      \end{scope}
      \end{tikzpicture}
    \big)
    ~\oplus~
    \Str\big(
      \begin{tikzpicture}[scale=0.75]
      \begin{scope}[minimum size=0, inner sep=1]
        % Vertices of string graph
        \draw (0,0) node[](v1){$ 1$};
        \draw (1,0) node[](v2){$ 1$};
        %
        % Arrows of string graph
        \draw[->] (v2) --node[pos=0.5,above](){$\alpha$} (v1);
      \end{scope}
      \end{tikzpicture}
    \big)
    ~\oplus~
    \Str\big(
      \begin{tikzpicture}[scale=0.75]
      \begin{scope}[minimum size=0, inner sep=1]
        % Vertices of string graph
        \draw (0,0) node[](v1){$ 2$};
        \draw (1,0) node[](v2){$ 2$};
        \draw (2,0) node[](v3){$ 2$};
        %
        % Arrows of string graph
        \draw[->] (v2) --node[pos=0.5,above](){$\delta$} (v1);
        \draw[->] (v3) --node[pos=0.5,above](){$\delta$} (v2);
      \end{scope}
      \end{tikzpicture}
    \big)
    $ using strips.} % In this example, the locations of the boundaries move outwards from columns $1$ and $6$ in the top row to columns $0$ and $7$ in the bottom row (where their existence is implied). New boundaries appear in columns $3$, $4$, $5$ (implied) and $6$. Altogether, these boundaries delimit the three direct summands of $\syzygy^1(\Str v_1) $.
  }[0.95\textwidth]{%
    \begin{tikzpicture}[scale=0.7]
      % Cells of 1st row (top row)
      \oddleftedgecell{(0,1)}{1 0}{};
      \evencell{(2,1)}{1 1}{};
      \oddcell{(4,1)}{1 2}{};
      \evencell{(6,1)}{1 3}{};
      \oddcell{(8,1)}{1 4}{};
      \evencell{(10,1)}{1 5}{};
      \oddcell{(12,1)}{1 6}{};
      \evencell{(14,1)}{1 7}{};
      \oddcell{(16,1)}{1 8}{};
      \evenrightedgecell{(18,1)}{1 9}{};
      
      % Cells of 2nd row (bottom row)
      \evenleftedgecell{(0,-1)}{2 0}{};
      \oddcell{(2,-1)}{2 1}{};
      \evencell{(4,-1)}{2 2}{};
      \oddcell{(6,-1)}{2 3}{};
      \evencell{(8,-1)}{2 4}{};
      \oddcell{(10,-1)}{2 5}{};
      \evencell{(12,-1)}{2 6}{};
      \oddcell{(14,-1)}{2 7}{};
      \evencell{(16,-1)}{2 8}{};
      \oddrightedgecell{(18,-1)}{2 9}{};
      
      \begin{scope}[inner sep=0, minimum size=0, node distance=7]
        % Entries of 1st row (top row)
        \oddtrivbdysyll{1 2}{}{$e$}{}{}{};
        \evenintsyll{1 3}{}{$\alpha\beta$}{};
        \oddintsyll{1 4}{}{$\delta$}{};
        \evenintsyll{1 5}{}{$\gamma$}{};
        \oddintsyll{1 6}{}{$\gamma\alpha$}{};
        \eventrivbdysyll{1 7}{}{$e'$}{}{}{};
        
        % Entries of the 2nd row (bottom row)
        \evenintsyll{2 2}{}{$\gamma\alpha$}{};
        \oddintsyll{2 3}{}{$ \gamma $}{};
        \evenbdysyll{2 4}{}{ $\delta^2$ }{}{}{};
        \oddbdysyll{2 5}{}{ $\alpha$}{}{}{};
        % (Nothing in column 6)
        \oddbdysyll{2 7}{}{$ \delta^2 $}{}{}{};
      \end{scope}
      
      % Write column indices
      \columnindexlabel{(1 1)}{0};
      \columnindexlabel{(1 2)}{1};
      \columnindexlabel{(1 3)}{2};
      \columnindexlabel{(1 4)}{3};
      \columnindexlabel{(1 5)}{4};
      \columnindexlabel{(1 6)}{5};
      \columnindexlabel{(1 7)}{6};
      \columnindexlabel{(1 8)}{7};
    \end{tikzpicture}
  }
  
  \vspace{4ex}
  
  \subcaptionbox{
    \label{subfig:syzygy-of-v5}
    \textit{Calculating $\syzygy^1(\Str v_5) = \Str\big(
      \begin{tikzpicture}[scale=0.75]
      \begin{scope}[minimum size=0, inner sep=1]
        % Vertices of string graph
        \draw (0,0) node[](v1){$ 1$};
        \draw (1,0) node[](v2){$ 1$};
        \draw (2,0) node[](v3){$ 2$};
        \draw (3,0) node[](v4){$ 1$};
        %
        % Arrows of string graph
        \draw[->] (v2) --node[pos=0.5,above](){$\alpha$} (v1);
        \draw[->] (v2) --node[pos=0.5,above](){$\beta$} (v3);
        \draw[->] (v3) --node[pos=0.5,above](){$\gamma$} (v4);
      \end{scope}
      \end{tikzpicture}
    \big)$ using strips.} % In this example, the location of the lefthand boundary in the lower row is in column $3$, one column inward from the location of the boundary of the upper row, namely column $2$. (The same phenomenon occurs with the location of the righthand boundary, moving from column $5$ to $4$.) This behaviour~-- for the boundary to move inwards~-- is exceptional, occurring only when the boundary syllable in the top row is a pin-boundary syllable. 
  }[ 0.95\textwidth
  ]{
    \begin{tikzpicture}[scale=0.7]
    % Cells in 1st row (top row)
    \evenleftedgecell{(0,1)}{1 0}{};
    \oddcell{(2,1)}{1 1}{};
    \evencell{(4,1)}{1 2}{};
    \oddcell{(6,1)}{1 3}{};
    \evencell{(8,1)}{1 4}{};
    \oddcell{(10,1)}{1 5}{};
    \evencell{(12,1)}{1 6}{};
    \oddcell{(14,1)}{1 7}{};
    \evencell{(16,1)}{1 8}{};
    \oddrightedgecell{(18,1)}{1 9}{};
    
    % Cells in 2nd row (bottom row)
    \oddleftedgecell{(0,-1)}{2 0}{};
    \evencell{(2,-1)}{2 1}{};
    \oddcell{(4,-1)}{2 2}{};
    \evencell{(6,-1)}{2 3}{};
    \oddcell{(8,-1)}{2 4}{};
    \evencell{(10,-1)}{2 5}{};
    \oddcell{(12,-1)}{2 6}{};
    \evencell{(14,-1)}{2 7}{};
    \oddcell{(16,-1)}{2 8}{};
    \evenrightedgecell{(18,-1)}{2 9}{};
    
    \begin{scope}[inner sep=0, minimum size=0, node distance=7]
      % Entries of 1st row (top row)
      \oddbdysyll{1 3}{}{$\alpha\beta$}{}{}{};
      \evenintsyll{1 4}{}{$\beta\gamma$}{};
      \oddintsyll{1 5}{}{$\alpha$}{};
      \evenbdysyll{1 6}{}{$\beta\gamma$}{}{}{};
      
      % Entries of 2nd row (bottom row)
      \oddbdysyll{2 4}{}{$\alpha$}{}{}{};
      \evenbdysyll{2 5}{}{$\beta\gamma$}{}{}{};
    \end{scope}
    
    % Write column indices
    \columnindexlabel{(1 1)}{0};
    \columnindexlabel{(1 2)}{1};
    \columnindexlabel{(1 3)}{2};
    \columnindexlabel{(1 4)}{3};
    \columnindexlabel{(1 5)}{4};
    \columnindexlabel{(1 6)}{5};
    \columnindexlabel{(1 7)}{6};
    \columnindexlabel{(1 8)}{7};
  \end{tikzpicture}
  }
  
  \caption[Movement of boundaries
  ]{
    \label{fig:movement-of-boundaries}
    \textit{Movement of boundaries.} Existing boundaries move one column inwards or outwards from one row to the next, the latter situation being more common by far. Additional boundaries may also appear when the lower strip ``branches'', representing a nontrivial decomposition of the syzygy (refer Figure \ref{fig:branching}). 
  }
  
\end{thesisfigure}

The strip representing $v_3$ in the top row of Figure \ref{subfig:syzygy-of-v3} has boundary syllables in columns 3 and 4, and neither of these is a pin-boundary syllable. Underneath it is a strip representing $\syzygy^1(\Str v_3)$. This actually has no boundary syllables, but we can infer \emph{where} stationary boundary syllables would round it off~-- specifically, columns 2 and 5~-- and \emph{which} stationary boundary syllables they would be. We see that, from one row to the next, the locations of the (potentially implied) boundaries move one column outward on each side.

The strip representing $v_2$ in the top row of Figure \ref{subfig:syzygy-of-v2} has boundary syllables in columns 1 and 6 but neither of these is a pin-boundary syllable. The second row has several syllables and several gaps. The boundaries present or implied in columns 3, 4, 5 and 6 are surplus to the present argument and we will revisit them later. Those aside we certainly infer boundary syllables in columns 0 and 7 of the bottom row, because there the blank cells appear in a peak with a nonblank cell. Once again, we notice that the locations of the (potentially implied) boundaries move one column outwards, from columns 1 and 6 to 0 and 7.

%Taken with the aforementioned boundary locations in consecutive pairs, they serve to demarcate the data of several strips. In other words, if we transplanted the column intervals $[0,3]$ or $[4,5]$ or $[6,7]$ of the bottom row into otherwise-blank rows, we would obtain strips representing the direct summands of $\syzygy^1( \Str v_1 )$. No such situation arose in Figure \ref{subfig:syzygy-of-v3} because $\syzygy^1( \Str v_3 )$ is indecomposable.

Lastly, turn to Figure \ref{subfig:syzygy-of-v5} which represents $v_5$ and its syzygy using strips and in which appear the dreaded pin-boundary syllables. Two things are remarkable about this figure. The first is that the boundary on each side moves one column inwards. The second is how this impacts columns 3 and 4. Taking column 3 as an illustrative example, we see that underneath the syllable $\big( \ilsyll{ \circ }{\beta\gamma}{0} \big)$ is $\big( \ilsyll{ \circ }{\alpha}{1} \big)$. The lower syllable is not $\big( \ilsyll{ \circ }{\beta\gamma}{0} \big)\descent=\big( \ilsyll{ \circ }{\alpha}{0} \big)$ but, rather, the boundary analogue. 

A general rule, obeyed in Figures \ref{subfig:syzygy-of-v3} and \ref{subfig:syzygy-of-v2} but violated here in \ref{subfig:syzygy-of-v5}, is that underneath a syllable $\bdp$ appears $\bdp\descent$. However, the cause of this misbehaviour is not the syllables in columns 3 and 4 of the top row. It is in fact the pin-boundary syllables $\big( \ilsyll{ \circ }{ \alpha\beta }{ 1 } \big)$ and $ \big( \ilsyll{ \circ }{ \beta\gamma }{ 1 } \big) $ that neighbour them.

We should think of pin-boundary syllables as interfering entities that cause boundaries to move one column inwards, therefore defying our otherwise-successful hopes of $\descent$-orbits for columns. Pin-boundary syllables $\big( \ilsyll{ \circ }{ p }{ 1 } \big)$ exist to represent a component $p\alpha$ of a commutativity relation of $A$, say $p\alpha-q$, where $\alpha$ is an appropriate arrow. It is thanks to their existence that our formalism calculates syzygies using a peakwise operation rather than a cellwise one, for their corrupting influence can extend to their valley neighbour.

The pinnacle of their meddlesome activity occurs when two pin-boundary syllables neighbour one another in a strip, a scenario arising iff the strip represents $P/\soc P$ for some pin module $P$. In this case the two boundaries on either side seem to pass through each other.

We mark such an event with special symbols, as witnessed in Figure \ref{subfig:syzygy-of-v4}. We call these \emph{virtual} syllables because, unlike other syllables, these reserved symbols do not represent any part of a string graph or string module; they are simply indicators of the presence of short exact sequences of the form
$
\begin{tikzcd}[sep=small, cramped]
  0 \ar[r]
    \& \soc P \ar[r]
      \& P \ar[r]
        \& P/\soc P \ar[r]
          \& 0\tcomma
\end{tikzcd}
$
for $P$ pin. These virtual syllables appear two-at-a-time in a valley underneath a peak of pin-boundary syllables. Each virtual syllable appears alone in its own peak. Directly below any virtual syllable belongs the appropriate stationary syllable.
\begin{thesisfigure}
  \centering
  
  \subcaptionbox{
    \label{subfig:syzygy-of-v4}
    \textit{Calculating $\syzygy^1( \Str v_4 ) = \Str\big(
      \raisebox{-0.3ex}{
      \begin{tikzpicture}[scale=0.75]
      \begin{scope}[minimum size=0, inner sep=1]
        % Vertices of string graph
        \draw (0,0) node[](v1){$ 1$};
      \end{scope}
      \end{tikzpicture}
      }
    \big)$ using strips.} The path $e$ is the stationary path at $\target\big( \ilsyll{ \circ }{ \alpha\beta }{ 1 } \big)$, the target of the lefthand pin-boundary syllable. Similarly $e'$.
  }[0.95\textwidth
  ]{
    \begin{tikzpicture}[scale=0.7]
      % Cells of the 1st row (top row)
      \oddleftedgecell{(0,2)}{1 0}{};
      \evencell{(2,2)}{1 1}{};
      \oddcell{(4,2)}{1 2}{};
      \evencell{(6,2)}{1 3}{};
      \oddcell{(8,2)}{1 4}{};
      \evencell{(10,2)}{1 5}{};
      \oddcell{(12,2)}{1 6}{};
      \evencell{(14,2)}{1 7}{};
      \oddcell{(16,2)}{1 8}{};
      \evenrightedgecell{(18,2)}{1 9}{};
      
      % Cells of the 2nd row
      \evenleftedgecell{(0,0)}{2 0}{};
      \oddcell{(2,0)}{2 1}{};
      \evencell{(4,0)}{2 2}{};
      \oddcell{(6,0)}{2 3}{};
      \evencell{(8,0)}{2 4}{};
      \oddcell{(10,0)}{2 5}{};
      \evencell{(12,0)}{2 6}{};
      \oddcell{(14,0)}{2 7}{};
      \evencell{(16,0)}{2 8}{};
      \oddrightedgecell{(18,0)}{2 9}{};
      
      % Cells of the 3rd row (bottom row)
      \oddleftedgecell{(0,-2)}{3 0}{};
      \evencell{(2,-2)}{3 1}{};
      \oddcell{(4,-2)}{3 2}{};
      \evencell{(6,-2)}{3 3}{};
      \oddcell{(8,-2)}{3 4}{};
      \evencell{(10,-2)}{3 5}{};
      \oddcell{(12,-2)}{3 6}{};
      \evencell{(14,-2)}{3 7}{};
      \oddcell{(16,-2)}{3 8}{};
      \evenrightedgecell{(18,-2)}{3 9}{};
      
      \begin{scope}[inner sep=0, minimum size=0, node distance=7]
        % Entries of the 1st row (top row)
        \oddbdysyll{1 4}{}{$\alpha\beta$}{}{}{};
        \evenbdysyll{1 5}{}{$\beta\gamma$}{}{}{};
        
        % Entries of the 2nd row
        \evenvirtualsyll{2 4}{}{$e$}{};
        \oddvirtualsyll{2 5}{}{$e'$}{};
        
        % Entries of the 3rd row (bottom row)
        \oddtrivbdysyll{3 4}{}{$e$}{}{}{};
        \eventrivbdysyll{3 5}{}{$e'$}{}{}{};
      \end{scope}
      
      % Write column indices
      \columnindexlabel{(1 1)}{0};
      \columnindexlabel{(1 2)}{1};
      \columnindexlabel{(1 3)}{2};
      \columnindexlabel{(1 4)}{3};
      \columnindexlabel{(1 5)}{4};
      \columnindexlabel{(1 6)}{5};
      \columnindexlabel{(1 7)}{6};
      \columnindexlabel{(1 8)}{7};
    \end{tikzpicture}
  }
  
  \caption[
    Syzygy of a socle-quotient of a pin module
  ]{
    \label{fig:syzygy-of-socle-quotient-of-pin}
    \textit{Syzygy of a socle-quotient of a pin module.}
  }
\end{thesisfigure}

The data of a virtual syllable $\big( \ilsyll{ \circ }{e}{0} \big)$ are a stationary path $e$ and a stability term $\ep=0$. This path $e$ in the overquiver represents the target vertex of a commutativity relation. The associated stationary syllable is then $\big( \ilsyll{ \circ }{e}{1} \big)$.

We recall that elsewhere a stability term $\ep=0$ denotes an interior syllable so it may be tempting to think of these virtual syllables as ``interior stationary'' syllables, in keeping with that precedent. We know of no benefit to be garnered from this interpretation and so we do not take it up in our exposition.
% Mercifully, when a door closes a window opens: Our prohibition from one purpose leaves these symbols available for another: signifying the syzygy of a socle-quotient of a pin module, as described above.

\nlpass{Branching} There is one final phenomenon to explore: when the syzygy of a string module has a nontrivial direct sum decomposition. This has already occurred under our very noses.

Figure \ref{subfig:syzygy-of-v2} depicts a syzygy calculation, with $v_2$ represented in the top row and $\syzygy^1( \Str v_2 )$ in the bottom row. Whereas in other examples the support has been an uninterrupted sequence of interior syllables bounded on either side by (potentially implied) boundary syllables, here the support is several such sequences, respectively delimited by the boundaries in columns 0 and 3, 4 and 5, and 6 and 7, and all flattened together onto the page in a single row.

Such flattening only muddies the waters. Our task is made much clearer by giving each of these individual strips their own row as in Figure \ref{fig:branching}. One consequence of this generosity is that we must generally index rows using the vertices of a rooted tree. This allows the nontrivial direct sum decomposition of the syzygy of a string to be reflected by \emph{branching} in the tree. 
\begin{thesisfigure}
  \centering
  
  \begin{tikzpicture}[scale=0.6]
      % Cells of 1st row (top row)
      \draw (0,1) node[](1){};
      
      \oddleftedgecell{(1) ++(0,0)}{1 0}{};
      \evencell{(1) ++(2,0)}{1 1}{};
      \oddcell{(1) ++(4,0)}{1 2}{};
      \evencell{(1) ++(6,0)}{1 3}{};
      \oddcell{(1) ++(8,0)}{1 4}{};
      \evencell{(1) ++(10,0)}{1 5}{};
      \oddcell{(1) ++(12,0)}{1 6}{};
      \evencell{(1) ++(14,0)}{1 7}{};
      \oddcell{(1) ++(16,0)}{1 8}{};
      \evenrightedgecell{(1) ++(18,0)}{1 9}{};
      
      % Entries of 1st row (top row)
      \begin{scope}[inner sep=0, minimum size=0, node distance=7]
        \oddtrivbdysyll{1 2}{}{$e$}{}{}{};
        \evenintsyll{1 3}{}{$\alpha\beta$}{};
        \oddintsyll{1 4}{}{$\delta$}{};
        \evenintsyll{1 5}{}{$\gamma$}{};
        \oddintsyll{1 6}{}{$\gamma\alpha$}{};
        \eventrivbdysyll{1 7}{}{$e'$}{}{}{};
      \end{scope}
      
      \newcommand{\rowdisplacement}{(0.8,-1.4)}
      \draw (0,-2) node[](2a){};
      \draw (2a) ++\rowdisplacement node[](2b){};
      \draw (2b) ++\rowdisplacement node[](2c){};
      
      % Cells of row 2a (backmost bottom row)
      \evenleftedgecell{(2a)++(0,0)}{2a 0}{};
      \oddcell{(2a)++(2,0)}{2a 1}{};
      \evencell{(2a)++(4,0)}{2a 2}{};
      \oddcell{(2a)++(6,0)}{2a 3}{};
      \evencell{(2a)++(8,0)}{2a 4}{};
      \oddcell{(2a)++(10,0)}{2a 5}{};
      \evencell{(2a)++(12,0)}{2a 6}{};
      \oddcell{(2a)++(14,0)}{2a 7}{};
      \evencell{(2a)++(16,0)}{2a 8}{};
      \oddrightedgecell{(2a)++(18,0)}{2a 9}{};
      
      % Entries of row 2a (backmost bottom row)
      \begin{scope}[inner sep=0, minimum size=0, node distance=7]
        \evenintsyll{2a 2}{}{$\gamma\alpha$}{};
        \oddintsyll{2a 3}{}{$ \gamma $}{};
        \evenbdysyll{2a 4}{}{ $\delta^2$ }{}{}{};
      \end{scope}
      
      % Cells of row 2b (middle bottom row)
      \evenleftedgecell{(2b)++(0,0)}{2b 0}{};
      \oddcell{(2b)++(2,0)}{2b 1}{};
      \evencell{(2b)++(4,0)}{2b 2}{};
      \oddcell{(2b)++(6,0)}{2b 3}{};
      \evencell{(2b)++(8,0)}{2b 4}{};
      \oddcell{(2b)++(10,0)}{2b 5}{};
      \evencell{(2b)++(12,0)}{2b 6}{};
      \oddcell{(2b)++(14,0)}{2b 7}{};
      \evencell{(2b)++(16,0)}{2b 8}{};
      \oddrightedgecell{(2b)++(18,0)}{2b 9}{};

      % Entries of row 2b (middle bottom row)
      \begin{scope}[inner sep=0, minimum size=0, node distance=7]
        \oddbdysyll{2b 5}{}{ $\alpha$}{}{}{};
      \end{scope}

      % Cells of row 2c (foremost bottom row)
      \evenleftedgecell{(2c)++(0,0)}{2c 0}{};
      \oddcell{(2c)++(2,0)}{2c 1}{};
      \evencell{(2c)++(4,0)}{2c 2}{};
      \oddcell{(2c)++(6,0)}{2c 3}{};
      \evencell{(2c)++(8,0)}{2c 4}{};
      \oddcell{(2c)++(10,0)}{2c 5}{};
      \evencell{(2c)++(12,0)}{2c 6}{};
      \oddcell{(2c)++(14,0)}{2c 7}{};
      \evencell{(2c)++(16,0)}{2c 8}{};
      \oddrightedgecell{(2c)++(18,0)}{2c 9}{};
      
      % Entries of row 2c (foremost bottom row)
      \begin{scope}[inner sep=0, minimum size=0, node distance=7]
        \oddbdysyll{2c 7}{}{$ \delta^2 $}{}{}{};
      \end{scope}
      
      % Write column indices
      \columnindexlabel{(1 1)}{0};
      \columnindexlabel{(1 2)}{1};
      \columnindexlabel{(1 3)}{2};
      \columnindexlabel{(1 4)}{3};
      \columnindexlabel{(1 5)}{4};
      \columnindexlabel{(1 6)}{5};
      \columnindexlabel{(1 7)}{6};
      \columnindexlabel{(1 8)}{7};
      
      % Draw row indices
      \newcommand{\rowindexdisplacement}{(-3,0)}
      
      \newcommand{\rowindex}[1]{
        \draw (#1) ++\rowindexdisplacement node[](#1 index){};
        \filldraw[lightgray] (#1 index) circle (0.15);
        % 1: label of row node
      }
      \newcommand{\rowindexconnector}[2]{
        \draw[lightgray] (#1 index) -- (#2 index);
        % 1: Label of row node for source index
        % 2: Label of row node for target index
      }
      
      \rowindex{1};
      \rowindex{2a};
      \rowindex{2b};
      \rowindex{2c};
      \rowindexconnector{1}{2a};
      \rowindexconnector{1}{2b};
      \rowindexconnector{1}{2c};
    \end{tikzpicture}
    
  \caption[
    Branching
  ]{
    \label{fig:branching}
    \textit{Branching.}  Rather than flattening the strips representing the direct summands of $\syzygy^1( \Str v_2 )$ onto a single row, as in Figure
    \ref{subfig:syzygy-of-v2}, we place them into rows of their own. Each of these rows is directly above the top row. In each, columns 0 to 7 inclusive are shown.
  }
\end{thesisfigure}

Another consequence is that 
% when one row has several rows immediately underneath it, has several valleys underneath it.
we must specify into which rows the cells of a valley under a given peak are placed. There is a general rule which we hope Figure \ref{fig:branching} is clear enough to illustrate.
%Namely, since the entries of a valley are either both interior or both not interior (which is nothing more than the manifestation in our formalism of the biseriality of $A$), the cells appear in the same row in the former case and distinct rows in the latter case, and obviously boundaries should appear in the same row as the sequence of interior syllables that they delimit.

\nlpass{Recap} We close this subsection by consolidating the examples we have explored.

We have seen strings represented as rows of symbols, rows that we call strips. A string graph can be represented by many strips, but any strip represents a unique string.

We have started to see how these symbols can be entered into an array, with rows indexed by the vertices of a rooted tree and columns indexed by the integers (which are also just the vertices of a graph:
$\begin{tikzcd}[sep=small, cramped]
  \cdots \ar[r, no head]
    \& 0 \ar[r, no head]
      \& 1 \ar[r, no head]
        \& 2 \ar[r, no head]
          \& 3 \ar[r, no head]
            \& \cdots
\end{tikzcd}
$), to represent syzygy-taking. At least, we have seen the first two levels of such an array, representing a string and its syzygy. The contents of one row determine the contents of the next rows according to a local rule: a peak above determines the valley below.

The symbols used are called syllables. A syllable comprises three pieces of information. Largely these syllables represent parts of a string module and whether that part is at the boundary or in the interior of the string.

A select few symbols are kept aside for when we represent socle-quotients of pin modules and their syzygies. Indeed, strings whose boundaries resemble the boundaries of such socle-quotients are the only ones that require the local rule to operate at the scope of peaks of syllables, rather than individual syllables. These errant examples aside, the location of the boundary moves one column outwards from one row to the next. New boundaries can also appear, indicating that a line of symbols are to be separated into different rows, with the resulting strips representing direct summands.

This is the core of our model for tabulating the syzygy information. What remains are the details of its construction and the patterns it illuminates.

\subsection{Strips and string modules}

In this subsection, we give the formal definition of strips and verify that they represent string modules.

\nlpass{Peak and valley compatibility} Let $(\bdp,\bdq)$ be a pair of (nonvirtual) syllables for $A$.
\begin{thesisfigure}
  \centering
  
  \subcaptionbox{%
    \textit{Peaks with exactly one blank entry.} The other entry may be a boundary or interior syllable.
  }[
    0.95\textwidth
  ]{
    \begin{tikzpicture}[scale=0.5]
      \peakrep{(-7.5, 0)}{pk1}{}{};
      \oddintsyll{pk1 l}{}{}{};
    
      \peakrep{(-2.5,0)}{pk2}{}{};
      \evenintsyll{pk2 r}{}{}{};
      
      \peakrep{(2.5,0)}{pk3}{}{};
      \oddbdysyll{pk3 l}{}{}{}{}{};
    
      \peakrep{(7.5, 0)}{pk4}{}{};
      \evenbdysyll{pk4 r}{}{}{}{}{};
    \end{tikzpicture}
  }
  
  \subcaptionbox{%
    \textit{Peaks with $\dagger$-partnered sources.} The syllables may be any combination of interior or boundary, provided their sources are $\dagger$-partners.
  }[
    0.95\textwidth
  ]{
    \begin{tikzpicture}[scale=0.5]
      \peakrep{(-7.5,0)}{pk1}{}{};
      \oddintsyll{pk1 l}{$i$}{}{};
      \evenintsyll{pk1 r}{$\idag$}{}{};
    
      \peakrep{(-2.5,0)}{pk2}{}{};
      \oddbdysyll{pk2 l}{$i$}{}{}{}{};
      \evenintsyll{pk2 r}{$\idag$}{}{};
      
      \peakrep{(2.5,0)}{pk3}{}{};
      \oddintsyll{pk3 l}{$i$}{}{};
      \evenbdysyll{pk3 r}{$\idag$}{}{}{}{};
    
      \peakrep{(7.5,0)}{pk4}{}{};
      \oddbdysyll{pk4 l}{$i$}{}{}{}{};
      \evenbdysyll{pk4 r}{$\idag$}{}{}{}{};
    \end{tikzpicture}
  }
  
  \caption[%
    Peak-compatible syllables
  ]{%
    \label{fig:peak-compatible-pairs}
    \textit{Peak-compatible pairs of syllables.}
  }
  
  \vspace{4ex}
  
  \subcaptionbox{%
    \textit{Blank valley.}
  }[%
    0.2\textwidth
  ]{
    \begin{tikzpicture}[scale=0.5]
      \valleyrep{(0,0)}{vl1}{}{}
    \end{tikzpicture}
  }  
  \subcaptionbox{%
    \textit{Boundary valleys.}
  }[%
    0.5\textwidth
  ]{
    \begin{tikzpicture}[scale=0.5]
      \valleyrep{(-2.5,0)}{vl1}{}{}
      \evenbdysyll{vl1 l}{}{}{}{}{};
      
      \valleyrep{(2.5,0)}{vl2}{}{};
      \oddbdysyll{vl2 r}{}{}{}{}{};
    \end{tikzpicture}
  }
  \subcaptionbox{%
    \textit{Interior valley.}
  }[%
    0.2\textwidth
  ]{
    \begin{tikzpicture}[scale=0.5]
      \valleyrep{(0,0)}{vl1}{}{}
      \evenintsyll{vl1 l}{}{}{$i$};
      \oddintsyll{vl1 r}{}{}{$\idag$};
    \end{tikzpicture}
  }
  
  \caption[
    Valley-compatible pairs
  ]{
    \label{fig:valley-compatible-pairs}
    \textit{Valley-compatible pairs of syllables.}
  }
\end{thesisfigure}
\begin{enumerate}
    \item We say $\bdp$ and $\bdq$ are \define{peak compatible} (see Figure \ref{fig:peak-compatible-pairs}) if
    \begin{enumerate}
        \item both are blank (this is called a \emph{blank peak}), or
        \item exactly one is blank (an \emph{implied peak}), or
        \item neither is blank and $\source(\bdp)^\dagger = \source(\bdq)$ and the orientations of $\bdp$ and $\bdq$ are respectively $ -1$ and $+1$.
    \end{enumerate} 
    
    \item We say $\bdp$ and $\bdq$ are \define{valley compatible} (see Figure \ref{fig:valley-compatible-pairs}) if
    \begin{enumerate}
        \item both are blank (this is called a \emph{blank valley}), or
        \item one is blank and the other is boundary (a \emph{boundary valley}), or
        \item both are interior and $\target(\bdp)^\dagger = \target(\bdq)$ and the orientations of $\bdp$ and $\bdq$ are respectively $+1$ and $-1$ (an \emph{interior valley}).
    \end{enumerate}
\end{enumerate}
These relations are illustrated in Figures \ref{fig:peak-compatible-pairs} and \ref{fig:valley-compatible-pairs}. No syllable is peak or valley compatible with itself except the blank syllable, which is both.

\nlpass{Peaks and valleys} A \define{peak} is a pair of peak-compatible syllables, which we draw as
\iltikzpicture{
  \peakrep{(0,0)}{pk1}{$\bdp$}{$\bdq$};
}, the orientations being implicit. (Note that $(\bdp^{-1},\bdq)$ is a peak-compatible pair iff $(\bdq^{-1},\bdp)$ is, by the above definition.) Similarly, we define \define{valleys} and denote them as
\iltikzpicture{
  \valleyrep{(0,0)}{vl1}{$\bdp$}{$\bdq$}
}. This notation makes defining the \define{reflection}
\iltikzpicture{
  \peakrep{(0,0)}{pk1}{$\bdq$}{$\bdp$};
}
of a peak or
\iltikzpicture{
  \valleyrep{(0,0)}{vl1}{$\bdq$}{$\bdp$}
} of a valley obvious.

\nrex Any two (horizontally) neighbouring cells in Figures \ref{subfig:stripifying-v1}, \ref{subfig:stripifying-v2}, \ref{subfig:stripifying-v3}, \ref{subfig:stripifying-v4}, \ref{subfig:stripifying-v5}, \ref{subfig:syzygy-of-v3}, \ref{subfig:syzygy-of-v2}, \ref{subfig:syzygy-of-v5} or \ref{fig:branching} are either a peak or valley (as appropriate).

\nlpass{Strips} A \define{strip $w$} is a concave, non-necessarily bounded juxtaposition of (nonvirtual) syllables alternately forming peaks and valleys; here, \emph{concave} means that no blank syllable is between two nonblank syllables.

Formally we consider the juxtaposition to be a single row of cells with columns indexed by $\integers$; thus, a strip $w$ is a function $w \colon \ilarrow{ \integers }{}{ \Syll(A) }$. Of course, the cells are in 1-to-1 correspondence with the columns. The \emph{entry} of cell $k$ is $w(k)$.

\nlpass{Support of a strip} It follows that a strip $w$ is blank on all but an interval subset of $\integers$, which we call its support $\supp w$

\nlpass{Neighbours} For any $k \in \supp w$, one out of $\{w(k),w(k+1)\}$ and $\{w(k-1), w(k)\}$ is a peak and one is a valley. The \define{peak neighbour} of $k$ is the cell forming a peak with $k$; the \define{valley neighbour} is defined similarly.

\nlpass{Reflection} By precomposing with the reflection $\ilmapsto{k}{}{-k}$ on $\integers$ and postcomposing with the orientation involution $\ilmapsto{\bdp^s}{}{\bdp^{-s}}$ on $\Syll(A)$, we obtain the \define{reflection} of $w$. This is clearly a well-defined operation since these involutions carry peaks to peaks and valleys to valleys. 

\nrex Examples of strips abounded in Figures \ref{subfig:stripifying-v1}, \ref{subfig:stripifying-v2}, \ref{subfig:stripifying-v3}, \ref{subfig:stripifying-v4}, \ref{subfig:stripifying-v5}, \ref{subfig:syzygy-of-v3}, \ref{subfig:syzygy-of-v2}, \ref{subfig:syzygy-of-v5} and \ref{fig:branching}. 

\npass Any strip represents a string module in an (hopefully) obvious way, which the following proposition formalises.

\begin{proposition}\label{prop:strips-give-string-graphs}
Any strip $w$ represents a well-defined string graph, hence a string module, and moreover $w$ and its reflection both represent the same string graph.

Conversely, any string graph can be represented by a strip.
\end{proposition}

\proof Towards proving the forward implication, write $p,q$ for the paths underlying syllables $\bdp,\bdq$ in $w$. Recall from Definition \ref{def:syllables} that these underlying paths $p,q$ satisfying $p,q \in N \setminus C$ and so, notably, represent $Q$-paths whose $A$-residue is linearly independent of all other $A$-paths.

We can view each $p,q$ as some linear subgraph of $\OQ$ of the form $\big(
\begin{tikzcd}[sep=small]
  \syllableplaceholder \ar[r]
    \& \syllableplaceholder \ar[r]
      \& \cdots \ar[r]
        \& \syllableplaceholder
\end{tikzcd}\big)
$, where each vertex and arrow canonically represents one of the ground quiver $Q$ of $A$. This representation gives rise to quiver homomorphisms $\ilarrow{p}{}{Q}$ or $\ilarrow{q}{}{Q}$. (In the case that $\bdp$ or $\bdq$ is blank, the associated homomorphism has empty domain.)

When
\iltikzpicture{\peakrep{(0,0)}{pk1}{ $\bdp$}{$\bdq$}}
is a peak, then the sources (if any) of $p, q$ are distinct but represent the same vertex of $Q$, and therefore have the same image in the quiver homomorphism. Therefore, in the disjoint union of the homomorphisms $\ilarrow{p}{}{Q}$ and $\ilarrow{q}{}{Q}$, we can identify the source vertices in the domains of paths $p,q$ coming from peak-neighbour syllables and obtain a well-defined graph homomorphism. In addition (still assuming
\iltikzpicture{\peakrep{(0,0)}{pk1}{ $\bdp$}{$\bdq$}}
is a peak), the first arrows of $p$ and $q$ are distinct and so represent distinct arrows of $Q$. It follows that, when we identify the source vertices as above, pairs of arrows in the domain that common source are mapped by this homomorphism to distinct $Q$-arrows.

The preceding paragraph remains true when we simultaneously replace: ``peak'' by ``valley''; ``peak-neighbour(ing)'' by ``valley-neighbour(ing)''; ``source'' by ``target'' and; ``source vertex'' by ``sink vertex''.

This act of identifying linear quivers at their source vertices or sink vertices yields a quiver $G$ whose underlying graph is linear also. Since the cells of $w$ are $\integers$-indexed, the underlying graph of our constructed quiver will be connected and either finite, unbounded in one direction or unbounded in two directions. As we have already commented, pairs of arrows whose common source is a source vertex or whose common target is a sink vertex are mapped to distinct $Q$-arrows. As also commented above, any maximal equioriented subquiver of $G$ (ie, a full subquiver that is sub-$1$-regular and that is specified by a set of vertices including one source vertex of $G$ and one sink vertex of $G$) represents an independent path; therefore an arbitrary equioriented subquiver of $G$ represents one too \emph{a fortiori}. These three conditions defining a string graph hold, so the first assertion is proven.

The second claim follows immediately, since the foregoing abstract construction is not affected by reflecting $w$.

For the converse assertion, let $v: \ilarrow{G}{}{Q}$ be a string graph, assumed without loss of generality to be indecomposable.

If $G$ contains no arrows, and so comprises a single vertex $x$ satisfying $v(x)=i \in Q_0$, then lift $i$ to an $\OQ$-vertex that we also denote $i$; we may represent $v$ using a function $w \colon \ilarrow {\integers}{}{ \Syll(A) }$ that takes blank values on all integers save on $0$, where $w(0)=\bde_i^1$. (This superscript $1$ is the orientation of the syllable.) The support of this function is concave in a trivial way, all potential peaks are blank except for one that is implied, and all potential valleys are blank except one that is boundary. Consequently $w$ truly is a strip. Otherwise, assume $G$ contains at least one arrow.

To each maximal equioriented subquiver $\Gamma \ceq \big(\begin{tikzcd}[sep=small]
    y \ar[r, "x_1"]
      \& \syllableplaceholder \ar[r, "x_2"]
        \& \cdots \ar[r, "x_\ell"]
          \& z
    \end{tikzcd}\big)$
of $G$ we may associate a syllable as follows, writing $p \ceq v(x_1 x_2 \cdots x_\ell)$ for convenience: if the sink vertex $z$ of this subquiver has indegree $1$ in $G$, we associate the boundary syllable $\widehat{\Gamma} \ceq \big(\ilsyll{\syllableplaceholder}{p}{1}\big)$; otherwise, we associate the interior syllable $\widehat{\Gamma} \ceq \big(\ilsyll{\syllableplaceholder}{p}{0}\big)$. (The orientations of these syllables is yet to be assigned.)

Choose one such maximal equioriented subquiver $\Gamma_0$ of $G$. Let $\Gamma_1$ be the maximal equioriented subquiver of $G$ sharing a sink vertex with $\Gamma_0$ (if any exists), let $\Gamma_2$ be that sharing a source vertex with $\Gamma_1$ (if any exists), let $\Gamma_2$ be that sharing a sink vertex with $\Gamma_1$ (if any exists), and so on, iteratively defining $\Gamma_k$ for $k \geq 0$. Similarly on the other side, let $\Gamma_{-1}$ be the maximal equioriented subquiver of $G$ sharing a source vertex with $\Gamma_0$ (if any exists), let $\Gamma_{-2}$ be that sharing a sink vertex with $\Gamma_{-1}$ (if any exists), and so on for $\Gamma_k$ with $k < 0$.

For any $k \in \integers$ for which $\Gamma_k$ exists, define $v(k) \ceq (\widehat{\Gamma}_k)^{(-1)^k}$. (This exponent $(-1)^k$ is the orientation we assign to the syllable $\widehat{\Gamma}_k$.) For remaining $k$, let $v(k)$ be blank.

Our iterative definition of the $\Gamma_k$, working outwards from $\Gamma_0$, ensures that $\supp v$ is concave. Potential peaks containing two nonblank syllables come from neighbouring subquivers $\Gamma_k, \Gamma_m$ whose common vertex is a source vertex of $G$. Because the two arrows incident on that source vertex have different images in $v$, we have that $\widehat{\Gamma}_k$ and $\widehat{\Gamma}_m$ are distinct and their sources are exchanged by $\dagger$. We deduce that these potential peaks truly are peaks.

A similar argument shows that all potential valleys are indeed valleys. The only technicality to note is that~-- when $\widehat{\Gamma}_{k}$ is nonblank~-- $\widehat{\Gamma}_k$ is an interior syllable iff the corresponding sink vertex of $G$ has indegree $2$ in $G$ iff $\widehat{\Gamma}_{k+(-1)^k}$ is defined. Thus, the potential valleys containing one nonblank syllable are boundary valleys and those containing two are interior valleys.

Since $v$ is concave, its peaks are peaks and its valleys valleys, we concluce that $v$ defines a string, as claimed.\qed

\subsection{Patches and projective covers}

In this subsection, we relate strips and patches in a fashion that emulates string modules and projective covers.

\begin{proposition}\label{prop:peaks-are-tops-of-patches}
Any peak \iltikzpicture{\peakrep{(0,0)}{pk1}{$\bdp$}{$\bdq$}} is the top row of exactly one patch.
\end{proposition}

\proof Lemma \ref{lem:patches-have-different-tops} implies that each peak is the top of at most one patch; what remains is the ``at least'' argument. If $\bdp$ and $\bdq$ are both blank then the peak is the top of the blank patch. If neither is blank then the assumption $\source(\bdp)^\dagger = \source(\bdq)$ ensures they are the top of some patch constructed in Subsection \ref{subsec:construction-of-patches}. If exactly one of $\bdp,\bdq$ is blank then the peak is the top row of the unique amended version of a patch constructioned in the previous sentence.\qed

\begin{corollary}\label{cor:patch-covers-of-strips-well-defined}
Any strip is the top row of a well-defined line of patches.
\end{corollary}

\proof We can view a strip as a line of peaks. Each of which is the top of some patch by the previous lemma, whence we obtain the line of patches.\qed

\begin{proposition}\label{prop:line-of-patches-is-projective-cover}
If the strip $w$ is the top row of a line of patches, then $\projcover(\Str w)$ is the direct sum of the associated projectives (using the association given in Definition \ref{def:projective-of-a-patch}).
\end{proposition}

\proof We know from Proposition \ref{prop:peaks-are-tops-of-patches} that to any nonblank peak of $w$ we may associate some nonblank, nonvirtual patch and hence some indecomposable projective $P_i$. We know from Proposition \ref{prop:strips-give-string-graphs} that the nonblank peaks of $w$ correspond to the source vertice in the corresponding string graph. These source vertices are a basis of $(\Str w)/ (\rad( \Str w))$, and each is fixed by exactly one primitive idempotent $e_j$. Clearly $i=j$ for each nonblank peak of $w$ and the result follows.\qed 

\nrex Let
$w \ceq \big( \begin{tikzpicture}[scale=0.75]
  \begin{scope}[minimum size=0, inner sep=1]
    % Vertices of string graph
    \draw (0,0) node[](v1){$ 1$};
    \draw (1,0) node[](v2){$ 1$};
    \draw (2,0) node[](v3){$ 2$};
    \draw (3,0) node[](v4){$ 2$};
    \draw (4,0) node[](v5){$ 1$};
    \draw (5,0) node[](v6){$ 1$};
    \draw (6,0) node[](v7){$ 2$};
    %
    % Arrows of string graph
    \draw[->] (v1) --node[pos=0.5,above](){$\alpha$} (v2);
    \draw[->] (v2) --node[pos=0.5,above](){$\beta$} (v3);
    \draw[->] (v4) --node[pos=0.5,above](){$\delta$} (v3);
    \draw[->] (v4) --node[pos=0.5,above](){$\gamma$} (v5);
    \draw[->] (v6) --node[pos=0.5,above](){$\alpha$} (v5);
    \draw[->] (v7) --node[pos=0.5,above](){$\gamma$} (v6);
  \end{scope}
\end{tikzpicture} \big)$. We find that $\projcover(\Str w)$ is given in Figure \ref{subfig:projective-cover-from-line-of-patches}. This module is represented by the line of patches in Figure \ref{subfig:line-of-patches-from-strip}, whose top row is a strip representing $w$.
\begin{thesisfigure}
  \centering
  
  \subcaptionbox{%
    \label{subfig:line-of-patches-from-strip}
    \textit{Line of patches from a strip.}
  }[%
    0.95\textwidth
  ]{%
    \begin{tikzpicture}[scale=0.7]
      % Cells of 1st row (top row)
      \oddleftedgecell{(0,1)}{1 0}{};
      \evencell{(2,1)}{1 1}{};
      \oddcell{(4,1)}{1 2}{};
      \evencell{(6,1)}{1 3}{};
      \oddcell{(8,1)}{1 4}{};
      \evencell{(10,1)}{1 5}{};
      \oddcell{(12,1)}{1 6}{};
      \evencell{(14,1)}{1 7}{};
      \oddcell{(16,1)}{1 8}{};
      \evenrightedgecell{(18,1)}{1 9}{};
      
      % Cells of 2nd row (bottom row)
      \evenleftedgecell{(0,-1)}{2 0}{};
      \oddcell{(2,-1)}{2 1}{};
      \evencell{(4,-1)}{2 2}{};
      \oddcell{(6,-1)}{2 3}{};
      \evencell{(8,-1)}{2 4}{};
      \oddcell{(10,-1)}{2 5}{};
      \evencell{(12,-1)}{2 6}{};
      \oddcell{(14,-1)}{2 7}{};
      \evencell{(16,-1)}{2 8}{};
      \oddrightedgecell{(18,-1)}{2 9}{};
      
      \begin{scope}[inner sep=0, minimum size=0, node distance=7]
        % Entries of 1st row (top row)
        \oddtrivbdysyll{1 2}{}{$e$}{}{}{};
        \evenintsyll{1 3}{}{$\alpha\beta$}{};
        \oddintsyll{1 4}{}{$\delta$}{};
        \evenintsyll{1 5}{}{$\gamma$}{};
        \oddintsyll{1 6}{}{$\gamma\alpha$}{};
        \eventrivbdysyll{1 7}{}{$e'$}{}{}{};
        
        % Entries of the 2nd row (bottom row)
        \evenintsyll{2 2}{}{$\gamma\alpha$}{};
        \oddintsyll{2 3}{}{$ \gamma $}{};
        \evenbdysyll{2 4}{}{ $\delta^2$ }{}{}{};
        \oddbdysyll{2 5}{}{ $\alpha$}{}{}{};
        % (Nothing in column 6)
        \oddbdysyll{2 7}{}{$ \delta^2 $}{}{}{};
      \end{scope}
      
    \end{tikzpicture}
  }
  
  \vspace{2ex}
  
  \subcaptionbox{%
    \label{subfig:projective-cover-from-line-of-patches}
    \textit{Corresponding projective module.}
  }[%
    0.95\textwidth
  ]{%
    \begin{tikzpicture}[scale=0.6]
      % Give "anchoring" nodes for nonzero summands
      \draw (-7,0) node[](1){};
      \draw (0,0) node[](2){};
      \draw (7,0) node[](3){};
      
      % Draw nodes of lefthand projective
      \draw (1) ++(2,6) node[](1e1){$e_1A$};
      \draw (1) ++(0,4) node[](1b){$\beta A$};
      \draw (1) ++(0,2) node[](1bc){$\beta\gamma A$};
      \draw (1) ++(2,0) node[](1bca){$\beta\gamma\alpha A$};
      \draw (1) ++(4,4) node[](1a){$\alpha A$};
      \draw (1) ++(4,2) node[](1ab){$\alpha\beta A$};
      
      % Draw nodes of middle projective
      \draw (2) ++(2,6) node[](2e2){$e_2 A$};
      \draw (2) ++(0,4) node[](2d){$\delta A$};
      \draw (2) ++(0,2) node[](2dd){$\delta^2 A$};
      \draw (2) ++(0,0) node[](2ddd){$\delta^3 A$};
      \draw (2) ++(4,4) node[](2c){$\gamma A$};
      \draw (2) ++(4,2) node[](2ca){$\gamma\alpha A$};
      
      % Draw nodes of righthand projective
      \draw (3) ++(2,6) node[](3e2){$e_2 A$};
      \draw (3) ++(0,4) node[](3c){$\gamma A$};
      \draw (3) ++(0,2) node[](3ca){$\gamma\alpha A$};
      \draw (3) ++(4,4) node[](3d){$\delta A$};
      \draw (3) ++(4,2) node[](3dd){$\delta^2 A$};
      \draw (3) ++(4,0) node[](3ddd){$\delta^3 A$};
      
      % Draw action lines on lefthand projective
      \draw (1e1) --node[pos=0.5,right](){$\beta$} (1b);
      \draw (1b) --node[pos=0.5,right](){$\gamma$} (1bc);
      \draw (1bc) --node[pos=0.5,right](){$\alpha$} (1bca);
      \draw (1e1) --node[pos=0.5,right](){$\alpha$} (1a);
      \draw (1a) --node[pos=0.5,right](){$\beta$} (1ab);
      \draw (1ab) --node[pos=0.5,right](){$\gamma$} (1bca);
      
      % Draw action lines on middle projective
      \draw (2e2) --node[pos=0.5,right](){$\delta$} (2d);
      \draw (2d) --node[pos=0.5,right](){$\delta$} (2dd);
      \draw (2dd) --node[pos=0.5,right](){$\delta$} (2ddd);
      \draw (2e2) --node[pos=0.5,right](){$\gamma$} (2c);
      \draw (2c) --node[pos=0.5,right](){$\alpha$} (2ca);
      
      % Draw action lines on righthand projective
      \draw (3e2) --node[pos=0.5,right](){$\delta$} (3d);
      \draw (3d) --node[pos=0.5,right](){$\delta$} (3dd);
      \draw (3dd) --node[pos=0.5,right](){$\delta$} (3ddd);
      \draw (3e2) --node[pos=0.5,right](){$\gamma$} (3c);
      \draw (3c) --node[pos=0.5,right](){$\alpha$} (3ca);
    \end{tikzpicture}
  }
  
  \caption[%
    Patch cover example
  ]{%
    \label{fig:patch-cover-example}
    \textit{An example of a patch cover.} The top row of (\subref{subfig:line-of-patches-from-strip}) is a strip representing the string module $\Str w$
    . This uniquely determines the line of patches in the rest of (\subref{subfig:line-of-patches-from-strip}), a line which models the projective module in (\subref{subfig:projective-cover-from-line-of-patches}).
  }
    
\end{thesisfigure}

\nlpass{Patch covers} We consequently define the \define{patch cover $\projcover w$} of a strip $w$ to be the corresponding line of patches, given in Corollary \ref{cor:patch-covers-of-strips-well-defined}. We emphasise that this line is ``infinite in both directions'', featuring blank patches on either end as necessary.

\subsection{Syzygy algorithm}

We know that any strip $w$ has a patch cover $\projcover w$, the bottom row of which either does or does not contain virtual syllables. These cases respectively correspond to whether $w$ does or does not represent the socle-quotient of a pin module. In the latter case, the bottom row of $\projcover w$ specifies the support of several strips $w_r$ (indexed by $r$) such that $ \syzygy^1(\Str w) = \bigoplus_r \Str( w_r )$. In the former case, the bottom row of $\projcover w$ does not specify a strip directly. Rather, it forms the top row of a second line of patches, the bottom row of which is a strip representing $\syzygy^1(\Str w)$.

%In this subsection we prove two statements: first, that the bottom row of $\projcover w$ specifies the supports of several strips $\{w_r\}_r$; second, that we have $\bigoplus_r \Str( w_r ) = \syzygy^1( \Str w)$.

We advise the reader that strips representing socle-quotients of pin modules are exceptional and often require special attention. They are either the subject of bespoke claims or exempted from general ones; this distinction will always be made clear when necessary.

We also comment that, by excluding socle-quotient strips $w$, we ensure that virtual syllables will not enter into discussion.

\centresubsec{Syzygy of a strip that does not represent the socle-quotient of a pin module}

\nlpass{Segments} Let $w$ be any strip (socle-quotients of pin modules allowed) and let $\projcover w$ its patch cover. Recall that the columns of $w$ and of $\projcover w$ are indexed by $\integers$. Any cell belongs to a single patch and any cell has a unique peak neighbour and a unique valley neighbour. Neighbours always lie in the same row.

Consider the equivalence relation \define{$\sim$} defined on cells of the bottom row of $\projcover w$ and generated by the following properties: \begin{enumerate*}
    \item a cell is always $\sim$-related to its peak neighbour; and
    \item a cell is $\sim$-related to its valley neighbour unless they belong to a patch associated to a projective string module.
\end{enumerate*}

We call the $\sim$-equivalence classes the \define{segments of $\projcover w$}.

\begin{proposition}\label{prop:syzygy-algorithm-for-all-but-soc-quots}
Suppose the strip $w$ does not represent the socle-quotient of a pin module.
\begin{enumerate}
    \item\label{subprop:segments-of-bottom-row-of-patch-cover-are-strips} The restriction of $\projcover w$ to any segment of $w$ gives the support of a valid strip $w_r$ .
    \item\label{subprop:segments-of-bottom-row-of-patch-cover-are-syzygies} Those strips $w_r$ (indexed by $r$) represent the direct summands of $\syzygy^1(\Str w)$; that is, $ \syzygy^1( \Str w ) \isom \bigoplus_r \Str(w_r) $.
\end{enumerate}
\end{proposition}

\proofof{\ref{subprop:segments-of-bottom-row-of-patch-cover-are-strips}} Fix a segment of $\projcover w$ and let $v$ be the corresponding restriction. To prove that $v$, extended by blank values where necessary, is a strip we must show that any putative peak is indeed a peak, that any putative valley is indeed a valley, and that $v$ satisfies the concavity property. Here, by \emph{putative} peak or valley, we mean a neighbouring pair of syllables with orientations
\iltikzpicture{
\peakrep{(0,0)}{pk1}{$-1$}{$+1$}
}
or
\iltikzpicture{
\valleyrep{(0,0)}{vl1}{$+1$}{$-1$}
}.

We first consider a putative peak
\iltikzpicture{
\peakrep{(0,0)}{pk1}{$\bdq$}{$\bdq'$}
}.
We certainly know that in the top row above it is a valley
\iltikzpicture{
\valleyrep{(0,0)}{vl1}{$\bdp$}{$\bdp'$}
} so there are three possibilities in the terminology of Figure \ref{fig:valley-compatible-pairs}: it is either an interior valley, a boundary valley or a blank valley. If it is an interior valley, then by Lemma \ref{lem:properties-of-patches}\ref{sublem:int-sylls-in-top-row-of-patches} we have $\bdq=\bdp\descent$ and $\bdq'=\bdp'\descent$ and by calculation $\source(\bdq)=\target(\bdp)$ and $\source(\bdq')=\target(\bdp')$, which together mean
$$
\source(\bdq)^\dagger=\target(\bdp)^\dagger=\target(\bdp')=\source(\bdq')\tcomma
$$
and so the putative peak truly is a peak, and an interior peak at that. If the putative peak is a boundary or blank valley then at least of one of $\bdp$ or $\bdp'$ is a blank syllable, which must lie in a blank patch by concavity of $w$; this forces at least one of $\bdq$ or $\bdq'$ to be blank and so, again, the syllables do form a peak.

Next we consider a potential valley
\iltikzpicture{
\valleyrep{(0,0)}{vl1}{$\bdq$}{$\bdq'$}
}
. Clearly this potential valley is the bottom row of a single patch; that patch must represent some projective module. If the patch is blank (and the projective module zero), its bottom row is blank and therefore a valley. If the patch represents a projective string module, then only one cell lies in the segment $v$. The entry (say $\bdq$) of that cell is $\bdp\descent$ for some syllable $\bdp$ satisfying $b_{\source(\bdp)}=1$, implying $\bdq$ is either blank or a boundary syllable and so the potential valley is either a blank or boundary valley. If the patch represents a pin module, then we deduce that validity of the valley from parts \ref{sublem:bottom-row-under-one-pin-boundary-is-boundary-valley} and \ref{sublem:bottom-row-of-generic-pin-patch-is-interior-valley} of Lemma \ref{lem:bottom-row-of-pin-patch}.

It only remains to prove concavity of $v$. Any positions outside the segment of $v$ have blank values so it suffices to consider those columns supported by the segment. Suppose for contradiction that there are columns indexed by integers $k<\ell<k'$ such that $v(k),v(k')$ are nonblank but $v(\ell)$ is blank.

We claim that, since $v(k)$ is nonblank, at least one of $w(k-1),w(k),w(k+1)$ is nonblank. If $w$ is blank in all three columns $k-1,k,k+1$ then the unique patch encompassing column $k$ is blank in the top row and so blank entirely, hence $v(k)$ is blank; our claim follows by contrapositive. In the same fashion we deduce that at least one of $w(k'-1),w(k'),w(k'+1)$ is nonblank.

We know that $k+1 \leq l \leq k'-1$ and so, by concavity of $w$, $w(\ell)$ is nonblank. But recall $v(\ell)$ is blank also. The only possibility is that $w(\ell)$ is some boundary syllable. This forces $w$ to be blank either on all columns to the left or to the right of $\ell$, by the concavity and valley properties of strips, which means that either the patch through column $k$ or that through column $k'$ is blank, hence at least of one of $v(k)$ or $v(k')$ must be blank. This is our desired contradiction.
\begin{thesisfigure}
  \centering
  \begin{tikzpicture}[scale=\standardscale]
    \patchrep{(0,0)}{$\bdp$}{$\bdp'$}{$\bdq$}{$\bdq'$};
  \end{tikzpicture}
  
  \caption[%
    Notation for syllables in a patch
  ]{%
    \label{fig:notation-for-syllables-in-patch} \textit{Notation for syllables in patch $X$.}
  }  
\end{thesisfigure}

\proofof{\ref{subprop:segments-of-bottom-row-of-patch-cover-are-syzygies}}
We know that the top row of $\projcover w$ is the strip $w$ and, thanks to part \ref{subprop:segments-of-bottom-row-of-patch-cover-are-strips}, that the bottom row also specifies a family of strips $w_r$ (indexed by $r$). Using Proposition \ref{prop:strips-give-string-graphs} we can turn these strips into string modules $\Str w$ and $\Str w_r$ with prescribed bases. The following argument will establish that there is a basis of the projective module $\projcover(\Str (w ))$ defined using the combinatorial data of the patch cover $\projcover w$ and admitting a partition into a top and bottom part, mirroring the top and bottom rows of $\projcover w$. The basis vectors of the bottom part correspond to the prescribed basis of $\bigoplus_r \Str w_r$ and span a submodule of $\projcover(\Str w )$, modulo which the basis vectors of the top part give rise to the prescribed basis of $\Str( w )$. The canonical inclusion and projection maps of this basis of $\projcover( \Str w )$ therefore yield a short exact sequence
$$
\begin{tikzcd}[sep=small,cramped]
0 \ar[r]
  \& \bigoplus\limits_r \Str w_r \ar[r]
    \& \projcover( \Str w ) \ar[r]
      \& \Str w \ar[r]
        \& 0\tcomma
\end{tikzcd}
$$
whence $\syzygy^1(\Str w) \isom \bigoplus_r \Str(w_r)$ as claimed. We advise the reader that the abstract argument pursued in our proof is perhaps more clearly understood through the concrete example that follows it.

For an initial basis of $\projcover( \Str w )$, let us take the disjoint union of the standard bases of its direct summands. We will alter this basis patch by patch.
\begin{thesisfigure}
  \centering
  
  \subcaptionbox{
    \label{subfig:syzygy-basis-of-proj-cover}
    Module $\projcover(\Str v_2)$ with basis $T \cup B$
  }[
    0.9\textwidth
  ]{
    \begin{tikzpicture}[scale=0.6]
      % Give "anchoring" nodes for nonzero summands
      \draw (-6,0) node[](1){};
      \draw (0,0) node[](2){};
      \draw (6,0) node[](3){};
      
      % Draw nodes of lefthand projective
      \draw (1) ++(2,6) node[](1e1){$e_1A$};
      \draw (1) ++(0,4) node[](1b){$\beta A$};
      \draw (1) ++(0,2) node[](1bc){$\beta\gamma A$};
      \draw (1) ++(2,0) node[](1bca){$\beta\gamma\alpha A$};
      \draw (1) ++(4,4) node[](1a){$\alpha A$};
      \draw (1) ++(4,2) node[](1ab){$\alpha\beta A$};
      
      % Draw nodes of middle projective
      \draw (2) ++(2,6) node[](2e2){$e_2 A$};
      \draw (2) ++(0,4) node[](2d){$\delta A$};
      \draw (2) ++(0,2) node[](2dd){$\delta^2 A$};
      \draw (2) ++(0,0) node[](2ddd){$\delta^3 A$};
      \draw (2) ++(4,4) node[](2c){$\gamma A$};
      \draw (2) ++(4,2) node[](2ca){$\gamma\alpha A$};
      
      % Draw nodes of righthand projective
      \draw (3) ++(2,6) node[](3e2){$e_2 A$};
      \draw (3) ++(0,4) node[](3c){$\gamma A$};
      \draw (3) ++(0,2) node[](3ca){$\gamma\alpha A$};
      \draw (3) ++(4,4) node[](3d){$\delta A$};
      \draw (3) ++(4,2) node[](3dd){$\delta^2 A$};
      \draw (3) ++(4,0) node[](3ddd){$\delta^3 A$};
      
      % Draw action lines on lefthand projective
      % Include node (lb) "left beta"
      \draw (1e1) --node[pos=0.5,right](){$\beta$} node[pos=0.5, inner sep=0, minimum size=0](lb){} (1b);
      \draw (1b) --node[pos=0.5,right](){$\gamma$} (1bc);
      \draw (1bc) --node[pos=0.5,right](){$\alpha$} (1bca);
      \draw (1e1) --node[pos=0.5,right](){$\alpha$} (1a);
      \draw (1a) --node[pos=0.5,right](){$\beta$} (1ab);
      \draw (1ab) --node[pos=0.5,right](){$\gamma$} (1bca);
      
      % Draw action lines on middle projective
      \draw (2e2) --node[pos=0.5,right](){$\delta$} (2d);
      \draw (2d) --node[pos=0.5,right](){$\delta$} (2dd);
      \draw (2dd) --node[pos=0.5,right](){$\delta$} (2ddd);
      \draw (2e2) --node[pos=0.5,right](){$\gamma$} (2c);
      \draw (2c) --node[pos=0.5,right](){$\alpha$} (2ca);
      
      % Draw action lines on righthand projective
      % Include node (rd) "right delta"
      \draw (3e2) --node[pos=0.5,right](){$\delta$} node[pos=0.5, inner sep=0, minimum size=0](rd){} (3d);
      \draw (3d) --node[pos=0.5,right](){$\delta$} (3dd);
      \draw (3dd) --node[pos=0.5,right](){$\delta$} (3ddd);
      \draw (3e2) --node[pos=0.5,right](){$\gamma$} (3c);
      \draw (3c) --node[pos=0.5,right](){$\alpha$} (3ca);
      
      \begin{scope}[inner sep=0, minimum size=0]
        \newcommand{\distance}{0.7}
        \draw (1ab)+(135:\distance) node[](t1){};
        \draw (1ab)+(225:\distance) node[](t2){};
        \draw (1ab)+(315:\distance) node[](t3){};
        \draw (2d)+(315:\distance) node[](t4){};
        \draw (2d)+(45:\distance) node[](t5){};
        \draw (2d)+(135:\distance) node[](t6){};
        
        \draw (2c)+(225:\distance) node[](t7){};
        \draw (3ca)+(225:\distance) node[](t8){};
        \draw (3ca)+(315:\distance) node[](t9){};
        \draw (3ca)+(45:\distance) node[](t10){};
        \draw (2c)+(45:\distance) node[](t11){};
        \draw (2c)+(135:\distance) node[](t12){};
        
        \draw[dotted] (t1) to[out=225,in=135] (t2) to[out=315,in=225] (t3) to[out=45,in=225] (t4) to[out=45,in=315] (t5) to[out=135,in=45] (t6) to[out=225, in=45] (t1);
        \draw[dotted] (t7) to[out=315,in=135] (t8) to[out=315,in=225] (t9) to[out=45,in=315] (t10) to[out=135,in=315] (t11) to[out=135,in=45] (t12) to[out=225,in=135] (t7);
        
        % Draw extra nodes (ll) "left" and (rr) "right"
        \draw (lb)+(-3,0.75) node[](ll){};
        \draw (rd)+(3,0.75) node[](rr){};
        
        \draw[dotted] (lb) to[out=315,in=135] (t1);
        \draw[dotted] (t5) to (t12);
        \draw[dotted] (t10) to (rd);
        \draw[dotted] (lb) to[out=135,in=0] (ll);
        \draw[dotted] (rd) to[out=45,in=180] (rr);
        
        \draw (ll)+(0.4,0.6) node[](){$T$};
        \draw (ll)+(0.4,-0.6) node[](){$B$};
      \end{scope}
    \end{tikzpicture}
  }

  \
  
  \subcaptionbox{
    \label{subfig:syzygy-basis-of-quotient}
    Module $\Str(v_2)$ with basis $T$
  }[
    0.9\textwidth
  ]{
    \begin{tikzpicture}[scale=\standardscale]
      % Draw nodes of module
      \draw (0,4) node[](v1){$e_1 A$};
      \draw (2,2) node[](v2){$\alpha A$};
      \draw (4,0) node[](v3){$\alpha \beta A$};
      \draw (6,2) node[](v4){$e_2 A$};
      \draw (8,0) node[](v5){$\gamma A$};
      \draw (10,2) node[](v6){$\gamma A$};
      \draw (12,4) node[](v7){$e_2 A$};
      
      % Draw action lines of module
      \draw (v1) --node[pos=0.5,right](){$\alpha$} (v2);
      \draw (v2) --node[pos=0.5,right](){$\beta$} (v3);
      \draw (v3) --node[pos=0.5,right](){$\delta$} (v4);
      \draw (v4) --node[pos=0.5,right](){$\gamma$} (v5);
      \draw (v5) --node[pos=0.5,right](){$\alpha$} (v6);
      \draw (v6) --node[pos=0.5,right](){$\gamma$} (v7);
    \end{tikzpicture}

  }

  \

  \subcaptionbox{
    \label{subfig:syzygy-basis-of-submodule}
    Module $\syzygy^1(\Str v_2)$ with basis $B$
  }[
    0.9\textwidth
  ]{
    \begin{tikzpicture}[scale=\standardscale]
      % Draw anchoring nodes for modules
      \draw (0,0) node[](l){};
      \draw (13,0) node[](c){};
      \draw (20,0) node[](r){};
      
      % Draw nodes of lefthand module
      \draw (l)+(0,6) node[](lv1){$\beta A$};
      \draw (l)+(2,4) node[](lv2){$\beta\gamma A$};
      \draw (l)+(4,2) node[](lv3){$\beta\gamma\alpha A$};
      \draw (l)+(6,4) node[](lv4){$(\alpha\beta A)-(\delta A)$};
      \draw (l)+(8,2) node[](lv5){$\delta^2 A$};
      \draw (l)+(10,0) node[](lv6){$\delta^3 A$};
      
      % Draw nodes of middle module
      \draw (c)+(0,2) node[](cv1){$\gamma\alpha A$};
      \draw (c)+(2,4) node[](cv2){$(\gamma A)-(\gamma\alpha A)$};
      
      % Draw nodes of righthand module
      \draw (r)+(0,6) node[](rv1){$\delta A$};
      \draw (r)+(2,4) node[](rv2){$\delta^2 A$};
      \draw (r)+(4,2) node[](rv3){$\delta^3 A$};
      
      % Draw action lines of lefthand module
      \draw (lv1) --node[pos=0.5,right](){$\gamma$} (lv2);
      \draw (lv2) --node[pos=0.5,right](){$\alpha$} (lv3);
      \draw (lv3) --node[pos=0.5,right](){$\gamma$} (lv4);
      \draw (lv4) --node[pos=0.5,right](){$\delta$} (lv5);
      \draw (lv5) --node[pos=0.5,right](){$\delta$} (lv6);
      
      % Draw action line of middle module
      \draw (cv1) --node[pos=0.5,right](){$\alpha$} (cv2);
      
      % Draw action line of righthand module
      \draw (rv1) --node[pos=0.5,right](){$\delta$} (rv2);
      \draw (rv2) --node[pos=0.5,right](){$\delta$} (rv3);
    \end{tikzpicture}
  }
  
  \caption[
    Example of syzygy-friendly basis
  ]{
    \label{fig:syzygy-friendly-basis}
    \textit{Example of syzygy-friendly basis.}
  }
\end{thesisfigure}

Let $X$ be a patch in $\projcover w$; if $X$ is an amended patch, then without loss of generality replace $X$ by the original patch from which it was amended. Write $P_X$ for the indecomposable projective associated to $X$. Denote the syllables of $X$ by $\bdp,\bdp',\bdq,\bdq'$ as in Figure \ref{fig:notation-for-syllables-in-patch}. Write $p$ for the \define{compression} $\big(\ilslpath{i}{\ell+\ep} \big)$ of the syllable $\bdp = \big(\ilsyll{i}{\ell}{\ep}\big)$ and, similarly, $p'$ for the compression of $\bdp'$. By definition these compressions are both nonstationary paths in $\OQ$.

The paths comparable with $p,p'$ in the prefix order represent basis vectors of $P_X$. We will divide the basis into two parts: an upper part and a lower part. In the upper part, we place all basis vectors represented by strict prefixes of $p$ and $p'$. In the lower part, we place all basic vectors represented by paths having $p$ or $p'$ as strict prefix. This leaves undetermined the fate of vectors represented by $p$ or $p'$. 

If $\bdp$ is boundary, then we place the vector represented by $p$ in the lower part; likewise with $p'$ if $\bdp'$ is boundary. If $\bdp$ is interior, then it necessarily appears in a valley \iltikzpicture{\valleyrep{(0,0)}{vl1}{$\bdu$}{$\bdp$}} next to another interior syllable $\bdu$, whose compression we denote $u$. We know $u$ represents some basis vector of a neighbouring indecomposable direct summand of $\projcover(\Str w)$. We replace $\{u,p\}$ in the basis of $\projcover(\Str w)$ by $\{u, u-p\}$, and then place $u$ in the upper part and $u-p$ in the lower part. Symmetrically, if $\bdp'$ is interior then it occurs in a valley \iltikzpicture{\valleyrep{(0,0)}{vl1}{$\bdp'$}{$\bdu'$}} neighbouring an interior syllable $\bdu'$, and we replace $\{p',u'\}$ in the basis of $\projcover(\Str w)$ by $\{p', p'-u'\}$, where $u'$ denotes the compression of $\bdu'$.

Performing these changes for all patches $X$ yields a well-defined, bipartitioned basis for $\projcover(\Str w)$.

The basis vectors of the lower part span a submodule of $\projcover (\Str w)$. Observe that the basis vectors represented by the compression of a boundary syllable, or by differences $u-p$ or $p'-u'$ of compressions of interior syllables, span the head, and each $\alpha \in Q_1$ annihilates at least one component of any such difference. Any vector not so annihilated is mapped to (some linear multiple of) another vector in the lower part. Remaining basis vectors in the lower part, represented by strict superpaths of the paths $p,p'$, are annihilated by all $A$-arrows save at most one. We deduce that the lower part is just the standard basis for $\bigoplus_r \Str w_r$, up to rescaling.

Modulo the lower basis vectors, compressions $p,u$ or $u,p'$ of neighbouring interior syllables are identified. We therefore retrieve the standard basis for $\Str w$. The result follows.\qed

\nrex We take up the calculation from Figures \ref{subfig:syzygy-of-v2}, \ref{fig:branching} and \ref{fig:patch-cover-example} that
$$
\syzygy^1(\Str v_2)
=
\Str\big(
  \begin{tikzpicture}[scale=0.75]
    \begin{scope}[minimum size=0, inner sep=1]
      % Vertices of string graph
      \draw (0,0) node[](v1){$ 2$};
      \draw (1,0) node[](v2){$ 1$};
      \draw (2,0) node[](v3){$ 1$};
      \draw (3,0) node[](v4){$ 2$};
      \draw (4,0) node[](v5){$ 2$};
      \draw (5,0) node[](v6){$ 2$};
      %
      % Arrows of string graph
      \draw[->] (v1) --node[pos=0.5,above](){$\gamma$} (v2);
      \draw[->] (v2) --node[pos=0.5,above](){$\alpha$} (v3);
      \draw[->] (v4) --node[pos=0.5,above](){$\gamma$} (v3);
      \draw[->] (v4) --node[pos=0.5,above](){$\delta$} (v5);
      \draw[->] (v5) --node[pos=0.5,above](){$\delta$} (v6);
    \end{scope}
  \end{tikzpicture}
\big)
~\oplus~
\Str\big(
  \begin{tikzpicture}[scale=0.75]
    \begin{scope}[minimum size=0, inner sep=1]
      % Vertices of string graph
      \draw (0,0) node[](v1){$ 1$};
      \draw (1,0) node[](v2){$ 1$};
      %
      % Arrows of string graph
      \draw[->] (v2) --node[pos=0.5,above](){$\alpha$} (v1);
    \end{scope}
  \end{tikzpicture}
\big)
~\oplus~
\Str\big(
  \begin{tikzpicture}[scale=0.75]
    \begin{scope}[minimum size=0, inner sep=1]
      % Vertices of string graph
      \draw (0,0) node[](v1){$ 2$};
      \draw (1,0) node[](v2){$ 2$};
      \draw (2,0) node[](v3){$ 2$};
      %
      % Arrows of string graph
      \draw[->] (v2) --node[pos=0.5,above](){$\delta$} (v1);
      \draw[->] (v3) --node[pos=0.5,above](){$\delta$} (v2);
    \end{scope}
  \end{tikzpicture}
\big)\tstop
$$
The patch cover $\projcover v_2$ designates a basis $T \cup B$ of $\projcover( \Str v_2 )$, as described above and illustrated in Figure \ref{subfig:syzygy-basis-of-proj-cover}. The basis $B$ spans the submodule $\syzygy^1(\Str v_2)$, while the residue of $T$ modulo $B$ is a basis of $\Str(v_2)$. These are seen in Figures \ref{subfig:syzygy-basis-of-quotient} and \ref{subfig:syzygy-basis-of-submodule}.

\nrmk Thanks to the careful phrasing of the preceding proposition and proof, the results hold for infinite-dimensional string modules too.
\begin{thesisfigure}
  \centering
  
  \begin{tikzpicture}[scale=0.60]
    % Top row of cells
    \evenleftedgecell{(-11,2)}{e -11 2};
    \oddcell{(-9,2)}{o -9 2}{}{};
    \evencell{(-7,2)}{e -7 2}{}{};
    \oddcell{(-5,2)}{o -5 2}{}{};
    \evencell{(-3,2)}{e -3 2}{}{};
    \oddcell{(-1,2)}{o -1 2}{}{};
    \evencell{(1,2)}{e 1 2}{}{};
    \oddcell{(3,2)}{o 3 2}{}{};
    \evencell{(5,2)}{e 5 2}{}{};
    \oddcell{(7,2)}{o 7 2}{}{};
    \evencell{(9,2)}{e 9 2}{}{};
    \oddrightedgecell{(11,2)}{o 11 2};
    
    % Bottom row of cells
    \oddleftedgecell{(-11,0)}{o -11 0};
    \evencell{(-9,0)}{e -9 0}{}{};
    \oddcell{(-7,0)}{o -7 0}{}{};
    \evencell{(-5,0)}{e -5 0}{}{};
    \oddcell{(-3,0)}{o -3 0}{}{};
    \evencell{(-1,0)}{e -1 0}{}{};
    \oddcell{(1,0)}{o 1 0}{}{};
    \evencell{(3,0)}{e 3 0}{}{};
    \oddcell{(5,0)}{o 5 0}{}{};
    \evencell{(7,0)}{e 7 0}{}{};
    \oddcell{(9,0)}{o 9 0}{}{};
    \evenrightedgecell{(11,0)}{e 11 0};
    
    \begin{scope}[node distance=8]
      % Syllables in the top row
      \oddintsyll{o -9 2}{}{$\gamma\alpha$}{};
      \evenintsyll{e -7 2}{}{$\delta^2$}{};
      \oddintsyll{o -5 2}{}{$\alpha\beta$}{};
      \evenintsyll{e -3 2}{}{$\beta\gamma$}{};
      \oddintsyll{o -1 2}{}{$\gamma\alpha$}{};
      \evenintsyll{e 1 2}{}{$\delta^2$}{};
      \oddintsyll{o 3 2}{}{$\alpha\beta$}{};
      \evenintsyll{e 5 2}{}{$\beta\gamma$}{};
      \oddintsyll{o 7 2}{}{$\gamma\alpha$}{};
      \evenintsyll{e 9 2}{}{$\delta^2$}{};
      
      % Syllables in the bottom row
      \evenbdysyll{e -9 0}{}{$\beta$}{}{}{};
      \oddbdysyll{o -7 0}{}{$\delta$}{}{}{};
      \evenintsyll{e -5 0}{}{$\gamma$}{};
      \oddintsyll{o -3 0}{}{$\alpha$}{};
      \evenbdysyll{e -1 0}{}{$\beta$}{}{}{};
      \oddbdysyll{o 1 0}{}{$\delta$}{}{}{};
      \evenintsyll{e 3 0}{}{$\gamma$}{};
      \oddintsyll{o 5 0}{}{$\alpha$}{};
      \evenbdysyll{e 7 0}{}{$\beta$}{}{}{};
      \oddbdysyll{o 9 0}{}{$\delta$}{}{}{};
    \end{scope}
  \end{tikzpicture}
  
  \caption[
    Syzygy of an infinite-dimensional string module
  ]{
    \label{fig:inf-dim-syz}
    \textit{Syzygy of an infinite-dimensional string module.}
  }
\end{thesisfigure}

\nrex Consider the infinite-dimensional string module
$$
\Str
\big(
  \begin{tikzpicture}[scale=0.75]
    \begin{scope}[minimum size=0, inner sep=1]
      % Vertices of string graph
      \draw (-1,0) node[](v0){$\cdots$};
      \draw (0,0) node[](v1){$ 1$};
      \draw (1,0) node[](v2){$ 1$};
      \draw (2,0) node[](v3){$ 2$};
      \draw (3,0) node[](v4){$ 2$};
      \draw (4,0) node[](v5){$ 2$};
      \draw (5,0) node[](v6){$ 1$};
      \draw (6,0) node[](v7){$ 1$};
      \draw (7,0) node[](v8){$ 2$};
      \draw (8,0) node[](v9){$ 1$};
      \draw (9,0) node[](v10){$ 1$};
      \draw (10,0) node[](v11){$ 2$};
      \draw (11,0) node[](v12){$ 2$};
      \draw (12,0) node[](v13){$ 2$};
      \draw (13,0) node[](v14){$ 1$};
      \draw (14,0) node[](v15){$ 1$};
      \draw (15,0) node[](v16){$ 2$};
      \draw (16,0) node[](v17){$ 1$};
      \draw (17,0) node[](v18){$ \cdots$};
      %
      % Arrows of string graph
      \newcommand{\stringarrow}[3]{
        \draw[->] (#1) --node[pos=0.5,above](){$ #3 $} (#2);
        }
      \stringarrow{v0}{v1}{}
      \stringarrow{v2}{v1}{\alpha}
      \stringarrow{v3}{v2}{\gamma}
      \stringarrow{v3}{v4}{\delta}
      \stringarrow{v4}{v5}{\delta}
      \stringarrow{v6}{v5}{\beta}
      \stringarrow{v7}{v6}{\alpha}
      \stringarrow{v7}{v8}{\beta}
      \stringarrow{v8}{v9}{\gamma}
      \stringarrow{v10}{v9}{\alpha}
      \stringarrow{v11}{v10}{\gamma}
      \stringarrow{v11}{v12}{\delta}
      \stringarrow{v12}{v13}{\delta}
      \stringarrow{v14}{v13}{\beta}
      \stringarrow{v15}{v14}{\alpha}
      \stringarrow{v15}{v16}{\beta}
      \stringarrow{v16}{v17}{\gamma}
      \stringarrow{v18}{v17}{}
    \end{scope}
  \end{tikzpicture}
\big)\tcomma
$$
whose string graph is an infinite repetition of $\big(\gamma\alpha\big)^{-1}\big(\delta^2\big)\big(\alpha\beta\big)^{-1}\big(\beta\gamma\big)$. Its syzygy is a direct sum of countably many copies of
$
\Str
\big(
  \begin{tikzpicture}[scale=0.75]
    \begin{scope}[minimum size=0, inner sep=1]
      % Vertices of string graph
      \draw (0,0) node[](v1){$ 2$};
      \draw (1,0) node[](v2){$ 2$};
      \draw (2,0) node[](v3){$ 1$};
      \draw (3,0) node[](v4){$ 1$};
      \draw (4,0) node[](v5){$ 2$};
      %
      % Arrows of string graph
      \newcommand{\stringarrow}[3]{
        \draw[->] (#1) --node[pos=0.5,above](){$ #3 $} (#2);
        }
      \stringarrow{v2}{v1}{\delta}
      \stringarrow{v2}{v3}{\gamma}
      \stringarrow{v4}{v3}{\alpha}
      \stringarrow{v4}{v5}{\beta}
    \end{scope}
  \end{tikzpicture}
\big)
$. This can be calculated using our formalism, as in the Figure \ref{fig:inf-dim-syz}.

\centresubsec{Syzygy of a strip that does represent the socle-quotient of a pin module}

\begin{proposition}\label{prop:syzygy-algorithm-for-soc-quots}
If the strip $w$ does represent the socle-quotient $P/\soc P$ of a pin module $P$, then the bottom row of $\projcover w$ has a single segment that is blank except for two virtual patches.

Moreover, the bottom row of $\projcover w$ is the top row of a new line of (virtual and blank) patches whose bottom row is a strip representing the simple module $\soc P = \syzygy^1(P/\soc P)$.
\end{proposition}

\proof This holds by construction. It is perhaps most easily articulated through Figure \ref{fig:what-happens-under-a-socle-quotient-of-pin}.\qed
\begin{thesisfigure}
  \centering
  
  \begin{tikzpicture}[scale=0.60]
    % Top row of cells
    \evenleftedgecell{(-11,2)}{e -11 2};
    \oddcell{(-9,2)}{o -9 2}{}{};
    \evencell{(-7,2)}{e -7 2}{}{};
    \oddcell{(-5,2)}{o -5 2}{}{};
    \evencell{(-3,2)}{e -3 2}{}{};
    \oddcell{(-1,2)}{o -1 2}{}{};
    \evencell{(1,2)}{e 1 2}{}{};
    \oddcell{(3,2)}{o 3 2}{}{};
    \evencell{(5,2)}{e 5 2}{}{};
    \oddcell{(7,2)}{o 7 2}{}{};
    \evencell{(9,2)}{e 9 2}{}{};
    \oddrightedgecell{(11,2)}{o 11 2};
    
    % Middle row of cells
    \oddleftedgecell{(-11,0)}{o -11 0};
    \evencell{(-9,0)}{e -9 0}{}{};
    \oddcell{(-7,0)}{o -7 0}{}{};
    \evencell{(-5,0)}{e -5 0}{}{};
    \oddcell{(-3,0)}{o -3 0}{}{};
    \evencell{(-1,0)}{e -1 0}{}{};
    \oddcell{(1,0)}{o 1 0}{}{};
    \evencell{(3,0)}{e 3 0}{}{};
    \oddcell{(5,0)}{o 5 0}{}{};
    \evencell{(7,0)}{e 7 0}{}{};
    \oddcell{(9,0)}{o 9 0}{}{};
    \evenrightedgecell{(11,0)}{e 11 0};
    
    % Bottom row of cells
    \evenleftedgecell{(-11,-2)}{e -11 -2};
    \oddcell{(-9,-2)}{o -9 -2}{}{};
    \evencell{(-7,-2)}{e -7 -2}{}{};
    \oddcell{(-5,-2)}{o -5 -2}{}{};
    \evencell{(-3,-2)}{e -3 -2}{}{};
    \oddcell{(-1,-2)}{o -1 -2}{}{};
    \evencell{(1,-2)}{e 1 -2}{}{};
    \oddcell{(3,-2)}{o 3 -2}{}{};
    \evencell{(5,-2)}{e 5 -2}{}{};
    \oddcell{(7,-2)}{o 7 -2}{}{};
    \evencell{(9,-2)}{e 9 -2}{}{};
    \oddrightedgecell{(11,-2)}{o 11 -2};
    
    % Pin boundary syllables in top row
    \begin{scope}[node distance=15]
      \oddbdysyll{o -1 2}{$i$}{}{}{}{};
      \evenbdysyll{e 1 2}{$\idag$}{}{}{}{};
    \end{scope}
    
    % Virtual syllables in the middle row
    \evenvirtualsyll{e -1 0}{$j$}{}{};
    \oddvirtualsyll{o 1 0}{$\jdag$}{}{};
    
    % Boundary stationary syllables in bottom row
    \begin{scope}[node distance=8]
      \oddbdysyll{o -1 -2}{$j$}{$0$}{}{}{};
      \evenbdysyll{e 1 -2}{$\jdag$}{$0$}{}{}{};
    \end{scope}
  \end{tikzpicture}
  
  ~
  
  \caption[%
    Syzygies of socle-quotients of pin modules
  ]{
    \label{fig:what-happens-under-a-socle-quotient-of-pin}
    \textit{Syzygies of socle-quotients of pin modules.} The top row is $w$, a strip representing a socle-quotient of a pin module. Its nonblank entries are pin-boundary syllables $\ilsyll{i}{a_i-1}{1}$ and $\ilsyll{\idag}{a_\idag-1}{1}$.
    
    The nonblank entries of the middle row are virtual syllables $\ilsyll{j}{0}{0}$ and $\ilsyll{\jdag}{0}{0}$, for $j=i-a_i$.
    
    The nonblank entries of the bottom row are stationary syllables $\ilsyll{j}{0}{1}$ and $\ilsyll{\jdag}{0}{1}$. The bottom row is a strip representing the simple module at a particular vertex, namely that represented by $j$ and $\jdag$.
    
    Each (appropriate) $2\times2$-region in the figure is a patch.
  }
\end{thesisfigure}

\centresubsec{Syzygy fabric}

\nlpass{Syzygy fabric} Fix some strip $w$ representing a string module. We can use the above algorithm to iteratively construct an array holding all the syzygy information about a given string module. We call this array the \define{syzygy fabric}.

The columns of the array are indexed by $\integers$. The rows are indexed by the vertices of a rooted tree $\mathcal{T}$. We construct the tree iteratively, specifying the contents of each row as we go.

The induction begins with the root vertex of the tree. We identify this row with the given strip $w$.

Iteratively, suppose that the $t$th row is populated with the strip $v$ (for some $t \in \mathcal{T}_0$). We stipulate that the child vertices of $t$ correspond to the segments of $\projcover v$, and that the entries of that child row equal the restriction of $\projcover v$ to the corresponding segment.

\nrex Let us consider the syzygy fabric of $\Str\big(
  \begin{tikzpicture}[scale=0.75]
  \begin{scope}[minimum size=0, inner sep=1]
    % Vertices of string graph
    \draw (0,0) node[](v1){$ 1$};
    \draw (1,0) node[](v2){$ 1$};
    \draw (2,0) node[](v3){$ 2$};
    \draw (3,0) node[](v4){$ 2$};
    \draw (4,0) node[](v5){$ 2$};
    \draw (5,0) node[](v6){$ 2$};
    %
    % Arrows of string graph
    \draw[->] (v1) --node[pos=0.5,above](){$\alpha$} (v2);
    \draw[->] (v2) --node[pos=0.5,above](){$\beta$} (v3);
    \draw[->] (v4) --node[pos=0.5,above](){$\delta$} (v3);
    \draw[->] (v5) --node[pos=0.5,above](){$\delta$} (v4);
    \draw[->] (v6) --node[pos=0.5,above](){$\delta$} (v5);
  \end{scope}
  \end{tikzpicture}
\big)$; this string module is also the indecomposable injective module $E_2$.
\begin{thesisfigure}
  \centering
  
  {
    \begin{tikzpicture}[scale=0.5]
      % Cells of 1st row (top row)
      \draw (0,1) node[](1){};
      
      \evenleftedgecell{(1) ++(0,0)}{1 0}{};
      \oddcell{(1) ++(2,0)}{1 1}{};
      \evencell{(1) ++(4,0)}{1 2}{};
      \oddcell{(1) ++(6,0)}{1 3}{};
      \evencell{(1) ++(8,0)}{1 4}{};
      \oddcell{(1) ++(10,0)}{1 5}{};
      \evencell{(1) ++(12,0)}{1 6}{};
      \oddcell{(1) ++(14,0)}{1 7}{};
      \evencell{(1) ++(16,0)}{1 8}{};
      \oddrightedgecell{(1) ++(18,0)}{1 9}{};
      
      % Entries of 1st row (top row)
      \begin{scope}[inner sep=0, minimum size=0, node distance=7]
        \oddtrivbdysyll{1 3}{}{$\source(\beta)$}{}{}{};
        \evenintsyll{1 4}{}{$\alpha\beta$}{};
        \oddintsyll{1 5}{}{$\delta^3$}{};
        \eventrivbdysyll{1 6}{}{$\source(\gamma)$}{}{}{};
      \end{scope}
      
      \newcommand{\rowdisplacement}{(0.8,-1.4)}
      \draw (0,-2) node[](2a){};
      \draw (2a) ++\rowdisplacement node[](2b){};
      \draw (2b) ++\rowdisplacement node[](2c){};
      
      % Cells of row 2a (backmost bottom row)
      \oddleftedgecell{(2a)++(0,0)}{2a 0}{};
      \evencell{(2a)++(2,0)}{2a 1}{};
      \oddcell{(2a)++(4,0)}{2a 2}{};
      \evencell{(2a)++(6,0)}{2a 3}{};
      \oddcell{(2a)++(8,0)}{2a 4}{};
      \evencell{(2a)++(10,0)}{2a 5}{};
      \oddcell{(2a)++(12,0)}{2a 6}{};
      \evencell{(2a)++(14,0)}{2a 7}{};
      \oddcell{(2a)++(16,0)}{2a 8}{};
      \evenrightedgecell{(2a)++(18,0)}{2a 9}{};
      
      % Entries of row 2a (backmost bottom row)
      \begin{scope}[inner sep=0, minimum size=0, node distance=7]
        \evenintsyll{2a 3}{}{$ \gamma\alpha $}{};
        \oddintsyll{2a 4}{}{ $\gamma$ }{};
        \eventrivbdysyll{2a 5}{}{$ \source(\delta) $}{}{}{};
      \end{scope}
      
      % Cells of row 2b (middle bottom row)
      \oddleftedgecell{(2b)++(0,0)}{2b 0}{};
      \evencell{(2b)++(2,0)}{2b 1}{};
      \oddcell{(2b)++(4,0)}{2b 2}{};
      \evencell{(2b)++(6,0)}{2b 3}{};
      \oddcell{(2b)++(8,0)}{2b 4}{};
      \evencell{(2b)++(10,0)}{2b 5}{};
      \oddcell{(2b)++(12,0)}{2b 6}{};
      \evencell{(2b)++(14,0)}{2b 7}{};
      \oddcell{(2b)++(16,0)}{2b 8}{};
      \evenrightedgecell{(2b)++(18,0)}{2b 9}{};

      % Entries of row 2b (middle bottom row)
      \begin{scope}[inner sep=0, minimum size=0, node distance=7]
        \oddbdysyll{2b 6}{}{ $\alpha$}{}{}{};
      \end{scope}
      
      % Draw row indices
      \newcommand{\rowindexdisplacement}{(-3,0)}
      
      \newcommand{\rowindex}[2]{
        \draw (#1) ++\rowindexdisplacement node[lightgray](#1 index){#2};
        % \filldraw[lightgray] (#1 index) circle (0.15);
        % 1: label of row node
      }
      \newcommand{\rowindexconnector}[2]{
        \draw[lightgray] (#1 index) -- (#2 index);
        % 1: Label of row node for source index
        % 2: Label of row node for target index
      }
      
      \rowindex{1}{(a)};
      \rowindex{2a}{(b)};
      \rowindex{2b}{(c)};
    %   \rowindex{2c};
      \rowindexconnector{1}{2a};
      \rowindexconnector{1}{2b};
    %   \rowindexconnector{1}{2c};
    \end{tikzpicture}
  }

  \vspace{4ex}

  {
    \begin{tikzpicture}[scale=0.5]
      % Cells of 1st row (top row)
      \draw (0,1) node[](1){};
      
      \evenleftedgecell{(1) ++(0,0)}{1 0}{};
      \oddcell{(1) ++(2,0)}{1 1}{};
      \evencell{(1) ++(4,0)}{1 2}{};
      \oddcell{(1) ++(6,0)}{1 3}{};
      \evencell{(1) ++(8,0)}{1 4}{};
      \oddcell{(1) ++(10,0)}{1 5}{};
      \evencell{(1) ++(12,0)}{1 6}{};
      \oddcell{(1) ++(14,0)}{1 7}{};
      \evencell{(1) ++(16,0)}{1 8}{};
      \oddrightedgecell{(1) ++(18,0)}{1 9}{};
      
      % Entries of 1st row (top row)
      \begin{scope}[inner sep=0, minimum size=0, node distance=7]
        % \oddtrivbdysyll{1 3}{}{$\source(\delta)$}{}{}{};
        \evenintsyll{1 4}{}{$\gamma\alpha$}{};
        \oddintsyll{1 5}{}{$\gamma$}{};
        \eventrivbdysyll{1 6}{}{$\source(\delta)$}{}{}{};
      \end{scope}
      
      \newcommand{\rowdisplacement}{(0.8,-1.4)}
      \draw (0,-2) node[](2a){};
      \draw (2a) ++\rowdisplacement node[](2b){};
      \draw (2b) ++\rowdisplacement node[](2c){};
      
      % Cells of row 2a (backmost bottom row)
      \oddleftedgecell{(2a)++(0,0)}{2a 0}{};
      \evencell{(2a)++(2,0)}{2a 1}{};
      \oddcell{(2a)++(4,0)}{2a 2}{};
      \evencell{(2a)++(6,0)}{2a 3}{};
      \oddcell{(2a)++(8,0)}{2a 4}{};
      \evencell{(2a)++(10,0)}{2a 5}{};
      \oddcell{(2a)++(12,0)}{2a 6}{};
      \evencell{(2a)++(14,0)}{2a 7}{};
      \oddcell{(2a)++(16,0)}{2a 8}{};
      \evenrightedgecell{(2a)++(18,0)}{2a 9}{};
      
      % Entries of row 2a (backmost bottom row)
      \begin{scope}[inner sep=0, minimum size=0, node distance=7]
        \evenbdysyll{2a 3}{}{$ \delta^2 $}{}{}{};
        % \oddintsyll{2a 4}{}{ $\gamma$ }{};
        % \eventrivbdysyll{2a 5}{}{$ \source(\delta) $}{}{}{};
      \end{scope}
      
      % Cells of row 2b (middle bottom row)
      \oddleftedgecell{(2b)++(0,0)}{2b 0}{};
      \evencell{(2b)++(2,0)}{2b 1}{};
      \oddcell{(2b)++(4,0)}{2b 2}{};
      \evencell{(2b)++(6,0)}{2b 3}{};
      \oddcell{(2b)++(8,0)}{2b 4}{};
      \evencell{(2b)++(10,0)}{2b 5}{};
      \oddcell{(2b)++(12,0)}{2b 6}{};
      \evencell{(2b)++(14,0)}{2b 7}{};
      \oddcell{(2b)++(16,0)}{2b 8}{};
      \evenrightedgecell{(2b)++(18,0)}{2b 9}{};

      % Entries of row 2b (middle bottom row)
      \begin{scope}[inner sep=0, minimum size=0, node distance=7]
        \oddtrivbdysyll{2b 4}{}{ $\source(\beta)$}{}{}{};
        \evenbdysyll{2b 5}{}{$ \alpha$}{}{}{};
      \end{scope}

      % Cells of row 2c (foremost bottom row)
      \oddleftedgecell{(2c)++(0,0)}{2c 0}{};
      \evencell{(2c)++(2,0)}{2c 1}{};
      \oddcell{(2c)++(4,0)}{2c 2}{};
      \evencell{(2c)++(6,0)}{2c 3}{};
      \oddcell{(2c)++(8,0)}{2c 4}{};
      \evencell{(2c)++(10,0)}{2c 5}{};
      \oddcell{(2c)++(12,0)}{2c 6}{};
      \evencell{(2c)++(14,0)}{2c 7}{};
      \oddcell{(2c)++(16,0)}{2c 8}{};
      \evenrightedgecell{(2c)++(18,0)}{2c 9}{};
      
      % Entries of row 2c (foremost bottom row)
      \begin{scope}[inner sep=0, minimum size=0, node distance=7]
        \oddbdysyll{2c 6}{}{$ \delta^2 $}{}{}{};
      \end{scope}
      
      % Draw row indices
      \newcommand{\rowindexdisplacement}{(-3,0)}
      
      \newcommand{\rowindex}[2]{
        \draw (#1) ++\rowindexdisplacement node[lightgray](#1 index){ #2};
        % \filldraw[lightgray] (#1 index) circle (0.15);
        % 1: label of row node
      }
      \newcommand{\rowindexconnector}[2]{
        \draw[lightgray] (#1 index) -- (#2 index);
        % 1: Label of row node for source index
        % 2: Label of row node for target index
      }
      
      \rowindex{1}{(b)};
      \rowindex{2a}{(d)};
      \rowindex{2b}{(c)};
      \rowindex{2c}{(d)};
      \rowindexconnector{1}{2a};
      \rowindexconnector{1}{2b};
      \rowindexconnector{1}{2c};
    \end{tikzpicture}
  }

  \vspace{4ex}

  {
    \begin{tikzpicture}[scale=0.5]
      % Cells of 1st row (top row)
      \draw (0,1) node[](1){};
      
      \oddleftedgecell{(1) ++(0,0)}{1 0}{};
      \evencell{(1) ++(2,0)}{1 1}{};
      \oddcell{(1) ++(4,0)}{1 2}{};
      \evencell{(1) ++(6,0)}{1 3}{};
      \oddcell{(1) ++(8,0)}{1 4}{};
      \evencell{(1) ++(10,0)}{1 5}{};
      \oddcell{(1) ++(12,0)}{1 6}{};
      \evencell{(1) ++(14,0)}{1 7}{};
      \oddcell{(1) ++(16,0)}{1 8}{};
      \evenrightedgecell{(1) ++(18,0)}{1 9}{};
      
      % Entries of 1st row (top row)
      \begin{scope}[inner sep=0, minimum size=0, node distance=7]
        \oddtrivbdysyll{1 4}{}{$\source(\beta)$}{}{}{};
        \evenbdysyll{1 5}{}{$\alpha$}{}{}{};
      \end{scope}
      
      \newcommand{\rowdisplacement}{(0.8,-1.4)}
      \draw (0,-2) node[](2a){};
      \draw (2a) ++\rowdisplacement node[](2b){};
    %   \draw (2b) ++\rowdisplacement node[](2c){};
      
      % Cells of row 2a (backmost bottom row)
      \evenleftedgecell{(2a)++(0,0)}{2a 0}{};
      \oddcell{(2a)++(2,0)}{2a 1}{};
      \evencell{(2a)++(4,0)}{2a 2}{};
      \oddcell{(2a)++(6,0)}{2a 3}{};
      \evencell{(2a)++(8,0)}{2a 4}{};
      \oddcell{(2a)++(10,0)}{2a 5}{};
      \evencell{(2a)++(12,0)}{2a 6}{};
      \oddcell{(2a)++(14,0)}{2a 7}{};
      \evencell{(2a)++(16,0)}{2a 8}{};
      \oddrightedgecell{(2a)++(18,0)}{2a 9}{};
      
      % Entries of row 2a (backmost bottom row)
      \begin{scope}[inner sep=0, minimum size=0, node distance=7]
        % \oddtrivbdysyll{2a 3}{}{$ \delta^2 $}{}{}{};
        \evenintsyll{2a 4}{}{ $\gamma\alpha$ }{};
        \oddintsyll{2a 5}{}{$ \gamma $}{};
      \end{scope}
      
    %   % Cells of row 2b (middle bottom row)
    %   \evenleftedgecell{(2b)++(0,0)}{2b 0}{};
    %   \oddcell{(2b)++(2,0)}{2b 1}{};
    %   \evencell{(2b)++(4,0)}{2b 2}{};
    %   \oddcell{(2b)++(6,0)}{2b 3}{};
    %   \evencell{(2b)++(8,0)}{2b 4}{};
    %   \oddcell{(2b)++(10,0)}{2b 5}{};
    %   \evencell{(2b)++(12,0)}{2b 6}{};
    %   \oddcell{(2b)++(14,0)}{2b 7}{};
    %   \evencell{(2b)++(16,0)}{2b 8}{};
    %   \oddrightedgecell{(2b)++(18,0)}{2b 9}{};

    %   % Entries of row 2b (middle bottom row)
    %   \begin{scope}[inner sep=0, minimum size=0, node distance=7]
    %     \evenbdysyll{2b 6}{}{ $\alpha$}{}{}{};
    %   \end{scope}
      
      % Draw row indices
      \newcommand{\rowindexdisplacement}{(-3,0)}
      
      \newcommand{\rowindex}[2]{
        \draw (#1) ++\rowindexdisplacement node[lightgray](#1 index){#2};
        % \filldraw[lightgray] (#1 index) circle (0.15);
        % 1: label of row node
      }
      \newcommand{\rowindexconnector}[2]{
        \draw[lightgray] (#1 index) -- (#2 index);
        % 1: Label of row node for source index
        % 2: Label of row node for target index
      }
      
      \rowindex{1}{(c)};
      \rowindex{2a}{(b)};
    %   \rowindex{2b};
    %   \rowindex{2c};
      \rowindexconnector{1}{2a};
    %   \rowindexconnector{1}{2b};
    %   \rowindexconnector{1}{2c};
    \end{tikzpicture}
  }

  \vspace{4ex}

  {
    \begin{tikzpicture}[scale=0.5]
      % Cells of 1st row (top row)
      \draw (0,1) node[](1){};
      
      \oddleftedgecell{(1) ++(0,0)}{1 0}{};
      \evencell{(1) ++(2,0)}{1 1}{};
      \oddcell{(1) ++(4,0)}{1 2}{};
      \evencell{(1) ++(6,0)}{1 3}{};
      \oddcell{(1) ++(8,0)}{1 4}{};
      \evencell{(1) ++(10,0)}{1 5}{};
      \oddcell{(1) ++(12,0)}{1 6}{};
      \evencell{(1) ++(14,0)}{1 7}{};
      \oddcell{(1) ++(16,0)}{1 8}{};
      \evenrightedgecell{(1) ++(18,0)}{1 9}{};
      
      % Entries of 1st row (top row)
      \begin{scope}[inner sep=0, minimum size=0, node distance=7]
        % \eventrivbdysyll{1 3}{}{$\source(\beta)$}{}{}{};
        % \oddtrivbdysyll{1 4}{}{$\source(\gamma)$}{}{}{};
        \evenbdysyll{1 5}{}{$\delta^2$}{}{}{};
        % \oddtrivbdysyll{1 6}{}{$\source(\gamma)$}{}{}{};
      \end{scope}
      
      \newcommand{\rowdisplacement}{(0.8,-1.4)}
      \draw (0,-2) node[](2a){};
      \draw (2a) ++\rowdisplacement node[](2b){};
      \draw (2b) ++\rowdisplacement node[](2c){};
      
      % Cells of row 2a (backmost bottom row)
      \evenleftedgecell{(2a)++(0,0)}{2a 0}{};
      \oddcell{(2a)++(2,0)}{2a 1}{};
      \evencell{(2a)++(4,0)}{2a 2}{};
      \oddcell{(2a)++(6,0)}{2a 3}{};
      \evencell{(2a)++(8,0)}{2a 4}{};
      \oddcell{(2a)++(10,0)}{2a 5}{};
      \evencell{(2a)++(12,0)}{2a 6}{};
      \oddcell{(2a)++(14,0)}{2a 7}{};
      \evencell{(2a)++(16,0)}{2a 8}{};
      \oddrightedgecell{(2a)++(18,0)}{2a 9}{};
      
      % Entries of row 2a (backmost bottom row)
      \begin{scope}[inner sep=0, minimum size=0, node distance=7]
        % \oddtrivbdysyll{2a 3}{}{$ \source(\beta) $}{}{}{};
        \evenbdysyll{2a 4}{}{ $\alpha$ }{}{}{};
      \end{scope}
      
      % Cells of row 2b (middle bottom row)
      \evenleftedgecell{(2b)++(0,0)}{2b 0}{};
      \oddcell{(2b)++(2,0)}{2b 1}{};
      \evencell{(2b)++(4,0)}{2b 2}{};
      \oddcell{(2b)++(6,0)}{2b 3}{};
      \evencell{(2b)++(8,0)}{2b 4}{};
      \oddcell{(2b)++(10,0)}{2b 5}{};
      \evencell{(2b)++(12,0)}{2b 6}{};
      \oddcell{(2b)++(14,0)}{2b 7}{};
      \evencell{(2b)++(16,0)}{2b 8}{};
      \oddrightedgecell{(2b)++(18,0)}{2b 9}{};

      % Entries of row 2b (middle bottom row)
      \begin{scope}[inner sep=0, minimum size=0, node distance=7]
        \oddtrivbdysyll{2b 5}{}{ $\source(\delta)$}{}{}{};
        % \eventrivbdysyll{2b 6}{}{ $\source(\gamma)$}{}{}{};
      \end{scope}
      
      % Draw row indices
      \newcommand{\rowindexdisplacement}{(-3,0)}
      
      \newcommand{\rowindex}[2]{
        \draw (#1) ++\rowindexdisplacement node[lightgray](#1 index){#2};
        % \filldraw[lightgray] (#1 index) circle (0.15);
        % 1: label of row node
      }
      \newcommand{\rowindexconnector}[2]{
        \draw[lightgray] (#1 index) -- (#2 index);
        % 1: Label of row node for source index
        % 2: Label of row node for target index
      }
      
      \rowindex{1}{(d)};
      \rowindex{2a}{(c)};
      \rowindex{2b}{(e)};
    %   \rowindex{2c};
      \rowindexconnector{1}{2a};
      \rowindexconnector{1}{2b};
    %   \rowindexconnector{1}{2c};
    \end{tikzpicture}
  }

  \vspace{4ex}

  { 
    \begin{tikzpicture}[scale=0.5]
      % Cells of 1st row (top row)
      \draw (0,1) node[](1){};
      
      \oddleftedgecell{(1) ++(0,0)}{1 0}{};
      \evencell{(1) ++(2,0)}{1 1}{};
      \oddcell{(1) ++(4,0)}{1 2}{};
      \evencell{(1) ++(6,0)}{1 3}{};
      \oddcell{(1) ++(8,0)}{1 4}{};
      \evencell{(1) ++(10,0)}{1 5}{};
      \oddcell{(1) ++(12,0)}{1 6}{};
      \evencell{(1) ++(14,0)}{1 7}{};
      \oddcell{(1) ++(16,0)}{1 8}{};
      \evenrightedgecell{(1) ++(18,0)}{1 9}{};
      
      % Entries of 1st row (top row)
      \begin{scope}[inner sep=0, minimum size=0, node distance=7]
        % \eventrivbdysyll{1 3}{}{$\source(\beta)$}{}{}{};
        \oddtrivbdysyll{1 4}{}{$\source(\delta)$}{}{}{};
        % \evenbdysyll{1 5}{}{$\delta^2$}{}{}{};
        % \oddtrivbdysyll{1 6}{}{$\source(\gamma)$}{}{}{};
      \end{scope}
      
      \newcommand{\rowdisplacement}{(0.8,-1.4)}
      \draw (0,-2) node[](2a){};
      \draw (2a) ++\rowdisplacement node[](2b){};
      \draw (2b) ++\rowdisplacement node[](2c){};
      
      % Cells of row 2a (backmost bottom row)
      \evenleftedgecell{(2a)++(0,0)}{2a 0}{};
      \oddcell{(2a)++(2,0)}{2a 1}{};
      \evencell{(2a)++(4,0)}{2a 2}{};
      \oddcell{(2a)++(6,0)}{2a 3}{};
      \evencell{(2a)++(8,0)}{2a 4}{};
      \oddcell{(2a)++(10,0)}{2a 5}{};
      \evencell{(2a)++(12,0)}{2a 6}{};
      \oddcell{(2a)++(14,0)}{2a 7}{};
      \evencell{(2a)++(16,0)}{2a 8}{};
      \oddrightedgecell{(2a)++(18,0)}{2a 9}{};
      
      % Entries of row 2a (backmost bottom row)
      \begin{scope}[inner sep=0, minimum size=0, node distance=7]
        % \oddtrivbdysyll{2a 3}{}{$ \source(\beta) $}{}{}{};
        \evenbdysyll{2a 4}{}{ $\delta^2$ }{}{}{};
      \end{scope}
      
      % Cells of row 2b (middle bottom row)
      \evenleftedgecell{(2b)++(0,0)}{2b 0}{};
      \oddcell{(2b)++(2,0)}{2b 1}{};
      \evencell{(2b)++(4,0)}{2b 2}{};
      \oddcell{(2b)++(6,0)}{2b 3}{};
      \evencell{(2b)++(8,0)}{2b 4}{};
      \oddcell{(2b)++(10,0)}{2b 5}{};
      \evencell{(2b)++(12,0)}{2b 6}{};
      \oddcell{(2b)++(14,0)}{2b 7}{};
      \evencell{(2b)++(16,0)}{2b 8}{};
      \oddrightedgecell{(2b)++(18,0)}{2b 9}{};

      % Entries of row 2b (middle bottom row)
      \begin{scope}[inner sep=0, minimum size=0, node distance=7]
        \oddbdysyll{2b 5}{}{ $\alpha$}{}{}{};
        % \eventrivbdysyll{2b 6}{}{ $\source(\gamma)$}{}{}{};
      \end{scope}
      
      % Draw row indices
      \newcommand{\rowindexdisplacement}{(-3,0)}
      \newcommand{\rowindex}[2]{
        \draw (#1) ++\rowindexdisplacement node[lightgray](#1 index){#2};
        % \filldraw[lightgray] (#1 index) circle (0.15);
        % 1: label of row node
      }
      \newcommand{\rowindexconnector}[2]{
        \draw[lightgray] (#1 index) -- (#2 index);
        % 1: Label of row node for source index
        % 2: Label of row node for target index
      }
      
      \rowindex{1}{(e)};
      \rowindex{2a}{(d)};
      \rowindex{2b}{(c)};
    %   \rowindex{2c};
      \rowindexconnector{1}{2a};
      \rowindexconnector{1}{2b};
    %   \rowindexconnector{1}{2c};
    \end{tikzpicture}
  }
  
  \caption[
    Parts of a syzygy fabric
  ]{
    \label{fig:syz-fab-of-rex-e2}
    \textit{Parts of the the syzygy fabric of $\Str\big(
      \begin{tikzpicture}[scale=0.75]
      \begin{scope}[minimum size=0, inner sep=1]
        % Vertices of string graph
        \draw (0,0) node[](v1){$ 1$};
        \draw (1,0) node[](v2){$ 1$};
        \draw (2,0) node[](v3){$ 2$};
        \draw (3,0) node[](v4){$ 2$};
        \draw (4,0) node[](v5){$ 2$};
        \draw (5,0) node[](v6){$ 2$};
        %
        % Arrows of string graph
        \draw[->] (v1) --node[pos=0.5,above](){$\alpha$} (v2);
        \draw[->] (v2) --node[pos=0.5,above](){$\beta$} (v3);
        \draw[->] (v4) --node[pos=0.5,above](){$\delta$} (v3);
        \draw[->] (v5) --node[pos=0.5,above](){$\delta$} (v4);
        \draw[->] (v6) --node[pos=0.5,above](){$\delta$} (v5);
      \end{scope}
      \end{tikzpicture}
    \big)=E_2$.}
  }
\end{thesisfigure}

First, let us describe the tree $\mathcal{T}$ indexing the rows of this array. This tree $\mathcal{T}$ is constructed as the limit of an increasing sequence $\mathcal{T}_{(0)} \subseteq \mathcal{T}_{(1 )}\subseteq \mathcal{T}_{(2)} \subseteq \cdots$ of trees. The initial tree $\mathcal{T}_{(0)}$ comprises a single vertex, labelled (a), and no arrows. Then, to any leaf vertex of the $k$th tree $\mathcal{T}_{(k)}$, we read the label of the leaf and add child vertices according to the following scheme; this process yields $T_{(k+1)}$.
\begin{center}$
\hspace{\fill}
\begin{tikzcd}[column sep=tiny]
    \& \text{(a)} \ar[dl] \ar[dr]
    \\
  \text{(b)}
    \&
      \& \text{(c)}
\end{tikzcd}
\hspace{\fill}
\begin{tikzcd}[column sep=tiny]
    \& \text{(b)} \ar[dl] \ar[d] \ar[dr]
    \\
  \text{(d)}
    \& \text{(c)}
      \& \text{(d)}
\end{tikzcd}
\hspace{\fill}
\begin{tikzcd}[column sep=tiny]
\text{(c)} \ar[d]
    \\
  \text{(b)}
\end{tikzcd}
\hspace{\fill}
\begin{tikzcd}[column sep=tiny]
    \& \text{(d)} \ar[dl] \ar[dr]
    \\
  \text{(c)}
    \&
      \& \text{(e)}
\end{tikzcd}
\hspace{\fill}
\begin{tikzcd}[column sep=tiny]
    \& \text{(e)} \ar[dl] \ar[dr]
    \\
  \text{(d)}
    \&
      \& \text{(c)}
\end{tikzcd}
\hspace{\fill}$
\end{center}
So, for example, any leaf of $T_{(k)}$ that is labelled (a) earns child vertices in $T_{(k+1)}$ respectively labelled (b) and (c), while a leaf labelled (b) earns child vertices labelled (d), (c) and (d), and so on.

How these rows are populated is illustrated in Figure \ref{fig:syz-fab-of-rex-e2}. By piecing together the fabric one row to next, we complete the array in the limit.

\chapter[Syzygy patterns and applications]{Syzygy patterns and applications to homological questions}\label{chap:patterns-and-applications}

% In this chapter, we prove some results about the syzygies of string modules

\section{Interlude on associated topics}

We hit pause on our string-and-syzygy focus to explore other pertinent areas of theory.

In this section, we prove several results that we will later need. Many of these appear stated in the literature without proof: a shortcoming we are happy to rectify. 

\subsection{Properties of colocalising subcategories}

To begin with, we demonstrate that colocalising subcategories often behave similarly to localising ones.

\nlpass{Products of complexes} Suppose that the members $m$ of some set index a collection of complexes
$
X_m^\bullet \ceq
  \big(
  \begin{tikzcd}[sep=small]
  \cdots \ar[r]
    \& X_m^{-1} \ar[rr, "d_m^{-1}"]
      \&\& X_m^{0} \ar[rr, "d_m^{0}"]
        \&\& X_m^1 \ar[rr, "d_m^1"]
          \&\& X_m^2 \ar[r]
            \& \cdots
  \end{tikzcd}
  \big)
$.
The obvious termwise product of these complexes, having for $k$th term the product $\prod_m X_m^k $ of $k$th terms and for $k$th boundary map the product $\prod_m d_m^k$ of $k$th boundary maps, is indeed the categorical product of the $X_m^\bullet$ in $\complexcat(A)$.

\nlpass{Inverse limits as cokernels}\label{psg:inv-lims-as-colims} A small diagram in $\complexcat(A)$, comprising a set of morphisms connecting a set $\{X_m^\bullet \colon m\}$ of objects, yields a morphism $F^\bullet \colon \ilarrow{  \prod_m X_m^\bullet  }{}{  \prod_m X_m^\bullet  }$ (a morphism which, to be clear, also lives in $\complexcat(A)$). One such small diagram is an inverse system of $\complexcat(A)$-morphisms:
$$
\Delta \ceq \big(\begin{tikzcd}
\cdots \ar[r]
    \& X_3^\bullet \ar[r,"f_3^\bullet"]
        \& X_2^\bullet \ar[r, "f_2^\bullet"]
        \& X_1^\bullet \ar[r, "f    _1^\bullet"]
                \& X_0^\bullet
\end{tikzcd}\big)\tstop
$$
(We mention that this inverse system in $\complexcat(A)$ gives rise to an inverse system in $\Modcat A$ for each $k \in \integers$, namely the system comprising the $k$th terms of each object and morphism.)

The inverse system $\Delta$ can be encoded using the morphism $F^\bullet \colon \ilarrow{ \prod_{m\geq0} X_m^\bullet }{}{ \prod_{m\geq0} X_m^\bullet  }$, whose nonzero components are the identity morphisms and the transition maps $-f_{m+1}^\bullet \colon \ilarrow{ X_{m+1}^\bullet }{}{ X_m^\bullet }$. We may draw $F^\bullet$ as follows, in a manner evoking the separated quiver construction for quivers.
\begin{equation}\label{diag:separated-inverse-system}
\begin{tikzcd}[row sep=small, column sep=large, inner sep=0.1ex]
\prod_m X_m^\bullet \ar[d, "F^\bullet"]
  \&
    \& \big( \cdots \ar[dr]
     \& X_3^\bullet \ar[d, equal] \ar[dr, near start, "-f_3^\bullet"]
       \& X_2^\bullet \ar[d, equal] \ar[dr, near start, "-f_2^\bullet"]
         \& X_1^\bullet \ar[d, equal] \ar[dr, near start, "-f_1^\bullet"]
           \& X_0^\bullet\big) \ar[d, equal]
\\
\prod_m X_m^\bullet
  \&
    \& \big( \cdots
     \& X_3^\bullet
       \& X_2^\bullet
         \& X_1^\bullet
           \& X_0^\bullet\big)
\end{tikzcd}
\end{equation}
The kernel and cokernel of $F^\bullet$ are noteworthy.

One can calculate that $\ker F^\bullet = \varprojlim \Delta$. Unfortunately, there is no guarantee in general that $\coker F^\bullet = 0$. Hence the exact sequence
\begin{equation}\label{seq:exact-seq-of-inv-limit}
\begin{tikzcd}[row sep=small]
0 \ar[r]
  \& \ker F^\bullet \ar[r]
    \& \prod_m X_m^\bullet \ar[r, "F^\bullet"]
      \& \prod_m X_m^\bullet \ar[r]
        \& \coker F^\bullet \ar[r]
          \& 0
            \\
  \& \varprojlim \Delta \ar[u, equals]
\end{tikzcd}
\end{equation}
is typically not a short one (confer \cite[\S3.5]{Wei94}).

This should be contrasted with a dual situation, whereby a direct limit of a direct system may be encoded as the cokernel of a morphism  $G^\bullet \colon \ilarrow{ \bigoplus_m X_m^\bullet }{}{ \bigoplus_m X_m^\bullet }$ on the coproduct. In other words, for any direct system $\Theta$ there exists $G^\bullet$ for which $\coker G^\bullet = \varinjlim \Theta$, constructed analogously to $F^\bullet$ above. One calculates in such a case that $\ker G^\bullet = 0$ always.

\nlpass{The Mittag-Leffler condition}\label{psg:mittag-leffler} Thankfully there are meaningful circumstances in which an $F^\bullet$ as above has $\coker F^\bullet = 0$. At least one such circumstance is captured by the following criterion \cite[\S3.5]{Wei94}

Let $\Delta$ be an inverse system of $\complexcat(A)$-morphisms $f_{m+1}^\bullet \colon \ilarrow{ X_{m+1}^\bullet }{}{ X_m^\bullet }$ ($m \geq 0$) that we encode as a morphism $F^\bullet$, as in passage \ref{psg:inv-lims-as-colims}.
% This system satisfies the \define{Mittag-Leffler condition} if the descending sequence
% $$
% \cdots \subseteq \im( f_{m+3}^\bullet f_{m+2}^\bullet f_{m+1}^\bullet ) \subseteq \im(  f_{m+2}^\bullet f_{m+1}^\bullet ) \subseteq \im( f_{m+1}^\bullet ) \subseteq X_m^\bullet
% $$
% stabilises.

We know that, for each degree $k \in \integers$, $\Delta$ induces an inverse system $\Delta^k$ of module homomorphisms $f^k_{m+1} \colon \ilarrow{ X_{m+1}^k }{}{ X_m^k }$. The inverse system $\Delta^k$ satisfies the \define{Mittag-Leffler condition} if the descending sequence
$$
\cdots \subseteq \im( f_{m+3}^k f_{m+2}^k f_{m+1}^k ) \subseteq \im(  f_{m+2}^k f_{m+1}^k ) \subseteq \im( f_{m+1}^k ) \subseteq X_m^k
$$
of images stabilises for each $m \geq 0$. This occurs when $X_m^k$ is finite dimensional for each $m$, for instance, or when each $f_{m+1}^k$ is surjective.

\begin{lemma}\label{lem:mittag-leffler-on-inv-systems-of-modules} \emph{\cite[Prop 3.5.7]{Wei94}}
Let $\Delta^k$ be as above. If $\Delta^k$ satisfies the Mittag-Leffler condition, then the $k$th term of $\coker F^\bullet$ is zero.
\end{lemma}

\begin{corollary}\label{cor:mittag-leffler-on-inv-systems-of-complexes}
Let $F^\bullet,\Delta,\Delta^k$ be as above. If $\Delta^k$ satisfies the Mittag-Leffler condition for each $k \in \integers$, then $\coker F^\bullet = 0$.
\end{corollary} 

\npass We are now able to give a list of properties for colocalising subcategories, roughly analogous to Proposition \ref{lem:properties-of-loc-subcats}.

\newcommand{\colocsubcat}{\categoryfont{C}}

\begin{proposition}\label{prop:properties-of-coloc-subcats}
Let $\colocsubcat$ be a colocalising subcategory of $\derivedcat(A)$.
\begin{enumerate}
    \item\label{subprop:concrete-coloc-subcat-properties}
    \begin{enumerate}
        \item\label{subprop:2-out-of-3-in-coloc} If
        $\big( \begin{tikzcd}[sep=small]
        0 \ar[r]
          \& X^\bullet \ar[r]
            \& Y^\bullet \ar[r]
              \& Z^\bullet \ar[r]
                \& 0
        \end{tikzcd} \big)$
        is a short exact sequence of complexes and two of the three objects $X^\bullet,Y^\bullet,Z^\bullet$ are in $\colocsubcat$, then so is the third.
        
        \item\label{subprop:shifts-in-coloc} If a complex $X^\bullet$ is in $\colocsubcat$, then so is $X^\bullet[r]$ for every $r\in\integers$.
        
        \item\label{subprop:qisomic-in-coloc} If $X^\bullet$ and $Y^\bullet$ are quasiisomorphic complexes and $X^\bullet$ is in $\colocsubcat$, then so is $Y^\bullet$.
        
        \item\label{subprop:products-in-coloc} If $\{X_i^\bullet \colon i \in \mathcal{I} \}$ is a set of objects of $\colocsubcat$, then $\prod_{i \in \mathcal{I}}X_i^\bullet$ is in $\colocsubcat$.
    \end{enumerate}
    
    \item\label{subprop:direct-factors-in-coloc} If $X^\bullet \prod Y^\bullet $ is in $\colocsubcat$, then so are $X^\bullet$ and $Y^\bullet$.
    
    \item\label{subprop:bounded-complexes-in-coloc} If $X^\bullet$ is a bounded complex, where the module $X^k$ is in $\colocsubcat$ for every $k \in \integers$, then $X^\bullet$ is in $\colocsubcat$.
    
    \item\label{subprop:inverse-system-in-coloc} If
    $
    \big(\begin{tikzcd}[sep=small]
    \cdots \ar[r]
      \& X_3^\bullet \ar[rr, "f_3^\bullet"]
        \&\& X_2^\bullet \ar[rr, "f_2^\bullet"]
          \&\& X_1^\bullet \ar[rr, "f_1^\bullet"]
            \&\& X_0^\bullet
    \end{tikzcd}\big)
    $ is an inverse system of $\complexcat(A)$-morphisms between complexes $X^\bullet_m$, if each resulting inverse system $
    \big(\begin{tikzcd}[sep=small]
    \cdots \ar[r]
      \& X_3^k \ar[rr, "f_3^k"]
        \&\& X_2^k \ar[rr, "f_2^k"]
          \&\& X_1^k \ar[rr, "f_1^k"]
            \&\& X_0^k
    \end{tikzcd}\big)
    $ of $k$th terms satisfies the Mittag-Leffler condition and if furthermore each complex $X_m^\bullet$ belongs to $\colocsubcat$, then $\varprojlim X_m^\bullet $ is in $\colocsubcat$.
    
    \item\label{subprop:bounded-below-complexes-in-coloc} If $X^\bullet$ is a bounded-below complex of modules $X^k$ ($k \in \integers$), each of which is in $\colocsubcat$, then $X^\bullet$ is in $\colocsubcat$.
\end{enumerate}
\end{proposition}

\proof Much of the proof proceeds in the same fashion as that of \cite[Prop 2.1]{Ric19}. Part \ref{subprop:concrete-coloc-subcat-properties} is a concrete description of a colocalising subcategory. Part \ref{subprop:direct-factors-in-coloc} follows from the same take on the Eilenberg swindle used by Rickard, with the additional comment that $X^\bullet \prod Y^\bullet = X^\bullet \oplus Y^\bullet$ since $X^\bullet \prod Y^\bullet$ has finitely many factors. Part \ref{subprop:bounded-complexes-in-coloc} results from the same induction on the length of the support of $X^\bullet$. To prove parts \ref{subprop:inverse-system-in-coloc} and \ref{subprop:bounded-below-complexes-in-coloc}, however, requires careful use of the Mittag-Leffler condition.

We address part \ref{subprop:inverse-system-in-coloc} first. As we did before in Passages \ref{psg:inv-lims-as-colims} and \ref{psg:mittag-leffler}, encode the inverse system as a $\complexcat(A)$-morphism
$
F^\bullet \colon \ilarrow{ \prod_{m\geq0} X_m^\bullet }{}{ \prod_{m\geq0} X_m^\bullet  }
$. The Mittag-Leffler assumption ensures that we may apply Corollary \ref{cor:mittag-leffler-on-inv-systems-of-complexes} and conclude that $\coker F^\bullet = 0$. This yields a short exact sequence
$$
\begin{tikzcd}
0 \ar[r]
  \& \varprojlim X_m^\bullet = \ker F^\bullet \ar[r]
    \& \prod_m X_m^\bullet \ar[r, "F^\bullet"]
      \& \prod_m X_m^\bullet \ar[r]
        \& 0
\end{tikzcd}\tstop
$$
whose middle term and righthand term belong to $\colocsubcat$ by assumption and part \ref{subprop:concrete-coloc-subcat-properties}\ref{subprop:products-in-coloc}, so the result follows by part \ref{subprop:concrete-coloc-subcat-properties}\ref{subprop:2-out-of-3-in-coloc}.

Now for part \ref{subprop:bounded-below-complexes-in-coloc}, assume $X^\bullet$ is supported in nonnegative degrees. (By part \ref{subprop:concrete-coloc-subcat-properties}\ref{subprop:shifts-in-coloc} this is no loss of generality.) For $m \geq 0$ let
$
\sigma^m( X^\bullet ) \ceq \big(
\begin{tikzcd}[sep=small]
\cdots \ar[r]
    \& 0 \ar[r]
        \& X^0 \ar[r]
            \& \cdots \ar[r]
                \& X^m \ar[r]
                    \& 0 \ar[r]
                        \& \cdots
\end{tikzcd}
\big)
$
be the brutal truncation of $X^\bullet$ below $m$. Each truncation is a bounded complex of members of $\colocsubcat$ by assumption, so in $\colocsubcat$ by part \ref{subprop:bounded-complexes-in-coloc}.

One finds that $X^\bullet$ is the inverse limit of its brutal truncations: this inverse system is shown vertically on the lefthand side of Figure \ref{fig:bounded-below-complex-as-inverse-limit-of-truncations}, and its contents are laid bare on the righthand side. As the caption to the figure makes clear, the inverse system of $k$th terms satisfies the Mittag-Leffler condition for each $k \in \integers$, therefore its inverse limit $X^\bullet$ belongs to $\colocsubcat$ by part \ref{subprop:inverse-system-in-coloc}. This establishes part \ref{subprop:bounded-below-complexes-in-coloc}, concluding the proof.\qed
\begin{thesisfigure}
\centering
$
\begin{tikzcd}
X^\bullet
  \&
    \& \cdots \ar[r]
      \& 0 \ar[r]
        \& X^0 \ar[r]
          \& X^1 \ar[r]
            \& X^2 \ar[r]
              \& X^3 \ar[r]
                \& \cdots
\\
\\
\vdots \ar[d]
  \&
    \& 
      \& \vdots\ar[d]
        \& \vdots \ar[d, equals]
          \& \vdots \ar[d, equals]
            \&\vdots \ar[d, equals]
              \& \vdots \ar[d]
                \&
                  \\
\sigma^{\leq 2}(X^\bullet) \ar[d]
  \&
    \& \cdots \ar[r]
      \& 0 \ar[r] \ar[d]
        \& X^0 \ar[r] \ar[d, equals]
          \& X^1 \ar[r] \ar[d, equals]
            \& X^2 \ar[r] \ar[d]
              \& 0 \ar[r] \ar[d]
                \& \cdots
                  \\
\sigma^{\leq 1}(X^\bullet) \ar[d]
  \&
    \& \cdots \ar[r]
      \& 0 \ar[r] \ar[d]
        \& X^0 \ar[r] \ar[d, equals]
          \& X^1 \ar[r] \ar[d]
            \& 0 \ar[r] \ar[d]
              \& 0 \ar[r] \ar[d]
                \& \cdots
                  \\
\sigma^{\leq 0}(X^\bullet )
  \&
    \& \cdots \ar[r]
      \& 0 \ar[r]
        \& X^0 \ar[r]
          \& 0 \ar[r]
            \& 0 \ar[r]
              \& 0 \ar[r]
                \& \cdots
\end{tikzcd}
$

\caption[%
  A bounded-below complex is an inverse limit of its brutal truncations 
]{%
  \label{fig:bounded-below-complex-as-inverse-limit-of-truncations}
  \textit{A bounded-below complex $X^\bullet$  is an inverse limit of its brutal truncations $\sigma^{\leq m}(X^\bullet)$.} The components of each morphism $\ilarrow{ \sigma^{\leq (m+1)}(X^\bullet) }{}{ \sigma^{\leq m}(X^\bullet) }$ are all surjections, being either maps onto $0$ or identity maps $\ilarrow{ X^k }{}{ X^k }$. Thus the inverse system of $k$th terms satisfies the Mittag-Leffler condition, for any $k \in \integers$.
}
\end{thesisfigure}

\subsection{Band modules}\label{subsec:band-modules}

Next, we take a detour to explore band modules and their syzygies. This is important for our investigation of syzygy-finiteness of an SB algebra later.

\nlpass{Finite powers of string graphs} Call a connected string graph $w$ with at least one arrow \define{powerable} if it has two sink vertices $i,j$, both with indegree 1 and such that $w(i) = w(j)$ but $w(\target^{-1}(i)) \neq w(\target^{-1}(j))$. 

The \define{$m$th power $w^m$} ($m\geq 1$) of any powerable string graph $w$ with sink vertices $i \neq j$ both having indegree 1, is the string graph obtained from $m$ disjoint copies of $w$ by identifying the $r$th copy of $j$ with the $(r+1)$th copy of $i$. Observe that $w^m$ is powerable too.

A string graph $v$ is \define{primitive} if it is not $w^m$ for any string graph $w$ and integer $m \geq 2$. Every powerable finite string graph is a power of some primitive powerable string graph.

\newcommand{\infpower}[1]{\widehat{#1}}

\nlpass{Biinfinite powers of string graphs} If $w$ is a powerable string graph with sink vertices $i \neq j$, both of degree 1, then we write \define{$\infpower{w}$} for the \define{biinfinite power} $\cdots www \cdots$ of $w$, formed from the $\integers$-indexed disjoint union of copies of $w$ by identifying the $r$th copy of $j$ with the $(r+1)$th copy of $i$.

\nrex Two examples of powerable string graphs for $A$ are $$\big(
  \begin{tikzpicture}[scale=0.75]
  \begin{scope}[minimum size=0, inner sep=1]
    % Vertices of string graph
    \draw (0,0) node[](v1){$ 1$};
    \draw (1,0) node[](v2){$ 1$};
    \draw (2,0) node[](v3){$ 2$};
    \draw (3,0) node[](v4){$ 1$};
    \draw (4,0) node[](v5){$ 1$};
    \draw (5,0) node[](v6){$ 2$};
    \draw (6,0) node[](v7){$ 1$};
    %
    % Arrows of string graph
    \draw[->] (v2) --node[pos=0.5,above](){$\alpha$} (v1);
    \draw[->] (v2) --node[pos=0.5,above](){$\beta$} (v3);
    \draw[->] (v3) --node[pos=0.5,above](){$\gamma$} (v4);
    \draw[->] (v5) --node[pos=0.5,above](){$\alpha$} (v4);
    \draw[->] (v5) --node[pos=0.5,above](){$\beta$} (v6);
    \draw[->] (v6) --node[pos=0.5,above](){$\gamma$} (v7);
  \end{scope}
  \end{tikzpicture}
\big)
\text{\hspace{0.5cm}and\hspace{0.5cm}}
\big(
  \begin{tikzpicture}[scale=0.75]
  \begin{scope}[minimum size=0, inner sep=1]
    % Vertices of string graph
    \draw (0,0) node[](v1){$ 1$};
    \draw (1,0) node[](v2){$ 1$};
    \draw (2,0) node[](v3){$ 2$};
    \draw (3,0) node[](v4){$ 2$};
    \draw (4,0) node[](v5){$ 2$};
    \draw (5,0) node[](v6){$ 1$};
    \draw (6,0) node[](v7){$ 1$};
    \draw (7,0) node[](v8){$ 2$};
    \draw (8,0) node[](v9){$ 1$};
    %
    % Arrows of string graph
    \draw[->] (v2) --node[pos=0.5,above](){$\alpha$} (v1);
    \draw[->] (v3) --node[pos=0.5,above](){$\gamma$} (v2);
    \draw[->] (v3) --node[pos=0.5,above](){$\delta$} (v4);
    \draw[->] (v4) --node[pos=0.5,above](){$\delta$} (v5);
    \draw[->] (v6) --node[pos=0.5,above](){$\beta$} (v5);
    \draw[->] (v7) --node[pos=0.5,above](){$\alpha$} (v6);
    \draw[->] (v7) --node[pos=0.5,above](){$\beta$} (v8);
    \draw[->] (v8) --node[pos=0.5,above](){$\gamma$} (v9);
  \end{scope}
  \end{tikzpicture}
\big)\tstop
$$
The former is $w^2$ for $w \ceq
\big(
  \begin{tikzpicture}[scale=0.75]
  \begin{scope}[minimum size=0, inner sep=1]
    % Vertices of string graph
    \draw (0,0) node[](v1){$ 1$};
    \draw (1,0) node[](v2){$ 1$};
    \draw (2,0) node[](v3){$ 2$};
    \draw (3,0) node[](v4){$ 1$};
    %
    % Arrows of string graph
    \draw[->] (v2) --node[pos=0.5,above](){$\alpha$} (v1);
    \draw[->] (v2) --node[pos=0.5,above](){$\beta$} (v3);
    \draw[->] (v3) --node[pos=0.5,above](){$\gamma$} (v4);
  \end{scope}
  \end{tikzpicture}
\big)
$. The latter is primitive.

\nlpass{Band graph} A \define{band graph (for $A$)} is a quiver homomorphism $w \colon \ilarrow{G}{}{Q} $ such that:
\begin{enumerate}
    \item the underlying (undirected) graph of $G$ is a connected cycle;
    \item
    \begin{enumerate}
        \item for any subgraph $\big(\begin{tikzcd}[sep=small,cramped]
    \syllableplaceholder \& i \ar[l,"x"'] \ar[r, "y"] \& \syllableplaceholder
    \end{tikzcd}\big)$ of $G$ featuring a source vertex $i$ and its two outgoing arrows, we have $w(x)\neq w(y)$;
        \item for any subgraph
    $\big(\begin{tikzcd}[sep=small,cramped]
    \syllableplaceholder \ar[r, "x"] \& i \& \ar[l,"y"'] \syllableplaceholder 
    \end{tikzcd}\big)$ of $G$ featuring a sink vertex $i$ and its two incoming arrows, we have $w(x)\neq w(y)$;
    \end{enumerate}
    \item for any path
    $\big(\begin{tikzcd}[sep=small]
    \syllableplaceholder \ar[r, "x_1"]
      \& \syllableplaceholder \ar[r, "x_2"]
        \& \cdots \ar[r, "x_\ell"]
          \& \syllableplaceholder
    \end{tikzcd}\big)$
    in $G$, the $A$-path $p\ceq w(x_1 x_2 \cdots x_\ell)$ is linearly independent of any other $A$-path.
\end{enumerate}

As with string graphs, we depict a band graph as a labelled graph, each vertex $v$ or arrow $x$ being labelled by $w(v)$ or $w(x)$.

% \nrex ``$\alpha^{-1}\beta\gamma$'' \textbf{[ADD EXAMPLE]}

\nlpass{Bands from string graphs} Let $w$ be any powerable string with sink vertices $i \neq j$ both having indegree $1$. Identifying $i$ and $j$ yields a band graph.

Every band graph $v$ is obtained in this fashion from some powerable string graph $w^m$.

\nlpass{Band modules} We follow the notation of \cite[\S1]{HZ16}.

Let $v$ be a primitive string graph with sink vertices $i \neq j$ both having indegree $1$, let $m \geq 1$ and let $\psi \colon \ilarrow{\kk^m }{}{ \kk^m }$ be an indecomposable vector-space automorphism with companion matrix as in Figure \ref{fig:companion-matrix}.

Let $i_r$ ($1 \leq r \leq m$) for the $m$ vertices of $v^m$ corresponding to $i$ and, abusing notation, also write $j$ for the sink vertex of $v^m$ corresponding to the $m$th copy of $j$. Each of these vertices is a basis vector of $\Str(v^r)$, and all of these particular basis vectors are fixed by a unique vertex idempotent $e_k$ of $A$. It follows that $j$ and the $i_r$ each span an isomorphic simple submodule of $\Str( v^m )$.
\begin{thesisfigure}
\centering
$
\begin{bmatrix}
  0
    & 0
      & 0
        & \cdots
          & 0
            & 0
              & \lambda_1
               \\
  1
    & 0
      & 0
        & \cdots
          & 0
            & 0
              & \lambda_2
                \\
  0
    & 1
      & 0
        & \cdots
          & 0
            & 0
              & \lambda_3
                \\
  \vdots
    & \vdots
      & \vdots
        & \ddots
          & \vdots
            & \vdots
              & \vdots
                \\
  0
    & 0
      & 0
        & \cdots
          & 0
            & 0
              & \lambda_{m-2}
                \\
  0
    & 0
      & 0
        & \cdots
          & 1
            & 0
              & \lambda_{m-1}
                \\
  0
    & 0
      & 0
        & \cdots
          & 0
            & 1
              & \lambda_m
\end{bmatrix}
$

\caption[
  Companion matrix
]{
  \label{fig:companion-matrix}
  \textit{Companion matrix}
}
\end{thesisfigure}

The \define{band module $\Bnd(v^m, \psi)$} is the quotient $\Str( v^m )/\langle j - \sum_{r=1}^m \lambda_r i_r \rangle$.

\npass Huisgen-Zimmermann states the following results in \cite{HZ16} but directs the interested reader to a forthcoming PhD thesis \cite{GalPHD} for the proof. We cannot find any trace of this thesis in the literature so, in its absence, we supply our own proofs.

\begin{proposition}\label{prop:syzygies-of-band-modules} Let $v$ be a primitive string graph and let $\infpower{v}$ be the corresponding biinfinite power. Moreover, let $m$ be any positive integer and $\psi$ an indecomposable automorphism of $\kk^m$ with companion matrix as in Figure \ref{fig:companion-matrix}. Then the following statements are equivalent.
\begin{enumerate}
    \item\label{subprop:syzygy-of-band-is-band} The syzygy $\syzygy^1(\Bnd(v^m, \psi))$ is a band module
    \item\label{subprop:syzygy-of-inf-string-is-indec} The syzygy $\syzygy^1(\Str \infpower{v} )$ is an indecomposable, infinite-dimensional string module.
    \item\label{subprop:projcover-all-pin} None of the indecomposable direct summands of $\projcover(\Str v)$ is a string module; that is, they are all pin.
\end{enumerate}
If these conditions fail to be satisfied, then $\syzygy^1(\Bnd(v^m, \psi))$ is a direct sum of string modules.
\end{proposition}

\proof The equivalence of \ref{subprop:syzygy-of-inf-string-is-indec} and \ref{subprop:projcover-all-pin} follows from our implementation of the syzygy algorithm in Propositions \ref{prop:syzygy-algorithm-for-all-but-soc-quots} and \ref{prop:syzygy-algorithm-for-soc-quots}. It suffices to recognise that direct summands of $\syzygy^1(\Str \infpower{v})$ correspond to segments, and distinct segments arise from string direct summands of $\projcover( \Str \infpower{v} )$.

Towards proving that \ref{subprop:syzygy-of-band-is-band} is equivalent to \ref{subprop:projcover-all-pin}, we will reconstruct $\Bnd(v^m, \psi)$ as a quotient of its projective cover $\projcover(\Bnd(v^m, \psi)) = \projcover(\Str v^m ) = \projcover(\Str v)^{\oplus m}$. Let $p_1,q_1,\dots,p_\ell,q_\ell$ be the $A$-paths such that (for subscripts $r$ taken modulo $\ell$) $p_r \neq q_r$ are the longest paths in the band graph $v$ whose source is a given source vertex while $q_r \neq p_{r+1}$ are the longest paths whose target is a given sink vertex. Let $i_r$ be the $Q$-vertex with $\source(p_r)=\source(q_r)=i_r$, so that $\projcover(\Bnd(v^m,\psi))=(P_{i_1} \oplus \cdots \oplus P_{i_\ell})^{\oplus m}$. Finally, let $z_{r,s}$ denote the generator $e_{i_r}$ of the $s$th copy of $P_{i_r}=e_{i_r} A$ in this direct sum (for $1 \leq r \leq \ell, 1 \leq s \leq m$).

In this notation, we have that $\Bnd(v^m, \psi)$ is the quotient of $ (P_{i_1} \oplus \cdots \oplus  P_{i_\ell})^{\oplus m} $ by the submodule $X$ generated by the following elements:
\begin{equation}\label{eqn:generators-of-syzygy}
\begin{array}{rl}
  z_{r,s}q_r + z_{r+1,s}p_{r+1}
    & (1 \leq r < \ell, 1 \leq s \leq m) \tcomma
      \\
  z_{\ell,s}q_\ell + z_{1,s+1}p_1
    & (1 \leq s < m) \tcomma
      \\
  z_{\ell,m}q_\ell - \sum\limits_{s=1}^m \lambda_s z_{1,s} p_1 \tstop
\end{array}
\end{equation}
By construction, we have that $X = \syzygy^1(\Bnd(v^m, \psi))$.

It is easy to identify $\rad X$.  Towards this goal, let $\alpha_r$ be the unique arrow such that $p_r\alpha_r$ is nonzero in $A$ if such an arrow exists, or $0$ otherwise. Similarly, let $\beta_r$ be the unique arrow such that $q_r \beta_r $ is nonzero in $A$, or $0$ if no such arrow exists. We know that generators of $\rad X = X\rad A$ are given by acting on the generators of $X$ by arrows $\gamma$. We find that, if we act on any generator in (\ref{eqn:generators-of-syzygy}) by any arrow $\gamma$, we obtain $0$ except in the following cases.
\begin{equation}\label{eqn:arrowed-gens-1}
\left.
\begin{array}{lcl}
(z_{r,s}q_r + z_{r+1,s}p_{r+1})\beta_r
  & =
    & z_{r,s}q_r\beta_r
      \\
(z_{r,s}q_r + z_{r+1,s}p_{r+1})\alpha_{r+1}
  & =
    & z_{r+1,s}p_{r+1}\alpha_{r+1}
\end{array}\right\}
(1 \leq r < \ell, 1 \leq s \leq m)
\end{equation}
\begin{equation}\label{eqn:arrowed-gens-2}
\left.
\begin{array}{lcl}
  (z_{\ell,s}q_\ell + z_{1,s+1}p_1) \beta_\ell
    & =
      & z_{\ell,s}q_{\ell}\beta_\ell
        \\
  (z_{\ell,s}q_\ell + z_{1,s+1}p_1) \alpha_1
    & =
      & z_{1,s+1}p_1  \alpha_1
\end{array}
\right\}
(1 \leq s < m)
\end{equation}
\begin{equation}\label{eqn:arrowed-gens-3}
\begin{array}{lcl}
(z_{\ell,m}q_\ell - \sum\limits_{s=1}^m \lambda_s z_{1,s} p_1) \beta_\ell
  & =
    & z_{\ell,m} q_\ell \beta_\ell
      \\
(z_{\ell,m}q_\ell - \sum\limits_{s=1}^m \lambda_s z_{1,s} p_1) \alpha_1
  & = 
    & - \sum\limits_{s=1}^m \lambda_s z_{1,s} p_1\alpha_1
\end{array}
\end{equation}
(It is possible that some of these quantities are $0$ also.) Consider the final element, $- \sum_{s=1}^m \lambda_s z_{1,s} p_1 \alpha$. This is a $\kk$-linear sum of elements $z_{1,s} p_1 \alpha_1$. Aside from the first such element $z_{1,1} p_1 \alpha_1 $, these all appear in (\ref{eqn:arrowed-gens-2}). Therefore the submodule of $ (P_{i_1} \oplus \cdots \oplus  P_{i_\ell})^{\oplus m} $ generated by the elements in (\ref{eqn:arrowed-gens-1}), (\ref{eqn:arrowed-gens-2}) and (\ref{eqn:arrowed-gens-3}) equals
$$\langle z_{r,s}p_r\alpha_r , z_{r,s}q_r\beta_r \colon 1 \leq r \leq \ell, 1 \leq s \leq m\rangle = \bigoplus\limits_{s=1}^m \langle z_{r,s}p_r\alpha_r , z_{r,s}q_r\beta_r \colon 1 \leq r \leq \ell\rangle\tstop
$$
The $s$th direct summand of the righthand expression is just the sum of the principal right ideals $p_r\alpha_r A$ and $q_r \beta_r A$ inside the $s$th copy of $P_{i_r}$. In a given copy, we have $p_r \alpha_r A \cap q_r \alpha_r A = 0$ iff $P_{i_r}$ is not pin. 

Consider the following two statements.
\begin{enumerate}
    \item For each generator $g$ of $X$ given in (\ref{eqn:generators-of-syzygy}), there are distinct arrows $\alpha,\beta$ with $g\alpha \neq 0 \neq g\beta$.
    \item For each simple direct summand $S$ of $\soc X$, there are generators $g,h$ of $X$ given in (\ref{eqn:generators-of-syzygy}) and distinct arrows $\alpha,\beta$ such that $S \leq \langle g\alpha, h\beta \rangle $.
\end{enumerate}
Our characterisation of $\rad X$ and Lemma \ref{subprop:rejection-lemma} imply that these statements are both simultaneously true iff all direct summands $P_{i_r}$ are pin. When they are both true, it is clear that $X$ is a band module. When either is false, we deduce that $X$ is string. This gives the desired result.\qed

\npass Recall from Passage \ref{psg:submods-and-quots} that $\soc X$ denotes the largest semisimple submodule of $X$, which is called the socle of $X$, and that $\head X$ denotes the largest semisimple quotient of $X$, which is called the head of $X$ and is isomorphic to $X/\rad X$.

\begin{corollary}
Suppose all indecomposable direct summands of $\projcover(\Bnd(v^m, \psi))$ are pin. Then $\soc(\Bnd(v^m,\psi)) \isom \head(\syzygy^1(\Bnd(v^m, \psi)) )$.
\end{corollary}

\begin{corollary}
Suppose all indecomposable direct summands of $\projcover(\Bnd(v^m, \psi))$ are pin. Then $\soc(\projcover (\Bnd(v^m,\psi))) \isom \head(\syzygy^2(\Bnd(v^m, \psi)) )$.
\end{corollary}

\begin{corollary}
Suppose all indecomposable direct summands of $\projcover(\Bnd(v^m, \psi))$ are pin and that $p_r,q_r$ ($1 \leq r \leq \ell$) are the longest paths in the band graph $v$, as in the proof of Proposition \ref{prop:syzygies-of-band-modules}. Suppose further that $u_r,v_r$ are $A$-paths such that $p_ru_r$ and $q_rv_r$ are components of commutativity relations.

Then $\syzygy^1(\Bnd(v^m, \psi)) = \Bnd(w^m, \varphi)$, where the longest paths in the band graph $w$ are $v_r$ and $u_{r+1}$, and where $\varphi$ is some linear automorphism $\ilarrow{\kk^m}{}{\kk^m}$.
\end{corollary}

\proof All of these follow from our construction of $\syzygy^1(\Bnd(v^m,\psi))$ above.\qed

\subsection{Width of strips}

We will be interested in results about syzygy-finiteness of certain string modules. This is equivalent to the existence of a bound on widths of indecomposable summands of string modules, as this subsection makes clear.

\nlpass{Width of a string graph} Our previous definition of width for strips coincides with the following definition for string graphs and string modules.

The \define{width $\width w$} of a finite, indecomposable string $w \colon \ilarrow{G}{}{Q}$ is
$$
\width w
\ceq
\big|
  \{
    v\in G_0
  \colon
    v \text{ is a sink or a source vertex}
  \}
\big|
-
1\tcomma
$$
and, in general, we define width of a string as the sum of the widths of its indecomposable parts. The width of an infinite string is $+\infty$.

\nrmks \ilitem{} We highlight that $\width(w)$ depends only on the (connected) source quiver $G$ rather than the image of any vertex or arrow.

\ilitem{} A string has width $0$ iff its domain has no arrows.

\ilitem{} A connected string has width $1$ iff its domain is an equioriented linear subgraph 
    $ \begin{tikzcd}[sep=small]
    \syllableplaceholder \ar[r]
      \& \syllableplaceholder \ar[r]
        \& \cdots \ar[r]
          \& \syllableplaceholder
            \end{tikzcd}$ (containing at least one arrow).

\ilitem{} Provided $G$ contains at least one arrow, $\width w$ is the number of distinct equioriented linear subgraphs of $G$ whose source and target are respectively source and sink vertices of $G$.

\nrex In our example algebra, $P_2$ is a string module $\Str\big(
\begin{tikzpicture}[scale=0.75]
  \begin{scope}[minimum size=0, inner sep=1]
    % Vertices of string graph
    \draw (0,0) node[](v1){$ 1$};
    \draw (1,0) node[](v2){$ 1$};
    \draw (2,0) node[](v3){$ 2$};
    \draw (3,0) node[](v4){$ 2$};
    \draw (4,0) node[](v5){$ 2$};
    \draw (5,0) node[](v6){$ 2$};
    %
    % Arrows of string graph
    \draw[->] (v2) --node[pos=0.5,above](){$\alpha$} (v1);
    \draw[->] (v3) --node[pos=0.5,above](){$\gamma$} (v2);
    \draw[->] (v3) --node[pos=0.5,above](){$\delta$} (v4);
    \draw[->] (v4) --node[pos=0.5,above](){$\delta$} (v5);
    \draw[->] (v5) --node[pos=0.5,above](){$\delta$} (v6);
  \end{scope}
  \end{tikzpicture}
\big)$ with width $2$, while $\Str\big(
      \begin{tikzpicture}[scale=0.75]
      \begin{scope}[minimum size=0, inner sep=1]
        % Vertices of string graph
        \draw (0,0) node[](v1){$ 2$};
        \draw (1,0) node[](v2){$ 1$};
        \draw (2,0) node[](v3){$ 1$};
        \draw (3,0) node[](v4){$ 2$};
        \draw (4,0) node[](v5){$ 2$};
        \draw (5,0) node[](v6){$ 2$};
        %
        % Arrows of string graph
        \draw[->] (v1) --node[pos=0.5,above](){$\gamma$} (v2);
        \draw[->] (v2) --node[pos=0.5,above](){$\alpha$} (v3);
        \draw[->] (v4) --node[pos=0.5,above](){$\gamma$} (v3);
        \draw[->] (v4) --node[pos=0.5,above](){$\delta$} (v5);
        \draw[->] (v5) --node[pos=0.5,above](){$\delta$} (v6);
      \end{scope}
      \end{tikzpicture}
    \big)$
  is a string module of width $3$. There are string modules of arbitrary even width having form $\Str\big(
      \begin{tikzpicture}[scale=0.75]
      \begin{scope}[minimum size=0, inner sep=1]
        % Vertices of string graph
        \draw (0,0) node[](v1){$ 1$};
        \draw (1,0) node[](v2){$ 1$};
        \draw (2,0) node[](v3){$ 2$};
        \draw (3,0) node[](v4){$ 1$};
        \draw (4,0) node[](v5){$ 1$};
        \draw (5,0) node[](v6){$ 2$};
        \draw (6,0) node[](v7){$ 1$};
        \draw (7,0) node[](v8){$ \cdots$};
        \draw (8,0) node[](v9){$ 1$};
        \draw (9,0) node[](v10){$ 1$};
        \draw (10,0) node[](v11){$ 2$};
        \draw (11,0) node[](v12){$ 1$};
        %
        % Arrows of string graph
        \draw[->] (v2) --node[pos=0.5,above](){$\alpha$} (v1);
        \draw[->] (v2) --node[pos=0.5,above](){$\beta$} (v3);
        \draw[->] (v3) --node[pos=0.5,above](){$\gamma$} (v4);
        \draw[->] (v5) --node[pos=0.5,above](){$\alpha$} (v4);
        \draw[->] (v5) --node[pos=0.5,above](){$\beta$} (v6);
        \draw[->] (v6) --node[pos=0.5,above](){$\gamma$} (v7);
        \draw[->] (v8) -- (v7);
        \draw[->] (v8) -- (v9);
        \draw[->] (v10) --node[pos=0.5,above](){$\alpha$} (v9);
        \draw[->] (v10) --node[pos=0.5,above](){$\beta$} (v11);
        \draw[->] (v11) --node[pos=0.5,above](){$\gamma$} (v12);
      \end{scope}
      \end{tikzpicture}
    \big)$. Emulating this construction but omitting the leftmost vertex and incident arrow, we construct string modules of arbitrary odd width.
    
\npass The following proposition and corollary establish that we can approximate the vector-space dimension of a string module using width.
    
\begin{proposition} Let $m\geq 0$ be the smallest integer such that $J^m=0$; here $J \ideal A$ is the Jacobson radical of $A$. If $\Str(w)$ is a nonzero string module for $A$, then
$$
\width( \Str w ) + 1
\leq
\kdim ( \Str w )
<
m\width( \Str w ) + 1
\tstop
$$
\end{proposition}

\proof If $w$ is empty, then the inequality becomes $0 \leq 0 < 1$, a triviality. Otherwise it suffices to assume that $w \colon \ilarrow{G}{}{Q}$ is indecomposable because width and vector-space dimension are additive with respective to indecomposables. If $G$ contains no arrows, the inequality becomes $1 \leq 1 < m+1$: another triviality. Hence assume that $G$ contains at least one arrow. Since $G$ is connected and linear, $|G_0| = |G_1|+1$.

For any equioriented linear subgraph
$ \begin{tikzcd}[sep=small]
    \syllableplaceholder \ar[r, "x_1"]
      \& \syllableplaceholder \ar[r, "x_2"]
        \& \cdots \ar[r, "x_\ell"]
          \& \syllableplaceholder
\end{tikzcd}$ featuring at least one arrow, condition \ref{defn:string-graph}\ref{subdefn:linear-subgraphs-linearly-independent} assures that $w(x_1 \cdots x_l)$ has nonzero residue in $A$; thus $1\leq l < m$. We know that $G$ is the combination of $\width w$ maximal such subgraphs, which implies $\width w \leq |G_1| < m \width(w) $. The result follows immediately since $\kdim(\Str w)=|G_0|=|G_1|+1$.\qed

\begin{corollary}\label{cor:width-approximates-kdim}
If $(w_k)_{k\geq 0}$ is a sequence of string graphs for $A$, then either of
$$\big( \width(\Str w_k) \big)_{k\geq 0}
\text{\hspace{1cm}and\hspace{1cm}}
\big(\kdim(\Str w_k) \big)_{k\geq 0}
$$ diverges to $+\infty$ iff the other does. Similarly, one sequence is bounded iff the other is.
\end{corollary}

\section{Main results}

We now return to our core focus. In this section, we explore the syzygy behaviour of string modules and the ramifications thereof in the derived category of an SB algebra. We also provide compelling evidence towards our conjectured result, that the injective (string) modules of an SB algebra have finite syzygy type (a term that we define in Passage \ref{def:syzygy-finite}).

\setcounter{subsection}{-1}
\subsection{Results in context}

To help the reader contextualise our results, we discuss in this subsection certain properties that the syzygies of a module may satisfy and, in particular, the impact of canonical classes of modules satisfying them.

Many syzygy patterns, generally to do with some notion of repetition, have been studied extensively before, often in the context of the finitistic dimension conjectures as in \cite{GHZ98}. Rickard picked up this line of research in his investigation of triangulated subcategories of $\derivedcat(A) $ \cite[\S7]{Ric19}, albeit using the dual notion of cosyzygies rather than syzygies. There, he pays particular attention to the patterns seen amongst cosyzygies of simple or projective modules.

We introduce several interrelated patterns below and explore circumstances in which certain modules exhibit them. None of the proofs appearing here are original~-- they all appear in Rickard's article \cite{Ric19} if not earlier~-- but the arguments serve as so strong a motivation for the work in this thesis that their omission would be a sin. Once we have an abstract grasp of the syzygy properties of interest, we then narrow our focus onto the simple and injective modules, motivating and explaining in broad terms the results that follow in this thesis.

\centresubsec{Syzygy patterns}

\nlpass{Syzygy repetition index} \cite[Def 2.1]{GHZ98} (but compare also \cite[Def 7.7]{Ric19}).\hspace{1ex} Let $X$ be a finite-dimensional $A$-module.

If there exists a $k \geq 0$ such that every indecomposable nonprojective direct summand of $\syzygy^k X$ occurs as a direct summand of $\syzygy^r X$ for infinitely many $r$, then the \define{syzygy repetition index} of $M$ is the least such $k$; if none exists, then this repetition index is $\infty$.

\nlpass{An equivalent formulation} It will prove useful in the sequel to reformulate the above definition. To this end, fix some finite-dimensional $A$-module $X$ as above. 

For $k \geq 0$, let \define{$\mathcal{A}_k$} denote the set of isomorphism classes of indecomposable, nonprojective modules $Y$ such that $Y$ is isomorphic to a direct summand of $\syzygy^r X$ for some $r \geq k$ and, similarly, let \define{$\mathcal{B}_k$} denote the set of isomorphism classes of indecomposable, nonprojective modules $Y$ for which
$Y$ is isomorphic to a direct summand of $\syzygy^r X$ for some $r \leq k$.

(We emphasise that the only difference between these two properties is the orientation of the inequality sign. Informally, we may remember that $\mathcal{A}_k$ contains the isoclasses seen at or $\mathcal{A}$fter the $k$th syzygy, while $\mathcal{B}_k$ contains those seen at or $\mathcal{B}$efore the $k$th syzygy.)

We evidently have that $\mathcal{A}_{k} \supseteq \mathcal{A}_{k+1}$ and $\mathcal{B}_k \subseteq \mathcal{B}_{k+1}$ for each $k$, that each successive set difference $\mathcal{A}_k \setminus \mathcal{A}_{k+1}$ and $\mathcal{B}_{k+1} \setminus \mathcal{B}_k $ is finite and that clearly $\mathcal{B}_0$ itself is finite. Furthermore, the $\mathcal{A}_k$ and $\mathcal{B}_k$ relate to one another by the following system of inclusions:
\begin{equation}\label{eqn:inclusions-of-aks-and-bks}
\mathcal{B}_0 \subseteq \mathcal{B}_1 \subseteq \mathcal{B}_2 \subseteq \cdots \subseteq \bigcup_{r \geq 0} \mathcal{B}_r = \mathcal{A}_0 \supseteq \mathcal{A}_1 \supseteq \mathcal{A}_2 \supseteq \cdots \supseteq \bigcap_{r \geq 0} \mathcal{A}_r\tstop
\end{equation}
(We know that $\mathcal{A}_0$, equalling a countable union $\mathcal{B}_0 \cup \big(\bigcup_{r > 0} (\mathcal{B}_{r+1}\setminus\mathcal{B}_r)\big)$ of finite sets, is countable. We deduce that all the $\mathcal{B}_k$ are countable too.)

Let $\mdefine{\mathcal{A}_{\omega}} \ceq \bigcap_{r \geq 0} \mathcal{A}_r$. We easily find that $\mathcal{A}_\omega$ comprises the isoclasses of those modules appearing as direct summands of $\syzygy^r X $ syzygies of $X$ for infinitely many indices $r$.

We find that $X$ has syzygy repetition index $k$ iff $\mathcal{A}_k = \mathcal{A}_\omega$ (and $k$ is minimal with this property).

\begin{proposition}\label{prop:fin-syz-rep-ind-implies-in-coloc}
\emph{(Dual in statement and proof to \cite[Prop 7.8]{Ric19}.)\hspace{1ex}}
If $X$ has syzygy repetition index $k < \infty$, then $X \in \Coloc(\Projcat A)$.
\end{proposition}

\proof Let $X_0, X_1, X_2, \dots$ be a countable (possibly finite) sequence of nonisomorphic indecomposable modules representing all the isomorphism classes in $\mathcal{A}_k = \mathcal{A}_\omega$. Write $\mathcal{X} \ceq \{X_0, X_1, X_2, \dots \}$
%, and for each $r \geq 0$ write
% $\mathcal{X}_r \ceq \{X' \in \mathcal{X} \colon X' \text{ is isomorphic to a direct summand of } \syzygy^s X \text{ for some } s \geq r \}$.
% These form a descending sequence: $\mathcal{X}_0 \supseteq \mathcal{X}_1 \supseteq \mathcal{X}_2 \supseteq \cdots $.
%
%We deduce that the aforementioned descending sequence actually stabilises immediately: $\mathcal{X}_r = \mathcal{X}_0$ for all $r \geq 0$.

Let $Z \ceq ( \prod_{t \geq 0}  X_t)^{\prod \omega}$ be the direct product of $\omega$ many copies of each indecomposable in $\mathcal{X}$. Evidently $ \syzygy^{k+1} X$ itself is a direct factor of $Z$.

By taking the direct product
$
\prod_{r \geq k+1} \big(
\begin{tikzcd}[sep=small, cramped]
  0 \ar[r]
    \& \syzygy(\syzygy^r X) \ar[r]
      \& \projcover( \syzygy^r X ) \ar[r]
        \& \syzygy^r X \ar[r]
          \& 0
\end{tikzcd}
\big)
$
of each of the syzygy short exact sequences of the $\syzygy^r X$, we obtain a short exact sequence whose leftmost and rightmost terms are each a direct product of countably many copies of each indecomposable in $\mathcal{X}$; that is, a short exact sequence having the form
$\big( \begin{tikzcd}[sep=small, cramped]
  0 \ar[r]
    \& Z \ar[r]
      \& P \ar[r]
        \& Z \ar[r]
          \& 0
\end{tikzcd} \big)$ where $P$ is projective.

We may splice this sequence together with itself to obtain a exact sequence
$$\big( \begin{tikzcd}[sep=small, cramped]
  0 \ar[r]
    \& Z \ar[r]
      \& P \ar[r]
        \& P \ar[r]
          \& P \ar[r]
            \& P \ar[r]
              \& \cdots
\end{tikzcd}
\big)
$$ in $\complexcat(\Modcat A)$, which yields a quasiisomorphism in $\complexcat(\Modcat A)$ between the complex containing only $Z$ in degree $0$ and the complex $\big( \begin{tikzcd}[sep=small, cramped]
  \cdots \ar[r]
    \& 0 \ar[r]
      \& 0 \ar[r]
        \& P \ar[r]
          \& P \ar[r]
            \& P \ar[r]
              \& \cdots
\end{tikzcd}
\big)$ supported in nonnegative degrees. The latter is a bounded-below complex of projectives, which belongs to $\Coloc(\Projcat A)$ by part \ref{subprop:bounded-below-complexes-in-coloc} of Proposition \ref{prop:properties-of-coloc-subcats}. The former, being quasiisomorphic to it, belongs there too by part \ref{subprop:qisomic-in-coloc} of the same result. Since $\syzygy^{k+1} X$ is a direct factor of $Z$, part \ref{subprop:direct-factors-in-coloc} of this proposition shows $X \in \Coloc( \Projcat A )$, as claimed.

The final step is to observe that the standard syzygy short exact sequences may be spliced together to give an exact sequence
$
\big(\begin{tikzcd}[cramped,sep=small]
  0 \ar[r]
    \& \syzygy^{k+1} X \ar[r]
      \& \projcover(\syzygy^k X) \ar[r]
        \& \projcover(\syzygy^{k-1} X) \ar[r]
          \& \cdots \ar[r]
            \& \projcover( \syzygy^1 X) \ar[r]
              \& \projcover X \ar[r]
                \& X \ar[r]
                 \& 0
\end{tikzcd}\big)
$.\linebreak
Similarly to the previous paragraph, this exact complex corresponds to a quasiisomorphism in $\complexcat(A)$ between the concentration in degree $0$ of $X$ and the complex
$$
\begin{tikzcd}[sep=small, cramped]
  \cdots \ar[r]
    \& 0 \ar[r]
      \& \syzygy^{k+1} X \ar[r]
        \& \projcover(\syzygy^k X) \ar[r]
          \& \projcover(\syzygy^{k-1} X ) \ar[r]
            \& \cdots \ar[r]
              \& \projcover(\syzygy^1 X) \ar[r]
                \& \projcover X \ar[r]
                  \& 0 \ar[r]
                    \& 0 \ar[r]
                      \& \cdots
\end{tikzcd}
$$
supported in degrees $-(k+1) \leq r \leq 0$. Since $X$ is quasiisomorphic to a bounded complex of members of $\Coloc(\Projcat A)$, we conclude that $X \in \Coloc(\Projcat A)$ as well. This completes the proof.
\qed

\nlpass{Weakly periodic} A finite-dimensional module $X$ is \define{weakly periodic} if $X$ is isomorphic to a direct summand of $\syzygy^k X$ for some $k \geq 1$.

\begin{proposition}\label{prop:weakly-periodic-in-coloc}
If $X$ is weakly periodic, then $X \in \Coloc(\Projcat A)$. 
\end{proposition}

\proof To fix notation, suppose that $k \geq 0$ and $Y \in \modcat A$ satisfy $\syzygy^k X \isom X \oplus Y$.

Observe that for $r\geq 0$ we have
$
\syzygy^{r+k} X \isom \syzygy^r(\syzygy^k X) \isom \syzygy^r( X \oplus Y ) \isom ( \syzygy^r X ) \oplus (\syzygy^r Y)
$, and so $\syzygy^r X$ is a direct summand of $\syzygy^{r+k} X$. Iteration of this argument shows that $\syzygy^r X$ is a direct summand of $\syzygy^{r+\ell k} X$ for any $\ell>0$. In particular, any indecomposable direct summand of any $\syzygy^r X$ is isomorphic to a direct summand of $\syzygy^s X$ for infinitely many $s \geq 0$.

We deduce that, in the notation of (\ref{eqn:inclusions-of-aks-and-bks}), $\mathcal{A}_0 = \mathcal{A}_\omega$. Therefore $X$ has syzygy repetition index $0 < \infty$. By Proposition \ref{prop:fin-syz-rep-ind-implies-in-coloc}, $X \in \Coloc(\Projcat A)$. \qed

\nlpass{Syzygy finiteness of a module}\label{def:syzygy-finite} A finite-dimensional module $X \in \modcat A$ is \define{syzygy finite} if there exists a module $Y \in \modcat A$
% and an integer $m \geq 0$
such that
$
\add\{\syzygy^r X \in \modcat A \colon r \geq 0\} = \add Y.
$

% In the notation of (\ref{eqn:inclusions-of-aks-and-bks}), $X$ is syzygy finite iff the ascending sequence of $\mathcal{B}_r$ stabilises iff $\mathcal{A}_0$ is finite.

\begin{proposition} \emph{(Dual in statement and proof to \cite[Prop 7.2]{Ric19}.)\hspace{1ex}}
If $X$ is syzygy finite, then $X \in \Coloc(\Projcat A)$.
\end{proposition}

\proof For $X$ to be syzygy finite means precisely that $\mathcal{A}_0$ is finite in the notation of (\ref{eqn:inclusions-of-aks-and-bks}). It follows that the descending sequence $(\mathcal{A}_r)_{r \geq 0}$ must stabilise at some index $k$ (assumed minimal), and so $X$ has syzygy repetition index $k < \infty$. By Proposition \ref{prop:fin-syz-rep-ind-implies-in-coloc}, $X$ belongs to $\Coloc(\Projcat A)$, as claimed.

\centresubsec{Simple modules and injective modules}

\npass First, let us dualise a proof of Rickard's to show that the only colocalising subcategory of $\derivedcat(A)$ to contain the injective modules is $\derivedcat(A)$ itself. The only part of Rickard's argument that needs any real adjustment is the appeal to the Mittag-Leffler condition.

\begin{proposition}\label{prop:injs-cogenerate}
\emph{(Dual in statement, and mostly dual in proof, to \cite[Prop 2.2]{Ric19}.)\hspace{0.5cm}} Let $\categoryfont{C} \subseteq \derivedcat(A)$ be a colocalising subcategory and assume that all injective modules belong to $\categoryfont{C}$. Then $\categoryfont{C}=\derivedcat(A)$.
\end{proposition}

\proof Let $X$ be any $A$-module. By standard results \cite[\S2.3]{Wei94}, $X$ has an injective resolution
$
\big(
\begin{tikzcd}[sep=small, cramped]
  0 \ar[r]
    \& M \ar[r]
      \& E^0 \ar[r]
        \& E^1 \ar[r]
          \& E^2 \ar[r]
            \& \cdots
\end{tikzcd}
\big)$ which, similarly to above, induces a quasiisomorphism in $\complexcat(A)$ between the concentration in degree $0$ of $M$ and the complex $\big(
\begin{tikzcd}[sep=small, cramped]
  \cdots \ar[r]
    \& 0 \ar[r]
      \& E^0 \ar[r]
        \& E^1 \ar[r]
          \& E^2 \ar[r]
           \& \cdots
\end{tikzcd}
\big)$
having term $E^0$ in degree $0$. The latter is a bounded-below complex all of whose terms belong to $\categoryfont{C}$, and so belongs it $\categoryfont{C}$ as well by Proposition \ref{prop:properties-of-coloc-subcats}\ref{subprop:bounded-below-complexes-in-coloc}. We deduce that $X$, being quasiisomorphic, lies in $\categoryfont{C}$ too.

Now let $X^\bullet \ceq \big(
\begin{tikzcd}[sep=small,cramped]
\cdots \ar[r]
  \& X^{-1} \ar[r, "d^{-1}"]
    \& X^0 \ar[r, "d^0"]
      \& X^1 \ar[r, "d^1"]
        \& X^2 \ar[r]
          \& \cdots
\end{tikzcd}
\big)$ be any complex. For any $r \leq 0$, let $$\mdefine{\sigma^{\geq r} X^\bullet} \ceq \big(
\begin{tikzcd}[cramped]
\cdots \ar[r]
  \& 0 \ar[r]
    \& \coker d^{r-1} \ar[r]
      \& X^{r+1} \ar[r, "d^{r+1}"]
        \& X^{r+2} \ar[r, "d^{r+2}"]
          \& X^{r+3} \ar[r]
            \& \cdots
\end{tikzcd}
\big)$$
be the \define{good truncation of $X^\bullet$ above $r$}. The $k$th term of this complex is: $X^k$ when $k > r$; $\coker d^{r-1}$ when $k = r$ and; $0$ otherwise. Aside from the $r$th boundary map $\ilarrow{\coker d^{r-1}=(X^{r}/\im d^{r-1})}{}{X^{r+1}}$ given by $ \ilmapsto{ x + \im{d^{r-1}} }{}{ xd^r }$ (which is well-defined, since $\im d^{r-1} \leq \ker d^r $), all other boundary morphisms are either $d^k$ or $0$. These are bounded-below truncations of modules that lie in $\categoryfont{C}$ so, by part \ref{subprop:bounded-below-complexes-in-coloc} of Proposition \ref{prop:properties-of-coloc-subcats}, $\categoryfont{C}$ contains these truncations.

For any $r \leq 0$ there are morphisms $\ilarrow{ \sigma^{r-1} X^\bullet}{}{\sigma^r X^\bullet}$ given as follows:
$$
\begin{tikzcd}[column sep=small,cramped]
\sigma^{r-1} X^\bullet \ar[d]
  \& \big(\cdots \ar[r]
    \& 0 \ar[r] \ar[d]
      \& \coker d^{r-2} \ar[r] \ar[d]
        \& X^r \ar[r] \ar[d, two heads]
          \& X^{r+1} \ar[r] \ar[d]
            \& X^{r+2} \ar[r] \ar[d]
              \& \cdots\big)
              \\
\sigma^r X^\bullet
  \& \big(\cdots \ar[r]
    \& 0 \ar[r]
      \& 0 \ar[r]
        \& \coker d^{r-1} \ar[r]
          \& X^{r+1} \ar[r]
            \& X^{r+2} \ar[r]
              \& \cdots\big)\tstop
\end{tikzcd}
$$
All components of this morphism of complexes are surjections, being variably identity maps $\ilarrow{X^k}{}{X^{k}}$, canonical quotients $\ilepi{X^r}{}{\coker d^{r-1} = (X^r/\im d^{r-1})}$ or maps onto $0$. We deduce that the resulting inverse system
$\big(
\begin{tikzcd}[sep=small, cramped]
\cdots \ar[r]
  \& \sigma^2 X^\bullet \ar[r]
    \& \sigma^1 X^\bullet \ar[r]
      \& \sigma^0 X^\bullet
\end{tikzcd}
\big)$ satisfies the Mittag-Leffler condition. Using part \ref{subprop:inverse-system-in-coloc} of Proposition \ref{prop:properties-of-coloc-subcats}, we find that the inverse limit of this system of truncations belongs to $\categoryfont{C}$. However this inverse limit is patently $X^\bullet$.

We deduce that $\categoryfont{C}$ contains every complex and so conclude that $\categoryfont{C} = \derivedcat(A)$, as claimed. \qed

\npass In the above argument, what is crucial is that an arbitrary module $X$ can demonstrably be ``built'' within the colocalising subcategory using known members. The problem of building an arbitrary module in the preceding proof was solved by appealing to the existence of injective resolutions; with $X$ found, no further mention of the injective modules was necessary. However, as Rickard explores, one could alternatively build $X$ from the simple modules by a familiar semisimple filtration.

% Rickard then notes that one can build the injective modules within a colocalising subcategory from the simple modules.

\begin{proposition}
\emph{(Dual in statement and proof to \cite[Lem 6.1]{Ric19}.)\hspace{0.5cm}} Let $\categoryfont{C} \subseteq \derivedcat(A)$ be a colocalising subcategory and assume that all simple modules belong to $\categoryfont{C}$. Then $\categoryfont{C} = \derivedcat(A)$.
\end{proposition}

\proof Any semisimple module $Y$ is a product of simple modules so $Y \in \Coloc(\Projcat A)$ by Proposition \ref{subprop:concrete-coloc-subcat-properties}\ref{subprop:products-in-coloc}. An arbitrary module $X$ can be built from semisimple modules by finitely many extensions. Part \ref{subprop:2-out-of-3-in-coloc} of this same proposition therefore shows $X \in \Coloc(\Projcat A)$. Now the proof of Proposition \ref{prop:injs-cogenerate} may be copied verbatim, beginning with the second paragraph, to establish the claim.\qed

\npass Equipped with this knowledge, Rickard provides examples of algebras where the simples have finite (co)syzygy type. His reasoning has recourse to several classic results \cite[\S7]{Ric19} and one ad-hoc argument \cite[\S6]{Ric19}.

It is not generally true that the simples for SB algebras even have finite syzygy repetition index~-- the simple module $S_1$ of our running example algebra being a straightforward witness to this transgression~-- and so any effort to prove their syzygy finiteness is certainly doomed.

However, countless calculations by the author found that the injective modules over SB algebras do have finite syzygy type. Of course, injective modules $E$ that are also projective satisfy $\syzygy^1 E = 0$ and so are syzygy finite in a trivial way. This eliminates the pin modules, reducing the matter only to the study of injective string modules.

The primary thrust of our research was consequently to prove that injective string modules for SB algebra are syzygy finite: a statement we conjecture to be true, but cannot yet prove in its entirety. What follows in the remaining subsections of this chapter are our partial results towards this conjecture. In reverse order, let us sketch a crude outline of our results.

Our final main result, in Subsection \ref{subsec:sb-algs-with-few-verts}, is that this conjecture is true when $A$ has at most $2$ isomorphism classes of simple modules. Preceding that discussion, in Subsection \ref{subsec:desc-and-injective-syllables} we show a technical result that deploys our formalism for syzygies in a fundamental way, namely to demonstrate that there is behaviour in the syzygy fabric suggestive of being syzygy-infinite and that the symbols used to describe injective string modules do not exhibit this behaviour. Lastly~-- which is to say initially~-- we provide in Subsection \ref{subsec:noncycles-in-pin-graphs} a subclass of SB algebras for which a more general result than our conjecture is true: that in this subclass all modules are syzygy-finite (in a uniform way).

The algebras in this subclass are those for which a particular associated quiver is acyclic. Acyclic components of that quiver provide some additional interest, since we can show the existence of syzygy patterns for simples on such components. Under mild conditions, we prove that those simples belong to $\Coloc(\Projcat A)$ and thus dually to $\Loc(\Injcat A)$. (Briefly, let us remark that, when all components of this quiver associated to $A$ are acyclic, our proofs therefore demonstrate independently from one another that the simples and the injectives belong to $\Coloc(\Projcat A)$ and to $\Loc(\Injcat A)$.)

\subsection{Noncycles in pin graphs}\label{subsec:noncycles-in-pin-graphs}

\npass{} Recall that $Q$ denotes the ordinary quiver of the special biserial algebra $A = \kk Q/\langle \rho \rangle$. Additionally recall that any syllable $\bdp$ is peak-compatible with a unique stationary syllable, namely $ \bde_{ (\source( \bdp ))^\dagger } $, and that the sidestep operation $\ilmapsto{ \bdp }{}{ \bde_{ (\source( \bdp ))^\dagger } }$ is denoted $\sidestep$.

% Further recall that the pin graph $\pingraph_A = \pingraph$ of $A$ has vertex set $Q_0$ and an arrow $\ilarrow{ i }{}{ j }$ iff there exists a pin module $P$ with $\head P \isom S_i $ and $\soc P \isom S_j $.

\nlpass{Pin graph} The \define{pin graph $\pingraph_A$} of a SB algebra $A$ is the quiver with vertex set $Q_0$ and with an arrow $\ilarrow{i}{}{j}$ iff there exists a pin module $P$ with $\head P \isom S_i$ and $\soc P \isom S_j$.

It is sub-$1$-regular but not necessarily connected. Evidently, $\pingraph_{A^\op}=(\pingraph_A)^\op$.

We know that $\pingraph_A$ is discrete (has no arrows) iff $A$ is monomial. If it is $1$-regular, then $A$ is selfinjective.

\nrex The pin graph of our running example algebra is $\big(\begin{tikzcd}[sep=small] 1\ar[loop left] \& 2 \end{tikzcd}\big)$.

\centresubsec{Simples not on cycles in $\pingraph_A$}

\npass Recall that a finite-dimensional module $X$ is \define{weakly periodic} if $X$ is isomorphic to a direct summand of $\syzygy^k X$ for some $k \geq 1$.

\nexample We will shortly prove a lemma about the existence of certain weakly periodic uniserial modules. This general argument is much more straightforward than the prose of the proof suggests, and is best illustrated with a concrete example. To this end, recall our running example algebra
$$
A \ceq \kk\big(
\begin{tikzcd}
  1 \ar[loop left, "\alpha"] \ar[r, "\beta", shift left]
    \& 2 \ar[l, "\gamma", shift left] \ar[loop right, "\delta"]
\end{tikzcd}
\big)
/
\big\langle
  \alpha^2,
  \beta\delta,
  \gamma\beta,
  \delta\gamma,
  \alpha\beta\gamma-\beta\gamma\alpha,
  \gamma\alpha\beta,
  \delta^4
\big\rangle\tstop
$$
Let $(\OQ, N, C, \dagger)$ be the unique permissible data for $A$. (In particular, $\OQ$ is the unique overquiver for $A$)

We see that $\gamma\beta, \beta\delta, \delta\gamma$ have zero residue in $A$. Taken in order, these arrows form the cycle
$$%\big(\raisebox{1ex}{
  \begin{tikzcd}[sep=small]
    2 \ar[r, "\gamma"]\ar[rr, phantom, ""{coordinate, name=C}]
      \& 1 \ar[r, "\beta" ]
        \& 2 \ar[ll,  rounded corners, "\delta", to path={
            -- ([xshift=1.5ex]\tikztostart.east)
            |- ([yshift=-2ex]C)[pos=1]\tikztonodes
            -| ([xshift=-1.5ex]\tikztotarget.west)
            -- (\tikztotarget.west)
            }]
\end{tikzcd}
%} \big)
% =
% \big(\raisebox{1ex}{
%   \begin{tikzcd}[sep=small]
%     2 \ar[r, "\delta"]\ar[rr, phantom, ""{coordinate, name=C}]
%       \& 2 \ar[r, "\gamma" ]
%         \& 1 \ar[ll,  rounded corners, "\beta", to path={
%             -- ([xshift=1.5ex]\tikztostart.east)
%             |- ([yshift=-2ex]C)[pos=1]\tikztonodes
%             -| ([xshift=-1.5ex]\tikztotarget.west)
%             -- (\tikztotarget.west)
%             }]
% \end{tikzcd}
% } \big)
$$
in the presenting quiver. They can be lifted to arrows in $\OQ$ that we also call $\gamma,\beta,\delta$, whose sources have stationary syllables $\bde_{\source(\gamma)}, \bde_{\source(\beta)}, \bde_{\source(\delta)}$ respectively.

The composition $\descent\sidestep$ permutes these stationary syllables:
$\begin{tikzcd}[sep=small]
  \bde_{\source(\gamma)} \ar[r, mapsto, "\descent\sidestep"]
    \& \bde_{\source(\beta)} \ar[r, mapsto, "\descent\sidestep"]
      \& \bde_{\source(\delta)} \ar[r, mapsto, "\descent\sidestep"]
        \& \bde_{\source(\gamma)}
\end{tikzcd}$. In particular this means that $\bde_{\source(\gamma)}$ is $(\descent\sidestep)$-periodic, but hence also that $\bde_{\source(\gamma)}\descent$ is $(\sidestep\descent)$-periodic. We observe that the direct summand $\gamma A$ of $\rad P_2$ can be represented by a strip whose only nonblank entry is the syllable $\bde_{\source(\gamma)}\descent=\big( \ilsyll{\circ}{\alpha}{1}\big)$.

The sequence of syllables, starting with $ \bde_{\source(\gamma)}\descent $ and alternately applying $\sidestep$ and $\descent$, can be seen in the syzygy fabric when calculating the syzygies of $\gamma A$. Starting from the $\bde_{\source(\gamma)}\descent=\big( \ilsyll{\circ}{\alpha}{1}\big)$ in the strip representation of $\gamma A$ in Figure \ref{fig:following-a-boundary}, we move to the right when applying $\sidestep$ and downwards when applying $\descent$.
\begin{thesisfigure}
  \centering
  
  \begin{tikzpicture}[scale=\standardscale]
    % Highlighting for cells
    \newcommand{\highlightcell}[1]{
      \draw[fill=promptColor!30] #1 ++(-1,1) rectangle ++(2,-2);
    }
    \highlightcell{(2,0)};
    \highlightcell{(4,0)};
    \highlightcell{(4,-2)};
    \highlightcell{(6,-2)};
    \highlightcell{(6,-4)};
    \highlightcell{(8,-4)};
    \highlightcell{(8,-6)};
    \highlightcell{(10,-6)};
    \highlightcell{(10,-8)};
    \highlightcell{(12,-8)};
    \highlightcell{(12,-10)};
    \highlightcell{(14,-10)};
  
    % Cells in 1st (top) row
    \oddleftedgecell{(-2,0)}{1 -1};
    \evencell{(0,0)}{1 0}{};
    \oddcell{(2,0)}{1 1}{};
    \evencell{(4,0)}{1 2}{};
    \oddcell{(6,0)}{1 3}{};
    \evenrightedgecell{(8,0)}{1 4};
    % Cells in 2nd row
    \oddleftedgecell{(-4,-2)}{2 -1};
    \evencell{(-2,-2)}{2 0}{};
    \oddcell{(0,-2)}{2 1}{};
    \evencell{(2,-2)}{2 2}{};
    \oddcell{(4,-2)}{2 3}{};
    \evencell{(6,-2)}{2 4}{};
    \oddcell{(8,-2)}{2 5}{};
    \evenrightedgecell{(10,-2)}{2 6};
    % Cells in 3rd row
    \oddleftedgecell{(2,-4)}{3 -1};
    \evencell{(4,-4)}{3 0}{};
    \oddcell{(6,-4)}{3 1}{};
    \evencell{(8,-4)}{3 2}{};
    \oddcell{(10,-4)}{3 3}{};
    \evenrightedgecell{(12,-4)}{3 4};
    % Cells in 4th row
    \oddleftedgecell{(4,-6)}{4 -1};
    \evencell{(6,-6)}{4 0}{};
    \oddcell{(8,-6)}{4 1}{};
    \evencell{(10,-6)}{4 2}{};
    \oddcell{(12,-6)}{4 3}{};
    \evenrightedgecell{(14,-6)}{4 4};
    % Cells in the 5th row
    \oddleftedgecell{(2,-8)}{5 -1};
    \evencell{(4,-8)}{5 0}{};
    \oddcell{(6,-8)}{5 1}{};
    \evencell{(8,-8)}{5 2}{};
    \oddcell{(10,-8)}{5 3}{};
    \evencell{(12,-8)}{5 4}{};
    \oddcell{(14,-8)}{5 5}{};
    \evenrightedgecell{(16,-8)}{5 6};
    % Cells in 6th (bottom) row
    \oddleftedgecell{(8,-10)}{6 -1};
    \evencell{(10,-10)}{6 0}{};
    \oddcell{(12,-10)}{6 1}{};
    \evencell{(14,-10)}{6 2}{};
    \oddcell{(16,-10)}{6 3}{};
    \evenrightedgecell{(18,-10)}{6 4};
    
    \begin{scope}[node distance=6]
      % Syllables in 1st (top) row
      \oddbdysyll{1 1}{}{ $ \alpha $ }{}{}{};
      \eventrivbdysyll{1 2}{}{ $ e_{\source(\beta)} $ }{}{}{};
      % Syllables in 2nd row
      \oddtrivbdysyll{2 1}{}{ }{}{}{};
      \evenintsyll{2 2}{}{$\beta \gamma$}{};
      \oddintsyll{2 3}{}{ $ \gamma \alpha $ }{};
      \eventrivbdysyll{2 4}{}{ $ e_{\source(\delta)} $ }{}{}{};
      % Syllables in 3rd row
      \oddbdysyll{3 1}{}{ $ \delta^2 $ }{}{}{};
      \eventrivbdysyll{3 2}{}{ $ e_{\source(\gamma)} $ }{}{}{};
      % Syllables in 4th row
      \oddbdysyll{4 1}{}{ $ \alpha  $ }{}{}{};
      \eventrivbdysyll{4 2}{}{ $ e_{\source(\beta)} $ }{}{}{};
      % Syllables in 5th row
      \oddtrivbdysyll{5 1}{}{ }{}{}{};
      \evenintsyll{5 2}{}{$\beta \gamma$}{};
      \oddintsyll{5 3}{}{ $ \gamma \alpha $ }{};
      \eventrivbdysyll{5 4}{}{ $ e_{\source(\delta)} $ }{}{}{};
      % Syllables in 6th (bottom) row
      \oddbdysyll{6 1}{}{ $ \delta^2 $ }{}{}{};
      \eventrivbdysyll{6 2}{}{ $ e_{\source(\gamma)} $ }{}{}{};
    \end{scope}
  \end{tikzpicture}
  
  \caption[%
    Following a boundary
  ]{%
    \label{fig:following-a-boundary}
    \textit{Following a boundary.} We represent $\gamma A$ as a strip module (1st row) and alternately step rightwards one cell and downwards one cell along the highlighted path, the entries we encounter being obtained by iteratively applying $\sidestep$ and $\descent$ with each step. We find that $\gamma A$ recurs (4th row), hence $\gamma A$ is weakly periodic. Similarly we find that $\delta A$ repeats (3rd and 6th row).
  }
\end{thesisfigure}

The stationary syllables we see along the way\footnote{ In the figure, all of the relevant stationary syllables are present. This is a white lie for the sake of illustration; in reality, the stationary syllables outside the 1st row would all be implied.} form a righthand boundary of successive strips. When $\bde_{\source(\gamma)}\descent=\big( \ilsyll{\circ}{\alpha}{1}\big)$ turns up for the second time, it is a lefthand boundary. The lefthand and righthand boundaries together delimit a strip, specifically one representing $\gamma A$. This concludes the calculation that $\gamma A$ is weakly periodic.

The same can be said for $\delta A$, the other direct summand of $\rad P_2$, since it also lies on a cycle whose pairwise products have zero residue in $A$. That cycle is the same one as before, and the relevant calculation also appears in (the bottom four rows of) Figure \ref{fig:following-a-boundary}.

\begin{lemma}\label{lem:weakly-periodic-uniserial-radical-summands} Recall that $\langle \rho \rangle$ denotes the defining ideal of $A = \kk Q/\langle \rho \rangle$.

Suppose that $i \in Q_0$ is a vertex lying on an $\ell$-cycle
$\big(\begin{tikzcd}[sep=small]
i \ar[r, "\alpha_1"]
  \& \cdot \ar[r, "\alpha_2"]
    \& \cdots \ar[r,"\alpha_\ell"]
      \& i
\end{tikzcd}\big)$ satisfying $\alpha_r \alpha_{r+1} \in \langle \rho \rangle$ for all $1 \leq r \leq \ell$ (subscripts understood modulo $\ell$), and suppose that $P_i$ is a string module. % Also write $\alpha \ceq \alpha_1$ for convenience.

Then $\alpha_1 A$ is a direct summand of $\syzygy^\ell(  \alpha_1 A)$.
\end{lemma}

\proof There exists a choice of permissible data $(\OQ,N,C,\dagger)$ for $A$ such that no path in $\OQ$ represents any $\alpha_r \alpha_{r+1}$ ($1 \leq r \leq \ell$, subscripts taken modulo $\ell$). This follows from the existence of the $\ell$-cycle of relations of length $2$.

Lift each $\alpha_r \in Q_1$ to an arrow of $\OQ$ that we also denote $\alpha_r$. Since each $\alpha_r \in \OQ_1$ does represent an arrow of the ground quiver (rather than an augmented arrow), we deduce that $\bde_{\source(\alpha_r)}\descent$ is defined for all $r$. From the existence of the cycle again, we deduce that $\bde_{\source(\alpha_1)} (\descent\sidestep)^r = \bde_{\source(\alpha_{r+1})}$.

In any patch where an $\bde_{\source(\alpha_r)}$ appears in the top row or is implied to be there, the cell underneath it contains either $\bde_{\source(\alpha_r)}\descent$ or the perturbed version thereof (and in the particular case of $\bde_{\source(\alpha_1)}$, the entry underneath must be $\bde_{\source(\alpha_1)}\descent$ because $P_i$ is a string module and so there is no perturbation). Because this pertubation only affects the target of the syllable rather than the source, both versions are taken by $\sidestep$ to $\bde_{\source(\alpha_r)}\descent\sidestep$

As $P_i$ is a string module, $\rad P_i$ is a direct sum of one or two uniserial modules. One of these direct summands is $\alpha_1 A$, which we can represent as a strip whose nonblank entries are the syllable $\bde_{\source(\alpha_1)}\descent$ in a peak along with the stationary syllable $\bde_{\source(\alpha_1)}\descent\sidestep = \bde_{\source(\alpha_2)}$. The argument of the previous paragraph establishes that, starting from $\bde_{\source(\alpha_2)}$ and successively moving to the peak neighbour of the supported child cell $r$ times as we did in Figure \ref{fig:following-a-boundary}, we reach cells whose present or implied entries are $ \bde_{\source(\alpha_{r+2})} $. 

After $\ell$ such diagonal moves, we find that $\bde_{\source(\alpha_{\ell+2})}$ occurs next to $\bde_{ \source(\alpha_1) }\descent$. Both of these are boundary syllables. Our $\ell$ moves have taken us down $\ell$ rows, none of which held virtual strips. We conclude that $\alpha_1 A$ is a direct summand of $\syzygy^\ell(\alpha_1 A)$, as claimed.\qed

\begin{theorem}\label{thm:acyclic-pin-simples-for-2-reg-quiv-in-subcats}
Suppose that the ordinary quiver of a SB algebra $A$ is $2$-regular. If $i \in Q_0$ does not lie on a cycle in the pin graph $\pingraph_A$ of $A$, then $S_i \in \Coloc( \Projcat A )$. Dually, $S_i \in \Loc( \Injcat A )$ also.
\end{theorem}

\proof Choose permissible data $(\OQ, N, C, \dagger)$ for $A$. For each $\alpha \in Q_1$, there exists a unique arrow $\beta \in Q_1$ with $\target( \alpha ) = \source( \beta )$ such that no $\OQ$-path represents $\alpha\beta$; this implies that $\alpha \beta $ has zero residue in $A$. This assignment $\ilmapsto{ \alpha }{}{ \beta }$ yields a permutation $ \pi \colon \ilarrow{ Q_1}{}{ Q_1 } $ of arrows. We know $Q_1$ is finite hence $\alpha\pi^\ell = \alpha$ for some $\ell > 0$.

Let $ m_i $ denote the maximal length of a path in $ \pingraph_A $ with source $ i $. We know $ m_i $ is finite because $ i $ does not lie on a cycle in $ \pingraph_A $. We proceed by induction on $ m_i $, starting with the basis case $ m_i=0 $ when $i$ is a sink vertex in $ \pingraph $ and so $ P_i $ is a string module.

Let $ \alpha_1, \alpha_2 \in Q_1$ be the distinct arrows with $\source( \alpha_1 )=\source( \alpha_2 )=i$, so that $ \rad P_i = \alpha_1 A \oplus \alpha_2 A $. In the canonical short exact sequence $
\big(
\begin{tikzcd}[sep=small]
  0 \ar[r]
    \& \rad P_i \ar[r]
      \& P_i \ar[r]
        \& \head P_i \ar[r]
          \& 0
\end{tikzcd}
\big)
$, the middle term clearly belongs to $\Coloc( \Projcat A )$. By the two-out-of-three property, the righthand term $S_i$ belongs there iff the lefthand term does iff both summands $\alpha_r A$ of the lefthand term do. 

For $r \in \{1,2\}$, let $\ell_r > 0$ be the minimal integer such that $\pi^{\ell_r}$ fixes $\alpha_r$. Taking in turn the arrows in the forward $\pi$-orbit of $\alpha_r$, we obtain a cycle
$\big(\begin{tikzcd}[]
i \ar[r, "\alpha_r"]
  \& \cdot \ar[r, "\alpha_r \pi"]
    \& \cdots \ar[r,"\alpha_r \pi^{\ell_r}"]
      \& i
\end{tikzcd}\big)$
whose arrows compose pairwise to have zero residue in $A$. We deduce $\alpha_r A$ is weakly periodic by Lemma
\ref{lem:weakly-periodic-uniserial-radical-summands} and so $\alpha_r A \in \Coloc( \Projcat A ) $ by Proposition \ref{prop:weakly-periodic-in-coloc}. This concludes the proof of the basis case.

For the inductive step, suppose $S_j \in \Coloc( \Projcat A )$ for all vertices $j$ with $m_j < m_i$. There exists an arrow $\ilarrow{i}{}{j}$ in $\pingraph_A$ for such a vertex $j$. By definition, $A$ has a pin module $P_i$ with $ \soc P_i = S_j$. By assumption, the lefthand term and middle term of the short exact sequence $
\big( \begin{tikzcd}[sep=small]
  0 \ar[r]
    \& \soc P_i \ar[r]
      \& P_i \ar[r]
        \& P_i/\soc P_i \ar[r]
          \& 0
\end{tikzcd} \big)
$ belong to $\Coloc( \Projcat A )$, thus its righthand term does.

We know that $\soc P_i$ is a $1$-dimensional subspace of $A$ being spanned by the residue of either component $p,q$ of the commutativity relation $p-q$ having source $i$. This is a two-sided ideal of $A$, since for any arrow $\gamma \in Q_1$ we have $\gamma p \equiv \gamma q \equiv 0 $ and $p \gamma \equiv q \gamma \equiv 0$ in $A$, while $e_i \soc P_i = \soc P_i = (\soc P_i)e_j$ and $e_k \soc P_i = 0 = (\soc P_i) e_k$ for all $k \neq i$. Hence consider the quotient algebra $A' \ceq A/\soc P_i$, writing $P'_k$ for the indecomposable projective $A'$-module at $k$. Evidently $P'_i \isom P_i/\soc P_i$ is a string module, while $P'_k \isom P_k$ for all other $k \neq i$. (These are isomorphisms of $A$-modules.) From the previous paragraph we deduce that $P'_k \in \Coloc( \Projcat A )$ for all vertices $k$ (including $i$).

The ground quiver of $A'$ equals that of $A$~-- in particular it is $2$-regular~-- therefore we can apply the argument of the inductive step to $\syzygy^1 ( \head P'_i)$. Because all of the projective $A'$-modules involved in that argument belong to $\Coloc( \Projcat A )$, this very calculation shows $\head P'_i \isom \head P_i \isom S_i \in \Coloc( \Projcat A )$ too. This concludes the inductive step, and hence the first statement.

By employing the preceding argument on the simple $A^\op$-module at $i$ and then applying the standard duality $\dual$, we similarly show that $S_i \in \Loc(\Injcat A)$, as claimed.\qed

\centresubsec{Syzygy finiteness}

\nlpass{Syzygy finiteness of an algebra}\label{def:syzygy-finite-alg}% A finite-dimensional module $X \in \modcat A$ is \define{syzygy finite} if there exists a module $Y \in \modcat A$
% and an integer $m \geq 0$
%such that
%$
%\add\{\syzygy^r X \in \modcat A \colon r \geq 0\} = \add Y.
%$
An algebra $A$ is \define{$m$-syzygy finite}, for $m \geq 0$, if there exists $Y \in \modcat A$ such that
$$\{\syzygy^r X \in \modcat A \colon r \geq m,  X \in \modcat A \} = \add Y
$$
and, more generally, it is \define{syzygy-finite} if it is $m$-syzygy finite for some $m \geq 0$. Otherwise, it is \define{syzygy-infinite}.

\nexamples\label{ex:syz-fin-algs} \ilitem{} Representation-finite algebras are exactly the $0$-syzygy finite algebras.

\ilitem{} Any radical-square zero algebra is $1$-syzygy finite, with additive generator $Y=S_1 \oplus \cdots S_n$.

\ilitem{} Monomial algebras $A$ are $2$-syzygy finite \cite{ZH91}, with additive generator a direct sum of principal right ideals of $A$ generated by paths.

\nrex Our running example algebra is syzygy-infinite, as can be found by calculating the syzygies of $S_1$.

\nlpass{Interior of a strip} Let $w$ be a strip for $A$. Define the \define{interior of $w$} as
$$
  \mdefine{ \interior w }
  \ceq
  \{
    k \in \integers
    \colon
    w(k) \text{ is an interior (nonvirtual) syllable}
  \}
$$
and the \define{interior width of $w$} as
$ \mdefine{ \interiorwidth w } \ceq | \interior w |$

\nrex Recall the strips $v_1,v_2,v_3,v_4,v_5$ of Subsection \ref{subsec:preview-of-strips}. We have
$
\interiorwidth(v_1) = \interiorwidth(v_5) = 2$,
$\interiorwidth(v_2) = 4$ and 
$
\interiorwidth(v_3) = \interiorwidth(v_4) = 0
$.

\nexample If $w$ represents a simple module, uniserial module or indecomposable projective string module, then $\interiorwidth w = 0$.

\nrmks \ilitem{} Note that ``rounding off'' does not affect $\interior w$ or $\interiorwidth w$. Thus all strips representing the same string graph have the same interior width and so interior width is a well-defined invariant of a string module.

\ilitem{} Also note that
$
\interiorwidth w \leq \width w \leq \interiorwidth(w) + 2
$ for any strip $w$.
Therefore if $(w_m)_{m \geq 0}$ is a sequence of strips, then $(\interiorwidth w_m)_{m \geq 0}$ diverges to $+\infty$ iff $(\width w_m)_{m \geq 0}$ does.

\ilitem{} Because $\interior w$ is a union of valleys, $\interiorwidth w$ is even (whenever finite).

\ilitem{} Clearly $\interior w$ is an interval subset of $\integers$. We have $\interior w = (L,R)$ if the (potentially implied) left and right boundaries of $w$ appear in columns $L$ and $R$ respectively.

\nlpass{Boundaries of strips}\label{psg:late-wide-strips-share-interior-column} Suppose $w_r$ ($r\in\{0,1\}$) is a strip with left- and righthand boundaries respectively in columns $L_r$ and $R_r$, and that $w_1$ is a syzygy strip of $w_0$.

We have previously seen that left and right boundaries move at most one column outwards with each syzygy; that is, $L_0 - 1 \leq L_1$ and $R_1 \leq R_0 + 1$. This behaviour carries over to the extrema of $\interior w$: ie,
$$
\min(\interior w_0) - 1
\leq
\min( \interior w_1 )
\text{\hspace{1cm}and\hspace{1cm}}
\max( \interior w_1 )
\leq
\max( \interior w_0 ) + 1
$$
A trivial manipulation of each inequality gives
$$
\min(\interior w_0)
\leq
\min( \interior w_1 ) +1
\text{\hspace{1cm}and\hspace{1cm}}
\max( \interior w_1 ) - 1
\leq
\max( \interior w_0 )\tcomma
$$
from which we deduce
\begin{equation}\label{eqn:nesting-interiors}
  \interior w_0
  \supseteq
  \big[
    \min(\interior w_1)+1,
    \max( \interior w_1 )-1
  \big]
  \tcomma
\end{equation}
unless $\interiorwidth w_1 =0$ (when this inclusion is meaningless because the extrema do not exist). We mention that this inclusion holds vacuously when $\interiorwidth w_1 = 2$, because then the righthand interval is empty. Therefore, the inclusion (\ref{eqn:nesting-interiors}) is only significant if $\interiorwidth w_1 > 2$, in which case it further asserts $\interiorwidth w_0 \geq \interiorwidth w_1 - 2$

Let us iterate this argument. Suppose for $m \geq 1$ that $w_r$ ($r \in \{0,1,\dots,m\}$) are strips, that each $w_{r+1}$ is a syzygy strip of $w_r$, and that $ \interiorwidth w_m > 2m $ so that our claims have substance. We find that
$$
  \begin{array}{rcl}
    \interior w_{m-r}
      & \supseteq
        & \big[
          \min( \interior w_m ) + r,
          \max( \interior w_m ) - r
         \big]
          \\
       & \supseteq
         & \big[
           \min( \interior w_m ) + m,
           \max( \interior w_m ) - m
          \big]
  \end{array}
$$
for each $r \in \{0, 1, \dots, m\}$. The first inclusion implies that $\interiorwidth w_{m-r} \geq \interiorwidth w_m - 2r $. The second implies that that there are (at least) two columns common to all interiors $\interior w_r$. This situation is depicted in Figure \ref{fig:interiors-of-syzygies}.
\begin{thesisfigure}
  \centering
  
  \begin{tikzpicture}[scale=\standardscale]
    % Draw cells in 1st (top) row
    \oddleftedgecell{(4,0)}{1 2};
    \evencell{(6,0)}{1 3}{};
    \oddcell{(8,0)}{1 4}{};
    \evenrightedgecell{(10,0)}{1 5};
    % Draw cells in 2nd row
    \oddleftedgecell{(2,-2)}{2 1};
    \foreach \x in {2, 4}
      {\evencell{(2*\x,-2)}{2 \x}{};}
    \foreach \x in {3, 5}
      {\oddcell{(2*\x,-2)}{2 \x}{};}
    \evenrightedgecell{(12,-2)}{2 6};
    % Draw cells in 3rd row
    \oddleftedgecell{(0,-4)}{3 0};
    \foreach \x in {2, 4, 6}
      {\oddcell{(2*\x,-4)}{3 \x}{};}
    \foreach \x in {1, 3, 5}
      {\evencell{(2*\x,-4)}{3 \x}{};}
    \evenrightedgecell{(14,-4)}{3 7};
    % Draw cells in 4th (bottom) row
    \oddleftedgecell{(-2,-6)}{4 -1};
    \foreach \x in {0, 2, 4, 6}
      {\evencell{(2*\x,-6)}{4 \x}{};}
    \foreach \x in {1, 3, 5, 7}
      {\oddcell{(2*\x,-6)}{4 \x}{};}
    \evenrightedgecell{(16,-6)}{4 8};
    
    % Draw syllables in 1st (top) row
    \evenintsyll{1 3}{}{}{};
    \oddintsyll{1 4}{}{}{};
    % Draw syllables in 2nd row
    \evenintsyll{2 2}{}{}{};
    \oddintsyll{2 3}{}{}{};
    \evenintsyll{2 4}{}{}{};
    \oddintsyll{2 5}{}{}{};
    % Draw syllables in 3rd row
    \evenintsyll{3 1}{}{}{};
    \oddintsyll{3 2}{}{}{};
    \evenintsyll{3 3}{}{}{};
    \oddintsyll{3 4}{}{}{};
    \evenintsyll{3 5}{}{}{};
    \oddintsyll{3 6}{}{}{};    
    % Draw syllables in 4th (bottom) row
    \evenintsyll{4 0}{}{}{};
    \oddintsyll{4 1}{}{}{};
    \evenintsyll{4 2}{}{}{};
    \oddintsyll{4 3}{}{}{};
    \evenintsyll{4 4}{}{}{};
    \oddintsyll{4 5}{}{}{};
    \evenintsyll{4 6}{}{}{};
    \oddintsyll{4 7}{}{}{};
  \end{tikzpicture}
  
  \caption[%
    Interiors of syzygies
  ]{%
    \label{fig:interiors-of-syzygies}
    \textit{Interiors of syzygies.} Reading from the bottom upwards one row at a time, the interior can narrow with each step by at most 2 (or rather by at most 1 on either side). When the bottom row is sufficiently wide, as occurs here, the interiors of all of the rows involved have at least one column in common.
  }
\end{thesisfigure}

That second implication is crucial, for if column $k$ is common to the interior of successive syzygy strips $w_0,\dots,w_m$, then the values taken there
give $m+1$ consecutive interior syllables in a $\descent$-orbit:
$$
\begin{tikzcd}[]
  w_0(k) \ar[r, mapsto, "\descent"]
    \& w_1(k) \ar[r, mapsto, "\descent"]
      \& \cdots \ar[r, mapsto, "\descent"]
        \& w_m(k)
\end{tikzcd}\tstop
$$

\begin{theorem}\label{thm:explicit-bounds-on-syzygy-type}
Suppose the longest path in the pin graph $\pingraph_A$ of $A$ has length $m < +\infty$.
\begin{enumerate}
    \item\label{subthm:string-syzygy-of-band-when-pin-graph-acyclic} If $X$ is a band module for $A$, then $\syzygy^{2m+1} X$ is a string module.
    \item\label{subthm:syzygy-finite-strings-when-pin-graph-acyclic} If $X$ is a (finite- or infinite-dimensional) string module for $A$, then $$\syzygy^{2m+ 2} X \in \add \{ \Str w \in \Modcat A \colon \width w \leq 4m + 6\}\tstop$$
    In particular, injective string modules for $A$ are syzygy finite.
\end{enumerate}
\end{theorem}

\proofof{\ref{subthm:string-syzygy-of-band-when-pin-graph-acyclic}} Clearly $2m+1 = 2(m+1) - 1$. If $\syzygy^{2m+1} X$ is a band module then so is $\syzygy^r X$ for all $r \in \{0,1,\dots,2m+1\}$. We deduce the existence of a $\pingraph_A$-path of length $m+1$, which is absurd.

\proofof{\ref{subthm:syzygy-finite-strings-when-pin-graph-acyclic}} Let $r > 0$. Suppose that the strip $w$ has $\interiorwidth w > 4r$ and that it represents an indecomposable direct summand of a $(2r)$th syzygy of a string module. By the discussion in passage \ref{psg:late-wide-strips-share-interior-column}, there exists a sequence of $2r+1$ consecutive interior syllables in a $\descent$-orbit for $A$ and, therefore, there exists a path of length $r$ in $\Phi_A$. 

We deduce that any $r$ satisfying those conditions further satisfies $r \leq m$. Hence, for any $s \geq m+1$, any direct summand $w$ of $\syzygy^{2s} X$  has $\interiorwidth w \leq 4(m+1)$ and so $\width w \leq 4(m+1) + 2 = 4m + 6$.\qed

\begin{corollary}\label{cor:pin-graph-acyclic-implies-syz-finite}
If the pin graph of $A$ is acyclic, then $A$ is syzygy finite.
\end{corollary}

\proof By Theorem \ref{thm:explicit-bounds-on-syzygy-type},
$$ \syzygy^{2m+2}(\syzygy^{2m+1} X) =  \syzygy^{4m+3} X \in \add \{ \Str w \in \Modcat A \colon \width w \leq 4m + 6\}
$$ for any $A$-module $X$, where $m$ is the length of the longest path in $\pingraph_A$.\qed

\begin{corollary}
If the pin graph of $A$ is acyclic, then injectives generate for $A$, and therefore the big and little finitistic dimension conjectures hold for $A$.
\end{corollary}

\proof We know $A$ is the regular $A$-module, and so $\dual A$ is an $A^\op$-module.

Since $\pingraph_A$ is acyclic, $\pingraph_{A^\op}$ is acyclic also. The preceding corollary implies that $\dual A$ has finite syzygy type. Applying the functor $\dual$, we find $A$ has finite cosyzygy type (see \cite[Def 7.1]{Ric19} for a definition). We deduce that $A \in \Loc(\Injcat A)$ by \cite[7.2]{Ric19} and so injectives generate for $A$ by \cite[Prop 2.2]{Ric19} (or, rather, the comment immediately following it).\qed

\subsection[Descent and injective syllables]{Descent $\descent$ and injective syllables}\label{subsec:desc-and-injective-syllables}

We are keenly interested in the syzygies of injective string modules $E$ of a SB algebra $A$. This is because if $E$ always has finite syzygy type (as we conjecture) then projectives cogenerate for $A$. Additionally, the calculation could be dualised to establish that projective string modules have finite cosyzygy type and hence that injectives generate for $A$. Both the generation and cogeneration statements relate to the big finitistic dimension conjecture, as discussed previously.

One purpose of the syzygy fabric formalism is to facilitate analysis of syzygies. In particular it lets us rigorously organise our syzygies one under the other. The columns are, largely, just orbits under some function $\descent$ and, largely, neighbouring columns do not interact. It is tragic to us that ``largely'' cannot be replaced by ``entirely'' here. Nonetheless, if we are to study injective string modules $\Str w$ as intended, then in the syzygy fabric of $w$ we should examine the columns containing them and descending from them. This is the focus of this subsection.

This subsection explores what happens ``directly underneath'' an injective string module $\Str w$ in the syzygy fabric, and it culminates in Proposition \ref{prop:descent-orbits-of-inj-sylls-do-not-have-interior-tails}. This rather technical result establishes that the columns descending from an injective strip are either blank after a finite number of rows or witness branching infinitely often. In the latter case, the number of rows between successive branchings is bounded, which gives a ``vertical'' constraint on the interiors of strips directly underneath $w$. We believe that from this vertical constraint it should be possible to deduce a ``horizontal'' constraint on widths of strips (as we have proved from in the case of cycle-free pin graphs) and so obtain results about the syzygy type of $\Str w$.

Proposition \ref{prop:descent-orbits-of-inj-sylls-do-not-have-interior-tails} is the first step in a journey whose conclusion would be a proof of finite syzygy type for injective string modules; in fact, we believe it to be the most technically difficult step. The next step would be to demonstrate the constraint of the previous paragraph. Two further steps would follow: how the interaction of columns affects our proposition, and what happens in the remainder of the syzygy fabric of $E$ besides the region directly underneath it. Towards these steps in the remainder of the journey, we offer the indicative worked examples at the end of the subsection and the results about ``small'' SB algebras in Appendix \ref{app:sb-algs-with-few-verts}.

% We have nothing more conclusive than conjectures founded on worked examples. Nonetheless, we offer our partial result in the firm belief that it solves the most difficult combinatorial aspect of an ultimate proof of finite syzygy type.

\nconv When we use (connected) $1$-regular quivers in the following, we will work with the acyclic cover
$\big( \begin{tikzcd}[sep=small]
\cdots \ar[r]
  \& \syllableplaceholder \ar[r]
    \& \syllableplaceholder \ar[r]
      \& \syllableplaceholder \ar[r]
        \& \cdots
\end{tikzcd} \big)$ rather than the finite quiver $\big(
\raisebox{1ex}
{\begin{tikzcd}[sep=small]
    \syllableplaceholder \ar[r]
      \& \syllableplaceholder \ar[r]\ar[rr, phantom, ""{coordinate, name=C}]
        \& \cdots \ar[r]
          \& \syllableplaceholder \ar[r]
            \& \syllableplaceholder \ar[llll, rounded corners, to path={
                    -- ([xshift=1.5ex]\tikztostart.east)
                    |- ([yshift=-2ex]C)[pos=1]\tikztonodes
                    -| ([xshift=-1.5ex]\tikztotarget.west)
                    -- (\tikztotarget.west)
                    }]
\end{tikzcd}}
\big)$ itself.

\centresubsec{Factorisable paths of an interior $\descent$-cycle}

\nlpass{Interior $\descent$-cycle} Suppose that $\{\bdu_r, \bdv_r \colon r \in \integers \}$ is a periodic family of interior syllables such that $\bdu_r \descent = \bdv_r$ and $\bdv_r \descent = \bdu_{r+1} $; by periodic here, we mean that $\bdu_r = \bdu_{r+m}$ and $\bdv_r = \bdv_{r+m}$ for some $m > 0$. Taking subscripts modulo $m$, we will call $(\bdu_1, \bdv_1, \dots, \bdu_m, \bdv_m)$ an \define{interior $\descent$-cycle}.

\nlpass{Elementary factors of an interior $\descent$-cycle} Write $u_r, v_r$ respectively for the underlying paths of the syllables $\bdu_r, \bdv_r$ on an interior $\descent$-cycle. We call $\{u_r, v_r \colon r \in \integers\}$ the \define{elementary factors} of the interior $\descent$-cycle. Since the $\bdu_r,\bdv_r$ are certainly not boundary syllables, the $u_r,v_r$ are certainly not stationary paths.

We directly find that $\target(u_r) = \source(v_r)$ and $\target(v_r)=\source(u_{r+1})$ for all $r$, and that the $\{u_r, v_r \colon r \in \integers\}$ form an antichain of paths in $\OQ$; more precisely, in a single component $\OC$ of $\OQ$. Our setup guarantees that pairwise products $u_r v_r$ and $v_r u_{r+1}$ of elementary factors represent components of commutativity relations.

\nrex Let $\OC$ be the component of the overquiver covering $\alpha,\beta,\gamma$. A part of the subpath poset $\Subpath(\OC)$ is shown in Figure \ref{fig:hasse-diagram-for-C}.

There are infinitely many $\OC$-paths lifting $\alpha$ and $\beta\gamma$. For $r\in\integers$, give the name $u_r$ to those $\OC$-paths lifting $\alpha$ and the name $v_r$ to the paths lifting $\beta\gamma$, in such a way that $\target(u_r)=\source(v_r)$ and $\target(v_r)=\source(u_{r+1})$ for all $r$. Define $\bdu_r \ceq \big( \ilsyll{\syllableplaceholder}{u_r}{0} \big)$ and $\bdv_r \ceq \big( \ilsyll{\syllableplaceholder}{v_r}{0} \big)$. We find that $\bdu_r \descent = \bdv_r$ and $\bdv_r \descent = \bdu_{r+1}$, and so in particular $\bdu_r = \bdu_{r+1}$ and $\bdv_r = \bdv_{r+1}$. Therefore $(\bdu_1, \bdv_1)$ is an interior $\descent$-cycle with elementary factors $u_r=\alpha,v_r=\beta\gamma$.

We verify that pairwise products $u_r v_r$ and $v_r u_{r+1}$ respectively represent $\alpha\beta\gamma$ and $\beta\gamma\alpha$, which are indeed components of a commutativity relation. %We can perhaps think of the elementary factors as representing interlocking halves of commutativity relations.
\begin{thesisfigure}
  \centering

  \begin{tikzpicture}[scale=0.65]
    \newcommand{\fadedcolour}{black!30}
    \newcommand{\nodeHspacing}[1]{(#1,0)}
    \newcommand{\nodeVspacing}[1]{(0,#1*1.5)}
    \newcommand{\basicnode}[4]{%
      \draw (o) ++\nodeHspacing{#1} ++\nodeVspacing{#2} node[\fadedcolour](#3){$ #4 $};
    }
    \newcommand{\nodeplaceholder}{\cdot}
    \newcommand{\facnodeplaceholder}{\circ}
    \newcommand{\facnodeframework}[5]{
      \draw (o) ++\nodeHspacing{#1} ++\nodeVspacing{#2} node[rectangle, rounded corners, inner sep=5, #5](#3){$ #4 $};
    }
    \newcommand{\eltfactornode}[4]{%
      \facnodeframework{#1}{#2}{#3}{#4}{black}
    }
    \newcommand{\facnode}[4]{%
      \facnodeframework{#1}{#2}{#3}{#4}{black}
    }
    %%  Place "lattice" of nodes
    % 0th (bottom) row
    \draw (0,0) node[](o){};
    \facnode{0}{0}{sa0}{\source(\alpha)};
    \facnode{2}{0}{sb0}{\source(\beta)};
    \basicnode{4}{0}{sc0}{\source(\gamma)};
    \facnode{6}{0}{sa1}{\source(\alpha)};
    \facnode{8}{0}{sb1}{\source(\beta)};
    \basicnode{10}{0}{sc1}{\source(\gamma)};
    \facnode{12}{0}{sa2}{\source(\alpha)};
    \facnode{14}{0}{sb2}{\source(\beta)};
    \basicnode{16}{0}{sc2}{\source(\gamma)};
    \facnode{18}{0}{sa3}{\source(\alpha)};
    % 1st row
    \eltfactornode{1}{1}{a0}{\alpha};
    \basicnode{3}{1}{b0}{\beta};
    \basicnode{5}{1}{c0}{\gamma};
    \eltfactornode{7}{1}{a1}{\alpha};
    \basicnode{9}{1}{b1}{\beta};
    \basicnode{11}{1}{c1}{\gamma};
    \eltfactornode{13}{1}{a2}{\alpha};
    \basicnode{15}{1}{b2}{\beta};
    \basicnode{17}{1}{c2}{\gamma};
    % 2nd row
    \basicnode{0}{2}{ca0}{\gamma\alpha};
    \basicnode{2}{2}{ab0}{\alpha\beta};
    \eltfactornode{4}{2}{bc0}{\beta\gamma};
    \basicnode{6}{2}{ca1}{\gamma\alpha};
    \basicnode{8}{2}{ab1}{\alpha\beta};
    \eltfactornode{10}{2}{bc1}{\beta\gamma};
    \basicnode{12}{2}{ca2}{\gamma\alpha};
    \basicnode{14}{2}{ab2}{\alpha\beta};
    \eltfactornode{16}{2}{bc2}{\beta\gamma};
    \basicnode{18}{2}{ca3}{\gamma\alpha};
    % 3rd row
    \basicnode{1}{3}{cab0}{\gamma\alpha\beta};
    \facnode{3}{3}{abc0}{\alpha\beta\gamma};
    \facnode{5}{3}{bca0}{\beta\gamma\alpha};
    \basicnode{7}{3}{cab1}{\gamma\alpha\beta};
    \facnode{9}{3}{abc1}{\alpha\beta\gamma};
    \facnode{11}{3}{bca1}{\beta\gamma\alpha};
    \basicnode{13}{3}{cab2}{\gamma\alpha\beta};
    \facnode{15}{3}{abc2}{\alpha\beta\gamma};
    \facnode{17}{3}{bca2}{\beta\gamma\alpha};
    % 4th row
    \basicnode{0}{4}{bcab0}{\nodeplaceholder};
    \basicnode{2}{4}{cabc0}{\nodeplaceholder};
    \facnode{4}{4}{abca0}{\nodeplaceholder};
    \basicnode{6}{4}{bcab1}{\nodeplaceholder};
    \basicnode{8}{4}{cabc1}{\nodeplaceholder};
    \facnode{10}{4}{abca1}{\nodeplaceholder};
    \basicnode{12}{4}{bcab2}{\nodeplaceholder};
    \basicnode{14}{4}{cabc2}{\nodeplaceholder};
    \facnode{16}{4}{abca2}{\nodeplaceholder};
    \basicnode{18}{4}{bcab3}{\nodeplaceholder};
    % 5th row
    \facnode{1}{5}{bcabc0}{\nodeplaceholder};
    \basicnode{3}{5}{cabca0}{\nodeplaceholder};
    \basicnode{5}{5}{abcab0}{\nodeplaceholder};
    \facnode{7}{5}{bcabc1}{\nodeplaceholder};
    \basicnode{9}{5}{cabca1}{\nodeplaceholder};
    \basicnode{11}{5}{abcab1}{\nodeplaceholder};
    \facnode{13}{5}{bcabc2}{\nodeplaceholder};
    \basicnode{15}{5}{cabca2}{\nodeplaceholder};
    \basicnode{17}{5}{abcab2}{\nodeplaceholder};
    %%  Draw connecting lines in Hasse diagram
    \newcommand{\hasseconnector}[2]{%
      \draw[\fadedcolour] (#1) -- (#2);
    }
    % 0th row to 1st row
    \hasseconnector{sa0}{a0};
    \hasseconnector{sb0}{b0};
    \hasseconnector{sc0}{c0};
    \hasseconnector{sa1}{a1};
    \hasseconnector{sb1}{b1};
    \hasseconnector{sc1}{c1};
    \hasseconnector{sa2}{a2};
    \hasseconnector{sb2}{b2};
    \hasseconnector{sc2}{c2};
    \hasseconnector{sb0}{a0};
    \hasseconnector{sc0}{b0};
    \hasseconnector{sa1}{c0};
    \hasseconnector{sb1}{a1};
    \hasseconnector{sc1}{b1};
    \hasseconnector{sa2}{c1};
    \hasseconnector{sb2}{a2};
    \hasseconnector{sc2}{b2};
    \hasseconnector{sa3}{c2};
    % 1st row to 2nd row
    \hasseconnector{a0}{ca0};
    \hasseconnector{b0}{ab0};
    \hasseconnector{c0}{bc0};
    \hasseconnector{a1}{ca1};
    \hasseconnector{b1}{ab1};
    \hasseconnector{c1}{bc1};
    \hasseconnector{a2}{ca2};
    \hasseconnector{b2}{ab2};
    \hasseconnector{c2}{bc2};
    \hasseconnector{a0}{ab0};
    \hasseconnector{b0}{bc0};
    \hasseconnector{c0}{ca1};
    \hasseconnector{a1}{ab1};
    \hasseconnector{b1}{bc1};
    \hasseconnector{c1}{ca2};
    \hasseconnector{a2}{ab2};
    \hasseconnector{b2}{bc2};
    \hasseconnector{c2}{ca3};
    % 2nd row to 3rd row
    \hasseconnector{ca0}{cab0};
    \hasseconnector{ab0}{abc0};
    \hasseconnector{bc0}{bca0};
    \hasseconnector{ca1}{cab1};
    \hasseconnector{ab1}{abc1};
    \hasseconnector{bc1}{bca1};
    \hasseconnector{ca2}{cab2};
    \hasseconnector{ab2}{abc2};
    \hasseconnector{bc2}{bca2};
    \hasseconnector{ab0}{cab0};
    \hasseconnector{bc0}{abc0};
    \hasseconnector{ca1}{bca0};
    \hasseconnector{ab1}{cab1};
    \hasseconnector{bc1}{abc1};
    \hasseconnector{ca2}{bca1};
    \hasseconnector{ab2}{cab2};
    \hasseconnector{bc2}{abc2};
    \hasseconnector{ca3}{bca2};
    % 3rd row to 4th row
    \hasseconnector{cab0}{bcab0};
    \hasseconnector{abc0}{cabc0};
    \hasseconnector{bca0}{abca0};
    \hasseconnector{cab1}{bcab1};
    \hasseconnector{abc1}{cabc1};
    \hasseconnector{bca1}{abca1};
    \hasseconnector{cab2}{bcab2};
    \hasseconnector{abc2}{cabc2};
    \hasseconnector{bca2}{abca2};
    \hasseconnector{cab0}{cabc0};
    \hasseconnector{abc0}{abca0};
    \hasseconnector{bca0}{bcab1};
    \hasseconnector{cab1}{cabc1};
    \hasseconnector{abc1}{abca1};
    \hasseconnector{bca1}{bcab2};
    \hasseconnector{cab2}{cabc2};
    \hasseconnector{abc2}{abca2};
    \hasseconnector{bca2}{bcab3};
    % 4th row to 5th row
    \hasseconnector{bcab0}{bcabc0};
    \hasseconnector{cabc0}{cabca0};
    \hasseconnector{abca0}{abcab0};
    \hasseconnector{bcab1}{bcabc1};
    \hasseconnector{cabc1}{cabca1};
    \hasseconnector{abca1}{abcab1};
    \hasseconnector{bcab2}{bcabc2};
    \hasseconnector{cabc2}{cabca2};
    \hasseconnector{abca2}{abcab2};
    \hasseconnector{cabc0}{bcabc0};
    \hasseconnector{abca0}{cabca0};
    \hasseconnector{bcab1}{abcab0};
    \hasseconnector{cabc1}{bcabc1};
    \hasseconnector{abca1}{cabca1};
    \hasseconnector{bcab2}{abcab1};
    \hasseconnector{cabc2}{bcabc2};
    \hasseconnector{abca2}{cabca2};
    \hasseconnector{bcab3}{abcab2};
  \end{tikzpicture}
  
  \caption[%
    Factorisable paths within the Hasse diagram of $\Subpath(\OC)$.
  ]{%
    \label{fig:hasse-diagram-for-C}
    \textit{Factorisable paths within the Hasse diagram of $\Subpath(\OC)$.} The $\OC$-paths are labelled by the $Q$-paths they represent. (Some labels are omitted for space concerns.) Factorisable paths are printed in black; the others, grey.
    
    Among the factorisable paths are the elementary factors $u_r \ceq \alpha$ and $v_r \ceq \beta \gamma$, where the index $r \in \integers$ increases rightwards in the diagram.
}
\end{thesisfigure}

\centresubsec{Grading}

\newcommand{\comppaths}{\categoryfont{C}}
%   Category of paths in the component \OC

\npass Henceforth fix an interior $\descent$-cycle $(\bdu_1, \bdv_1, \dots, \bdu_m, \bdv_m)$ with elementary factors $\{u_r, v_r \colon r \in \integers\}$ all lying on a component $\OC$ of the overquiver. The following definitions will depend on the cycle; we simply make this dependence implicit.

Additionally, write $\comppaths$ for the path category of $\OC$.

\newcommand{\factpaths}{\categoryfont{F}}
%   Category of factorisable paths, which we may identify with category of paths in a quiver $\OF$ 
\newcommand{\OF}{\mathcal{F}}

\nlpass{Factorisable paths and grading} A \define{factorisable path} is a path that can be written as a product of elementary factors $u_r,v_r$. We include in this definition the endpoints $\target(u_r)=\source(v_r)$ and $\target(v_r)=\source(u_{r+1})$ of the elementary factors, which we consider to be $0$-fold products and which we call \define{factorisable vertices}. A factorisable path $p$ has \define{grade} $g$, written $\mdefine{\grade p} = g$, when it is a $g$-fold product of elementary factors.

The factorisable paths give a full subcategory $\factpaths \subseteq \comppaths$, namely the full subcategory on the factorisable vertices. We can also identify $\factpaths$ with the path category of the quiver
$$
\OF \ceq
\big(
  \begin{tikzcd}
    \cdots \ar[r]
      \& \source(u_r) \ar[r, "u_r"]
        \& \source(v_r) \ar[r, "v_r"]
          \& \source(u_{r+1}) \ar[r, "u_{r+1}"]
            \& \source(v_{r+1}) \ar[r]
              \& \cdots
  \end{tikzcd}
\big)\tstop
$$
In this light, $\grade p$ is simply the length of $p$ viewed as a path in $\OF$. We deduce that $\grade$ is additive; that is, $\grade(pp')=\grade p + \grade p'$.

\nrex We have elementary factors $u_r \ceq \alpha $ and $v_r \ceq \beta\gamma$ and the factorisable vertices are $\source(\alpha)$ and $\source(\beta)$. We can identify these with the vertices and arrows of the quiver
$$
  \OF \ceq \big(
  \begin{tikzcd}
    \cdots \ar[r]
      \& \source(\alpha) \ar[r, "\alpha"]
        \& \source(\beta) \ar[r, "\beta\gamma"]
          \& \source(\alpha) \ar[r, "\alpha"]
            \& \source(\beta) \ar[r]
              \& \cdots
  \end{tikzcd}
\big)\tstop
$$

The $2$-fold products, namely the paths representing $\alpha\beta\gamma$ and $\beta\gamma\alpha$, have grade $2$. The $3$-fold products, namely the paths representing $\alpha\beta\gamma\alpha$ and $\beta\gamma\alpha\beta\gamma$, have grade 3. And so on.

\nlpass{Grading $\comppaths$ using a functor $f \colon \ilarrow{ \comppaths }{}{ \factpaths }$}\label{def:grading-paths-and-f} We give a functor $f \colon \ilarrow{ \comppaths }{}{ \factpaths }$ that assigns a factorisable path to an arbitrary path in $\OC$. This functor is far easier to picture and use than the ensuing prose description may suggest, so we encourage the reader to consult the examples and illustration which follow in the next passage and Figure \ref{fig:example-of-grading}.
\begin{thesisfigure}
  \centering
  
  \begin{tikzpicture}[scale=0.65]
    \newcommand{\nodeHspacing}[1]{(#1,0)}
    \newcommand{\nodeVspacing}[1]{(0,#1*1.5)}
    \newcommand{\basicnode}[4]{%
      \draw (o) ++\nodeHspacing{#1} ++\nodeVspacing{#2} node[](#3){$ #4 $};
    }
    \newcommand{\nodeplaceholder}{\cdot}
    \newcommand{\facnodeplaceholder}{\circ}
    %%  Place "lattice" of nodes
    % 0th (bottom) row
    \newcommand{\boundingrectangle}{
      (-0.5,0) rectangle (18,7.7)
    }
    \newcommand{\faccolour}{promptColor}
    \newcommand{\evenfaccolour}{\faccolour!5}
    \newcommand{\oddfaccolour}{\faccolour!15}
    \newcommand{\fadedcolour}{black!30}
     
  \begin{scope}
     \clip \boundingrectangle;
     
    \draw (0,0) node[](o){};
    % (-1)st row (needed for some highlighted regions)
    \basicnode{5}{-1}{z0}{ };
    \basicnode{11}{-1}{z1}{ };
    \basicnode{17}{-1}{z2}{ };
    % 0th row
    \basicnode{0}{0}{sa0}{\source(\alpha)};
    \basicnode{2}{0}{sb0}{\source(\beta)};
    \basicnode{4}{0}{sc0}{\source(\gamma)};
    \basicnode{6}{0}{sa1}{\source(\alpha)};
    \basicnode{8}{0}{sb1}{\source(\beta)};
    \basicnode{10}{0}{sc1}{\source(\gamma)};
    \basicnode{12}{0}{sa2}{\source(\alpha)};
    \basicnode{14}{0}{sb2}{\source(\beta)};
    \basicnode{16}{0}{sc2}{\source(\gamma)};
    \basicnode{18}{0}{sa3}{\source(\alpha)};
    % 1st row
    \basicnode{-1}{1}{c-1}{};
    \basicnode{1}{1}{a0}{\alpha};
    \basicnode{3}{1}{b0}{\beta};
    \basicnode{5}{1}{c0}{\gamma};
    \basicnode{7}{1}{a1}{\alpha};
    \basicnode{9}{1}{b1}{\beta};
    \basicnode{11}{1}{c1}{\gamma};
    \basicnode{13}{1}{a2}{\alpha};
    \basicnode{15}{1}{b2}{\beta};
    \basicnode{17}{1}{c2}{\gamma};
    \basicnode{19}{1}{a3}{};
    % 2nd row
    \basicnode{0}{2}{ca0}{\gamma\alpha};
    \basicnode{2}{2}{ab0}{\alpha\beta};
    \basicnode{4}{2}{bc0}{\beta\gamma};
    \basicnode{6}{2}{ca1}{\gamma\alpha};
    \basicnode{8}{2}{ab1}{\alpha\beta};
    \basicnode{10}{2}{bc1}{\beta\gamma};
    \basicnode{12}{2}{ca2}{\gamma\alpha};
    \basicnode{14}{2}{ab2}{\alpha\beta};
    \basicnode{16}{2}{bc2}{\beta\gamma};
    \basicnode{18}{2}{ca3}{};
    % 3rd row
    \basicnode{-1}{3}{bca-1}{};
    \basicnode{1}{3}{cab0}{\gamma\alpha\beta};
    \basicnode{3}{3}{abc0}{\alpha\beta\gamma};
    \basicnode{5}{3}{bca0}{\beta\gamma\alpha};
    \basicnode{7}{3}{cab1}{\gamma\alpha\beta};
    \basicnode{9}{3}{abc1}{\alpha\beta\gamma};
    \basicnode{11}{3}{bca1}{\beta\gamma\alpha};
    \basicnode{13}{3}{cab2}{\gamma\alpha\beta};
    \basicnode{15}{3}{abc2}{\alpha\beta\gamma};
    \basicnode{17}{3}{bca2}{\beta\gamma\alpha};
    \basicnode{19}{3}{cab3}{};
    % 4th row
    \basicnode{0}{4}{bcab0}{};
    \basicnode{2}{4}{cabc0}{};
    \basicnode{4}{4}{abca0}{\alpha\beta\gamma\alpha};
    \basicnode{6}{4}{bcab1}{};
    \basicnode{8}{4}{cabc1}{};
    \basicnode{10}{4}{abca1}{\alpha\beta\gamma\alpha};
    \basicnode{12}{4}{bcab2}{};
    \basicnode{14}{4}{cabc2}{};
    \basicnode{16}{4}{abca2}{\alpha\beta\gamma\alpha};
    \basicnode{18}{4}{bcab3}{};
    % 5th row
    \basicnode{-1}{5}{abcab-1}{};
    \basicnode{1}{5}{bcabc0}{\beta\gamma\alpha\beta\gamma};
    \basicnode{3}{5}{cabca0}{};
    \basicnode{5}{5}{abcab0}{};
    \basicnode{7}{5}{bcabc1}{\beta\gamma\alpha\beta\gamma};
    \basicnode{9}{5}{cabca1}{};
    \basicnode{11}{5}{abcab1}{};
    \basicnode{13}{5}{bcabc2}{\beta\gamma\alpha\beta\gamma};
    \basicnode{15}{5}{cabca2}{};
    \basicnode{17}{5}{abcab2}{};
    \basicnode{19}{5}{bcabc3}{};
    % 6th row (needed for some highlighted regions)
    \basicnode{4}{6}{y0}{};
    \basicnode{6}{6}{t0}{};
    \basicnode{10}{6}{y1}{};
    \basicnode{12}{6}{t1}{};
    \basicnode{16}{6}{y2}{};
    \basicnode{18}{6}{t2}{};

    %% Draw highlighting regions
    % Leading factorisable path in the 0th row
    \fill[\evenfaccolour, rounded corners] ($(sa0)+0.5*(1,0)$) -- ($(c-1)-0.25*(1,0)+0.75*(0,1.5)$) -- ($(c-1)-1*(1,0)$) -- ($(sa0)-0.25*(1,0)-0.75*(0,1.5)$) -- cycle;
    \fill[\evenfaccolour, rounded corners] ($(sb0)+0.5*(1,0)$) -- ($(sb0)-0.25*(1,0)+0.75*(0,1.5)$) -- ($(sb0)-1*(1,0)$) -- ($(sb0)-0.25*(1,0)-0.75*(0,1.5)$) -- cycle;
    \fill[\evenfaccolour, rounded corners] ($(sa1)+0.5*(1,0)$) -- ($(c0)-0.25*(1,0)+0.75*(0,1.5)$) -- ($(sc0)-1*(1,0)$) -- ($(z0)-0.25*(1,0)-0.75*(0,1.5)$) -- cycle;
    \fill[\evenfaccolour, rounded corners] ($(sb1)+0.5*(1,0)$) -- ($(sb1)-0.25*(1,0)+0.75*(0,1.5)$) -- ($(sb1)-1*(1,0)$) -- ($(sb1)-0.25*(1,0)-0.75*(0,1.5)$) -- cycle;
    \fill[\evenfaccolour, rounded corners] ($(sa2)+0.5*(1,0)$) -- ($(c1)-0.25*(1,0)+0.75*(0,1.5)$) -- ($(sc1)-1*(1,0)$) -- ($(z1)-0.25*(1,0)-0.75*(0,1.5)$) -- cycle;
    \fill[\evenfaccolour, rounded corners] ($(sb2)+0.5*(1,0)$) -- ($(sb2)-0.25*(1,0)+0.75*(0,1.5)$) -- ($(sb2)-1*(1,0)$) -- ($(sb2)-0.25*(1,0)-0.75*(0,1.5)$) -- cycle;
    \fill[\evenfaccolour, rounded corners] ($(sa3)+0.5*(1,0)$) -- ($(c2)-0.25*(1,0)+0.75*(0,1.5)$) -- ($(sc2)-1*(1,0)$) -- ($(z2)-0.25*(1,0)-0.75*(0,1.5)$) -- cycle;
    % Leading factorisable path in the 1st row
    \fill[\oddfaccolour, rounded corners] ($(a0)+0.5*(1,0)$) -- ($(ca0)-0.25*(1,0)+0.75*(0,1.5)$) -- ($(ca0)-1*(1,0)$) -- ($(a0)-0.25*(1,0)-0.75*(0,1.5)$) -- cycle;
    \fill[\oddfaccolour, rounded corners] ($(a1)+0.5*(1,0)$) -- ($(ca1)-0.25*(1,0)+0.75*(0,1.5)$) -- ($(ca1)-1*(1,0)$) -- ($(a1)-0.25*(1,0)-0.75*(0,1.5)$) -- cycle;
    \fill[\oddfaccolour, rounded corners] ($(a2)+0.5*(1,0)$) -- ($(ca2)-0.25*(1,0)+0.75*(0,1.5)$) -- ($(ca2)-1*(1,0)$) -- ($(a2)-0.25*(1,0)-0.75*(0,1.5)$) -- cycle;
    \fill[\oddfaccolour, rounded corners] ($(a3)+0.5*(1,0)$) -- ($(ca3)-0.25*(1,0)+0.75*(0,1.5)$) -- ($(ca3)-1*(1,0)$) -- ($(a3)-0.25*(1,0)-0.75*(0,1.5)$) -- cycle;
    % Leading factorisable path in the 2nd row
    \fill[\oddfaccolour, rounded corners] ($(bc0)+0.5*(1,0)$) -- ($(bc0)-0.25*(1,0)+0.75*(0,1.5)$) -- ($(b0)-1*(1,0)$) -- ($(b0)-0.25*(1,0)-0.75*(0,1.5)$) -- cycle;
    \fill[\oddfaccolour, rounded corners] ($(bc1)+0.5*(1,0)$) -- ($(bc1)-0.25*(1,0)+0.75*(0,1.5)$) -- ($(b1)-1*(1,0)$) -- ($(b1)-0.25*(1,0)-0.75*(0,1.5)$) -- cycle;
    \fill[\oddfaccolour, rounded corners] ($(bc2)+0.5*(1,0)$) -- ($(bc2)-0.25*(1,0)+0.75*(0,1.5)$) -- ($(b2)-1*(1,0)$) -- ($(b2)-0.25*(1,0)-0.75*(0,1.5)$) -- cycle;
    % Leading factorisable path in the 3rd row
    \fill[\evenfaccolour, rounded corners] ($(abc0)+0.5*(1,0)$) -- ($(cabc0)-0.25*(1,0)+0.75*(0,1.5)$) -- ($(cab0)-1*(1,0)$) -- ($(ab0)-0.25*(1,0)-0.75*(0,1.5)$) -- cycle;
    \fill[\evenfaccolour, rounded corners] ($(bca0)+0.5*(1,0)$) -- ($(bca0)-0.25*(1,0)+0.75*(0,1.5)$) -- ($(bca0)-1*(1,0)$) -- ($(bca0)-0.25*(1,0)-0.75*(0,1.5)$) -- cycle;
    \fill[\evenfaccolour, rounded corners] ($(abc1)+0.5*(1,0)$) -- ($(cabc1)-0.25*(1,0)+0.75*(0,1.5)$) -- ($(cab1)-1*(1,0)$) -- ($(ab1)-0.25*(1,0)-0.75*(0,1.5)$) -- cycle;
    \fill[\evenfaccolour, rounded corners] ($(bca1)+0.5*(1,0)$) -- ($(bca1)-0.25*(1,0)+0.75*(0,1.5)$) -- ($(bca1)-1*(1,0)$) -- ($(bca1)-0.25*(1,0)-0.75*(0,1.5)$) -- cycle;
    \fill[\evenfaccolour, rounded corners] ($(abc2)+0.5*(1,0)$) -- ($(cabc2)-0.25*(1,0)+0.75*(0,1.5)$) -- ($(cab2)-1*(1,0)$) -- ($(ab2)-0.25*(1,0)-0.75*(0,1.5)$) -- cycle;
    \fill[\evenfaccolour, rounded corners] ($(bca2)+0.5*(1,0)$) -- ($(bca2)-0.25*(1,0)+0.75*(0,1.5)$) -- ($(bca2)-1*(1,0)$) -- ($(bca2)-0.25*(1,0)-0.75*(0,1.5)$) -- cycle; 
    \fill[\evenfaccolour, rounded corners] ($(cab3)+0.5*(1,0)$) -- ($(cab3)-0.25*(1,0)+0.75*(0,1.5)$) -- ($(cab3)-1*(1,0)$) -- ($(cab3)-0.25*(1,0)-0.75*(0,1.5)$) -- cycle; 
    % Leading factorisable paths in the 4th row
    \fill[\oddfaccolour, rounded corners] ($(abca0)+0.5*(1,0)$) -- ($(cabca0)-0.25*(1,0)+0.75*(0,1.5)$) -- ($(cabca0)-1*(1,0)$) -- ($(abca0)-0.25*(1,0)-0.75*(0,1.5)$) -- cycle;
    \fill[\oddfaccolour, rounded corners] ($(abca1)+0.5*(1,0)$) -- ($(cabca1)-0.25*(1,0)+0.75*(0,1.5)$) -- ($(cabca1)-1*(1,0)$) -- ($(abca1)-0.25*(1,0)-0.75*(0,1.5)$) -- cycle;
    \fill[\oddfaccolour, rounded corners] ($(abca2)+0.5*(1,0)$) -- ($(cabca2)-0.25*(1,0)+0.75*(0,1.5)$) -- ($(cabca2)-1*(1,0)$) -- ($(abca2)-0.25*(1,0)-0.75*(0,1.5)$) -- cycle;
    % Leading factorisable paths in the 5th row
    \fill[\oddfaccolour, rounded corners] ($(bcabc0)+0.5*(1,0)$) -- ($(bcabc0)-0.25*(1,0)+0.75*(0,1.5)$) -- ($(bcab0)-1*(1,0)$) -- ($(bcab0)-0.25*(1,0)-0.75*(0,1.5)$) -- cycle;
    \fill[\oddfaccolour, rounded corners] ($(bcabc1)+0.5*(1,0)$) -- ($(bcabc1)-0.25*(1,0)+0.75*(0,1.5)$) -- ($(bcab1)-1*(1,0)$) -- ($(bcab1)-0.25*(1,0)-0.75*(0,1.5)$) -- cycle;
    \fill[\oddfaccolour, rounded corners] ($(bcabc2)+0.5*(1,0)$) -- ($(bcabc2)-0.25*(1,0)+0.75*(0,1.5)$) -- ($(bcab2)-1*(1,0)$) -- ($(bcab2)-0.25*(1,0)-0.75*(0,1.5)$) -- cycle;
    \fill[\oddfaccolour, rounded corners] ($(bcabc3)+0.5*(1,0)$) -- ($(bcabc3)-0.25*(1,0)+0.75*(0,1.5)$) -- ($(bcab3)-1*(1,0)$) -- ($(bcab3)-0.25*(1,0)-0.75*(0,1.5)$) -- cycle;
    
    \fill[\evenfaccolour, rounded corners] ($(t0)+0.5*(1,0)$) -- ($(t0)-0.25*(1,0)+0.75*(0,1.5)$) -- ($(y0)-1*(1,0)$) -- ($(abcab0)-0.25*(1,0)-0.75*(0,1.5)$) -- cycle;
    \fill[\evenfaccolour, rounded corners] ($(t1)+0.5*(1,0)$) -- ($(t1)-0.25*(1,0)+0.75*(0,1.5)$) -- ($(y1)-1*(1,0)$) -- ($(abcab1)-0.25*(1,0)-0.75*(0,1.5)$) -- cycle;
    \fill[\evenfaccolour, rounded corners] ($(t2)+0.5*(1,0)$) -- ($(t2)-0.25*(1,0)+0.75*(0,1.5)$) -- ($(y2)-1*(1,0)$) -- ($(abcab2)-0.25*(1,0)-0.75*(0,1.5)$) -- cycle;
    
    %%  Draw connecting lines in Hasse diagram
    \newcommand{\hasseconnector}[2]{%
      \draw[\fadedcolour] (#1) -- (#2);
    }
    % 0th row to 1st row
    \hasseconnector{sa0}{a0};
    \hasseconnector{sb0}{b0};
    \hasseconnector{sc0}{c0};
    \hasseconnector{sa1}{a1};
    \hasseconnector{sb1}{b1};
    \hasseconnector{sc1}{c1};
    \hasseconnector{sa2}{a2};
    \hasseconnector{sb2}{b2};
    \hasseconnector{sc2}{c2};
    \hasseconnector{sb0}{a0};
    \hasseconnector{sc0}{b0};
    \hasseconnector{sa1}{c0};
    \hasseconnector{sb1}{a1};
    \hasseconnector{sc1}{b1};
    \hasseconnector{sa2}{c1};
    \hasseconnector{sb2}{a2};
    \hasseconnector{sc2}{b2};
    \hasseconnector{sa3}{c2};
    % 1st row to 2nd row
    \hasseconnector{a0}{ca0};
    \hasseconnector{b0}{ab0};
    \hasseconnector{c0}{bc0};
    \hasseconnector{a1}{ca1};
    \hasseconnector{b1}{ab1};
    \hasseconnector{c1}{bc1};
    \hasseconnector{a2}{ca2};
    \hasseconnector{b2}{ab2};
    \hasseconnector{c2}{bc2};
    \hasseconnector{a0}{ab0};
    \hasseconnector{b0}{bc0};
    \hasseconnector{c0}{ca1};
    \hasseconnector{a1}{ab1};
    \hasseconnector{b1}{bc1};
    \hasseconnector{c1}{ca2};
    \hasseconnector{a2}{ab2};
    \hasseconnector{b2}{bc2};
    \hasseconnector{c2}{ca3};
    % 2nd row to 3rd row
    \hasseconnector{ca0}{cab0};
    \hasseconnector{ab0}{abc0};
    \hasseconnector{bc0}{bca0};
    \hasseconnector{ca1}{cab1};
    \hasseconnector{ab1}{abc1};
    \hasseconnector{bc1}{bca1};
    \hasseconnector{ca2}{cab2};
    \hasseconnector{ab2}{abc2};
    \hasseconnector{bc2}{bca2};
    \hasseconnector{ab0}{cab0};
    \hasseconnector{bc0}{abc0};
    \hasseconnector{ca1}{bca0};
    \hasseconnector{ab1}{cab1};
    \hasseconnector{bc1}{abc1};
    \hasseconnector{ca2}{bca1};
    \hasseconnector{ab2}{cab2};
    \hasseconnector{bc2}{abc2};
    \hasseconnector{ca3}{bca2};
    % 3rd row to 4th row
    \hasseconnector{cab0}{bcab0};
    \hasseconnector{abc0}{cabc0};
    \hasseconnector{bca0}{abca0};
    \hasseconnector{cab1}{bcab1};
    \hasseconnector{abc1}{cabc1};
    \hasseconnector{bca1}{abca1};
    \hasseconnector{cab2}{bcab2};
    \hasseconnector{abc2}{cabc2};
    \hasseconnector{bca2}{abca2};
    \hasseconnector{cab0}{cabc0};
    \hasseconnector{abc0}{abca0};
    \hasseconnector{bca0}{bcab1};
    \hasseconnector{cab1}{cabc1};
    \hasseconnector{abc1}{abca1};
    \hasseconnector{bca1}{bcab2};
    \hasseconnector{cab2}{cabc2};
    \hasseconnector{abc2}{abca2};
    \hasseconnector{bca2}{bcab3};
    % 4th row to 5th row
    \hasseconnector{bcab0}{bcabc0};
    \hasseconnector{cabc0}{cabca0};
    \hasseconnector{abca0}{abcab0};
    \hasseconnector{bcab1}{bcabc1};
    \hasseconnector{cabc1}{cabca1};
    \hasseconnector{abca1}{abcab1};
    \hasseconnector{bcab2}{bcabc2};
    \hasseconnector{cabc2}{cabca2};
    \hasseconnector{abca2}{abcab2};
    \hasseconnector{cabc0}{bcabc0};
    \hasseconnector{abca0}{cabca0};
    \hasseconnector{bcab1}{abcab0};
    \hasseconnector{cabc1}{bcabc1};
    \hasseconnector{abca1}{cabca1};
    \hasseconnector{bcab2}{abcab1};
    \hasseconnector{cabc2}{bcabc2};
    \hasseconnector{abca2}{cabca2};
    \hasseconnector{bcab3}{abcab2};
    
  \end{scope}
  
  %%    Add labels to "catchment regions"
  % Factorisable paths in 0th row
  \draw (sa0) node[](){$\source(\alpha)$};
  \draw (sb0) node[](){$\source(\beta)$};
  \draw (sa1) node[](){$\source(\alpha)$};
  \draw (sb1) node[](){$\source(\beta)$};
  \draw (sa2) node[](){$\source(\alpha)$};
  \draw (sb2) node[](){$\source(\beta)$};
  \draw (sa3) node[](){$\source(\alpha)$};
  % Factorisable paths in 1st row
  \draw (a0) node[](){$\alpha$};
  \draw (a1) node[](){$\alpha$};
  \draw (a2) node[](){$\alpha$};
  % Factorisable paths in 2nd row
  \draw (bc0) node[](){$\beta\gamma$};
  \draw (bc1) node[](){$\beta\gamma$};
  \draw (bc2) node[](){$\beta\gamma$};
  % Factorisable paths in 3rd row
  \draw (abc0) node[](){$\alpha\beta\gamma$};
  \draw (bca0) node[](){$\beta\gamma\alpha$};
  \draw (abc1) node[](){$\alpha\beta\gamma$};
  \draw (bca1) node[](){$\beta\gamma\alpha$};
  \draw (abc2) node[](){$\alpha\beta\gamma$};
  \draw (bca2) node[](){$\beta\gamma\alpha$};
  % Factorisable paths in 4th row
  \draw (abca0) node[](){$\alpha\beta\gamma\alpha$};
  \draw (abca1) node[](){$\alpha\beta\gamma\alpha$};
  \draw (abca2) node[](){$\alpha\beta\gamma\alpha$};
  % Factorisable paths in 5th row
  \draw (bcabc0) node[](){$\beta\gamma\alpha\beta\gamma$};
  \draw (bcabc1) node[](){$\beta\gamma\alpha\beta\gamma$};
  \draw (bcabc2) node[](){$\beta\gamma\alpha\beta\gamma$};
  
  \begin{scope}[\fadedcolour]
    %%  Nonfactorisable paths
    % Length 0
    \draw (sc0) node[](){$\source(\gamma)$};
    \draw (sc1) node[](){$\source(\gamma)$};
    \draw (sc2) node[](){$\source(\gamma)$};
    % Length 1
    \draw (b0) node[](){$\beta$};
    \draw (c0) node[](){$\gamma$};
    \draw (b1) node[](){$\beta$};
    \draw (c1) node[](){$\gamma$};
    \draw (b2) node[](){$\beta$};
    \draw (c2) node[](){$\gamma$};
    % Length 2
    \draw (ca0) node[](){$\gamma\alpha$};
    \draw (ab0) node[](){$\alpha\beta$};
    \draw (ca1) node[](){$\gamma\alpha$};
    \draw (ab1) node[](){$\alpha\beta$};
    \draw (ca2) node[](){$\gamma\alpha$};
    \draw (ab2) node[](){$\alpha\beta$};
    % Length 3
    \draw (cab0) node[](){$\gamma\alpha\beta$};
    \draw (cab1) node[](){$\gamma\alpha\beta$};
    \draw (cab2) node[](){$\gamma\alpha\beta$};
    % Length 4
    \draw (bcab0) node[](){$\nodeplaceholder$};
    \draw (cabc0) node[](){$\nodeplaceholder$};
    \draw (bcab1) node[](){$\nodeplaceholder$};
    \draw (cabc1) node[](){$\nodeplaceholder$};
    \draw (bcab2) node[](){$\nodeplaceholder$};
    \draw (cabc2) node[](){$\nodeplaceholder$};
    % Length 5
    \draw (cabca0) node[](){$\nodeplaceholder$};
    \draw (abcab0) node[](){$\nodeplaceholder$};
    \draw (cabca1) node[](){$\nodeplaceholder$};
    \draw (abcab1) node[](){$\nodeplaceholder$};
    \draw (cabca2) node[](){$\nodeplaceholder$};
    \draw (abcab2) node[](){$\nodeplaceholder$};
  \end{scope}
  \end{tikzpicture}
  
  \caption[%
    An example of grading
  ]{%
    \label{fig:example-of-grading}
    \textit{An example of grading.} As before, the factorisable paths are in black and the others in grey. Paths with even grade lie in lighter-coloured regions. Paths with odd grade lie in darker regions. Grade increases up the figure.
    
    At the rightmost corner of any region is a factorisable path. Conversely, any factorisable path is at the rightmost corner of its own region. All paths in a given region are mapped to that factorisable path by $f$.
  }
\end{thesisfigure}

Let $i \in \OC_0$. (These are the objects of $\comppaths$.) Amongst paths having source $i$, there is a unique such $q$ for which $\target(q)$ is the only factorisable vertex through which $q$ passes. We define $f \colon \ilmapsto{ i }{}{ \target(q) }$.

For an arbitrary $\OC$-path, $p \ceq \big( \ilsltpath{\source(p)}{}{\target(p)} \big)$, we then define $f(p) \ceq \big( \ilsltpath{ f(\source(p)) }{}{ f(\target(p)) } \big)$. Towards a second, equivalent, definition in words, observe that there is a unique shortest path $p'$, having $p$ as a prefix, such that $\target(p')$ is factorisable. Amongst the set of suffixes of $p'$ that are factorisable, $f(p)$ is the longest. A third, less formal but perhaps easier to remember, definition is to add as many arrows as possible to the target of $p$, and then remove as many arrows as possible from the source of the resulting path, with the proviso that neither $\source(\alpha)$ nor $\target(\beta)$ is a factorisable vertex for any added arrow $\alpha$ or removed arrow $\beta$.

We mention that $f$ fixes factorisable paths.

That $f$ respects identities is clear. That $f$ respects composition follows from the following straightforward calculation, which uses the identity $\target(p)=\source(p')$ for composable paths $p,p'$. (For clarity, we underline that $f(p)$ denotes the image of $p$ in the functor $f$, and that $f(p)f(p')$ denotes the composition ``$f(p)$ then $f(p')$'' in $\factpaths$.)
$$
\begin{array}{rcl}
f(p)f(p')
  & =
    & \big(
    \ilsltpath{f(\source(p))}{}{f(\target(p))}
  \big)
  \big(
    \ilsltpath{f(\source(p'))}{}{f(\target(p'))}
  \big)
    \\
  & =
    & \big( \ilsltpath{f(\source(p))}{}{f(\target(p'))} \big)
      \\
  & =
    & \big( \ilsltpath{f(\source(pp'))}{}{f(\target(pp'))} \big)
      \\
  & =
    & f(pp')
\end{array}
$$

For an arbitrary $\OC$-path $p$, we define $\mdefine{\grade p} \ceq \grade f(p) $. This is consistent with the definition of factorisable paths since they are fixpoints of $f$.

\nexample In Figure \ref{fig:hasse-diagram-for-C} we drew the Hasse diagram for $\Subpath(\OC)$ and printed the factorisable paths in black, to distinguish them among the other paths in grey. This colouring is recreated in Figure \ref{fig:example-of-grading}.

Additionally in Figure \ref{fig:example-of-grading}, we highlight the ``catchment regions''~-- or perhaps more formally the preimage sets~-- of each factorisable path. We spot for instance that $\source(\gamma)$, $\gamma$ and $\source(\alpha)$ all lie in a single region, reflecting the fact that $f$ takes them all to the common image $\source(\alpha)$. It follows that $\grade \source(\gamma) = \grade \gamma = \grade \source(\alpha)=0 $. As a second example, we find that $f(\beta)=f(\beta\gamma)$ and so $\grade \beta = \grade \beta\gamma = 1$.

By construction, $f$ and therefore $\grade$ are constant on each region. We have colour-coded each region in Figure \ref{fig:example-of-grading} by the parity of grade (ie, whether $\grade$ is odd or even there). The importance of parity will become clear later.

\npass There are some immediate properties of factorisable paths, $f$ and $\grade$ that we collect in the following lemmas and corollary.

\begin{lemma}\label{lem:f-is-monotone} Both $f$ and $\grade$ are weakly monotone with respect to the subpath order $\leq$. This means that if $p \leq q$ then $f(p) \leq f(q)$ and $\grade p \leq \grade q $.
\end{lemma}

\proof If $y,z$ are paths satisfying $ypz = q$, then we calculate $f(q)=f(ypz)=f(y) f(p) f(z) \geq f(p)$ and $\grade q = \grade(ypz) = \grade y + \grade p + \grade z \geq \grade p$.\qed

\begin{lemma}\label{lem:grade-comparison} Write $\leq$ for the subpath order on $\Subpath(\OC)$. Suppose the $\OC$-path $p$ has grade $g$.  
\begin{enumerate}
    \item\label{sublem:grade-g+1-subpath}
    \begin{enumerate}
        \item\label{sublem:existence-grade-g+1-subpath} There exists a factorisable path $y$ of grade $g+1$ with $p<y$ (strictly).
        \item\label{sublem:uniqueness-grade-g+1-subpath} If $\source(p)$ is not factorisable, then there is a unique such $y$. 
    \end{enumerate}
    
    \item\label{sublem:grade-g-1-superpath}
    \begin{enumerate}
        \item\label{sublem:existence-grade-g-1-superpath} There exists a factorisable path $z$ of grade $g-1$ with $z < p$ (strictly), provided that $g > 0$.
        \item\label{sublem:uniqueness-grade-g-1-superpath} If $\target(p)$ is not factorisable, then there is a unique such $z$.
    \end{enumerate}
\end{enumerate}
\end{lemma}

\proof This proof is much clearer when visualised, as in Figure \ref{fig:factorisable-subpaths-and-superpaths}, than when written in prose, as below.
\begin{thesisfigure}
  \centering
  
  \begin{tikzpicture}[scale=0.65]
    \newcommand{\nodeHspacing}[1]{(#1,0)}
    \newcommand{\nodeVspacing}[1]{(0,#1*1.5)}
    \newcommand{\basicnode}[4]{%
      \draw (o) ++\nodeHspacing{#1} ++\nodeVspacing{#2} node[](#3){$ #4 $};
    }
    \newcommand{\nodeplaceholder}{\cdot}
    \newcommand{\facnodeplaceholder}{\circ}
    %%  Place "lattice" of nodes
    % 0th (bottom) row
    \newcommand{\boundingrectangle}{
      (-0.2,-0.2*1.5) rectangle (12.2,7.2*1.5)
    }
    \newcommand{\faccolour}{promptColor}
    \newcommand{\evenfaccolour}{\faccolour!5}
    \newcommand{\oddfaccolour}{\faccolour!15}
    \newcommand{\fadedcolour}{black!30}
    \newcommand{\hasseconnector}[2]{%
      \draw[\fadedcolour] (#1) -- (#2);
    }
    
    \begin{scope}
    
    \clip \boundingrectangle;
    
    \draw (0,0) node[](o){};

    \basicnode{0}{0}{0 0}{\nodeplaceholder};
    \basicnode{2}{0}{2 0}{\nodeplaceholder};
    \basicnode{4}{0}{4 0}{\nodeplaceholder};
    \basicnode{6}{0}{6 0}{x_1 \cdots x_{g-1}}
    \basicnode{8}{0}{8 0}{\nodeplaceholder};
    \basicnode{10}{0}{10 0}{\nodeplaceholder};
    \basicnode{12}{0}{12 0}{\nodeplaceholder};
    
    \basicnode{1}{1}{1 1}{\nodeplaceholder};
    \basicnode{3}{1}{3 1}{\nodeplaceholder};
    \basicnode{5}{1}{5 1}{\nodeplaceholder};
    \basicnode{7}{1}{7 1}{\nodeplaceholder};
    \basicnode{9}{1}{9 1}{\nodeplaceholder};
    \basicnode{11}{1}{11 1}{\nodeplaceholder};
    
    \basicnode{2}{2}{2 2}{\nodeplaceholder};
    \basicnode{4}{2}{4 2}{\nodeplaceholder};
    \basicnode{6}{2}{6 2}{\nodeplaceholder};
    \basicnode{8}{2}{8 2}{p''};
    \basicnode{10}{2}{10 2}{\nodeplaceholder};
    
    \basicnode{3}{3}{3 3}{\nodeplaceholder};
    \basicnode{5}{3}{5 3}{\nodeplaceholder};
    \basicnode{7}{3}{7 3}{\nodeplaceholder};
    \basicnode{9}{3}{9 3}{\nodeplaceholder};    

    \basicnode{4}{4}{4 4}{\nodeplaceholder};
    \basicnode{6}{4}{6 4}{p};
    \basicnode{8}{4}{8 4}{\nodeplaceholder};
    \basicnode{10}{4}{10 4}{\nodeplaceholder};

    \basicnode{3}{5}{3 5}{\nodeplaceholder};
    \basicnode{5}{5}{5 5}{\nodeplaceholder};
    \basicnode{7}{5}{7 5}{p'};
    \basicnode{9}{5}{9 5}{\nodeplaceholder};
    \basicnode{11}{5}{11 5}{\nodeplaceholder};

    \basicnode{2}{6}{2 6}{\nodeplaceholder};
    \basicnode{4}{6}{4 6}{\nodeplaceholder};
    \basicnode{6}{6}{6 6}{\nodeplaceholder};
    \basicnode{8}{6}{8 6}{\nodeplaceholder};
    \basicnode{10}{6}{10 6}{\nodeplaceholder};
    
    \basicnode{1}{7}{1 7}{\nodeplaceholder};
    \basicnode{3}{7}{3 7}{\nodeplaceholder};
    \basicnode{5}{7}{5 7}{\nodeplaceholder};
    \basicnode{7}{7}{7 7}{\nodeplaceholder};
    \basicnode{9}{7}{9 7}{\nodeplaceholder};
    
    \basicnode{0}{8}{0 8}{\nodeplaceholder};
    \basicnode{2}{8}{2 8}{\nodeplaceholder};
    \basicnode{4}{8}{4 8}{\nodeplaceholder};
    \basicnode{6}{8}{6 8}{\nodeplaceholder};
    \basicnode{8}{8}{8 8}{\nodeplaceholder};
    
    % Draw highlighted regions
    \fill[\evenfaccolour, rounded corners] ($(6 0)+0.5*(1,0)$) -- ($(3 3)-0.25*(1,0)+0.75*(0,1.5)$) -- ($(0 0)-1*(1,0)$) -- ($(3, -3*1.5)-0.25*(1,0)-0.75*(0,1.5)$) -- cycle;
    \fill[\evenfaccolour, rounded corners] ($(12 0)+0.5*(1,0)$) -- ($(10 2)-0.25*(1,0)+0.75*(0,1.5)$) -- ($(8 0)-1*(1,0)$) -- ($(10, -2*1.5)-0.25*(1,0)-0.75*(0,1.5)$) -- cycle;
    \fill[\evenfaccolour, rounded corners] ($(11 5)+0.5*(1,0)$) -- ($(8 8)-0.25*(1,0)+0.75*(0,1.5)$) -- ($(7 7)-1*(1,0)$) -- ($(10 4)-0.25*(1,0)-0.75*(0,1.5)$)
    -- cycle;
    \fill[\evenfaccolour, rounded corners] ($(5 7)+0.5*(1,0)$) -- ($(2, 10*1.5)-0.25*(1,0)+0.75*(0,1.5)$) -- ($(0 8)-1*(1,0)$) -- ($(3 5)-0.25*(1,0)-0.75*(0,1.5)$)
    -- cycle;
    
    \fill[\oddfaccolour, rounded corners] ($(9 3)+0.5*(1,0)$) -- ($(6 6)-0.25*(1,0)+0.75*(0,1.5)$) -- ($(4 4)-1*(1,0)$) -- ($(7 1)-0.25*(1,0)-0.75*(0,1.5)$)
    -- cycle;
    
    %%  Draw Hasse connectors
    % 0th to 1st, /
    \hasseconnector{0 0}{1 1};
    \hasseconnector{2 0}{3 1};
    \hasseconnector{4 0}{5 1};
    \hasseconnector{6 0}{7 1};
    \hasseconnector{8 0}{9 1};
    \hasseconnector{10 0}{11 1};
    % 0th to 1st, \
    \hasseconnector{2 0}{1 1};
    \hasseconnector{4 0}{3 1};
    \hasseconnector{6 0}{5 1};
    \hasseconnector{8 0}{7 1};
    \hasseconnector{10 0}{9 1};
    \hasseconnector{12 0}{11 1};
    % 1st to 2nd, /
    \hasseconnector{1 1}{2 2};
    \hasseconnector{3 1}{4 2};
    \hasseconnector{5 1}{6 2};
    \hasseconnector{7 1}{8 2};
    \hasseconnector{9 1}{10 2};
    % 1st to 2nd, \
    \hasseconnector{3 1}{2 2};
    \hasseconnector{5 1}{4 2};
    \hasseconnector{7 1}{6 2};
    \hasseconnector{9 1}{8 2};
    \hasseconnector{11 1}{10 2};
    % 2nd to 3rd, /
    \hasseconnector{2 2}{3 3};
    \hasseconnector{4 2}{5 3};
    \hasseconnector{6 2}{7 3};
    \hasseconnector{8 2}{9 3};
    % 2nd to 3rd, \
    \hasseconnector{4 2}{3 3};
    \hasseconnector{6 2}{5 3};
    \hasseconnector{8 2}{7 3};
    \hasseconnector{10 2}{9 3};
    % 3rd to 4th, /
    \hasseconnector{3 3}{4 4};
    \hasseconnector{5 3}{6 4};
    \hasseconnector{7 3}{8 4};
    \hasseconnector{9 3}{10 4};
    % 3rd to 4th, \
    \hasseconnector{5 3}{4 4};
    \hasseconnector{7 3}{6 4};
    \hasseconnector{9 3}{8 4};
    % 4th to 5th, /
    \hasseconnector{4 4}{5 5};
    \hasseconnector{6 4}{7 5};
    \hasseconnector{8 4}{9 5};
    \hasseconnector{10 4}{11 5};
    % 4th to 5th, \
    \hasseconnector{4 4}{3 5};
    \hasseconnector{6 4}{5 5};
    \hasseconnector{8 4}{7 5};
    \hasseconnector{10 4}{9 5};
    % 5th to 6th, /
    \hasseconnector{3 5}{4 6};
    \hasseconnector{5 5}{6 6};
    \hasseconnector{7 5}{8 6};
    \hasseconnector{9 5}{10 6};
    % 5th to 6th, \
    \hasseconnector{3 5}{2 6};
    \hasseconnector{5 5}{4 6};
    \hasseconnector{7 5}{6 6};
    \hasseconnector{9 5}{8 6};
    \hasseconnector{11 5}{10 6};
    % 6th to 7th, /
    \hasseconnector{2 6}{3 7};
    \hasseconnector{4 6}{5 7};
    \hasseconnector{6 6}{7 7};
    \hasseconnector{8 6}{9 7};
    % 6th to 7th, \
    \hasseconnector{2 6}{1 7};
    \hasseconnector{4 6}{3 7};
    \hasseconnector{6 6}{5 7};
    \hasseconnector{8 6}{7 7};
    \hasseconnector{10 6}{9 7};
    
    %%  Redraw nonfactorisable paths
    \begin{scope}[\fadedcolour]
      \draw (0 0) node[](){$\nodeplaceholder$};
      \draw (2 0) node[](){$\nodeplaceholder$};
      \draw (4 0) node[](){$\nodeplaceholder$};
      \draw (8 0) node[](){$\nodeplaceholder$};
      \draw (10 0) node[](){$\nodeplaceholder$};
      
      \draw (1 1) node[](){$\nodeplaceholder$};
      \draw (3 1) node[](){$\nodeplaceholder$};
      \draw (5 1) node[](){$\nodeplaceholder$};
      \draw (7 1) node[](){$\nodeplaceholder$};
      \draw (9 1) node[](){$\nodeplaceholder$};
      \draw (11 1) node[](){$\nodeplaceholder$};
      
      \draw (2 2) node[](){$\nodeplaceholder$};
      \draw (4 2) node[](){$\nodeplaceholder$};
      \draw (6 2) node[](){$\nodeplaceholder$};
    %   \draw (8 2) node[](){$\nodeplaceholder$};
      \draw (10 2) node[](){$\nodeplaceholder$};
      
      \draw (3 3) node[](){$\nodeplaceholder$};
      \draw (5 3) node[](){$\nodeplaceholder$};
      \draw (7 3) node[](){$\nodeplaceholder$};
    %   \draw (9 3) node[](){$\nodeplaceholder$};
      
      \draw (4 4) node[](){$\nodeplaceholder$};
    %   \draw (6 4) node[](){$\nodeplaceholder$};
      \draw (8 4) node[](){$\nodeplaceholder$};
      \draw (10 4) node[](){$\nodeplaceholder$};
      
      \draw (3 5) node[](){$\nodeplaceholder$};
      \draw (5 5) node[](){$\nodeplaceholder$};
    %   \draw (7 5) node[](){$\nodeplaceholder$};
      \draw (9 5) node[](){$\nodeplaceholder$};
      
      \draw (2 6) node[](){$\nodeplaceholder$};
      \draw (4 6) node[](){$\nodeplaceholder$};
      \draw (6 6) node[](){$\nodeplaceholder$};
      \draw (8 6) node[](){$\nodeplaceholder$};
      \draw (10 6) node[](){$\nodeplaceholder$};
      
      \draw (1 7) node[](){$\nodeplaceholder$};
      \draw (3 7) node[](){$\nodeplaceholder$};
      \draw (7 7) node[](){$\nodeplaceholder$};
      \draw (9 7) node[](){$\nodeplaceholder$};
    \end{scope}
    
    \end{scope}
    
    \begin{scope}
      \draw (6 0) node[](){$x_1 \cdots x_{g-1}$};
      \draw (12 0) node[](){$x_2 \cdots x_g$};
      \draw (8 2) node[](){$p''$};
      \draw (9 3) node[](){$x_1 \cdots x_g$};
      \draw (6 4) node[](){$p$};
      \draw (7 5) node[](){$p'$};
      \draw (5 7) node[](){$x_0 \cdots x_g$};
      \draw (11 5) node[](){$ x_1 \cdots x_{g+1} $};
    \end{scope}
  \end{tikzpicture}
  
  \caption[%
    Factorisable subpaths and superpaths
  ]{%
    \label{fig:factorisable-subpaths-and-superpaths}
    \textit{Factorisable subpaths and superpaths.} Like all paths in the central (darker) region, $p$ is a strict superpath of $x_1 \cdots x_{g-1}$) and a strict subpath of $x_0 \cdots x_g$. 
    
    Related to $p$ are the paths $p'$ and $p''$. Notice that $p'$ is also a strict superpath of $x_2 \cdots x_g$ and $p''$ is a strict subpath of $x_1 \cdots x_{g+1}$. 
  }
\end{thesisfigure}

Each part of the lemma comprises a statement about a strict relation $<$, which implies the existence statement, and also a uniqueness statement. However, the strictness of the relations $z < p < y$ comes automatically from the assumed (strict) inequalities $\grade z < \grade p < \grade y$, so it suffices to prove a statement about the weaker relation $\leq$.

Let us set up notation. Suppose $y \ceq x_0 x_1 \cdots x_g $ is a factorisable path of grade $g+1$, having elementary factors $x_r$, such that $f(p) = x_1 x_2 \cdots x_g$ if $g>0$ or $f(p)=\target(x_0)$ if $g=0$. If $g=1$, then define $z \ceq \source(x_1)$. If $g>1$, then define $z = x_1 x_2 \cdots x_{g-1}$. 

We proceed to show that the $y$ and $z$ just constructed have the claimed properties, starting with the claims about $y$ in part \ref{sublem:grade-g+1-subpath}\ref{sublem:existence-grade-g+1-subpath}. By the (second) definition of $f(p)$, amongst the strict suffixes of $y $ is a path $p'$ having $p$ as a prefix and $f(p)$ as a suffix. We know that $y$ also has $f(p)$ as a suffix, and that the set of paths having $f(p)$ as a suffix is a chain. If $y \leq p'$ then $g+1 = \grade y \leq \grade p' = g$, which is absurd, so $y > p'$. We deduce that $y < p' \leq p$ as claimed. As for part \ref{sublem:grade-g-1-superpath}\ref{sublem:existence-grade-g-1-superpath}, assume $g > 0$. The above definition of $z$ ensures it is the factorisable prefix of $f(p)$ having grade $g-1$. Amongst the suffixes of $p$ is a prefix $p''$ of $x_1 x_2 \cdots x_g$ also having grade $g$. The set of prefixes of $f(p)$ is a chain which includes both $p''$ and $z$, but if $p'' \leq z$ then we obtain the contradiction $g = \grade p'' \leq \grade z = g-1$ by monotonicity, thus $z < p''$ and so $z < p$.

This brings us to the uniqueness statements. We prove part \ref{sublem:grade-g+1-subpath}\ref{sublem:uniqueness-grade-g+1-subpath}. Let $y' \ceq x_s x_{s+1} \cdots x_{g+s}$ be another factorisable path with grade $g+1$, for $0 \leq s \leq g$. (This constraint on $s$ is only to ensure $y$ and $y'$ have common subpaths.) Further let $y'' \ceq x_s x_{s+1} \cdots x_g = \inf\{y,y'\}$. Note that $\grade y'' = g+1-s$.

If $p< y$ and $p < y'$ then $p \leq y''$, whence $g = \grade p \leq \grade y'' = g+1-s$. When $s>1$ we have a contradiction and when $s=0$ we have the tautology $y=y'$. This leaves only the case $s=1$, when $y' = x_1 x_2 \cdots x_{g+1}$. In order to have $p<y'$ and $f(p)=y$, we must have that $f(p)$ is a prefix of $y$, which implies $\source(y)$ is a factorisable vertex. The result follows.

We observe for part \ref{sublem:grade-g-1-superpath}\ref{sublem:uniqueness-grade-g-1-superpath} that all paths have a supremum and so the assumption $0 \leq s \leq g$ is unnecessary, but the remainder of the proof dualises naturally.\qed

\begin{corollary}\label{cor:grade-and-residue} Let $p$ be a $\OC$-path.
\begin{enumerate}
    \item\label{subcor:grade-below-2-means-nonzero-residue} If $ \grade p < 2$, then $p$ has nonzero $A$-residue.
    \item\label{subcor:grade-above-2-means-zero-residue} If $ \grade p > 2$, then $p$ has zero $A$-residue.
\end{enumerate}
\end{corollary}

\proof Using Lemma \ref{lem:grade-comparison}, we can compare $p$ with a factorisable path $x_1 x_2$ having grade 2. We know that $ x_1 x_2 $ represents a component of a commutativity relation. The strict subpaths of $x_1 x_2$ have nonzero $A$-residue (since they represent uniserial $A$-modules) while the strict superpaths of $x_1 x_2$ have zero $A$-residue. The result follows. \qed

\centresubsec{Odd and even grade}

\nlpass{Compression and $\descent$}%To a syllable $\bdp \ceq \big( \ilsyll{ i }{ \ell }{ \ep } \big)$, we associate the $\OQ$-path $ p \ceq \big( \ilsltpath{ i }{ \ell + \ep }{ \syllableplaceholder } \big)$, which we call the \define{compression} of $\bdp$.
Recall the $\OQ$-path $ p \ceq \big( \ilsltpath{ i }{ \ell + \ep }{ \syllableplaceholder } \big)$ is the \define{compression} of the syllable $\bdp \ceq \big( \ilsyll{ i }{ \ell }{ \ep } \big)$.

When $\ep=0$, the compression of $\bdp$ equals the underlying path of $p$. When $\ep = 1$, the compression of $\bdp$ covers the underlying path of $\bdp $ in the prefix order.

We can define $\descent$ on $\OQ$-paths compatibly with its definition on syllables. If $p \ceq \big( \ilslpath{i}{\ell} \big)$ satisfies $\len p > 0 $ and $\ell < a_i + b_i$, then we define $p \descent \ceq \big( \ilsltpath{ i-\ell }{m}{ \syllableplaceholder} \big)$ for $m \ceq a_i + b_i - \ell$. It is straightforward to verify that $\bdp \descent$ is defined iff the compression $p \descent$ is defined, in which case compression and $\descent$ commute on $\bdp$.

\nlpass{Shortest dependent paths and $\descent$} Recall from passage \ref{def:source-encoding} the source encoding of the permissible data of $A$ as integer sequence $(a_i)_{i \in \OQ_0}$ and bit sequence $(b_i)_{i \in \OQ_0}$.

For $i \in \OC_0$, let $\mdefine{r_i} \ceq \big( \ilsltpath{ i }{ a_i + b_i }{ \syllableplaceholder } \big)$. When $b_i=0$ (which means exactly that $i$ represents the source of a component of a commutativity relation), $r_i$ represents that component. When $b_i = 1$, $r_i$ is the shortest path with source $i$ that has zero $A$-residue. Either way, $r_i$ is the shortest path with source $i$ whose residue  dependent linearly on another $A$-path (perhaps for trivial reasons). For this reason, we call $r_i$ the \define{shortest dependent path} with source $i$.

The importance of shortest dependent paths becomes clear when we examine how $ \descent $ interacts with compression, for if $ p \in \supp \descent $ is a path then $(p)(p \descent) = r_{\source(p)}$.

\begin{lemma}
Any shortest dependent path $r_i$ satisfies $\grade r_i = 2$.
\end{lemma}

\proof For convenience, let $g \ceq \grade r_i$.

Either $r_i$ represents a component of a commutativity relation or it does not. We treat these cases in turn.

If $r_i$ does represent a component of a commutativity relation, then it has nonzero $A$-residue and therefore $g \leq 2$, by the contrapositive to Corollary \ref{cor:grade-and-residue}\ref{subcor:grade-above-2-means-zero-residue}. If $g \leq 1 $ then, by Lemma \ref{lem:grade-comparison}\ref{sublem:grade-g+1-subpath}\ref{sublem:existence-grade-g+1-subpath}, $r_i$ is a strict subpath of some factorisable path $x_1 x_2$, which also represents a component of a commutativity relation, which is absurd. Thus $g \geq 2$ which gives $g = 2$.

Otherwise, $r_i$ represents a monomial relation. We know in this case that $r_i$ is not factorisable because its source is certainly not a factorisable vertex. We also know that $r_i$ has zero $A$-residue, and so $\grade r_i \geq 2$ by the contrapositive to Corollary \ref{cor:grade-and-residue}\ref{subcor:grade-below-2-means-nonzero-residue}. Write  $x_1 x_2 \cdots x_g \ceq f(p)$ for elementary factors $x_r$. %Suppose for contradiction that $\grade r_i \geq 3$.

From Lemma \ref{lem:grade-comparison}\ref{sublem:grade-g-1-superpath} we deduce that $r_i$ is the strict superpath of a unique factorisable path $x_1 x_2 \cdots x_{g-1}$ of grade $g-1$. Let $p,q$ be the paths such that $r_i = p x_1 x_2 \cdots x_{g-1} q$. We know that $p$ has positive length, because $\source(r_i)$ is not a factorisable vertex, and that $q$ has positive length, since otherwise $f(r_i)=x_1 x_2 \cdots x_{g-1}$ which is an absurdity.

Towards a contradiction assume that $g \geq 3$, which implies that $p x_1 x_2$ is a strict prefix of $r_i$. Since $p x_1 x_2$ is a strict superpath of $x_1 x_2$ (which represents a component of a commutativity relation), the $A$-residue of the former is zero and so (trivially) depends linearly on another $A$-path. Since $r_i$ is the minimal such path with source $\source(r_i)=\source(p x_1 x_2)$ having this property, we find $r_i \leq p x_1 x_2$. But $p x_1 x_2 < r_i$, which is our desired contradiction. Thus $2 \leq g \leq 2$, from which equality ensues.\qed

\begin{corollary}\label{lem:descent-respects-parity}
If $p$ is a $\OC$-path and $p\descent$ is defined, then $\grade( p \descent) = 2 - \grade( p)$. 
\end{corollary}

\proof Previous results give $ 2 = \grade r_{\source(p)} = \grade( (p) (p\descent) ) = \grade p + \grade(p \descent)$.\qed

\npass We have just proven that, where defined, $\descent$ respects the parity of $\grade$. Now it remains to prove that injective syllables have even grade and elementary factors have odd grade. The second statement is trivial.

\begin{lemma}\label{lem:elementary-factors-have-grade-1}
For all elementary factors $u_r, v_r$ we have $\grade u_r = \grade v_r = 1$.
\end{lemma}

\proof This is immediate from the definition of $\grade$.

\nlpass{Injective syllables} An \define{injective syllable} is a syllable of the form $\big( \injsyll{c_i}{i}{\ep} \big)$ for some $\ep \in \{0,1\}$, where $i \in \OQ_0$ satisfies $d_i=1$.

By definition of $c_i$, $\big( \ilsltpath{\syllableplaceholder}{c_i}{i} \big)$ is the longest path with target $i$ having nonzero $A$-residue.

\nrex The only vertices $i \in \OC_0$ for which $d_i = 1$ are the $\source(\gamma)$. These satisfy $c_{\source(\gamma)}=2$. The only injective syllables on this component are therefore $\big( \injsyll{\alpha\beta}{\source(\gamma)}{0} \big)$ and  $\big( \injsyll{\alpha\beta}{\source(\gamma)}{1} \big)$. These have compressions $\alpha\beta$ and $\alpha\beta\gamma$ respectively. It is easy to verify that $\alpha\beta$ has nonzero $A$-residue as claimed.

\begin{lemma}\label{lem:injective-syllables-have-grade-2}
If $p$ is the compression of an injective syllable $\bdp \ceq \big( \injsyll{c_i}{i}{\ep} \big)$, then $\grade p = 2$.
\end{lemma}

\proof First we claim that it does not matter whether $\ep=0$ or $\ep =1 $. Recall from Passage \ref{def:grading-paths-and-f} the definition of $f$ on the objects of $\OC$ (ie, on vertices). Since $d_i=0$ by assumption, $i$ does not represent the target of any commutativity relation and so certainly $i$ is not an factorisable vertex. This implies the path $q$ in the definition of $f(i)$ has positive length, and so $f(i)=f(i-1)$, from which our claim follows. We therefore assume $\ep=0$ to ensure $\target(p)=i$.

Let $x_1 x_2 x_3$ be the factorisable path (with elementary factors $x_r$) such that $x_3$ passes through $i$, and let $q$ be the strict prefix of $x_3$ with $\target(q) = i$. Necessarily, $\len q > 0$ and $\grade q = 1$.

Consider $x_2 q$, whose target is $i$. We know $x_2 q$ has nonzero $A$-residue because it is a strict subpath of $x_2 x_3$. We deduce $p \geq x_2 q$ and hence $\grade p \geq \grade( x_2 q) = \grade(x_2) + \grade(q) = 2$ by monotonicity. But $\grade p \leq 2$ because $p$ has nonzero $A$-residue. We thus confirm $\grade p = 2$.\qed

\nexample The compressions of injective syllables for our running example algebra $A$ are $\alpha\beta$ and $\alpha\beta\gamma$. We verify in consultation with Figure \ref{fig:example-of-grading} that these both do have grade $2$.

\begin{proposition}\label{prop:descent-orbits-of-inj-sylls-do-not-have-interior-tails}
Let $\bdp$ be an injective syllable. If $\bdp \descent^r$ is defined for all $r \geq 0$, then the forward orbit $(\bdp \descent^r)_{r \geq 0}$ has no tail $(\bdp\descent^r)_{r \gg 0}$ containing only interior syllables.
\end{proposition}

\proof If there exists no $\descent$-cycle $
\begin{tikzcd}[sep=small]
  \bdu_r \arrow[r, mapsto, "\descent"]
    \& \bdv_r \arrow[r, mapsto, "\descent"]
      \& \bdu_{r+1}
\end{tikzcd}
$ ($r \in \integers$) of interior syllables lying on the same $\OQ$-component as $\bdp$ then the claim follows automatically. Hence assume there does exist such a cycle.

There are only finitely many syllables (up to periodicity) and so there certainly exists a (periodic) interior $\descent$-cycle $(\bdu_1, \bdv_1, \dots, \bdu_m, \bdv_m)$, with underlying paths $u_1,v_1,\cdots u_m, v_m$. Associated to this cycle is a grading $\grade$ that we can carry over to syllables by stipulating that $\grade( \bdq)$ equals the grade of the compression of $\bdq$. Corollary \ref{lem:descent-respects-parity} and Lemmas \ref{lem:elementary-factors-have-grade-1} and \ref{lem:injective-syllables-have-grade-2} combine to prove that (for any $r \in \integers$)
$$\big( \grade(\bdu_r\descent^s) \big)_{s \geq 0} = \big( \grade(\bdv_r \descent^s) \big)_{s \geq 0}$$
is consistently odd while $\big( \grade(\bdp\descent^s) \big)_{s \geq 0}$ is consistently even.

We deduce that no $\bdu_r$ or $\bdv_r$ appears in $(\bdp \descent^s)_{s \geq 0}$, from which the result follows.

\nlpass{Example: injective syllables and perturbation} The previous lemma addresses the behaviour of an injective syllable $\bdp$ when only $\descent$ is applied, which is the typical behaviour of a column in the syzygy fabric. We know from previous discussion that the impact, if any, of a neighbouring column is to perturb syllables; that is, to change the stability term $\ep$ of a syllable $\big( \ilsyll{i}{\ell}{\ep} \big)$ from $0$ to $1$. It is sensible to wonder what happens if some combination of $\descent$ and perturbation is applied to an injective syllable $\bdp$.

Recall our running example algebra $A$. The injective syllables of interest%
\footnote{\ \newcommand{\OD}{\mathcal{D}}
The algebra $A$ has two other injective syllables $\big(\injsyll{\delta^3}{\source(\delta)}{0}\big)$ and $\big(\injsyll{\delta^3}{\source(\delta)}{1}\big)$, both of which lie on a different component $\OD$ of $\OQ$ from the other syllables discussed. Because all vertices $i\in \OD_0$ have $b_i=1$, whenever $\bdq \descent$ is defined for some syllable $\bdq$ on $\OD$,  $\bdq\descent$ is a boundary syllable, and therefore these injective syllables on $\OD$ satisfy the strengthened claims of this passage for more trivial reasons.}
of $A$ are $\big( \ilsyll{\source(\alpha)}{\alpha\beta}{0} \big)$ and $\big( \ilsyll{\source(\alpha)}{\alpha\beta}{1} \big)$, the latter being the perturbation of the former. We calculate
$$
\begin{tikzcd}[sep=small]
\big( \source(\alpha) \patharrow{rr}{\alpha\beta}
  \&
    \& \syllableplaceholder \patharrow{r}{0}
      \& \syllableplaceholder\big) \ar[r, mapsto, "\descent"]
        \& \big( \source(\gamma) \patharrow{rr}{\gamma}
          \&
            \& \syllableplaceholder \patharrow{r}{0}
              \& \syllableplaceholder\big) \ar[r, mapsto, shift left, "\descent"]
                \& \big( \source(\alpha) \patharrow{rr}{\alpha} \ar[l, mapsto, shift left]
                  \&
                    \& \syllableplaceholder \patharrow{r}{1}
                      \& \syllableplaceholder\big)\tstop
\end{tikzcd}
$$
The leftmost and middle of these syllables are interior, so can be perturbed to $\big( \ilsyll{\source(\alpha)}{\alpha\beta}{1} \big) \notin \supp \descent $ and $\big( \ilsyll{\source(\gamma)}{\gamma}{1} \big)$ respectively. The forward $\descent$-orbit of the latter is as follows.
\begin{equation}\label{eqn:descent-orbit-of-gamma-1}
\begin{tikzcd}[sep=small]
\big( \source(\gamma) \patharrow{rr}{\gamma}
  \&
    \& \syllableplaceholder \patharrow{r}{1}
      \& \syllableplaceholder\big) \ar[r, mapsto, "\descent"]
        \& \big( \source(\beta) \patharrow{rr}{e_{\source(\beta)}}
          \&
            \& \syllableplaceholder \patharrow{r}{1}
              \& \syllableplaceholder\big) \ar[r, mapsto, shift left, "\descent"]
                \& \big( \source(\gamma) \patharrow{rr}{\gamma\alpha} \ar[l, mapsto, shift left]
                  \&
                    \& \syllableplaceholder \patharrow{r}{0}
                      \& \syllableplaceholder\big)\tstop
\end{tikzcd}
\end{equation}
Of these, only the rightmost is interior. Its perturbation is $\big( \ilsyll{\source(\gamma)}{\gamma\alpha}{1} \big) \notin \supp \descent$.

The preceding calculation exhausts all possible combinations of $\descent$ and perturbation, starting from an injective syllable of interest. We highlight that nowhere in this calculation do $\big(\ilsyll{\source(\alpha)}{\alpha}{0}\big)$ and $\big(\ilsyll{\source(\beta)}{\beta\gamma}{0}\big)$ arise, these being the syllables on the only interior $\descent$-cycle for $A$.

We find in this example that a stronger statement than the previous lemma holds. Not only does the $\descent$-orbit of an injective syllable have no tails of interior syllables, as shown above, but moreover no amount of perturbation of such an orbit ever yields one with such a tail.

This stronger statement holds in all SB algebras we have ever considered, but if it is to be proven then some refinement of our grading will be necessary. To see this, consider the three syllables in (\ref{eqn:descent-orbit-of-gamma-1}). Their compressions are respectively $\gamma\alpha$, $\beta$ and $\gamma\alpha$. These all have grade $1$, as shown in Figure \ref{fig:example-of-grading}, and so $\grade$ alone cannot distinguish them from the syllables of an interior $\descent$-cycle.

% \begin{enumerate}
%     \item the first option goes here
%     \item \begin{enumerate}
%         \item the second option goes here
%         \item \begin{enumerate}
%             \item the third option goes here
%             \item \begin{enumerate}
%                 \item the fourth option goes here
%                 \item the fourth$'$ option goes here
%             \end{enumerate}
%         \end{enumerate}
%     \end{enumerate}
% \end{enumerate}

\subsection{SB algebras with few vertices}\label{subsec:sb-algs-with-few-verts}

We conjecture that injective modules over an SB algebra have finite syzygy type. A general proof eludes us but, using the tools developed in this thesis, we can address the situation for small SB algebras. This is the focus of this subsection.

In the following, as ever, $A = \kk Q / \langle \rho \rangle $ is a basic, connected SB algebra with ordinary quiver $Q$. For convenience, we have $Q=\{1,2,\dots,n\}$.

We also introduce some new notation regarding cycles; that is, paths whose source and target coincide. For any cycle $\alpha_1\alpha_2\cdots \alpha_m$ with source and target $i$, we let \define{$[\alpha_1\alpha_2\cdots \alpha_m]$} denote the equivalence class
$$
\mdefine{[\alpha_1 \alpha_2 \cdots \alpha_m]}
\ceq
\{
  \alpha_1 \alpha_2 \alpha_3 \cdots \alpha_{m-1} \alpha_m,~
  \alpha_2 \alpha_3 \alpha_4 \cdots \alpha_m \alpha_1,~
  \dots,~
  \alpha_m \alpha_1 \alpha_2 \cdots \alpha_{m-2} \alpha_{m-1}
\}
$$
of its rotates. Additionally, we will stipulate that $(\alpha_1 \alpha_2 \cdots \alpha_m)^0 \ceq e_i$.

\begin{lemma}\label{lem:inj-have-fin-syz-type-over-2-vertices}
If $Q$ has $2$ vertices, then the indecomposable injective string modules for $A$ have finite syzygy type and, dually, the indecomposable projective string modules have finite cosyzygy type.
\end{lemma}

\proof The dual result follows immediately from the first, since we may apply the first result to $A^\op$ and then use the vector-space duality functor $\dual$.

The rest of the proof is by brute force, requiring two stages. The first stage is to enumerate all possible permissible data for SB algebras on $2$ vertices. We will do this over the coming paragraphs, using previous results to cut down our sample space. Once the list of cases has been scoped out, we then have to actually perform the calculations. This is the second stage. It only comprises direct computations of the sort that our \sbstrips\ package performs easily. We relegate the outcomes of these straightforward efforts to Appendix \ref{app:sb-algs-with-few-verts}.

Towards enumerating all possible permissible data, let us first consider the pin graph $\pingraph_A$ of $A$. This is necessarily a sub-$1$-regular quiver on $2$ vertices and so has $0$, $1$ or $2$ arrows. If $\pingraph_A$ has $2$ arrows, then $A$ is necessarily selfinjective and the result follows. Similarly if $\pingraph_A$ has $0$ arrows, then $A$ is necessarily monomial and the result follows. If it has $1$ arrow that is not a loop, then $\pingraph_A$ is acyclic and so $A$ is syzygy-finite by Theorem \ref{thm:explicit-bounds-on-syzygy-type}. This only leaves the case that $\pingraph_A$ comprises a single loop and a single isolated vertex; without loss of generality, we may assume $\pingraph_A = \big( \begin{tikzcd}[] 1 \ar[loop right] \& 2\end{tikzcd}
\big)$.

Next, we find an overquiver $\OQ$ and vertex exchange map $\dagger$ compatible with $\pingraph_A$. Obviously $\OQ$ is a $1$-regular quiver on $4$ vertices. Let $1,1'$ and $2,2'$ respectively be the two $\OQ$-vertices representing $1,2 \in Q_0$. Compatibility with $\pingraph_A$ requires that there be $\OQ$-paths $p,q$ satisfying $\{\source(p),\source(q)\} = \{1,1'\} = \{\target(p),\target(q)\}$ and satisfying $\len p, \len q \geq 2 $. This implies the existence of $\OQ$-cycles at $1,1'$, the shortest of which we will respectively denote by $\sigma_1$ and $\sigma_{1'}$; by exchanging $1$ and $1'$ where necessary, we may assume $\len \sigma_1 \leq \len \sigma_{1'}$. Since cycles of distinct length are disjoint and disjoint cycles have no arrows in common, we deduce that
$$
\len \sigma_1 \neq \len \sigma_{1'}
\Longrightarrow
[\sigma_{1}] \neq [\sigma_{1'}]
\Longrightarrow
(\len \sigma_1) + (\len \sigma_{1'}) \leq 4
\tstop
$$
Consequently, we have $(\len \sigma_1, \len \sigma_{1'}) \in \big\{ (1,1), (1,2), (1,3), (2,2), (3,3), (4,4) \big\}$.

If there were loops at $1$ and $1'$, or if there were a $2$-cycle $\big(\begin{tikzcd}[sep=small]
1 \ar[r] \& 1' \ar[r] \& 1
\end{tikzcd}\big)$, then this would force $Q$ to be disconnected. Therefore we can entirely eliminate the case $(\len \sigma_1, \len \sigma_{1'})=(1,1)$ and we can further constrain the case $(\len \sigma_1, \len \sigma_{1'})=(2,2)$.

Up to exchanging $\ilcorresp{1}{}{1'}$ and $\ilcorresp{2}{}{2'}$ if need be, the remaining possibilities for $\OQ$ are shown exhaustively in Figure \ref{fig:candidate-overquivers}. Aside from Figures \ref{subfig:untwisted-4-4-overquiver} and \ref{subfig:twisted-4-4-overquiver} where $(\len \sigma_1, \len \sigma_{1'})=(4,4)$, there is a unique possible overquiver for each potential cycle-length profile.
\begin{thesisfigure}
\centering

\subcaptionbox{%
  \label{subfig:1-2-overquiver}
  $(1,2)$
}[%
  0.3\textwidth
]{
  \begin{tikzcd}
  1 \ar[loop left, "\alpha"]
    \& 2 \arrow[dl, shift left, "\gamma"]
      \\
  1' \arrow[ur, shift left, "\beta"]
    \& 2' \arrow[loop right, "\delta"]
  \end{tikzcd}
}
\subcaptionbox{%
  \label{subfig:1-3-overquiver}
  $(1,3)$
}[%
  0.3\textwidth
]{
  \begin{tikzcd}
  1 \ar[loop left, "\alpha"]
    \& 2 \arrow[d, "\gamma"]
      \\
  1' \arrow[ur, "\beta"]
    \& 2' \arrow[l, "\delta"]
  \end{tikzcd}
}
\subcaptionbox{%
  \label{subfig:2-2-overquiver}
  $(2,2)$
}[%
  0.3\textwidth
]{
  \begin{tikzcd}
  1 \ar[r, shift left, "\alpha"]
    \& 2 \ar[l, shift left, "\beta"]
      \\
  1' \ar[r, shift left, "\gamma"]
    \& 2' \ar[l, shift left, "\delta"]
  \end{tikzcd}
}

\vspace{4ex}

\subcaptionbox{%
  \label{subfig:3-3-overquiver}
  $(3,3)$
}[%
  0.3\textwidth
]{
  \begin{tikzcd}
  1 \ar[d, "\alpha"']
    \& 2 \ar[l, "\gamma"']
      \\
  1' \ar[ur, "\beta"']
    \& 2' \ar[loop right, "\delta"]
  \end{tikzcd}
}
\subcaptionbox{%
  \label{subfig:untwisted-4-4-overquiver}
  $(4,4)$
}[%
  0.3\textwidth
]{
  \begin{tikzcd}
  1 \ar[r, "\alpha"]
    \& 2 \ar[d, "\beta"]
      \\
  1' \ar[u, "\delta"]
    \& 2' \ar[l, "\gamma"]
  \end{tikzcd}
}
\subcaptionbox{%
  \label{subfig:twisted-4-4-overquiver}
  $(4,4)$
}[%
  0.3\textwidth
]{
  \begin{tikzcd}
  1 \ar[r, "\alpha"] \ar[from=dr, "\delta" very near end]
    \& 2 \ar[dl, crossing over, "\beta" very near start]
      \\
  1' \ar[r, "\gamma"']
    \& 2'
  \end{tikzcd}
}

\caption[%
  Candidate overquivers
]{%
  \label{fig:candidate-overquivers}
  \textit{Candidate overquivers.} Each is shown with its cycle-length profile $(\len \sigma_1, \len \sigma_{1'})$
}
\end{thesisfigure}

Henceforth fix some overquiver $\OQ$ from among the candidates. Our next job is to determine the component set $C$; in essence, this means finding the $p$ and $q$ mentioned above. These paths must be incomparable in $\Subpath(\OQ)$. If $1$ and $1'$ lie on distinct connected components (as in the top row of Figure \ref{fig:candidate-overquivers}), then we may set simply $p \ceq \sigma_1^r$ and $q \ceq \sigma_{1'}^s$ (for exponents $r,s \geq 0$ independently chosen to ensure $\len p, \len q \geq 2$) and we will obtain the incomparability condition automatically. If $1$ and $1'$ lie on the same connected component (as in the bottom row), then there are two cases to consider. Let $u_1$ and $u_{1'}$ respectively denote the shortest path $\big(\ilsltpath{1}{}{1'}\big)$ and  $\big(\ilsltpath{1'}{}{1}\big)$. For some $r \geq 0$, our components are
$$
\text{either }
p \ceq (u_1 u_{1'})^r u_1
\text{ and }
q \ceq (u_{1'} u_1)^r u_{1'}
\text{\hspace{1cm}or }
p \ceq (u_1 u_{1'})^r
\text{ and }
q \ceq (u_{1'} u_1)^r\tstop
$$
(In the first case, the total number of elementary factors $u_i$ in each component is odd; in the second case, the number is even. In both cases, the exponent $r$ must again be chosen to satisfy the length requirement.)

Once $C$ is chosen, the final step is to compatibly choose the set $N$ of nonzero paths (or, equivalently, its complement). This choice may be performed independently for each connected component of $\OQ$, and is subject to our usual convention that stationary paths are always nonzero. Besides that, compatibility with $C$ requires that subpaths of component paths belong to $N$ but strict superpaths of component paths do not. On connected components $\OC$ of $\OQ$ that do feature a component path, this leaves only a finite number of choices for $\OC \cap N$, since all sufficiently long paths are superpaths of a component path. On connected components that do not, our choice of $\OC \cap N$ is unconstrained. The only such unconstrained connected components $\OC$ occur in the overquivers of Figures \ref{subfig:1-2-overquiver} and \ref{subfig:3-3-overquiver}, where they comprise a single arrow $\delta$, therefore choosing $\OC \cap N$ amounts to choosing some smallest positive power $t$ of $\delta$ which lies outside $N$.

Taken altogether, this gives permissible data $(\OQ,N,C,\dagger)$ for $A$.

We know that there are finitely many options for $\OQ$. As we have seen above and will see in greater detail in Appendix \ref{app:sb-algs-with-few-verts}, for each of these options the component paths $C$ as well as minimal elements of $N \setminus \Subpath(\OQ)$ can be parameterised by integer exponents $r,s,t \geq 0$. The parameter space can be partitioned into finitely many regions, in each of which the behaviour can be described uniformly. This reduces the calculation to a finite effort, which we undertake in the appendix.\qed

\begin{theorem}\label{thm:ig-and-pc-if-A-has-2-simples}
If $Q$ has at most $2$ vertices, then the indecomposable injective string modules for $A$ have finite syzygy type and, dually, the indecomposable projective string modules have finite cosyzygy type.
\end{theorem}

\proof If $Q$ has $1$ vertex, then it is selfinjective or monomial and the claim holds. The $2$-vertex case was settled in the previous lemma.\qed

\begin{corollary}
If $Q$ has at most $2$ vertices, then injectives generate and projectives cogenerate for $A$.

Consequently $A$ satisfies the big and little finitistic dimension conjectures under this hypothesis.
\end{corollary}

\nrmk Our Lemma \ref{lem:inj-have-fin-syz-type-over-2-vertices} is not technically a new result. Differing terminology aside, Salmer{\'o}n proved in \cite{Sal94} that the regular $A$-module has finite cosyzygy type, and thereby deduced that SB algebras on at most $2$ vertices satisfy the little finitistic dimension conjecture.

Nonetheless, we include Lemma \ref{lem:inj-have-fin-syz-type-over-2-vertices} and its proof for four reasons, the simplest reason being an aesthetic one: we believe that this proof beautifully uses the discrete information of a SB algebra that we have laboured to quantify. Each part of the permissible tuple $(\OQ,N,C,\dagger)$ is given its turn in the spotlight, as is the pin graph $\pingraph_A$.

The second reason is the obscurity of Salmer{\'o}n's article. It features in a Mexican journal whose back-catalogue is not available online and which we could not access using the standard channels. We only know the contents of the paper because we contacted Salmer{\'o}n directly and he generously replied with photographs of a printed copy. Indeed, we only know of the existence of the paper thanks to a passing reference in an article by Erdmann \etal\ \cite{EHIS04}, in which an entirely different approach is used to show all SB algebras satisfy the little finitistic dimension conjecture. It appears that the ensuing focus in the literature on Erdmann \etal's finitistic conclusion has come at the expense of (ours and) Salmer{\'o}n's syzygy-finiteness criterion. This is a wrong we are happy to right by recording a proof in an open-access thesis.

The third reason is a methodological one. Salmer{\'o}n's approach was to represent string modules by string graphs in the conventional way and then determine their syzygies in the fashion that Liu and Morin would later formally codify \cite{LM04}, which is to say he wrote string graphs simply as juxtapositions of independent paths and their formal inverses but did not mark those syllables at the boundary of a string graph as different from those in the interior. Accordingly, he makes no acknowledgement of the pin-boundary phenomenon over which we have agonised greatly. He describes performing the requisite calculations by hand but does not mention any effort to implement these calculations in a computer. We venture that our approach, which does highlight this important phenomenon and which has been realised computationally, is merit-worthy by being transparent and, thanks to the existence of \sbstrips, practical.

The final, and most crucial, reason is independence. Our work arose independently of Salmer{\'o}n's. It is true that we both explore the syzygy-finiteness of injectives and that we both reduce the amount of brute-force calculation required by using auxiliary results. Our auxiliary results (Corollary \ref{cor:pin-graph-acyclic-implies-syz-finite}) and Salmer{\'o}n's auxiliary result (\cite[Prop 2.1]{Sal94}) both find cases in which an SB algebra is (co)syzygy-finite. But his requires that no vertex of the ground quiver have total degree $4$, which forces there to be no $\pingraph_A$-paths of length $2$, which in particular implies that $\pingraph_A$ is acyclic. This is the hypothesis of our auxiliary result, hence our reduction step generalises his. Not only did we prove the result independently of Salmer{\'o}n, but we did so more economically than he did.

% We identified the key criterion (namely the syzygy-finiteness of injectives) independently, well before we learned anyone else had successfully exploited it. 

% Salmer{\'o}n's conclusion that small SB algebras satisfy the little finitistic dimension conjecture was generalised to all SB algebras in \cite{EHIS04}, but the authors of the latter paper exploit a different condition, with no (immediate) connection to the criterion. We only know of Salmer{\'o}n's paper thanks to a passing reference in \cite{EHIS04}, which mentions only the finitistic conclusion and not the criterion of syzygy-finiteness. Indeed, we only know the contents of Salmer{\'o}n's paper because we contacted him directly and he replied with photographs of a printed copy. His paper is printed in a small Mexcian journal whose back issues we could not access online nor through the University of Bristol library. This obscurity is the second reason. This thesis will eventually be open-access

\appendix
\chapter{SB algebras with 2 vertices}\label{app:sb-algs-with-few-verts}

Lemma \ref{lem:inj-have-fin-syz-type-over-2-vertices} claims that the indecomposable injective string modules over a $2$-vertex SB algebra $A$ have finite syzygy type. Every SB algebra admits a pin graph $\pingraph_A$ and permissible data $(\OQ, N, C, \dagger)$; the argument contained in Section \ref{subsec:sb-algs-with-few-verts} constraints what $\pingraph_A$ and $(\OQ,N,C,\dagger)$ can be.

One constraint is that $\pingraph_A = \big( \begin{tikzcd}[] 1 \ar[loop right] \& 2\end{tikzcd}
\big)$, which forces the indecomposable injective $E_1$ at $1$ is pin. This leaves $E_2$ as the only (indecomposable) injective string module up to isomorphism, reducing the necessary calculations to just the syzygies of this one module.

We have seen that only the overquivers in Figure \ref{fig:candidate-overquivers} are possible; recall that there these overquivers were labelled (a) through (f). We consider each overquiver $\OQ$ in order, further subdividing each case using the finitely many families of choices of $C$ compatible with a given $\OQ$ and, then, the finitely many families of choices of $N$ compatible with a given $C$.

We determine $N$ by determining its complement $\Subpath(\OQ) \setminus N $. In addition to the strict superpaths of component paths, there may be other generators of the complement; by this, we mean additional $\OQ$-paths $p$ that are incomparable to any component path, that satisfy $p \leq q \Longrightarrow q \notin N$ and that are minimal with respect to this property.

In each instance, the members of $C$ as well as the additional generators of $\Subpath(\OQ) \setminus N$ are parameterised by some nonnegative integers $t,s,r$ each satisfying a situational lower bound. The author performed explicit calculations using the \sbstrips\ package for small, explicit choices of parameters~-- roughly, those parameter choices in $\{ (t,s,r) \in \naturals^3 \colon  0 \leq t + s + r \leq 10\}$. From this sample patterns became clear and it was possible to partition parameter space into finitely many regions, with the behaviour in any given region able to written as a function of the parameters. Since \sbstrips\ is unable to perform general symbolic calculation, the author subsequently had to verify by hand that, within any given region, the formulas given for syzygy behaviour were accurate.

As an example partition of parameter space, let us consider case (\subref{subfig:2-2-overquiver})\ref{item:fin-many-sporadic-values}. There are two regions of parameter space that are finite (each is in fact a singleton) while the remaining region is unbounded. By contrast, in case (\subref{subfig:1-2-overquiver})\ref{item:example-to-discuss} there are five regions: two regions that are singletons; two unbounded regions where $s$ and $r$ may independently be arbitrarily large but $t$ is fixed; and one region where each of $t$, $s$ and $r$ may independently take arbitrarily large values.

This exemplifies a general truth: in every case of parameter space, there is exactly one region where all parameters may independently take on arbitrarily large values. In regions other than that single ``largest'' unbounded region, at least one parameter forces the existence of a component of length $2$ or a generator of $\Subpath(\OQ)\setminus N$ of length $1$ or $2$. 

% Using the \sbstrips\ package, the author ran sample calculations on those algebras in each case for which the parameter values were small~-- roughly speaking, those parameter values $(t,s,r)$ satisfying $\{ (t,s,r) \in \naturals^3 \colon  0 \leq t + s + r \leq 10\}$ (where parameters were assigned value $0$ when absent). 

% When $t=1$, the imposed relation is $\delta$ (which has length $1$). 

% The list of possibilities are as follows. The first level of the list dictates $\OQ$, the second level specifies $C$ and the third determines $N$ by mentioning any additional generators of its complement. 
% \begin{enumerate}
%     \item \begin{enumerate}
%         \item $C = \{\alpha^r, (\beta\gamma)^s\}$ ($r \geq 2$, $s \geq 1$)
%         \begin{enumerate}
%             \item $\delta^t$ ($t \geq 1$)
%             \item $\delta^t$ ($t \geq 1$) and $(\gamma\beta)^s$
%         \end{enumerate}
%     \end{enumerate}
%     \item \begin{enumerate}
%         \item $C = \{\alpha^r, (\beta\gamma\delta)^s\}$ ($r \geq 2$, $s \geq 1$)
%         \begin{enumerate}
%             \item \textit{No further generators}
%             \item $(\gamma\delta\beta)^s\gamma$
%             \item $ $
%         \end{enumerate}
%     \end{enumerate}
% \end{enumerate}

To record the indecomposable direct summands of syzygies of $E_2$ in each case, we will draw its \define{syzygy quiver}.\footnote{\ Cibils originated this concept \cite[\S5]{Cib93}, but we prefer the variant of Howard \cite[\S3]{How15} that permits multiple arrows between vertices since we believe it is more transparent.} The vertices of the syzygy quiver in each case are the (isomorphism classes of) indecomposable modules occurring as direct summands of $\syzygy^k E_2$ for some $k \geq 0$. We draw exactly $m$ arrows $\ilarrow{X}{}{Y}$ iff $\syzygy^1 X$ has exactly $m$ direct summands isomorphic to $Y$. That only finitely many vertices appear in each quiver implies that $E_2$ is syzygy finite. The length of the longest simple path in this quiver will give the index at which the additive class of syzygies of $E_2$ stabilise, and the vertices on cycles in this quiver will give the direct summands of the additive generator of this stable class.

For the sake of space, we will denote string modules just by their string graphs, and we will write strings graphs as words of paths and their formal inverses. Thus if
$
v_2 \ceq \big(
  \begin{tikzpicture}[scale=0.75]
  \begin{scope}[minimum size=0, inner sep=1]
    % Vertices of string graph
    \draw (0,0) node[](v1){$ 1$};
    \draw (1,0) node[](v2){$ 1$};
    \draw (2,0) node[](v3){$ 2$};
    \draw (3,0) node[](v4){$ 2$};
    \draw (4,0) node[](v5){$ 1$};
    \draw (5,0) node[](v6){$ 1$};
    \draw (6,0) node[](v7){$ 2$};
    %
    % Arrows of string graph
    \draw[->] (v1) --node[pos=0.5,above](){$\alpha$} (v2);
    \draw[->] (v2) --node[pos=0.5,above](){$\beta$} (v3);
    \draw[->] (v4) --node[pos=0.5,above](){$\delta$} (v3);
    \draw[->] (v4) --node[pos=0.5,above](){$\gamma$} (v5);
    \draw[->] (v6) --node[pos=0.5,above](){$\alpha$} (v5);
    \draw[->] (v7) --node[pos=0.5,above](){$\gamma$} (v6);
  \end{scope}
  \end{tikzpicture}
\big)\tcomma
$
then $\Str v_2$ will be simply reported as $\big(\alpha\beta\big)\big(\delta\big)^{-1}\big(\gamma\big)\big(\gamma\alpha\big)^{-1}$. 

\setlist[enumerate,1]{label=(\roman*)}
\setlist[enumerate,2]{label=(\Alph*)}
\setlist[enumerate,3]{label=(\Roman*)}
\setlist[enumerate,4]{label=(\alph*)}

\section*{Overquiver (\subref{subfig:1-2-overquiver})}
Let $\OQ$ be the overquiver in Figure \ref{subfig:1-2-overquiver}. The vertices $1$ and $1'$ lie on distinct connected components, therefore any choice of components is of the form $C=\{\alpha^r, (\beta\gamma)^s\}$ for some $r \geq 2$ and $s \geq 1$. We also must add the generator $\delta^t$ of the complement, for some $t \geq 1$.

\begin{enumerate}
    \item Suppose there are no additional generators of $\Subpath(\OQ) \setminus N$.
    
    \begin{framed}
    \begin{tabular}[t]{p{3cm}p{3.5cm}p{6cm}}
      \textbf{Parameter values}
        & \textbf{Syzygy quiver}
          & \textbf{Legend}
            \\
      $t = 1$, $s \geq 1$, $r \geq 2$
        & \fbox{$\begin{tikzcd}[sep=small]
            E_2
          \end{tikzcd}$}
          & 
            \\ & & \\
      $t = 2$, $s \geq 1$, $r \geq 2$
        & \raisebox{-5ex}{\fbox{$\begin{tikzcd}[sep=small]
          E_2 \ar[r, shift left] \ar[r, shift right]
            \& S_1 \ar[loop above] \ar[d]
              \& Z \ar[l, shift left] \ar[l, shift right]
                \\
            \& X \ar[r]
              \& Y \ar[u]
        \end{tikzcd}$}}
          & $X \ceq \big((\beta\gamma)^{s-1}\beta\big)$\newline
          $Y \ceq \big( \alpha^{r-1} \big)$\newline  
          $Z \ceq \big( (\gamma\beta)^{s-1} \gamma \big)$
            \\ & & \\
      any other $(t,s,r)$
        & \raisebox{-4ex}{\fbox{$\begin{tikzcd}[row sep=small]
            E_2 \ar[r] \ar[d] \ar[dr]
              \& Y \ar[dl, crossing over] \ar[d, shift left]
                \\
            X \ar[dr]
              \& S_1 \ar[u, shift left] \ar[l]
                \& U \ar[l] \ar[ul]
                  \\
              \& Z \ar[ur]
            \end{tikzcd}$}}
          & $X \ceq \big((\beta\gamma)^{s-1}\beta\big)$\newline
    $Y \ceq \big(\delta\big)$\newline
    $Z \ceq \big(\alpha^{r-1}\big)$\newline
    $U \ceq \big((\gamma\beta)^{s-1}\big)$
    \\
    \end{tabular}
    \end{framed}
    
    \item\label{item:example-to-discuss} Suppose there is the additional generator $(\gamma\beta)^s$.

    \begin{framed}    
    \begin{tabular}[t]{p{3cm}p{3.5cm}p{6cm}}
    \textbf{Parameter values}
      & \textbf{Syzygy quiver}
        & \textbf{Legend}
          \\
    $t=1$, $s \geq 1$, $r \geq 2$
      & \fbox{$\begin{tikzcd}[row sep=small]
        E_2 \ar[r]
          \& X \ar[r]
          \& P_2
      \end{tikzcd}$}
        & $X \ceq \big(\alpha^{r-1}\big)$
          \\
            & & \\
    $t = 2$, $s = 1$, $r = 2$
      & \raisebox{-2ex}{\fbox{$\begin{tikzcd}[row sep=small]
      E_2 \ar[r] \ar[d]
        \& X \ar[loop right] \ar[d]
          \\
      S_1 \ar[ur]
        \& S_2 \ar[l] \ar[loop right]
      \end{tikzcd}$}}
        & $X \ceq \big( \alpha \big)\big( \gamma \big)^{-1}$
    \\ & & \\
    $t=2$,\newline any other $(s,r)$
      & \raisebox{-4ex}{\fbox{$\begin{tikzcd}[row sep=small]
      E_2 \ar[r] \ar[d]
        \& X \ar[dl]
          \\
      Y \ar[r] \ar[d, shift left]
        \& S_1 \ar[u] \ar[loop right]
          \\
      Z \ar[u, shift left] \ar[ur]
      \end{tikzcd}$}}
        & $X \ceq \big((\beta\gamma)^{s-1}\big)$,\newline
      $Y \ceq \big(\alpha^r\big) \big(\gamma\big)^{-1}$,\newline
      $Z \ceq \big((\gamma\beta)^{s-1}\gamma\big) \big( \alpha^{r-1}\big)^{-1} \big( (\beta\gamma)^{s-1} \big)$
      \\ & & \\
    $t=3$, $s=1$, $r=2$
      & \raisebox{-4ex}{\fbox{$\begin{tikzcd}[row sep=small]
      E_2 \ar[r] \ar[d]
        \& Y \ar[loop right] \ar[d]
          \\
      S_1 \ar[ur]
        \& X \ar[l] \ar[dl, shift left]
          \\
      S_2 \ar[u] \ar[ur, shift left]
      \end{tikzcd}$}}
        & $X \ceq \big( \delta^{t-2}\big)$,\newline
        $ Y \ceq \big( \alpha \big)\big(\gamma\big)^{-1}$
          \\ & & \\
    any other $(t,s,r)$
      & \raisebox{-4ex}{\fbox{$\begin{tikzcd}[row sep=small]
        E_2 \ar[r] \ar[d]
          \& Z  \ar[r, shift left] \ar[d]
            \& U \ar[l, shift left] \ar[dl, shift left]
              \\
        X \ar[ur]
          \& Y \ar[l] \ar[dl, shift left]
            \\
        S_2 \ar[u] \ar[ur, shift left]
      \end{tikzcd}$}}
        & $X \ceq \big( (\beta\gamma)^{s-1} \big)$,\newline
        $Y \ceq \big( \delta^{t-2} \big)$,\newline
        $Z \ceq \big( \alpha^{r-1} \big)\big( \gamma \big)^{-1}$,\newline
        $U \ceq \big((\gamma\beta)^{s-1}\gamma\big)\big( \alpha \big)^{-1}\big((\beta\gamma)^{s-1}\big)$
          \\
    \end{tabular}
    \end{framed}
\end{enumerate}

\section*{Overquiver (\subref{subfig:1-3-overquiver})}
Let $\OQ$ be the overquiver in Figure \ref{subfig:1-3-overquiver}. The vertices $1$ and $1'$ lie on distinct connected components, therefore any choice of components is of the form $C=\{\alpha^r, (\beta\gamma\delta)^s\}$ for some $r \geq 2$ and $s \geq 1$.
\begin{enumerate}
    \item Suppose there are no additional generators of $\Subpath(\OQ) \setminus N$.
    
    \begin{framed}
    \begin{tabular}[t]{p{3cm}p{3.5cm}p{6cm}}
    \textbf{Parameter values}
      & \textbf{Syzygy quiver}
        & \textbf{Legend}
          \\
      any $(s,r)$
        & \raisebox{-3ex}{\fbox{$\begin{tikzcd}[row sep=small]
          E_2 \ar[r, shift left] \ar[r, shift left, bend left] \ar[d]
            \& X \ar[r]
              \& U \ar[dl]
                \\
          Y \ar[ur, shift left] \ar[ur, shift right]
            \& Z \ar[l, shift left] \ar[u]
        \end{tikzcd}$}}
          & $X \ceq \big((\beta\gamma\delta)^{s-1}\beta\gamma\big)$, 
          $Y \ceq \big(\gamma\big)$,\newline
          $Z \ceq \big( (\gamma\delta\beta)^{s-1}\gamma\delta \big)$, $U \ceq \big( \alpha^{r-1} \big)$
    \end{tabular}
    \end{framed}
    
    \item Suppose there is the additional generator $(\gamma\delta\beta)^s\gamma$.
    
    \begin{framed}
    \begin{tabular}[t]{p{3cm}p{3.5cm}p{6cm}}
    \textbf{Parameter values}
      & \textbf{Syzygy quiver}
        & \textbf{Legend}
          \\
    any $(s,r)$
      & \raisebox{-3ex}{\fbox{$\begin{tikzcd}[row sep=small]
        E_2 \ar[r] \ar[d] \ar[dr, shift left]
          \& X \ar[r]
            \& U \ar[d]
              \\
        Y \ar[loop below] \ar[r, shift left]
          \& S_1 \ar[u] \ar[l, shift left]
            \& Z \ar[l] \ar[ul]
      \end{tikzcd}$}}
        & $X \ceq \big((\beta\gamma\delta)^{s-1}\beta\gamma\big)$,\newline
        $Y \ceq \big((\delta\beta\gamma)^{s-1}\delta\beta\big)$,\newline
        $Z \ceq \big((\gamma\delta\beta)^{s-1}\gamma\delta\big)$, $U \ceq \big( \alpha^{r-1} \big)$    
    \end{tabular}
    \end{framed}
    
    \item Suppose there is the additional generator $(\gamma\delta\beta)^s$.
    
    \begin{framed}
    \begin{tabular}[t]{p{1.8cm}p{4.7cm}p{6cm}}
    \textbf{Parameter values}
      & \textbf{Syzygy quiver}
        & \textbf{Legend}
          \\
      $s=1$, $r=2$
        & \raisebox{-6ex}{\fbox{$\begin{tikzcd}[row sep=small]
          E_2 \ar[d, shift left] \ar[d, shift right] \ar[dr]
            \& S_1 \ar[r]
              \& W \ar[loop above] \ar[d]
                \&
                  \\
          Y \ar[r] \ar[loop below]
            \& Z \ar[u] \ar[r]
              \& X\ar[dr]
                \&
                  \\
            \&\& U\ar[u]
                \& V \ar[l]
        \end{tikzcd}$}}
          & $X \ceq \big(\beta\gamma\big)$, $Y \ceq \big(\delta\big)$,\newline
          $Z \ceq \big(\gamma\big)$, $U \ceq \big( \gamma\delta \big)$,\newline
          $V \ceq \big( \alpha \big)$, $W \ceq \big( \alpha \big)\big( \gamma\delta \big)^{-1}$
            \\ & & \\
    any other\newline $(s,r)$
      & \raisebox{-4ex}{\fbox{$\begin{tikzcd}[row sep=small]
        E_2 \ar[d, shift left] \ar[d, shift left] \ar[dr]
          \& X \ar[r]
            \& F \ar[r, shift left] \ar[d]
              \& G \ar[l, shift left] \ar[dl]
                \\
        Z \ar[loop below] \ar[r]
          \& U \ar[r] \ar[u]
            \& Y \ar[dr]
              \&
                \\
          \& \& V \ar[u]
               \& W \ar[l]
      \end{tikzcd}$}}
        & $X \ceq \big( (\beta\gamma\delta)^{s-1} \big)$, $Y \ceq \big( (\beta\gamma\delta)^{s-1}\beta\gamma \big)$,\newline
        $Z \ceq \big( (\delta\beta\gamma)^{s-1}\delta \big)$, $U \ceq \big( \gamma \big)$,\newline
        $V \ceq \big( (\gamma\delta\beta)^{s-1}\gamma\delta \big)$, $W \ceq \big(\alpha^{r-1}\big)$,\newline
        $F \ceq \big(\alpha^{r-1}\big)\big(\gamma\delta\big)^{-1}$,\newline
        $G \ceq \big((\gamma\delta\beta)^{s-1}\gamma\delta\big)\big( \alpha \big)^{-1}\big(\beta\gamma\delta\big)^{s-1}$      
    \end{tabular}
    \end{framed}
    \item Suppose there is the additional generator $(\delta\beta\gamma)^s$.
    
    \begin{framed}
    \begin{tabular}[t]{p{1.8cm}p{4.7cm}p{6cm}}
    \textbf{Parameter values}
      & \textbf{Syzygy quiver}
        & \textbf{Legend}
          \\
      $s=1$, $r=2$
        & \raisebox{-2ex}{\fbox{$\begin{tikzcd}[row sep=small]
          E_2 \ar[r] \ar[d]
            \& Z \ar[r] \ar[d]
              \& Y \ar[loop right]
                \\
          X \ar[ur]
            \& V \ar[l] \ar[r]
              \& U \ar[u] \ar[loop right]
        \end{tikzcd}$}}
          & $X \ceq \big(\beta\big)$, $Y \ceq \big( \delta\beta \big)$,\newline
          $Z \ceq \big( \delta \big)\big(\alpha\big)^{-1}$, $U \ceq  \big(\delta\big)\big(\alpha\big)^{-1} \big( \beta \big)$,\newline
          $V \ceq \big(\beta \big) \big(\alpha\big)^{-1} \big(\gamma\delta\big)^{-1}$
    \end{tabular}
    \end{framed}
    \begin{framed}
    \begin{tabular}[t]{p{1.8cm}p{4.7cm}p{6cm}}
    \textbf{Parameter values}
      & \textbf{Syzygy quiver}
        & \textbf{Legend}
          \\
      any other\newline $(s,r)$
        & \raisebox{-5ex}{\fbox{$\begin{tikzcd}[row sep=small]
          E_2 \ar[d] \ar[r]
            \& Z \ar[r] \ar[d]
              \& Y \ar[loop above]
                \\
          X \ar[ur]
            \& W \ar[l] \ar[r]
              \& U \ar[u] \ar[r, shift left]
                \& V \ar[l, shift left] \ar[ul, shift right]
        \end{tikzcd}$}}
          & $X \ceq \big( (\beta\gamma\delta)^{s-1}\beta \big)$,\newline
          $Y \ceq \big((\delta\beta\gamma)^{s-1}\delta\beta\big)$, $Z \ceq \big(\delta\big)\big(\alpha\big)^{-1}$,\newline
          $U \ceq \big( \delta \big)\big(\alpha\big)^{-1}\big( \beta \big)$,\newline
          $V \ceq \big((\beta\gamma\delta)^{s-1}\beta\big)\big(\alpha\big)^{-1}\big((\delta\beta\gamma)^{s-1}\delta\big)$,\newline
          $W \ceq \big( (\beta\gamma\delta)^{s-1} \beta \big)\big(\alpha\big)^{-1}\big((\gamma\delta\beta)^{s-1}\gamma\delta\big)$      
    \end{tabular}
    \end{framed}
    
    \item Suppose there are the additional generators $(\gamma\delta\beta)^s$ and $(\delta\beta\gamma)^s$
    
    \begin{framed}
    \begin{tabular}[t]{p{3cm}p{3.5cm}p{6cm}}
    \textbf{Parameter values}
      & \textbf{Syzygy quiver}
        & \textbf{Legend}
          \\
      any $(s,r)$
        & \raisebox{-3ex}{\fbox{$\begin{tikzcd}[row sep=small]
          E_2 \ar[d] \ar[r]
            \& Z \ar[dl, shift right] \ar[r, shift left]
              \& U \ar[l, shift left] \ar[d]
                \\
          Y \ar[loop below] \ar[r, shift left]
            \& S_1 \ar[l, shift left] \ar[r]
              \& X \ar[ul, shift left]
        \end{tikzcd}$}}
          & $ X \ceq \big((\beta\gamma\delta)^{s-1}\beta\big)$, $Y \ceq \big( (\delta\beta\gamma)^{s-1} \delta \big)$,\newline
          $Z \ceq \big( \delta \big) \big( \alpha^{r-1} \big)^{-1} $, \newline
          $U \ceq \big((\beta\gamma\delta)^{s-1}\beta\alpha\big)\big( (\gamma\delta\beta)^{-1}\gamma\delta\big)^{-1}$
    \end{tabular}
    \end{framed}

    \item Suppose there is the additional generator $(\delta\beta\gamma)^{s-1}\delta\beta$.
    
    \begin{framed}
    \begin{tabular}[t]{p{3cm}p{3.5cm}p{6cm}}
    \textbf{Parameter values}
      & \textbf{Syzygy quiver}
        & \textbf{Legend}
          \\
      any $(s,r)$
        & \raisebox{-2ex}{\fbox{$\begin{tikzcd}[row sep=small]
          E_2 \ar[d]
            \& Y \ar[d, shift left] \ar[r, shift left]
              \& X \ar[l, shift left]
                \\
          U \ar[r, shift left] \ar[r, shift right]
            \& Z \ar[u, shift left] \ar[ur, shift right]
        \end{tikzcd}$}}
          & $ X \ceq \big((\beta\gamma\delta)^{s-1}\big)$, $Y \ceq \big(\alpha^r\big)\big( \gamma\delta \big)^{-1}$,\newline
          $ Z \ceq \big((\gamma\delta\beta)^{s-1}\gamma\delta\big)\big(\alpha\big)^{-1}\big((\beta\gamma\delta)^{s-1}\big)$,\newline
          $U \ceq \big( \alpha^{r-1} \big)\big(  \gamma\delta \big)^{-1}\big(\delta\big)(\alpha^{r-1})^{-1}$
    \end{tabular}
    \end{framed}
\end{enumerate}

\section*{Overquiver (\subref{subfig:2-2-overquiver})}\label{app-sec:fin-many-sporadic-values}

Let $\OQ$ be the overquiver in Figure \ref{subfig:2-2-overquiver}. The vertices $1$ and $1'$ lie on distinct connected components, therefore any choice of components is of the form $C=\{(\alpha\beta)^r, (\gamma\delta)^s\}$ for some $r \geq 1$ and $s \geq 1$.

By simultaneously exchanging $\ilcorresp{1}{}{1'}$, $\ilcorresp{2}{}{2'}$, $\ilcorresp{\alpha}{}{\gamma}$ and $\ilcorresp{\beta}{}{\delta}$, we may assume that $1 \leq r \leq s$.

\begin{enumerate}
    \item Suppose there are no additional generators of $\Subpath(\OQ) \setminus N$.
    
    \begin{framed}
    \begin{tabular}[t]{p{3cm}p{3.5cm}p{6cm}}
    \textbf{Parameter values}
      & \textbf{Syzygy quiver}
        & \textbf{Legend}
          \\
      any $(s,r)$
        & \raisebox{-4ex}{\fbox{$\begin{tikzcd}[row sep=small]
          E_2 \ar[d] \ar[dr] \ar[r]
            \& U \ar[r, shift left]
              \& X \ar[l, shift left] \ar[dl]
                \\
          Y \ar[d, shift left]
            \& S_1 \ar[u] \ar[l]
              \\
          Z \ar[u, shift left] \ar[ur]
        \end{tikzcd}$}}
          & $X \ceq \big((\delta\gamma)^{s-1}\delta\big)$,\newline $Y \ceq \big((\gamma\delta)^{s-1}\gamma\big)$,\newline
          $Z \ceq \big( (\beta\alpha)^{r-1}\beta \big)$,\newline $U \ceq \big( (\alpha\beta)^{r-1}\alpha \big)$
    \end{tabular}
    \end{framed}
    
    \item Suppose there is the additional generator $(\beta\alpha)^r$.
    
    \begin{framed}
    \begin{tabular}[t]{p{2.8cm}p{3.5cm}p{6.2cm}}
    \textbf{Parameter values}
      & \textbf{Syzygy quiver}
        & \textbf{Legend}
          \\
    $s=1$, $r=1$
      & \raisebox{-2ex}{\fbox{$\begin{tikzcd}[row sep=small]
      E_2 \ar[d] \ar[dr, shift left]
        \& X \ar[r, shift left]
          \& Y \ar[l, shift left]
            \\
      S_1 \ar[r, shift left]
        \& Z \ar[l, shift left] \ar[u, shift left] \ar[u, shift right]
      \end{tikzcd}$}}
      & $X \ceq \big(\gamma\big)$, $Y \ceq \big( \beta \big)$,
      $Z \ceq \big(\alpha\big)\big(\delta\gamma\big)^{-1}$
        \\
    \end{tabular}
    \end{framed}
    
    \begin{framed}
    \begin{tabular}[t]{p{2.8cm}p{3.5cm}p{6.2cm}}
    \textbf{Parameter values}
      & \textbf{Syzygy quiver}
        & \textbf{Legend}
          \\
    $s=2$, $r=1$
      & \raisebox{-4ex}{\fbox{$\begin{tikzcd}[row sep=small]
        E_2 \ar[d] \ar[r, shift left]
          \& U \ar[r] \ar[d] \ar[dl, shift left]
            \& X \ar[d]
              \\
        S_1 \ar[ur, shift left]
          \&Y \ar[d, shift left]
            \& V \ar[l] \ar[loop below]
              \\
          \& Z \ar[u, shift left]
      \end{tikzcd}$}}
        & $X \ceq \big(\gamma\big)$, $Y \ceq \big(\gamma\delta\gamma\big)$, $Z \ceq \big(\beta\big)$,\newline
        $U \ceq \big( \delta\gamma\delta \big)\big(\beta\big)^{-1}$, $V \ceq \big(\gamma\delta\big)\big(\beta\big)^{-1}$
    \\  &   &   \\
    any other $(s,r)$
      & \raisebox{-4ex}{\fbox{$\begin{tikzcd}[row sep=small]
        E_2 \ar[r, shift left] \ar[d]
          \& U \ar[r] \ar[d] \ar[dl, shift left]
            \& W \ar[d]
              \\
        Z \ar[ur, shift left]
          \& X \ar[d, shift left]
            \& V \ar[l] \ar[d, shift left]
              \\
          \& Y \ar[u, shift left]
            \& F \ar[ul] \ar[u, shift left]
      \end{tikzcd}$}}
       & $X \ceq \big( (\gamma\delta)^{s-1}\gamma \big)$, $Y \ceq \big( (\beta\alpha)^{r-1}\beta \big)$,\newline
       $Z \ceq \big((\alpha\beta)^{r-1}\big)$,\newline
       $U \ceq \big( (\delta\gamma)^{s-1}\delta \big)\big(\beta\big)^{-1}$,\newline
       $V \ceq \big((\gamma\delta)^{r=1}\big)\big(\beta\big)^{-1} $, 
       $W \ceq \big( (\alpha\beta)^{r-1} \big)$,\newline
       $F \ceq \big( (\beta\alpha)^{r-1}\beta \big)\big(\gamma\delta\big)^{-1}\big((\alpha\beta)^{r-1}\big)$
    \end{tabular}
    \end{framed}
    
    \item \label{item:fin-many-sporadic-values} Suppose there is the additional generator $(\delta\gamma)^s$
    \begin{framed}
    \begin{tabular}[t]{p{3cm}p{3.5cm}p{6cm}}
    \textbf{Parameter values}
      & \textbf{Syzygy quiver}
        & \textbf{Legend}
          \\
      $s=1$, $r=1$
        & \raisebox{-2ex}{\fbox{$\begin{tikzcd}[row sep=small]
          E_2 \ar[d] \ar[dr, shift left]
            \& Y \ar[r, shift left] \ar[d, shift left]
              \& X \ar[l, shift left]
                \\
          S_1 \ar[r, shift left]
            \& Z \ar[l, shift left] \ar[u, shift left]
        \end{tikzcd}$}}
          & $X \ceq \big(\delta\big)$, $Y \ceq \big( \alpha \big)$, $Z \ceq \big( \beta\alpha \big)\big( \gamma \big)^{-1}$
    \\  &   &   \\
    $s=2$, $r=1$
      & \raisebox{-4ex}{\fbox{$\begin{tikzcd}[row sep=small]
        E_2 \ar[r, shift left] \ar[d]
          \& V \ar[r] \ar[dl, shift left] \ar[d]
            \& W \ar[d]
              \\
        Z \ar[ur, shift left]
          \& U \ar[d, shift left]
            \& X \ar[loop below] \ar[l]
              \\
          \& Y \ar[u, shift left]
      \end{tikzcd}$}}
        & $X \ceq \big( \delta \big)$, $Y \ceq \big( \delta\gamma\delta \big)$, $Z \ceq \big( \gamma\delta \big)$,\newline
        $U \ceq \big( \alpha \big)$, $V \ceq \big( \delta \big)\big(\beta\big)^{-1}$, \newline
        $W \ceq \big( \gamma\delta \big)^{-1}\big(\alpha\big)$
            \\  & & \\
    other $(s,r)$
      & \raisebox{-4ex}{\fbox{$\begin{tikzcd}[row sep=small]
        E_2 \ar[r, shift left] \ar[d]
          \& V \ar[r] \ar[d] \ar[dl, shift left]
            \& W \ar[d]
              \\
        Y \ar[ur, shift left]
          \& Z \ar[d, shift left]
            \& U \ar[l] \ar[d, shift left]
              \\
          \& X \ar[u, shift left]
            \& F \ar[ul] \ar[u, shift left]
      \end{tikzcd}$}}
        & $X \ceq \big((\delta\gamma)^{s-1}\delta\big)$, $Y \ceq \big( (\gamma\delta)^{s-1} \big)$,\newline
        $Z \ceq \big( (\alpha\beta)^{r-1}\alpha \big)$,\newline
        $U \ceq \big( \delta \big)\big((\alpha\beta)^{r-1}\big)^{-1}$,\newline
        $V \ceq \big( \delta \big)\big( (\beta\alpha)^{r-1}\beta \big)^{-1}$,\newline $W \ceq \big((\gamma\delta)^{r-1}\big)$,\newline
        $F \ceq \big( (\gamma\delta)^{r-1} \big)\big( \alpha\beta \big)^{-1} \big( (\delta\gamma)^{s-1} \delta \big)$
    \end{tabular}
    \end{framed}
    
    \item Suppose there exist additional generators $(\beta\alpha)^r$ and $(\gamma\delta)^s$.
    
    \begin{framed}
    \begin{tabular}[t]{p{1.8cm}p{5cm}p{5.7cm}}
    \textbf{Parameter values}
      & \textbf{Syzygy quiver}
        & \textbf{Legend}
          \\
      $s=1$, $r=1$
        & \raisebox{-2ex}{\fbox{$\begin{tikzcd}[row sep=small]
          E_2 \ar[r]
            \& Y \ar[r, bend left, shift left]\ar[r, shift left=1] \ar[r, shift right=1] \ar[r, bend right, shift right]
              \& S_1 \ar[r, shift left, bend left]
                \& X \ar[l, shift right] \ar[l, shift left] \ar[l, shift left, bend left]
        \end{tikzcd}$}}
          & $X \ceq \big( \beta \big)\big( \delta \big)^{-1}$,\newline
          $Y \ceq \big(\delta\big)\big(\beta\big)^{-1}\big(\delta\big)\big(\beta\big)^{-1}$
        \\  &   &   \\
      other $(s,r)$
        & \raisebox{-4ex}{\fbox{$\begin{tikzcd}[row sep=small]
          E_2 \ar[r, shift left]
            \& V \ar[dr, shift left, bend left] \ar[dr, shift left] \ar[dl] \ar[dl, bend right]
              \\
          X \ar[r, shift left]
            \& Z \ar[r, shift left] \ar[l, shift left] \ar[l, shift left, bend left]
              \& Y \ar[dl, bend left]
                \\
            \& U \ar[ul, bend left] \ar[ur] \ar[ur, bend left]
        \end{tikzcd}$}}
          & $X \ceq \big( (\gamma\delta)^{s-1} \big)$, $Y \ceq \big( (\alpha\beta)^{r-1} \big)$,\newline
          $ Z \ceq \big(\delta\big)\big( (\beta\alpha)^{r-1}\beta \big)^{-1} $,\newline
          $U \ceq \big( (\delta\gamma)^{s-1}\delta \big)\big( \beta \big)^{-1}$,\newline
          $V \ceq \big((\delta\gamma)^{s-1}\delta\big)\big(\beta \big)^{-1} \big( \delta \big) \big((\beta\alpha)^{r-1}\beta\big)^{-1} $
    \end{tabular}
    \end{framed}    
\end{enumerate}

\section*{Overquiver (\subref{subfig:3-3-overquiver})}

Let $\OQ$ be the overquiver in Figure \ref{subfig:3-3-overquiver}. The vertices $1$ and $1'$ lie on the same connected component. The shortest paths between them are $\big(\ilsltpath{1}{\alpha}{1'}\big)$ and $\big(\ilsltpath{1'}{\beta\gamma}{1}\big)$. The components can either be factorised into an even number of these paths or an odd number. In the even case, the component set is of the form $C=\{(\alpha\beta\gamma)^r, (\gamma\alpha\beta)^r\}$ for some $r \geq 1$. In the odd case, it has the form $C =\{ (\alpha\beta\gamma)^r\alpha, (\beta\gamma\alpha)^r\beta\gamma \}$ for some $r \geq 1$. In either case, we must also add a generator $\delta^s$ of the complement, for some $s \geq 1$.

\begin{enumerate}
    \item In the even case suppose there are no additional generators.
    
    \begin{framed}
    \begin{tabular}[t]{p{3cm}p{3.5cm}p{6cm}}
    \textbf{Parameter values}
      & \textbf{Syzygy quiver}
        & \textbf{Legend}
          \\
      $s = 1$, $r \geq 1$
        & \raisebox{-0ex}{\fbox{$\begin{tikzcd}[row sep=small]
          E_2
        \end{tikzcd}$}}
          & \\  & & \\
      $s=2$, $r \geq 2$
        & \raisebox{-2ex}{\fbox{$\begin{tikzcd}[row sep=small]
          E_2 \ar[d] \ar[r]
            \& S_2 \ar[dl]
              \\
          X \ar[r]
            \& Y \ar[u]
        \end{tikzcd}$}}
          & $X \ceq \big((\alpha\beta\gamma)^{r-1}\alpha\beta\big)$,\newline
          $Y \ceq \big( (\gamma\alpha\beta)^{r-1}\gamma\alpha \big)$
            \\ & & \\
    any other $(s,r)$
      & \raisebox{-4ex}{\fbox{$\begin{tikzcd}[row sep=small]
          E_2 \ar[rr] \ar[dr, shift left] \ar[dd]
           \&
             \& X \ar[dd] \ar[from=ddll, bend right]
               \\
            \& Y \ar[dl, shift left] \ar[ur] \ar[from=dr, crossing over]
              \\
          S_2 \ar[ur, shift left]
           \&
             \& Z \ar[ll, shift left]
        \end{tikzcd}$}}
          & $X \ceq \big( (\alpha\beta\gamma)^{r-1}\alpha\beta\big)$,\newline
          $ Y \ceq \big( \delta^{s-2} \big) $,\newline
          $Z \ceq \big( (\gamma\alpha\beta)^{r-1}\gamma\alpha \big)$
    \end{tabular}
    \end{framed}
    
    \item In the even case, suppose there is the additional generator $(\gamma\alpha\beta)^r$.
    
    \begin{framed}
    \begin{tabular}[t]{p{3cm}p{3.5cm}p{6cm}}
    \textbf{Parameter values}
      & \textbf{Syzygy quiver}
        & \textbf{Legend}
          \\
      $s = 1$, $r \geq 1$
        & \raisebox{-0ex}{\fbox{$\begin{tikzcd}[row sep=small]
          E_2 \ar[r]
            \& X
        \end{tikzcd}$}}
          & $X \ceq \big( (\beta\gamma\alpha)^{r-1}\beta\gamma \big)$
        \\  &   &   \\
    $s = 2$, $r \geq 2$
      & \raisebox{-2ex}{\fbox{$\begin{tikzcd}[row sep=small]
      E_2 \ar[r, shift left] \ar[r, shift right] \ar[d]
        \& S_2 \ar[dl] \ar[loop right]
          \\
      X \ar[r]
        \& Y \ar[u, shift left] \ar[u, shift right]
    \end{tikzcd}$}}
      & $X \ceq \big( (\gamma\alpha\beta)^{r-1}\gamma\alpha \big)$,\newline
      $Y \ceq \big( (\beta\gamma\alpha)^{r-1}\beta\gamma \big)$
        \\ & & \\
    any other $(s,r)$
      & \raisebox{-4ex}{\fbox{$\begin{tikzcd}[row sep=small]
      E_2 \ar[dd] \ar[rr, shift left]
        \&
          \& X \ar[dl]  \ar[from=ddll, bend left] \ar[ddll, bend left]
            \\
        \& Y \ar[dr, crossing over]  \ar[from=ul, crossing over]
          \\
      S_2 \ar[ur]
        \&
          \& Z \ar[ll, shift left] \ar[uu]
    \end{tikzcd}$}}
      & $X \ceq \big( \delta^{s-2} \big)$,\newline
      $Y \ceq \big( (\gamma\alpha\beta)^{r-1}\gamma\alpha \big)$,\newline
      $Z \ceq \big( (\beta\gamma\alpha)^{r-1}\beta\gamma \big)$
    \end{tabular}
    \end{framed}
    
    \item In the odd case, suppose there are no additional generators.
    
    \begin{framed}
    \begin{tabular}[t]{p{3cm}p{3.5cm}p{6cm}}
    \textbf{Parameter values}
      & \textbf{Syzygy quiver}
        & \textbf{Legend}
          \\
      $r=1$, $s=1$
        & \raisebox{-2ex}{\fbox{$\begin{tikzcd}[row sep=small]
          E_2 \ar[d]
           \\
           X \ar[r, shift left] \ar[r, shift right]
              \& Y \ar[loop above]
        \end{tikzcd}$}}
          & $X \ceq \big( \gamma \big)$, $Y \ceq \big(\beta\gamma\big)$
            \\  & & \\
        $r=1$, $s=2$
        & \raisebox{-2ex}{\fbox{$\begin{tikzcd}[row sep=small]
          E_2 \ar[d] \ar[r]
            \& S_2 \ar[dl]
              \\
          X \ar[r]
            \& Y \ar[u] \ar[from=ul, crossing over]
        \end{tikzcd}$}}
          & $X \ceq \big( (\gamma\alpha\beta)\gamma \big)$,\newline
          $Y \ceq \big(  \beta\gamma \big)$
    \\
    \end{tabular}
    \end{framed}
    \begin{framed}
    \begin{tabular}[t]{p{3cm}p{3.5cm}p{6cm}}
    \textbf{Parameter values}
      & \textbf{Syzygy quiver}
        & \textbf{Legend}
          \\
    $r=1$, $s \geq 3$
      & \raisebox{-4ex}{\fbox{$\begin{tikzcd}[row sep=small]
      E_2 \ar[d] \ar[r]
        \& 3 \ar[d] \ar[from=ddl]
          \\
      X \ar[d, shift left] \ar[ur]
        \& 4 \ar[l, crossing over] \ar[from=ul, crossing over] \ar[loop right]
          \\
      S_2 \ar[u, shift left] 
      \end{tikzcd}$}}
        & $X \ceq \big( \delta^{s-2} \big)$,\newline
        $Y \ceq \big( \gamma\alpha\beta\gamma \big)$,\newline
        $Z \ceq \big( \beta\gamma \big)$
            \\  & & \\
      $r \geq 2$, $s = 1$
        & \raisebox{-2ex}{\fbox{$\begin{tikzcd}[row sep=small]
          E_2 \ar[d]
            \& Y \ar[d, shift left]
              \\
          X \ar[r, shift left] \ar[r, shift right]
              \& Z \ar[u, shift left]
        \end{tikzcd}$}}
          & $X \ceq \big(\gamma\big)$, $Y \ceq \big( \beta\gamma \big)$,\newline
          $Z \ceq \big( (\beta\gamma\alpha)^{r-1}\beta\gamma \big)$
            \\ & & \\
            $r \geq 2$, $s = 2$
        & \raisebox{-4ex}{\fbox{$\begin{tikzcd}[row sep=small]
          E_2 \ar[rr] \ar[dr] \ar[dd]
            \&
              \& X \ar[dd]
                \\
            \& S_2 \ar[ur] \ar[loop above]
              \\
          Y \ar[ur] \ar[rr, shift left]
            \&
              \& Z \ar[ul] \ar[ll, shift left]
        \end{tikzcd}$}}
          & $X \ceq \big( (\gamma\alpha\beta)^{r}\gamma \big)$, $Y \ceq \big( \beta\gamma \big)$,\newline
          $Z \ceq \big( (\beta\gamma\alpha)^{r-1}\beta\gamma \big)$
            \\ & & \\
      any other $(r,s)$
        & \raisebox{-6ex}{\fbox{$\begin{tikzcd}[row sep=small]
          E_2\ar[rr] \ar[dr] \ar[dd]
            \&
              \& Z \ar[dl] \ar[dd, shift left]
                \\
            \& X \ar[dl] \ar[dd, shift left] \ar[from=dd, shift left]
              \\
          Y \ar[rr, crossing over]
            \&
              \& U \ar[ul] \ar[uu, shift left]
                \\
            \& S_2 \ar[ul]
        \end{tikzcd}$}}
          & $X \ceq \big( \delta^{s-1} \big)$,\newline
          $Y \ceq \big( (\gamma\alpha\beta)^s\gamma \big)$,\newline
          $ Z \ceq \big( \beta\gamma \big) $,\newline
          $U \ceq \big( (\beta\gamma\alpha)^{r-1}\beta\gamma \big)$
    \end{tabular}
    \end{framed}
\end{enumerate}

\section*{Overquiver (\subref{subfig:untwisted-4-4-overquiver})}
Let $\OQ$ be the overquiver in Figure \ref{subfig:untwisted-4-4-overquiver}. The vertices $1$ and $1'$ lie on the same connected component. The shortest paths between them are $\big(\ilsltpath{1}{\alpha\beta\gamma}{1'}\big)$ and $\big(\ilsltpath{1'}{\delta}{1}\big)$. The components can either be factorised into an even number of these paths or an odd number. In the even case, the component set is of the form $C=\{(\alpha\beta\gamma\delta)^r, (\delta\alpha\beta\gamma)^r\}$ for some $r \geq 1$. In the odd case, it has the form $C =\{ (\alpha\beta\gamma\delta)^r\alpha\beta\gamma, (\delta\alpha\beta\gamma)^r\delta \}$ for some $r \geq 0$.

\begin{enumerate}
    \item In the even case, suppose there are no additional generators.
    
    \begin{framed}
    \begin{tabular}[t]{p{3cm}p{3.5cm}p{6cm}}
    \textbf{Parameter values}
      & \textbf{Syzygy quiver}
        & \textbf{Legend}
          \\
      any $r$
        & \raisebox{-2ex}{\fbox{$\begin{tikzcd}[row sep=small]
          E_2 \ar[r, shift left] \ar[r, shift right] \ar[d]
            \& Z \ar[d, shift left]
              \\
          X \ar[ur, shift left] \ar[ur, shift right]
            \& Y \ar[l] \ar[u, shift left]
        \end{tikzcd}$}}
          & $ X \ceq \big( \beta \big)$, $Y \ceq \big( (\beta\gamma\delta\alpha)^{r-1}\beta\gamma\delta\big)$,\newline
          $Z \ceq \big( (\delta\alpha\beta\gamma)^{r-1} \delta\alpha\beta \big)$
    \end{tabular}
    \end{framed}
    
    \item In the even case, suppose there is the additional generator $(\beta\gamma\delta\alpha)^r\beta$.
    
    \begin{framed}
    \begin{tabular}[t]{p{3cm}p{3.5cm}p{6cm}}
    \textbf{Parameter values}
      & \textbf{Syzygy quiver}
        & \textbf{Legend}
          \\
      any $r$
        & \raisebox{-2ex}{\fbox{$\begin{tikzcd}[row sep=small]
          E_2 \ar[r] \ar[dr] \ar[d]
            \& Z \ar[dr, shift left]
              \\
          X \ar[r, shift left] \ar[loop below]
            \& S_2 \ar[l, shift left] \ar[u]
              \& Y \ar[l, shift left] \ar[ul, shift left]
        \end{tikzcd}$}}
          & $X \ceq \big( (\gamma\delta\alpha\beta)^{r-1}\gamma\delta\alpha \big)$,\newline
          $ Y \ceq \big( (\beta\gamma\delta\alpha)^{r-1}\beta\gamma\delta \big) $,\newline
          $ Z \ceq \big( (\delta\alpha\beta\gamma)^{r-1} \delta\alpha\beta \big) $
    \end{tabular}
    \end{framed}
    
    \item In the even case, suppose there is the additional generator $(\beta\gamma\delta\alpha)^r$.
    
    \begin{framed}
    \begin{tabular}[t]{p{3cm}p{3.5cm}p{6cm}}
    \textbf{Parameter values}
      & \textbf{Syzygy quiver}
        & \textbf{Legend}
          \\
      any $r$
        & \raisebox{-4ex}{\fbox{$\begin{tikzcd}[row sep=small]
          E_2 \ar[r] \ar[d, shift left] \ar[d, shift right]
            \& Y \ar[r] \ar[d]
              \& U \ar[d, shift left]
                \\
          X \ar[loop below] \ar[ur]
            \& V \ar[d, shift left]
              \& W \ar[u, shift left] \ar[l, shift left] \ar[l, shift left]
                \\
            \& Z \ar[u, shift left]
        \end{tikzcd}$}}
          & $X \ceq \big((\gamma\delta\alpha\beta)^{r-1}\gamma\delta\big)$, $Y \ceq \big( \beta \big)$,\newline
          $ Z \ceq \big( (\beta\gamma\delta\alpha)^{r-1}\beta\gamma\delta \big)$,\newline
          $ U \ceq \big( (\delta\alpha\beta\gamma)^{r-1} \delta \big)$,\newline
          $ V \ceq \big( (\delta\alpha\beta\gamma)^{r-1}\delta\alpha\beta \big)$,\newline
          $ W \ceq \big( (\beta\gamma\delta\alpha)^{r-1}\gamma\delta\alpha \big)\big( \beta\gamma \big)^{-1}$
            \\
    \end{tabular}
    \end{framed}
    
    \item In the even case, suppose that there is the additional generator $(\gamma\delta\alpha\beta)^r$.
    
    \begin{framed}
    \begin{tabular}[t]{p{3cm}p{3.5cm}p{6cm}}
    \textbf{Parameter values}
      & \textbf{Syzygy quiver}
        & \textbf{Legend}
          \\
      any $r$
        & \raisebox{-3ex}{\fbox{$\begin{tikzcd}[row sep=small]
          E_2 \ar[r, shift left] \ar[d]
            \& U \ar[dl, shift left] \ar[d] \ar[r]
              \& V \ar[d, shift left]
                \\
          Y \ar[ur, shift left]
            \& X \ar[loop below]
              \& Z \ar[l, shift left] \ar[l, shift right] \ar[u, shift left]
        \end{tikzcd}$}}
          & $X \ceq \big( (\gamma\delta\alpha\beta)^{r-1} \gamma\delta\alpha \big)$,\newline
          $ Y \ceq \big( (\delta\alpha\beta\gamma)^{r-1} \delta\alpha \big)$,\newline
          $Z \ceq \big(\gamma\big)\big( (\gamma\delta\alpha\beta)^{r-1}\gamma\delta \big)^{-1}$,\newline
          $U \ceq \big(\gamma\big)\big((\beta\gamma\delta\alpha)^{r-1}\beta\gamma\delta\big)^{-1}$,\newline
          $V \ceq \big((\delta\alpha\beta\gamma)^{r-1}\big)$
            \\
    \end{tabular}
    \end{framed}
    
    \item In the even case, suppose that there are the additional generators $(\beta\gamma\delta\alpha)^r$ and $(\gamma\delta\alpha\beta)^r$.
    
    \begin{framed}
    \begin{tabular}[t]{p{3cm}p{3.5cm}p{6cm}}
    \textbf{Parameter values}
      & \textbf{Syzygy quiver}
        & \textbf{Legend}
          \\
      any $r \geq 0$
        & \raisebox{-2ex}{\fbox{$\begin{tikzcd}[row sep=small]
          E_2 \ar[r] \ar[dr, shift left]
            \& Z \ar[d] \ar[r, shift left] \ar[r, shift right]
              \& Y \ar[l, shift left, bend left]
                \\
            \& X \ar[loop left] \ar[r, shift left]
              \& S_2 \ar[l] \ar[u]
        \end{tikzcd}$}}
          & $X \ceq \big( (\gamma\delta\alpha\beta)^{r-1}\gamma\delta \big)$,\newline
          $Y \ceq \big( (\delta\alpha\beta\gamma)^{r-1} \delta\alpha \big)$,\newline
          $Z \ceq \big( (\delta\gamma\delta\alpha)^{r-1}\beta\gamma\delta \big)\big(\gamma\big)^{-1}$
            \\
    \end{tabular}
    \end{framed}
    
    \item In the even case, suppose that there is the additional generator $(\gamma\delta\alpha\beta)^{r-1}\gamma\delta\alpha$.
    
    \begin{framed}
    \begin{tabular}[t]{p{1.4cm}p{3.2cm}p{7.9cm}}
    \textbf{Parameter values}
      & \textbf{Syzygy quiver}
        & \textbf{Legend}
          \\
      any $r$
        & \raisebox{-2ex}{\fbox{$\begin{tikzcd}[row sep=small]
          E_2 \ar[d] \\
            Z \ar[r, bend left, shift left] \ar[r, shift left] \ar[r, shift right] \ar[r, shift right, bend right]
              \& X \ar[r, shift left, bend left]
               \& Y \ar[l, shift right] \ar[l, shift left] \ar[l, shift left, bend left]
        \end{tikzcd}$}}
          & $X \ceq \big( (\delta\alpha\beta\gamma)^{r-1}\delta \big)$,\newline
          $ Y \ceq \big((\beta\gamma\delta\alpha)^{r-1}\beta\gamma\delta\big)\big( \beta\gamma \big)^{-1} $,\newline
          $Z \ceq \big((\beta\gamma\delta\alpha)^{r-1}\beta\gamma\delta\big)\big(\beta\gamma\big)^{-1}\big(\gamma\big)\big((\beta\gamma\delta\alpha)^{r-1}\beta\gamma\delta \big)^{-1} $
    \end{tabular}
    \end{framed}
    
    \item In the odd case, suppose that there are no additional generators.
    
    \begin{framed}
    \begin{tabular}[t]{p{3cm}p{2.5cm}p{7cm}}
    \textbf{Parameter values}
      & \textbf{Syzygy quiver}
        & \textbf{Legend}
          \\
      $r \geq 1$
        & \raisebox{-2ex}{\fbox{$\begin{tikzcd}[row sep=small]
          E_2 \ar[d] \\
          Y \ar[r, shift left] \ar[r, shift right]
            \& X \ar[loop above]
        \end{tikzcd}$}}
          & $ X \ceq \big( (\alpha\beta\gamma\delta)^r \big) \big( \alpha\beta\gamma \big)^{-1} \big((\delta\alpha\beta\gamma)^r\big) $,\newline
          $Y \ceq \big((\alpha\beta\gamma\delta)^r\big)\big( \beta\gamma \big)^{-1} \big(\gamma\big) \big( (\alpha\beta\gamma\delta)^r \big)^{-1}$
    \end{tabular}
    \end{framed}
\end{enumerate}

\section*{Overquiver (\subref{subfig:twisted-4-4-overquiver})}

Let $\OQ$ be the overquiver in Figure \ref{subfig:twisted-4-4-overquiver}. The vertices $1$ and $1'$ lie on the same connected component. The shortest paths between them are $\big(\ilsltpath{1}{\alpha\beta}{1'}\big)$ and $\big(\ilsltpath{1'}{\gamma\delta}{1}\big)$. The components can either be factorised into an even number of these paths or an odd number. In the even case, the component set is of the form $C=\{(\alpha\beta\gamma\delta)^r, (\gamma\delta\alpha\beta)^r\}$ for some $r \geq 1$. In the odd case, it has the form $C =\{ (\alpha\beta\gamma\delta)^r\alpha\beta, (\gamma\delta\alpha\beta)^r\gamma\delta \}$ for some $r \geq 0$.

Simultaneously exchanging $\ilcorresp{1}{}{1'}$, $\ilcorresp{2}{}{2'}$, $\ilcorresp{\alpha}{}{\gamma}$, $\ilcorresp{\beta}{}{\delta}$ yields an automorphism of $\OQ$ which fixes $C$ and which respects $\pingraph$. Therefore, we need only consider choices of $N$ up to this symmetry.

\begin{enumerate}
    \item In the even case, suppose that there are no additional generators.
    
    \begin{framed}
    \begin{tabular}[t]{p{3cm}p{3.5cm}p{6cm}}
    \textbf{Parameter values}
      & \textbf{Syzygy quiver}
        & \textbf{Legend}
          \\
      any $r$
        & \raisebox{-5ex}{\fbox{$\begin{tikzcd}[row sep=small]
          E_2 \ar[d] \ar[dr, shift left] \ar[drr, shift left=2]
            \\
          Y \ar[d]
            \& S_2 \ar[l, shift right] \ar[r, shift left]
              \& U \ar[dl, shift left]
                \\
          X \ar[ur] \ar[urr]
            \& Z \ar[ul, crossing over] \ar[u, crossing over]
        \end{tikzcd}$}}
          & $ X \ceq \big( (\delta\alpha\beta\gamma)^{r-1}\delta\alpha\beta \big)$,\newline
          $Y \ceq \big( (\alpha\beta\gamma\delta)^{r-1} \alpha\beta\gamma \big)$,\newline
          $ Z \ceq \big( (\beta\gamma\delta\alpha)^{r-1}\beta\gamma\delta \big)$,\newline
          $U \ceq \big( (\gamma\delta\alpha\beta)^{r-1}\gamma\delta\alpha \big)\big((\delta\alpha\beta\gamma)^r\big)^{-1}$
            \\
    \end{tabular}
    \end{framed}
    
    \item In the even case, suppose that there is the additional generator $(\beta\gamma\delta\alpha)^r$.
    
    \begin{framed}
    \begin{tabular}[t]{p{3cm}p{3.5cm}p{6cm}}
    \textbf{Parameter values}
      & \textbf{Syzygy quiver}
        & \textbf{Legend}
          \\
      any $r$
        & \raisebox{-4ex}{\fbox{$\begin{tikzcd}[row sep=small]
          E_2 \ar[r, shift left] \ar[d]
            \& Z \ar[dl, shift left]
              \\
          U \ar[ur, shift left] \ar[dr, shift left] \ar[dr, shift right]
            \& S_2 \ar[u] \ar[d]
              \& X \ar[l] \ar[ul]
                \\
          \& Y \ar[ur]
        \end{tikzcd}$}}
          & $X \ceq \big( (\delta\alpha\beta\gamma)^{r-1}\delta\alpha\beta \big)$,\newline
          $ Y \ceq \big( (\alpha\beta\gamma\delta)^{r-1}\alpha\beta\gamma \big) $,\newline
          $ Z \ceq \big((\gamma\delta\alpha\beta)^{r-1} \gamma\delta \big) $,\newline
          $U \ceq \big( (\beta\gamma\delta\alpha)^{r-1}\beta\gamma\delta \big)\big(\beta\big)^{-1}$
            \\
        \end{tabular}
    \end{framed}
    
    \item In the even case, suppose there are the additional generators $(\beta\gamma\delta\alpha)^r$ and $(\delta\alpha\beta\gamma)^r$.
    
    \begin{framed}
    \begin{tabular}[t]{p{1.5cm}p{3.5cm}p{7.5cm}}
    \textbf{Parameter values}
      & \textbf{Syzygy quiver}
        & \textbf{Legend}
          \\
      any $r$
        & \raisebox{-4ex}{\fbox{$\begin{tikzcd}[row sep=small]
          E_2 \ar[r, shift left]
            \& V \ar[dr, shift left] \ar[dr, shift right] \ar[dl, shift left] \ar[dl, shift right]
              \\
          X \ar[r, shift left] \ar[r, shift right]
            \& Z \ar[r, shift left] \ar[r, shift right]
              \& Y \ar[dl, shift left, bend left]
                \\
            \& U \ar[ur, shift right] \ar[ul, shift left] \ar[ul, shift left, bend left]
        \end{tikzcd}$}}
          & $ X \ceq \big( (\alpha\beta\gamma\delta)^{r-1} \alpha\beta \big)$,\newline
          $ Y \ceq \big((\gamma\delta\alpha\beta)^{r-1}\gamma\delta\big) $,\newline
          $ Z \ceq \big((\delta\alpha\beta\gamma)^{r-1} \delta\alpha\beta \big)\big(\delta\big)^{-1}$,\newline
          $U \ceq \big( (\beta\gamma\delta\alpha)^{r-1}\beta\gamma\delta\big) \big(\beta\big)^{-1} \big)$,\newline
          $V \ceq \big( (\beta\gamma\delta\alpha)^{r-1}\beta\gamma\delta \big)\big(\beta\big)^{-1} \big(\delta\big)\big((\delta\alpha\beta\gamma)^{r-1} \delta\alpha\beta\big)$
    \end{tabular}
    \end{framed}
    
    \item In the odd case, suppose there are no additional generators.
    \begin{framed}
    \begin{tabular}[t]{p{3cm}p{3.5cm}p{6cm}}
    \textbf{Parameter values}
      & \textbf{Syzygy quiver}
        & \textbf{Legend}
          \\
      any $r$
        & \raisebox{-5ex}{\fbox{$\begin{tikzcd}[row sep=small]
          E_2 \ar[d] \ar[dr] \ar[drr, shift left]
            \\
          Y \ar[drr, bend right]
            \& S_2 \ar[l] \ar[r]
              \& U \ar[dll, bend left, crossing over]
                \\
          Z \ar[u] \ar[ur, crossing over]
            \&\& X\ar[u] \ar[ul, crossing over]
        \end{tikzcd}$}}
          & $X \ceq \big((\delta\alpha\beta\gamma)^r \delta\big)$,\newline
          $Y \ceq \big((\alpha\beta\gamma\delta)^r \alpha\big)$,\newline
          $Z \ceq \big( (\beta\gamma\delta\alpha)^r\beta\big)$,\newline
          $U \ceq \big((\gamma\delta\alpha\beta)^r\gamma\big)$
    \end{tabular}
    \end{framed}
    
    \item In the odd case, suppose there is the additional generator $(\beta\gamma\delta\alpha)^r\beta\gamma$.
    
    \begin{framed}
    \begin{tabular}[t]{p{3cm}p{3.5cm}p{6cm}}
    \textbf{Parameter values}
      & \textbf{Syzygy quiver}
        & \textbf{Legend}
          \\
      $r \geq 0$
        & \raisebox{-4ex}{\fbox{$\begin{tikzcd}[row sep=small]
          E_2 \ar[d] \ar[r]
            \& U \ar[dl, shift left] \ar[dr, shift left] \ar[dr, shift right]
              \\
          Z \ar[ur, shift left]
            \& S_1 \ar[l] \ar[r]
              \& Y \ar[dl]
                \\
            \& X \ar[ul] \ar[u]
        \end{tikzcd}$}}
          & $X \ceq \big( (\delta\alpha\beta\gamma)^r \delta \big)$,
          $ Y \ceq \big( (\alpha\beta\gamma\delta)^r \alpha \big) $,\newline
          $ Z \ceq \big( (\gamma\delta\alpha\beta)^{r-1} \big) $,\newline
          $ U \ceq \big( (\beta\gamma\delta\alpha)^{r-1}\beta \big)\big( \delta \big)^{-1} $
    \end{tabular}
    \end{framed}

    \newpage    
    \item In the odd case, suppose there are the additional generators $(\beta\gamma\delta\alpha)^r\beta\gamma$ and $(\delta\alpha\beta\gamma)^r\delta\alpha$.
    
    \begin{framed}
    \begin{tabular}[t]{p{2.7cm}p{3.5cm}p{6.3cm}}
    \textbf{Parameter values}
      & \textbf{Syzygy quiver}
        & \textbf{Legend}
          \\
      $r=0$
        & \raisebox{-2ex}{\fbox{$\begin{tikzcd}[row sep=small]
          E_2 \ar[d]
            \\
          Y \ar[r, bend left, shift left] \ar[r, shift left] \ar[r, shift right] \ar[r, shift right, bend right]
            \& S_1 \ar[r, shift left=2]
              \& X \ar[l] \ar[l, shift left = 2]
        \end{tikzcd}$}}
          & $X \ceq (\beta)(\delta)^{-1}$,\newline
          $ Y \ceq \big( \delta \big)\big(\beta\big)^{-1}\big(\delta\big)\big(\beta\big)^{-1}$
            \\ & & \\
      any other $r$
        & \raisebox{-3ex}{\fbox{$\begin{tikzcd}[row sep=small]
          E_2 \ar[r, shift left = 1]
            \& Y \ar[r, shift left=2] \ar[from=dr, shift left] \ar[from=dr, shift right]
              \& Z \ar[l]
                \\
          V \ar[r, shift right=1] \ar[r, shift right=3] \ar[ur, shift left] \ar[ur, shift right]
            \& X \ar[from=ur, shift left,crossing over] \ar[from=ur, shift right, crossing over] \ar[r, shift right=1]
              \& U \ar[l, shift left = 3]
        \end{tikzcd}$}}
          & $X \ceq \big((\alpha\beta\gamma\delta)^r\big)$,\newline
          $Y \ceq \big((\gamma\delta\alpha\beta)^r\big)$,\newline
          $Z \ceq \big((\beta\gamma\delta\alpha)^r\beta\big)\big(\delta\big)^{-1}$,\newline
          $U \ceq \big( (\delta\alpha\beta\gamma)^r \delta \big)\big(\beta\big)^{-1}$,\newline
          $V \ceq \big( (\delta\alpha\beta\gamma)^r\delta \big)\big( \beta \big)^{-1}\big(\delta\big)\big( (\beta\gamma\delta\alpha)^r \beta \big)^{-1}$
    \end{tabular}
    \end{framed}
\end{enumerate}

\chapter{\sbstrips\ documentation}\label{app:sbstrips-documentation}

Here we include the documentation of the \sbstrips\ package.

\includepdf[fitpaper=true,frame=true,openright=false,pages=2-,scale=0.8]{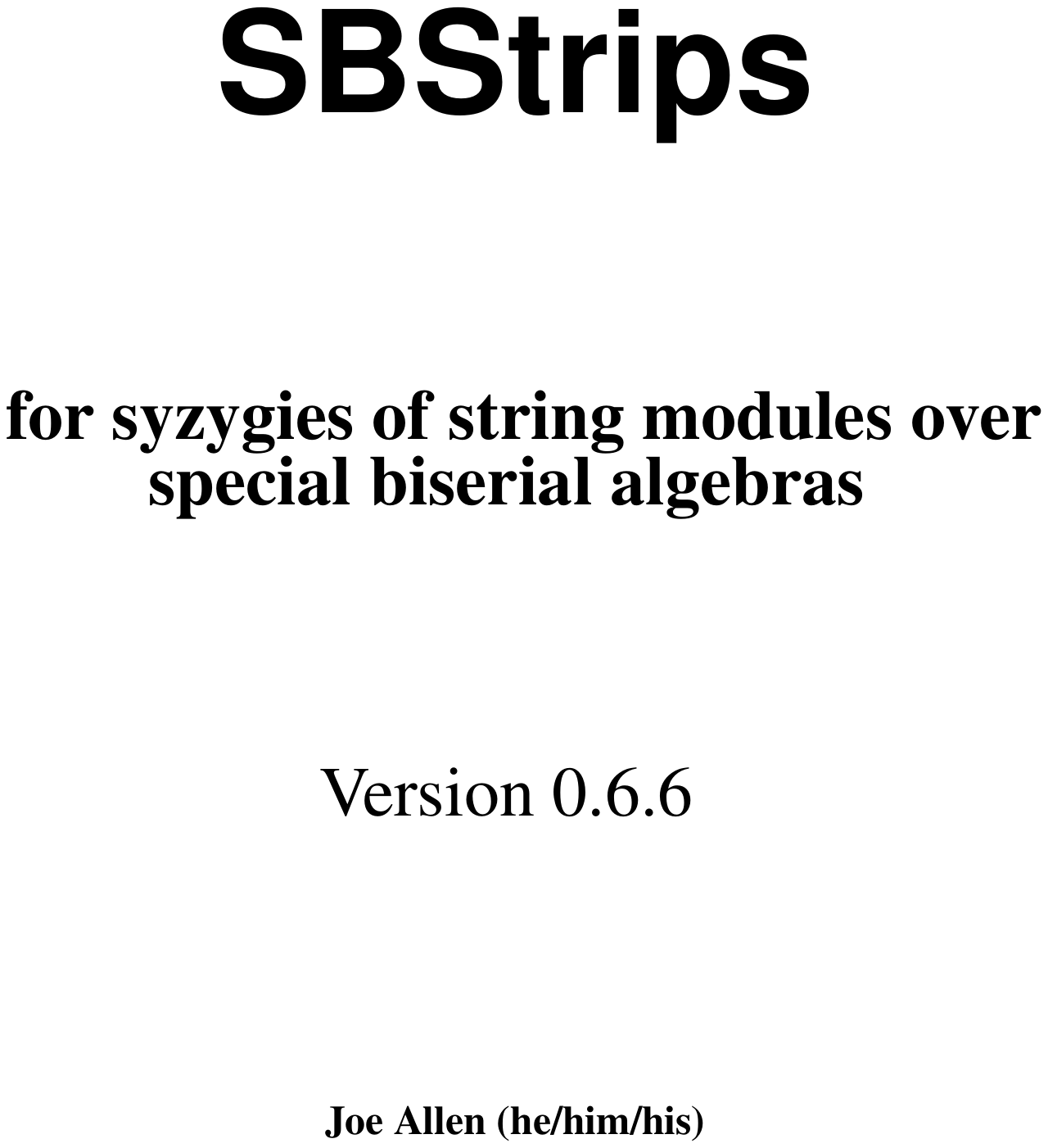}

%%%%%%% BIBLIOGRAPHY

\bibliographystyle{alpha}
\bibliography{main}

\end{document}